\documentclass[12pt]{article}

\usepackage{comment,fancyvrb,xcolor}
\includecomment{comment}\specialcomment{skip}{\begingroup \color{purple}{\it Skip: }}{\endgroup}
 \excludecomment{skip} %remove "%" to exclude remarks

\usepackage{amsmath,enumerate,amsbsy,amsfonts,amssymb,mathabx,amscd,graphicx,algorithm,amstext,enumitem,rotating,caption,subcaption,setspace,comment}

\usepackage{pdflscape}
\RequirePackage[colorlinks,citecolor=blue,urlcolor=blue]{hyperref}
\usepackage[round,authoryear]{natbib}
\usepackage{authblk}
\usepackage{mathrsfs}
\usepackage{appendix}
\usepackage{chngcntr}
\usepackage{epsfig}
\newtheorem{definition}{Definition}
\newtheorem{theorem}{Theorem}

\newtheorem{lemma}{Lemma}

\newtheorem{proposition}{Proposition}
\newtheorem{condition}{Condition}
\newtheorem{remark}{Remark} 

\usepackage{multirow,makecell,colortbl}

\newcommand{\bTheta}{\boldsymbol \Theta}
\newcommand{\bvarsigma}{\boldsymbol \varsigma}
\newcommand{\bvarrho}{\boldsymbol \varrho}

\newcommand{\bXi}{\boldsymbol \Xi}

\newcommand{\bpsi}{\boldsymbol \psi}

\newcommand{\bphi}{\boldsymbol \phi}
\newcommand{\bPhi}{\boldsymbol \Phi}

\newcommand{\btheta}{\boldsymbol \theta}

\newcommand{\bbeta}{\boldsymbol \beta}

\newcommand{\bvarphi}{\boldsymbol \varphi}

\def\T{{\scriptscriptstyle \top} }

\newcommand{\be}{{\mathbf e}}

\newcommand{\bx}{{\mathbf x}}
\newcommand{\by}{{\mathbf y}}
\newcommand{\bu}{{\mathbf u}}
\newcommand{\bv}{{\mathbf v}}

\newcommand{\bg}{{\mathbf g}}

\newcommand{\bb}{{\mathbf b}}
\newcommand{\bz}{{\mathbf z}}

\newcommand{\bB}{{\bf B}}

\newcommand{\bP}{{\bf P}}

\newcommand{\bX}{{\bf X}}
\newcommand{\bY}{{\bf Y}}

\newcommand{\bR}{{\bf R}}

\newcommand{\bW}{{\bf W}}
\newcommand{\bH}{{\bf H}}

\newcommand{\btc}{{\rm c}}

\newcommand{\var}{\text{Var}}

\makeatletter
\DeclareRobustCommand\widecheck[1]{{\mathpalette\@widecheck{#1}}}
\def\@widecheck#1#2{%
	\setbox\z@\hbox{\m@th$#1#2$}%
	\setbox\tw@\hbox{\m@th$#1%
		\widehat{%
			\vrule\@width\z@\@height\ht\z@
			\vrule\@height\z@\@width\wd\z@}$}%
	\dp\tw@-\ht\z@
	\@tempdima\ht\z@ \advance\@tempdima2\ht\tw@ \divide\@tempdima\thr@@
	\setbox\tw@\hbox{%
		\raise\@tempdima\hbox{\scalebox{1}[-1]{\lower\@tempdima\box
				\tw@}}}%
	{\ooalign{\box\tw@ \cr \box\z@}}}
\makeatother

\allowdisplaybreaks[3]

\begin{document}

%\title{\bf \Large Function-on-Function Linear Lagged Regression in High Dimensions, \textcolor{red}{this title needs to be changed}}
\title{\bf \Large Autoregressive Networks with Dependent Edges}

\author[1,2]{\small Jinyuan Chang}
\author[3]{\small Qin Fang}
\author[4,*]{\small Eric D. Kolaczyk}
\author[5]{\small Peter W. MacDonald}
\author[6]{\small Qiwei Yao}
\affil[1]{\it \footnotesize Joint Laboratory of Data Science and Business Intelligence, Institute of Statistical Interdisciplinary Research, Southwestern University of Finance and Economics, Chengdu,  China
}
\affil[2]{\it \footnotesize State Key Laboratory of Mathematical Sciences, Academy of
Mathematics and Systems Science, Chinese Academy of Sciences, Beijing, China}
\affil[3]{\it \footnotesize Business School, The University of Sydney, Sydney, Australia}
\affil[4]{\it \footnotesize Department of Mathematics and Statistics,
McGill University,
Montreal, Canada}
\affil[5]{\it \footnotesize Department of Statistics and Actuarial Science,
University of Waterloo,
Waterloo, Canada}
\affil[6]{\it \footnotesize Department of Statistics, The London School of Economics and Political Science, London,  U.K.}
\affil[*]{\it \footnotesize Address for correspondence: Eric D. Kolaczyk, Department of Mathematics and Statistics, McGill University, Montreal, Canada. Email: eric.kolaczyk@mcgill.ca
}
\setcounter{Maxaffil}{0}
\renewcommand\Affilfont{\itshape\small}
\date{}
\maketitle
\vspace{-1cm}

\begin{abstract}
We propose an autoregressive framework for modelling dynamic networks with dependent edges. It encompasses models that accommodate, for example, transitivity, degree heterogenenity, and other stylized features often observed in real network data. By assuming the edges of networks at each time are independent conditionally on their lagged values, the models, which exhibit a close connection with temporal ERGMs, facilitate both simulation and the maximum likelihood estimation in a straightforward manner. Due to the possibly large number of parameters in the models, the natural MLEs may suffer from slow convergence rates. An improved estimator for each component parameter is proposed based on an iteration employing projection, which mitigates the impact of the other parameters (\citeauthor{ChangChenTangWu2021}, \citeyear{ChangChenTangWu2021}; \citeauthor{ChangShiZhang2023}, \citeyear{ChangShiZhang2023}).
% \citep{ChangChenTangWu2021,ChangShiZhang2023}. 
Leveraging a martingale difference structure, the asymptotic distribution of the improved estimator is derived without the assumption of stationarity. The limiting distribution is not normal in general, although it reduces to normal when the underlying process satisfies some mixing conditions. Illustration with a transitivity model was carried out in both simulation and a real network data set.
\end{abstract}

\noindent {\small{\it Key words}: conditional independence, dynamic networks, maximum likelihood estimation, stylized features of network data, transitivity.}

%\newpage
\renewcommand{\theequation}{\arabic{equation}}
\section{Introduction}
\label{sec.intro}

Dynamic network modelling with dependent edges is  practically important and relevant but technically challenging.
It is natural that the sequence of edges observed over time between two nodes will depend on the behavior of other edges in the network. This dependence represents known social phenomena such as reciprocity, transitivity, degree heterogeneity (popularity), which are thought to govern the evolution of real-world networks in diverse applications.
%Without dependence across different edges, it is impossible to incorporate into the models some stylized features often observed in real network data such as transitivity or degree heterogeneity. 
On the other hand, dependent edges make the dynamic structures of network processes complex, and statistical inference challenging -- not only the computation but also the accompanying theory.

In the growing research area of dynamic network modelling, work continues in three interrelated directions: model specification, computation, and theory. 
%Model specification has progressed rapidly, with frameworks like (separable) temporal exponential random graph models \citep[(S)TERGM,][]{hanneke2010, krivitsky2014} and stochastic actor-oriented models \citep[SAOM,][]{snijders2017stochastic,koskinen2023multilevel} having great flexibility to specify interpretable, mechanistic models with edge dependence. 
%Computation and implementation have followed, with the aforementioned frameworks now having well-developed software packages \citep{statnet_web,ripley2025manual}, which scale to varying extents. 
Our work contributes model specification, computation, and substantial and highly non-trivial supporting theory.
Specifically, we develop a class of dependent-edge autoregressive (dep. edge AR) models for network data observed in discrete time. This class is amenable to theoretical analysis, with formal asymptotic analysis for the maximum likelihood estimators.
Our theoretical results in turn motivate a statistically and computationally efficient approach to parameter estimation, which we implement in an open-source software package.

We now briefly review some existing models for dynamic networks. We provide a more thorough review in Section~\ref{app:extendedlit} of the supplementary material.
Under our discrete-time view, we assume that we observe time-indexed network {\em snapshots} at regular intervals, or can coerce relational event data into network snapshots sampled at a constant time resolution.
When relational event data with continuous time stamps are available, classes like stochastic actor-oriented models (SAOM) for edge formation and dissolution events \citep{snijders2017stochastic}, and relational event models (network-structured multivariate point process models) for instantaneous relational events \citep[]{PerryWolfe2013,matias2018semiparametric,butts2023relational}
exist, and are well-developed, flexible, and computationally robust.
%From our discrete-time viewpoint, with the exception of SAOMs, which have been used to analyze discretely observed network snapshots, we do not compare further to these approaches here.
The exponential random graph model \citep[ERGM,][]{robins2007} has been extended to the dynamic setting in seminal work on (separable) temporal ERGMs \citep[(S)TERGM,][]{hanneke2010,krivitsky2014}. 
(S)TERGMs, along with a related class of logistic network regression models \citep[LNR,][]{almquist2014logistic} allow for flexible model specification through parametric models which specify the transition
distributions of network snapshots.
Previous works on AR network models instead specify transition distributions for individual edge variables, assuming a form of {\em edgewise dependence}: the only dependence between snapshots is between edge variables with the same (unordered) indices $\{i,j\}$ \citep{Jiang2023, Jiang2025}.
Finally, there is a large body of work which extends existing network models like the stochastic block model (SBM) and latent space model (LSM) to the dynamic setting based on a smooth or stochastically evolving sequence of latent
variables \citep{sewell2015latent,Durante2016,matias2017,pensky2019spectral,gallagher2021spectral,suveges2023,zhang2024change}.
While latent variable models do induce edge dependence both over time and within each snapshot, this dependence
is implicit and cannot be configured easily to accommodate stylized features of real network
data (e.g. transitivity).

Our new modelling framework directly extends the AR network models of \cite{Jiang2023} and \cite{Jiang2025}.
We specify the transition probabilities of forming a new edge or dissolving an existing edge between each pair of nodes explicitly depending on its history, and on the histories of other edge processes. This enlarged form of the transition probabilities makes the model flexible enough to accommodate stylized features such as transitivity and degree heterogeneity. The resulting network processes have dependent edges, which is radically different from previous AR network models.
%considered in \citeauthor{Jiang2023} \citeyearpar{Jiang2023} and \citeauthor{Jiang2025} \citeyearpar{Jiang2025}.
%Similar to \cite{hanneke2010} and \cite{almquist2014logistic}, we assume that the edges are conditionally independent given their joint histories. This makes not only statistical inference more transparent but also the corresponding theoretical analysis feasible. 
%This conditional independence avoids the difficulties caused by the normalized constants in ERGMs; see \cite{hanneke2010} and \cite{Leifeld2018}.
Nevertheless, based on a conditional independence assumption similar to \cite{hanneke2010} and \cite{almquist2014logistic}, we can build up a martingale difference structure which facilitates asymptotic analysis for the maximum likelihood estimators (see Section \ref{sec:estimation}).
Additionally, we allow the number of parameters in the model to diverge together with the network size to accommodate, for example, heterogeneity in forming and dissolving edges modelled with node-specific parameters.
% , even when the number of nodes in the network diverges at a faster rate than the number of observed networks.
%To restate our contribution, we have developed a modelling framework that substantially pushes the boundaries of statistical theory regarding parameter estimation and inference for dynamic network models with edge dependence.
%Our modelling framework greatly expands on previous AR network model classes developed by \citeauthor{Jiang2023} \citeyearpar{Jiang2023} and \citeauthor{Jiang2025} \citeyearpar{Jiang2025}. 
To the best of our knowledge, high-dimensional asymptotic analysis of this kind has not been done for a discrete-time dynamic network model class as flexible as the one we specify in this work. 
%Nor does similar analysis exist for more flexible (S)TERGMs or SAOMs. 
We visualize the relative scope of this advancement in Table~\ref{tab:model_classes}, which is explained in full detail in Section~\ref{app:extendedlit} of the supplementary material.

\begin{table}[h] 
\begin{tabular}{ccccc}
\multirow{4}{*}{\makecell{{\bf Between} \\ {\bf snapshot}}} & \multicolumn{1}{c|}{\textit{Full dep.}}        & \cellcolor[gray]{0.85}{\bf dep. edge AR}, LNR &                     & (S)TERGM, SAOM  \\
                                  & \multicolumn{1}{c|}{\textit{Local dep.}}       & \cellcolor[gray]{0.85}                                          &                     &                    \\
                                  & \multicolumn{1}{c|}{\textit{Edgewise dep.}}    & \cellcolor[gray]{0.85} AR  &                     &                    \\
                                  & \multicolumn{1}{c|}{\textit{Static/Full ind.}} & \cellcolor[gray]{0.85} ER                                  & LD-ERGM   & ERGM           \\ \cline{3-5}
                                  &                                                & \textit{Full ind.}                        & \textit{Local dep.} & \textit{Full dep.} \\
                                  &                                                & \multicolumn{3}{c}{{\bf Within snapshot}}
\end{tabular}
\caption{Classification of marginal edge dependence structures in some existing static and dynamic network models (see Section~\ref{app:extendedlit}). Between snapshot dependence structures are contained in one another from low to high. Within snapshot dependence structures are contained in one another from left to right. Our model class, ``{\bf dep. edge AR}'', can accommodate dependence structures in the entire shaded region. \label{tab:model_classes}}
\end{table}

%The rest of the paper is organized as follows. Section~\ref{sec:lit_review} provides a detailed (albeit necessarily still incomplete) review of the literature.  
Section \ref{sec:sec2} presents the general AR network framework with dependent edges. 
%The stationarity and ergodicity for AR(1) models are explored based on the fundamental limittheorem for homogeneous Markov chains. 
We also discuss the relationship between the proposed AR models, and TERGMs in Section \ref{sec:sec2}. Section \ref{sec:3examples} contains three concrete AR models which are designed to model, respectively, degree heterogeneity, persistence, and transitivity -- those are among the stylized features often observed in real network data. Section \ref{sec:estimation} presents the estimation procedures based on the maximum likelihood principle for the parameters in the AR models and the associated asymptotic theory. In so doing, we introduce the concepts of local parameters and global parameters, which need to be identified and estimated separately and may also entertain different convergence rates. 
%The initial maximum likelihood estimator (MLE) suffers from slower convergence rates due to the diverging number of local parameters. 
An improved estimator for each component parameter is obtained by projecting the score function onto the corresponding direction, which mitigates the impact of the other parameters (\citeauthor{ChangChenTangWu2021}, \citeyear{ChangChenTangWu2021}; \citeauthor{ChangShiZhang2023}, \citeyear{ChangShiZhang2023}).
%Based on a martingale difference structure, the asymptotic distribution of the improved estimator is derived without the assumption of stationarity. 
The limiting distribution of the improved estimator is derived without the assumption of stationarity. It is not normal in general, but it reduces to normal when the underlying process satisfies some mixing conditions which hold for many stationary processes.
Numerical illustration with a real dynamic network data set is reported in Section \ref{sec:real}.
An online supplement contains simulation studies for the proposed transitivity model, additional real data results and all the technical proofs. 
  
{\it Notation.}  For any positive integer $r$, write $[r] = \{1,\dots, r\}$ and $\mathbb{R}_+^r=\{(x_1,\ldots,x_r)^{\T}:x_i>0~\textrm{for any}~i\in[r]\}$. For any $x, y \in \mathbb{R},$ we write $x \vee y = \max(x,y).$  For two positive sequences $\{a_n\}$ and $\{b_n\}$, we write $a_n\ll b_n$ if  $\lim \sup_{n \to \infty} a_n/b_n = 0$, and 
 $a_n \lesssim b_n$ if $\lim \sup_{n \to \infty} a_n/b_n <\infty$. 
For any vector $\bb = (b_1, \dots, b_r)^{\T} \in \mathbb{R}^r,$ we let $\bb_{-l}$ denote the sub-vector of $\bb$ by removing the $l$-th component $b_l$. Given an index set $\mathcal{M} \subset [r] $, we let $\bb_{\mathcal{M}}$ denote the sub-vector of $\bb$ that consists of the components of $\bb$ with the indices in $\mathcal{M}$.
For any $r_1 \times r_2$ real matrix $\bB$, denote by $\bB^{\T}$ its transpose. %When $r_1 = r_2$, 
We use $\lambda_{\min}(\bB)$ to denote  the smallest eigenvalue of a square matrix $\bB$. For any countable set
$\mathcal{U}$, $|\mathcal{U}|$ denotes its cardinality.

\section{AR($m$) network framework} \label{sec:sec2}
\subsection{Model}
\label{sec:models}
Consider a dynamic network process defined on  $p$ nodes denoted by $1, \ldots, p$. Let $\bX_t \equiv (X_{i,j}^t)_{p\times p}$ be the adjacency matrix at time $t$, where $X_{i,j}^t =1$ denotes the existence of an edge between nodes $i$ and $j$ at time $t$, and $X_{i,j}^t=0$ otherwise. For simplicity, we only consider undirected networks without self-loops, i.e., $X_{i,i}^t \equiv 0$ and $X_{i,j}^t= X_{j,i}^t$. The main idea can be applied to directed networks.

\begin{definition}[AR($m$) networks] \label{def_AR}
 Conditionally on $\{\bX_{s}\}_{s\leq t-1}$,
the edges $\{X_{i,j}^t\}_{1\le i <j\le p}$ are mutually independent with 
\begin{align} \label{alphaij}
\alpha_{i,j}^{t-1} & \equiv \mathbb{P}(X_{i,j}^t=1\,|\,X_{i,j}^{t-1} =0,  \bX_{t-1}\setminus X_{i,j}^{t-1}, \bX_{t-k}~{\rm for}~k\ge 2)\nonumber\\ 
&= \mathbb{P}(X_{i,j}^t=1\,|\,X_{i,j}^{t-1} =0,  \bX_{t-1}\setminus X_{i,j}^{t-1}, \bX_{t-2}, \ldots, \bX_{t-m})\,,\\
    \label{betaij}
\beta_{i,j}^{t-1} &\equiv \mathbb{P}(X_{i,j}^t=0\,|\,X_{i,j}^{t-1} =1,  \bX_{t-1}\setminus X_{i,j}^{t-1}, \bX_{t-k}~{\rm for}~k\ge 2)\nonumber\\  &
= \mathbb{P}(X_{i,j}^t=0\,|\,X_{i,j}^{t-1} =1,  \bX_{t-1}\setminus X_{i,j}^{t-1},\bX_{t-2}, \ldots, \bX_{t-m})\,,
\end{align}
where $m\ge 1$ is an integer. %, and $\bZ_{i,j}$ are some exogenous variables.
\end{definition}

An AR($m$) network process defined above is a Markov chain with
order $m$. Based on \eqref{alphaij} and \eqref{betaij}, 
it holds that for $x\in\{0, 1\}$,
\[
\mathbb{P}(X_{i,j}^t=x\,|\,\bX_{t-1}, \ldots, \bX_{t-m}) =
(1 - \gamma_{i,j}^{t-1})^{1-x}( \gamma_{i,j}^{t-1})^x\,,
\]
where
\begin{equation} \label{def1}
\gamma_{i,j}^{t-1} = \alpha_{i,j}^{t-1} +X_{i,j}^{t-1} (1 - \alpha_{i,j}^{t-1} -\beta_{i,j}^{t-1})
=\mathbb{P}(X_{i,j}^t=1\,|\,\bX_{t-1}, \ldots, \bX_{t-m})\,,
\end{equation}
% &\mathbb{P}(X_{i,j}^t=1\,|\,\bX_{t-1}, \ldots, \bX_{t-m}) =
%\alpha_{i,j}^{t-1} +X_{i,j}^{t-1} (1 - \alpha_{i,j}^{t-1} -\beta_{i,j}^{t-1}) \equiv \gamma_{i,j}^{t-1}\,,
%\\
%& = 1 - P(X_{i,j}^t=0|\bZ_{i,j}, \bX^{t-1}, \cdots, \bX^{t-m}).
%\end{align}
i.e.,
$
X_{i,j}^t\,|\,\bX_{t-1},\ldots,\bX_{t-m}\sim {\rm Bernoulli}(\gamma_{i,j}^{t-1})\,,\; 1\leq i<j\leq p\,.
$
Clearly % $\{\bX_t\}_{t\geq1}$ is not a homogeneous Markov chain, and 
edges
$X_{i,j}^t$, for different $(i,j)$, are not necessarily independent of each other in practice. We may impose various forms for the conditional probabilities $\alpha_{i,j}^{t-1}$ and $\beta_{i,j}^{t-1}$ to reflect different stylized features of network data that capture such dependency. Put 
\begin{equation}\label{eq:alphaijbetaij}
\begin{split}
\alpha_{i,j}^{t-1} = f_{i,j}( \bX_{t-1}\setminus X_{i,j}^{t-1}, \bX_{t-2},\ldots, \bX_{t-m}; \btheta_0)\,,\\
\beta_{i,j}^{t-1} = g_{i,j} (\bX_{t-1}\setminus X_{i,j}^{t-1}, \bX_{t-2}, \ldots, \bX_{t-m}; \btheta_0)\,,
\end{split}
\end{equation}
where $f_{i,j}$'s and $g_{i,j}$'s are known functions, and $\btheta_0\in\bTheta\subset\mathbb{R}^q$ is a $q$-dimensional unknown true parameter vector. For any $\btheta\in\bTheta$, write 
\begin{equation}\nonumber 
% \label{eq:alphaijbetaij:theta}
\begin{split}
\alpha_{i,j}^{t-1}(\btheta) = f_{i,j}(\bX_{t-1}\setminus X_{i,j}^{t-1}, \bX_{t-2}, \ldots, \bX_{t-m}; \btheta)\,,\\
\beta_{i,j}^{t-1}(\btheta) = g_{i,j} (\bX_{t-1}\setminus X_{i,j}^{t-1}, \bX_{t-2}, \ldots, \bX_{t-m}; \btheta)\,.
\end{split}
\end{equation}
Then $\alpha_{i,j}^{t-1}=\alpha_{i,j}^{t-1}(\btheta_0)$ and $\beta_{i,j}^{t-1}=\beta_{i,j}^{t-1}(\btheta_0)$.

 Modelling dynamic networks by Markov or/and AR models is not new. See, for example, \cite{snijder2005}, \cite{ludkin18}, \cite{yang2011}, \cite{yudovina2015}, and \cite{Jiang2023}. However, most available Markov models are designed for  Erd\"os-Renyi networks with independent edges. Our setting provides a general framework to accommodate various dependency structures across different edges. Some practical network models satisfying this general framework are introduced in Section \ref{sec:3examples}.

%\subsection{Stationarity of AR(1) network processes} %\label{sec:stationarity}

For the special AR(1) processes (i.e., $m=1$), if both $f_{i,j}$ and $g_{i,j}$ in (\ref{eq:alphaijbetaij}) are always positive and smaller than 1 for all $1\le i < j\le p$, $\{\bX_t\}_{t\geq1}$ is an irreducible homogeneous Markov chain with $2^{p(p-1)/2}$ states. Therefore when $p$ is  fixed,  (i) there exists a unique stationary distribution, and (ii)  if $\bX_0$ is activated according to this stationary distribution, the process $\{ \bX_t\}_{t\geq 1}$ is strictly stationary and ergodic. See Theorems 3.1, 3.3 and 4.1 in Chapter 3 of \cite{bremaud1998}.  
Hence the degree heterogeneity model introduced below in Section \ref{sec:density} 
and the transitivity model introduced in Section \ref{sec:transitivity} are strictly stationary for any fixed constant $p$ if all the transition probability functions $\alpha_{i,j}^{t-1}$ and $\beta_{i,j}^{t-1}$ are strictly between 0 and 1. 
It is worth pointing out that the ergodicity only holds for any fixed constant $p$. Hence we cannot take for granted that the sample means of $\bX_t$ and/or its summary statistics converge when $p$ diverges together with the sample size, even when $\bX_t$ is stationary. This causes the major theoretical challenges. Note that  stationarity is  not an asymptotic
property while ergodicity~is. 
%{\color{red}add some discussion on the technical assumption of other papers. whether these papers require stationarity and endogenity?} 

\subsection{Relationship to  TERGMs}
\label{sec: relation2ERGM}
Our model assumes that the edges are conditionally independent given their lagged values. 
While general TERGMs can specify models with full within-snapshot dependence, a similar conditional independence assumption is made by \cite{hanneke2010}, to ensure non-degeneracy of the network transition distribution.
Compared to this restricted class of TERGMs studied by \cite{hanneke2010}, instead of requiring the transitions to follow an exponential family distribution, we may flexibly define the probability for forming a new edge in (\ref{alphaij}), and dissolving an existing edge in (\ref{betaij}). Those two probability functions can be in any desired forms as presented in (\ref{eq:alphaijbetaij}), provided that their values lie between $0$ and $1$.
%In other words, we may impose closed-form parametric functions for $\alpha_{i,j}^{t-1}$ and $\beta_{i,j}^{t-1}$, and those functions need only be between 0 and 1. %Hence, the likelihood functions are explicitly available, 
%This %not only makes both estimation and simulation based on our models directly feasible, but also 
%which allows us to depict more explicitly in our models some stylized features often observed in real network data. 
See Section \ref{sec:3examples} for examples. The numerical analyses 
in Section~\ref{sec:real} and Section~\ref{sec:simulation} of the supplementary material indicate that the proposed AR models are capable of simulating and reflecting some observed interesting dynamic network phenomena.

A TERGM \citep{hanneke2010, krivitsky2014} is specified by
\begin{equation*}
  \mathbb{P}({\bf X}_t \,|\, {\bf X}_{t-1}, \ldots, {\bf X}_{t-m} ; \btheta)\, \propto \exp\{ \bvarsigma(\btheta)^{\T} \bvarrho({\bf X}_t, {\bf X}_{t-1}, \ldots, {\bf X}_{t-m})\}\,,
\end{equation*}
where $\bvarsigma$ maps the underlying unknown parameters to a natural parameter, and $\bvarrho$ maps the present and past network configurations to a low-dimensional vector of sufficient statistics.
Under the additional assumption of conditional edge independence, $\bvarrho$ will decompose into a sum over the edges, and any TERGM can be expressed in our AR($m$) framework as
\begin{equation} \label{tergm_form}
\begin{split}
\alpha_{i,j}^{t-1}(\btheta) =&~ \frac{\exp\{\bphi(\btheta)^{\T} \bu_{i,j}({\bf X}_{t-1}
\setminus X_{i,j}^{t-1}, {\bf X}_{t-2}, \ldots, {\bf X}_{t-m} )\} }{
1+ \exp\{\bphi(\btheta)^{\T} \bu_{i,j}({\bf X}_{t-1} \setminus X_{i,j}^{t-1}, {\bf
X}_{t-2}, \ldots, {\bf X}_{t-m} )\}}\,,  \\
\beta_{i,j}^{t-1}(\btheta) =&~ \frac{\exp\{\bpsi(\btheta)^{\T} \bv_{i,j}({\bf X}_{t-1} \setminus X_{i,j}^{t-1}, {\bf X}_{t-2}, \ldots, {\bf X}_{t-m} )\}}
{1+ \exp\{\bpsi(\btheta)^{\T} \bv_{i,j}({\bf X}_{t-1} \setminus X_{i,j}^{t-1}, {\bf X}_{t-2}, \ldots, {\bf X}_{t-m} )\}}\,, 
\end{split}
\end{equation}
where $\bu_{i,j}(\cdot)$ and $\bv_{i,j}(\cdot)$ have closed-form expressions depending on the sufficient statistics of the original TERGM. 
Notice that the possible forms of the transition probabilities in \eqref{tergm_form} are restricted: one must take a linear combination of the natural parameter and sufficient statistics, and pass through an inverse logistic link function.

If $\bphi(\btheta)$ and $\bpsi(\btheta)$ in \eqref{tergm_form} are replaced by, respectively, 
$\bphi(\btheta_\alpha)$ and $\bpsi(\btheta_\beta)$, where $\btheta_\alpha$
and $\btheta_\beta$ are two sets of different parameters, $\bX_t$ follows a  STERGM with conditional independent
edges given the past networks \citep{krivitsky2014}. See
 Section \ref{app:relation2ERGM} of the supplementary material for detailed discussion and derivations regarding the relationship of the proposed dependent-edge AR models and TERGMs. 

\vspace{-0.2cm}
\section{Some interesting AR network models} \label{sec:3examples}
\vspace{-0.2cm}
To illustrate the usefulness of the AR($m$) framework proposed above, we state three AR($m$) network models which
reflect various stylized features in real network data.  In all three models, the parameters $\{ \xi_i\}_{i=1}^p$ and
$\{\eta_i\}_{i=1}^p$ reflect node heterogeneity in, respectively, forming a new edge and dissolving an existing edge. Specifically, the larger $\xi_i$ is, the more likely node $i$ will form new edges with other nodes, and the larger $\eta_i$ is, the more likely the existing edges between node $i$ and the others will be dissolved.
Instances of these three models can be simulated using our development R package \texttt{arnetworks}. 
Our package provides an implementation of our maximum-likelihood-based parameter estimation and inference procedure for the transitivity model (Section~\ref{sec:transitivity}), as described in Section~\ref{sec:estimation} and Section~\ref{sec:simulation} of the supplementary material,  and a general method-of-moments-based estimation approach introduced in Section~\ref{sec:estimation2} of the supplementary material, which can be applied to a broad class of dynamic network models, including but not limited to the three models discussed below. More details are available at \url{https://github.com/peterwmacd/arnetworks}.

\vspace{-0.3cm}
\subsection{Degree heterogeneity model}  \label{sec:density}
\vspace{-0.1cm}
For any $i\neq j$, let 
\begin{align*}
\vartheta_{i,j}^{t-1} &= \exp\{ a_0 D_{-i,-j}^{t-1}
    + a_1 (D_{i}^{t-1} +  D_{j}^{t-1})\}\,,\\
\varpi_{i,j}^{t-1} &=  \exp\{ b_0 (1-D_{-i,-j}^{t-1})
    + b_1 (2-D_{i}^{t-1}- D_{j}^{t-1})\}\,,
    \end{align*}
with
\begin{align*} 
D_{-i,-j}^{t-1}= \frac{1}{(p-2)(p-3)}\sum_{k,\ell:\, k,\ell \neq i,j,\, k \neq \ell} X_{k, \ell}^{t-1}\,,\qquad D_i^{t-1}= \frac{1}{p-1}\sum_{\ell:\,\ell \ne i} X_{i, \ell}^{t-1}\,,
\end{align*}
where $D_i^{t-1}$ and
$D_j^{t-1}$ are, respectively, the (normalized) degrees of node $i$ and node $j$ at time $t-1$, and $D_{-i,-j}^{t-1}$ is the (normalized) average degree excluding nodes $i$ and $j$ at time $t-1$.
We specify the transition probabilities % in (\ref{alphaij}) and (\ref{betaij}) as $\alpha_{i,j}^{t-1}=\alpha_{i,j}^{t-1}(\btheta_0)$ and $\beta_{i,j}^{t-1}=\beta_{i,j}^{t-1}(\btheta_0)$ with
as follows:
\begin{align} \label{densityalphaij}
    \alpha_{i,j}^{t-1}(\btheta) = ~ { \xi_i \xi_j \vartheta_{i,j}^{t-1}
   \over 1 +\vartheta_{i,j}^{t-1} + \varpi_{i,j}^{t-1}} \,,
    \qquad
    \beta_{i,j}^{t-1}(\btheta) =~   { \eta_i \eta_j \varpi_{i,j}^{t-1}
   \over 1 +\vartheta_{i,j}^{t-1} + \varpi_{i,j}^{t-1}}\,.
     \end{align} 
This is an AR(1) model with parameter vector $\btheta= (a_0, a_1,  b_0, b_1,  \xi_1, \ldots, \xi_p, \eta_1, \ldots, \eta_p)^{\T}\in\bTheta\subset\mathbb{R}_+^{2p+4}$. %, where $a_0$, $a_1$,  $ b_0$, $ b_1$,  $\xi_1, \ldots, \xi_p$, $\eta_1, \ldots, \eta_p$ are positive constants.
This model is able to capture a ``rich get richer'' effect in network evolution, whereby nodes with higher degree at time $t-1$ are more likely to form new connections at time $t$, and vice versa.
In particular, the propensity to form a new edge between nodes $i$ and $j$ at time $t$ is positively impacted by 
$D_{-i,-j}^{t-1}$, $D_i^{t-1}$ and $D_j^{t-1}$, and the propensity to dissolve an existing edge between nodes $i$ and $j$ at time $t$ is negatively impacted by these three quantities. 

\cite{hanneke2010} proposed a TERGM including a density statistic (equivalent to network average degree). In (\ref{densityalphaij}), we explicitly specify the impact from the previous snapshot's average degree on forming a new edge and dissolving an existing edge, while the model defined in Section 2.1 of \cite{hanneke2010} depends on the current snapshot's average degree, and does not differentiate the representations for these two types of impact. 
Within the STERGM framework, the edge counts model of \cite{krivitsky2014} assumes that the collection of all newly formed edges is conditionally independent of the collection of all newly dissolved edges given their history, and the two conditional distributions are controlled by different parameters.

\subsection{Persistence model} \label{sec:persistence}
We define the transition probabilities %in (\ref{alphaij}) and (\ref{betaij}) as $\alpha_{i,j}^{t-1}=\alpha_{i,j}^{t-1}(\btheta_0)$ and $\beta_{i,j}^{t-1}=\beta_{i,j}^{t-1}(\btheta_0)$ with
%\begin{align*}
 %   \alpha_{i,j}^{t-1}=& \theta_i \theta_j  \exp [- a\{ 1+ (1-X_{i,j}^{t-2}) +(1-X_{i,j}^{t-2})(1-X_{i,j}^{t-3})\}  ],\\[1ex]  
 %   \beta_{i,j}^{t-1} = & \eta_i \eta_j \exp \{- b( 1 + X_{i,j}^{t-2}+X_{i,j}^{t-2}  X_{i,j}^{t-3})\} .
%\end{align*}
%{\color{green} An alternative form:
\begin{equation} \label{eq:persistence}
\begin{split}
    \alpha_{i,j}^{t-1}(\btheta)=&~ \xi_i \xi_j  \exp [-  1 - a\{ (1-X_{i,j}^{t-2}) +(1-X_{i,j}^{t-2})(1-X_{i,j}^{t-3}) \} ]\,,\\  
    \beta_{i,j}^{t-1}(\btheta) = &~ \eta_i \eta_j \exp \{-  1 - b( X_{i,j}^{t-2}+X_{i,j}^{t-2}  X_{i,j}^{t-3})\} \,. 
    \end{split}
\end{equation}
%Now we can let $\xi_i, \eta_i \ge 0$ only ...}\\
This is an AR(3) model with parameter vector $\btheta = (a, b, \xi_1, \ldots, \xi_p,
\eta_1, \ldots, \eta_p)^{\T}\in\bTheta\subset\mathbb{R}_+^{2p+2}$. %, where $a, b,$ $\xi_1, \ldots, \xi_p, \eta_1, \ldots, \eta_p$ are positive constants.
The probability to form a new edge between nodes $i$ and $j$ at time $t$ is reduced if $X_{i,j}^{t-2} =0$, and it is reduced further if, in addition, $X_{i,j}^{t-3}=0$. 
The probability to dissolve an existing edge is reduced if
$X_{i,j}^{t-2} =1$, and it is reduced further if, in addition, $X_{i,j}^{t-3}=1$. Hence if the edge status between two nodes is unchanged for 2 or 3 time periods, the probability for it remaining unchanged next time is larger than that otherwise. 

Model (\ref{eq:persistence}) defines an AR(3) network process $\bX_t = (X_{i,j}^t)_{p\times p}$ with $p(p-1)/2$ independent edge processes. Although the conclusion on the AR(1) stationarity in the last paragraph of Section \ref{sec:models} does not apply directly, this AR(3) network process is also strictly stationary, which is implied by the fact that $\{X_{i,j}^t\}_{t\ge 1}$ is strictly stationary for each $1\le i<j\le p$. Formally, for given $(i,j)$ such that $1\leq i<j\leq p$, let
$\bY_t = (X_{i,j}^t, X_{i,j}^{t-1}, X_{i,j}^{t-2})^\T$.
Then $\{\bY_t\}_{t\geq1}$ is a homogeneous Markov chain with $2^3 =8$ states. Let $\bP$ denote
the transition probability matrix of $\{\bY_t\}_{t\geq1}$. Then 
$\bP$ is a $8\times 8$ matrix with only 2 positive elements in each row and each column, provided that $\xi_i \xi_j, \, \eta_i \eta_j \in (0, e)$. It is straightforward to check that each row or column of $\bP^2 $ has only 4 positive elements, and, more importantly, all the elements of $\bP^3$ is positive. Hence, the Markov chain $\{\bY_t\}_{t\geq1}$ is irreducible. By Theorems 3.1 and 3.3 in Chapter 3 of \cite{bremaud1998},  the process $\{\bY_t\}_{t\geq1}$ is strictly stationary, and so is $\{X_{i,j}^t\}_{t\geq1}$.

Persistent connectivity or non-connectivity is widely observed in, for example, brain networks, gene connections and social networks. A related TERGM including  a stability statistic defined in \cite{hanneke2010}, does not consider lags of order 2 and 3, and does not differentiate between the propensity for retaining an existing edge and that for retaining a no-edge status.

\subsection{Transitivity model} \label{sec:transitivity}

We propose an AR(1) model to reflect the feature of transitivity which
refers to the phenomenon that nodes are more likely to link if they share
links in common (i.e., ``the friend of my friend is also my friend'').
 To this end, we specify the transition probabilities as follows:
 %we let $\alpha_{i,j}^{t-1}=\alpha_{i,j}^{t-1}(\btheta_0)$ and $\beta_{i,j}^{t-1}=\beta_{i,j}^{t-1}(\btheta_0)$ in (\ref{alphaij}) and (\ref{betaij}) with 
\begin{equation} \label{trans-model}
\begin{split}
    \alpha_{i,j}^{t-1}(\btheta) &= \frac{\xi_i \xi_j\exp(a U_{i,j}^{t-1}) }{1+ \exp(a U_{i,j}^{t-1})+ \exp(b V_{i,j}^{t-1})}\,,\\
    \beta_{i,j}^{t-1}(\btheta) &=\frac{ \eta_i \eta_j \exp(b V_{i,j}^{t-1})}  {1+ \exp(a U_{i,j}^{t-1})+ \exp(b V_{i,j}^{t-1})}\,,
    \end{split}
\end{equation}
where $\btheta = (a, b, \xi_1, \ldots, \xi_p, \eta_1, \ldots, \eta_p)^\top\in\mathbb{R}_+^{2p+2}$, and
\begin{equation}
    \label{UVij}
    \begin{split}
U_{i,j}^{t-1} =&~ {1\over p-2}\sum_{k:\,k\neq i,j} X_{i,k}^{t-1} X_{j,k}^{t-1}\,,\\ V_{i,j}^{t-1} =&~{1 \over p-2} \sum_{k:\,k\neq i, j}\{X_{i,k}^{t-1}(1-X_{j,k}^{t-1})+ (1-X_{i,k}^{t-1})X_{j,k}^{t-1}\}\,.
\end{split}
\end{equation}
%{\color{red} (I have changed the divisor from $(p-1)$  to $(p-2)$. Both $p$ or $(p-2)$ are more reasonable than $(p-1)$ -- QY.)}\\
 The pair $(U_{i,j}^{t-1}, V_{i,j}^{t-1})$ characterizes the  number of
nodes with which both nodes $i$ and $j$ are connected,  and the  number of nodes
with which only one of $i$ and $j$ is connected at  time $t-1$.  The larger $U_{i,j}^{t-1}$ is (i.e., the more common friends $i$ and $j$ share at time $t-1$), the more likely $X_{i,j}^{t} = 1$. The larger $V_{i,j}^{t-1}$ is, the more likely $X_{i,j}^{t} = 0$. This reflects the transitivity of the networks. High levels of transitivity are found in various networks including friendship networks, industrial supply-chains, international trade flows, and alliances across firms and nations.   Note that the quantity $U_{i,j}^{t-1}$, used in \cite{graham2016} to define the edge status of $X_{i,j}^t$, reflects
the information based on which companies such as Facebook and LinkedIn have recommended new links to their customers. Also see the TERGM including a transitivity statistic defined in \cite{hanneke2010}.

We may use different parameters $a$ and $b$ in defining $\alpha_{i,j}^{t-1}(\btheta)$ and $\beta_{i,j}^{t-1}(\btheta)$ in \eqref{trans-model}. We do not pursue this more general form as (i) using different $\xi_i$ and $\eta_i$  reflects already the differences in the propensity between forming a new edge and dissolving an existing edge, and, perhaps more importantly, (ii) since most practical networks are sparse, the effective sample size for estimating the transition probability from the state of an existing edge is small. Therefore estimating the parameters only occurring in $\beta_{i,j}^{t-1}(\btheta)$ will be harder than those in $\alpha_{i,j}^{t-1}(\btheta)$. Using the same $a$ and $b$ in both $\alpha_{i,j}^{t-1}(\btheta)$ and $\beta_{i,j}^{t-1}(\btheta)$ improves the estimation by pulling the information together. See also the relevant simulation results in Section \ref{sec:more.gen.model} of the supplementary material.

\section{Estimation} \label{sec:estimation}

%{\color{blue}
In this section, we present estimation procedures based on the maximum likelihood principle, and show they are valid if the number of nodes $p$ and the sample size $n$ satisfy the restriction $p\ll n^{\delta}$ for some constant $\delta>0$, which allows the number of nodes $p$ to either be fixed or diverge together with $n$.  
When $p$ diverges with $n$, both the ergodicity and the central limit theorem for stationary Markov chains no longer apply even when $\bX_t$ is stationary  (see the last paragraph in Section \ref{sec:models}). This causes the major theoretical challenges faced here. Based on the conditional independence in our setting,  we can construct an appropriate martingale difference sequence, which facilitates %appropriate partial sums of which are amendable to 
the required asymptotic analysis without the stationarity assumption, regardless of whether $p$ is fixed or diverges together with $n$. 
%An alternative and computationally more efficient estimation procedure based on the method of moments principle will be presented in Section \ref{sec:estimation2}.}

\subsection{General approach}\label{sec:gen.appr}

The natural units of observation in our model are the $X^t_{i,j}$, indicating presence or absence of an edge between nodes $i$ and $j$ at time $t$.  Intuitively, the extent to which these observations can contribute useful information to the estimation of a given element of $\btheta_0$ depends in turn on the extent to which that element plays a consistent role over time $t$ in the corresponding probabilities
$$
\gamma_{i,j}^{t-1}(\btheta) =\alpha_{i,j}^{t-1}(\btheta) +X_{i,j}^{t-1}\{1 - \alpha_{i,j}^{t-1}(\btheta) -\beta_{i,j}^{t-1}(\btheta)\}\,.
$$ 
By (\ref{def1}), we have $\gamma_{i,j}^{t-1}=\gamma_{i,j}^{t-1}(\btheta_0)$. We formalize the above intuition as follows.
\begin{definition}[Global/local parameters] 
Write $\btheta=(\theta_1,\ldots,\theta_q)^{\T}$, where $q\ge 1$ is the total number of parameters.  Let 
\begin{align*}
\mathcal{G}=\big\{l\in[q]:\gamma_{i,j}^{t-1}(\btheta)~{\rm involves}~\theta_l~{\rm for~ all}~1\le i< j\le p~{\rm and}~t\in[n]\setminus[m]\big\}\,.
\end{align*}
We call $\btheta_{\mathcal{G}}$ and $\btheta_{\mathcal{G}^{\btc}}$,
respectively, the global parameter vector and the local parameter vector. 
\end{definition}
% In all three models presented in Section \ref{sec:3examples}, $\{\xi_i\}_{i=1}^p$ and $\{ \eta_i \}_{i=1}^p$ are local parameters. 

In the degree heterogeneity model described in Section~{\rm\ref{sec:density}}, the global parameter vector takes the form $\btheta_{\mathcal{G}} = (a_0, a_1, b_0, b_1)^{\T}$. For the persistence and transitivity models introduced in Sections~{\rm\ref{sec:persistence}} and {\rm\ref{sec:transitivity}}, we have $\btheta_{\mathcal{G}} = (a, b)^{\T}$. In these three models, the local parameters are represented by $\btheta_{\mathcal{G}^{\btc}} = (\xi_1, \dots, \xi_p, \eta_1, \dots, \eta_p)^{\T}$. Recall that $m$ denotes the order of the AR network process. We have $m = 1$ for both the degree heterogeneity model and transitivity model, and $m = 3$ for the persistence model. As we will  discuss in Section \ref{sec:globaliden}, the global parameter vector $\btheta_{0,\mathcal{G}}$ and the local parameter vector $\btheta_{0,\mathcal{G}^{\btc}}$ need to be treated differently, which we accomplish via partial likelihoods. The resulting estimators may also entertain different convergence rates.

We develop the estimation theory for our models in three stages below.  Sufficient conditions for identification of $\btheta_0$ are established with respect to an expected partial log-likelihood 
$\ell^{(l)}_{n,p}(\btheta)$, defined in \eqref{eq:idelocf} below.  An initial estimator $\tilde{\btheta}$ results from maximizing the corresponding partial log-likelihoods $\hat \ell^{(l)}_{n,p}(\btheta)$ defined in \eqref{logL:local} below, for each $l\in[q]$.   
  Finally, because of the potential high-dimensionality of our models (number of local parameters increasing with number of nodes), these estimators suffer from slow rates of convergence.  We offer estimators with improved rate of convergence, derived as a refinement of the initial estimator via the notion of projected score functions. To study the convergence properties of the proposed estimation procedures, a key step is to investigate  the convergence rate of
$\sup_{\btheta \in \bTheta}|\hat{\ell}_{n,p}^{(l)}(\btheta) - \ell_{n,p}^{(l)}(\btheta)|$. As standard techniques are not applicable due to the temporal dependence inherent in our dependent-edge AR network models, we construct  a martingale difference structure for the network models and use such structure to address the theoretical challenge rigorously. See Lemma~\ref{la:uniformcov} and its proof in the supplementary material for further details.

\subsection{Identification of $\btheta_{0}$}\label{sec:globaliden}

Let $\mathcal{F}_t$ be the $\sigma$-field generated by $\{\bX_1,\ldots,\bX_t\}$. For any $l\in[q]$, define
\begin{equation}\nonumber
\mathcal{S}_l=\big\{(i,j):1\leq i< j\leq p~\textrm{and}~\gamma_{i,j}^{t-1}(\btheta)~\textrm{involves}~\theta_l~\textrm{for any}~t\in[n]\setminus[m]\big\}\,.
\end{equation}
If $\theta_l$ is a global parameter, $\mathcal{S}_l=\{(i,j): 1\le i<j \le p\}$.  
For estimating $\theta_l$ for $l\in [q],$
put
\begin{align} 
%&\ell_{n,p}^{(l)}(\btheta)=\frac{1}{(n-m)|\mathcal{S}_l|}\sum_{t=m+1}^n\sum_{  (i,j) \in \mathcal{S}_l}\mathbb{E}_{\mathcal{F}_{t-1}}\bigg[\log \{1- \gamma_{i,j}^{t-1}(\btheta)\}\notag\\
%&~~~~~~~~~~~~~~~~~~~~~~~~~~~~~~~~~~~~~~~~~~~~~~~~~~~~~~~~~~~~~~~~~~~+  X_{i,j}^t\log \bigg\{\frac{\gamma_{i,j}^{t-1}(\btheta)}{1- \gamma_{i,j}^{t-1}(\btheta)}\bigg\}\bigg]\,,\\
\ell_{n,p}^{(l)}(\btheta)=\frac{1}{(n-m)|\mathcal{S}_l|}\sum_{t=m+1}^n\sum_{  (i,j) \in \mathcal{S}_l}\mathbb{E}_{\mathcal{F}_{t-1}}\big\{\log\big[\{\gamma_{i,j}^{t-1}(\btheta)\}^{X_{i,j}^t}\{1-\gamma_{i,j}^{t-1}(\btheta)\}^{1-X_{i,j}^t}\big]\big\}\,,
\label{eq:idelocf}
\end{align} 
%\begin{align}\label{eq:barelltheta}
%&\ell_{n,p}(\btheta)=\frac{1}{(n-m)p(p-1)}\sum_{t=m+1}^n\sum_{i,j:\, i\neq j}\mathbb{E}_{\mathcal{F}_{t-1},\btheta_0}\bigg[\log \{1- \gamma_{i,j}^{t-1}(\btheta)\}\notag\\
%&~~~~~~~~~~~~~~~~~~~~~~~~~~~~~~~~~~~~~~~~~~~~~~~~~~~~~~~~~~~~~~~~~~~+  X_{i,j}^t\log \bigg\{\frac{\gamma_{i,j}^{t-1}(\btheta)}{1- \gamma_{i,j}^{t-1}(\btheta)}\bigg\}\bigg]\,,
%\end{align}
where $\mathbb{E}_{\mathcal{F}_{t-1}}(\cdot)$ denotes the conditional expectation given $\mathcal{F}_{t-1}$ with the unknown true  parameter vector $\btheta_0$. %Notice that 
%\[
%\log \{1- \gamma_{i,j}^{t-1}(\btheta)\} +  X_{i,j}^t\log \bigg\{\frac{\gamma_{i,j}^{t-1}(\btheta)}{1- \gamma_{i,j}^{t-1}(\btheta)}\bigg\}=\log\big[\{\gamma_{i,j}^{t-1}(\btheta)\}^{X_{i,j}^t}\{1-\gamma_{i,j}^{t-1}(\btheta)\}^{1-X_{i,j}^t}\big]\,.
%\]
For any $t\in[n]\setminus[m]$ and $1\leq i<j\leq p$, due to $\log x\leq x-1$ for any $x>0$, we have
\begin{align*}
&\mathbb{E}_{\mathcal{F}_{t-1}}\big\{\log\big[\{\gamma_{i,j}^{t-1}(\btheta)\}^{X_{i,j}^t}\{1-\gamma_{i,j}^{t-1}(\btheta)\}^{1-X_{i,j}^t}\big]\big\}\\
&~~~~~~~~~~~~-\mathbb{E}_{\mathcal{F}_{t-1}}\big\{\log\big[\{\gamma_{i,j}^{t-1}(\btheta_0)\}^{X_{i,j}^t}\{1-\gamma_{i,j}^{t-1}(\btheta_0)\}^{1-X_{i,j}^t}\big]\big\}\\
%=\,&\mathbb{E}_{\mathcal{F}_{t-1},\btheta_0}\Bigg[\log\frac{\{\gamma_{i,j}^{t-1}(\btheta)\}^{X_{i,j}^t}\{1-\gamma_{i,j}^{t-1}(\btheta)\}^{1-X_{i,j}^t}}{\{\gamma_{i,j}^{t-1}(\btheta_0)\}^{X_{i,j}^t}\{1-\gamma_{i,j}^{t-1}(\btheta_0)\}^{1-X_{i,j}^t}}\Bigg]\\
&~~~~~\leq\mathbb{E}_{\mathcal{F}_{t-1}}\Bigg[\frac{\{\gamma_{i,j}^{t-1}(\btheta)\}^{X_{i,j}^t}\{1-\gamma_{i,j}^{t-1}(\btheta)\}^{1-X_{i,j}^t}}{\{\gamma_{i,j}^{t-1}(\btheta_0)\}^{X_{i,j}^t}\{1-\gamma_{i,j}^{t-1}(\btheta_0)\}^{1-X_{i,j}^t}}\Bigg]-1\\
&~~~~~=\frac{\gamma_{i,j}^{t-1}(\btheta)}{\gamma_{i,j}^{t-1}(\btheta_0)}\cdot\gamma_{i,j}^{t-1}(\btheta_0)+\frac{1-\gamma_{i,j}^{t-1}(\btheta)}{1-\gamma_{i,j}^{t-1}(\btheta_0)}\cdot\{1-\gamma_{i,j}^{t-1}(\btheta_0)\}-1=0\,,
\end{align*}
which implies
$
\ell_{n,p}^{(l)}(\btheta)\leq\ell_{n,p}^{(l)}(\btheta_0)$ for any $\btheta\in\bTheta$. Notice that 
\begin{align*}
&\mathbb{E}_{\mathcal{F}_{t-1}}\big\{\log\big[\{\gamma_{i,j}^{t-1}(\btheta)\}^{X_{i,j}^t}\{1-\gamma_{i,j}^{t-1}(\btheta)\}^{1-X_{i,j}^t}\big]\big\}\\
&~~~~~~~~~~~~-\mathbb{E}_{\mathcal{F}_{t-1}}\big\{\log\big[\{\gamma_{i,j}^{t-1}(\btheta_0)\}^{X_{i,j}^t}\{1-\gamma_{i,j}^{t-1}(\btheta_0)\}^{1-X_{i,j}^t}\big]\big\}=0
\end{align*}
if and only if 
\begin{align}\label{eq:iden}
\frac{\{\gamma_{i,j}^{t-1}(\btheta)\}^{X_{i,j}^t}\{1-\gamma_{i,j}^{t-1}(\btheta)\}^{1-X_{i,j}^t}}{\{\gamma_{i,j}^{t-1}(\btheta_0)\}^{X_{i,j}^t}\{1-\gamma_{i,j}^{t-1}(\btheta_0)\}^{1-X_{i,j}^t}}\equiv1\,,
\end{align}
where \eqref{eq:iden} is equivalent to $\gamma_{i,j}^{t-1}(\btheta)=\gamma_{i,j}^{t-1}(\btheta_0)$. Hence, for any $\btheta\in\bTheta\setminus\{\btheta_0\}$, $\ell_{n,p}^{(l)}(\btheta)=\ell_{n,p}^{(l)}(\btheta_0)$ if and only if $\gamma_{i,j}^{t-1}(\btheta)=\gamma_{i,j}^{t-1}(\btheta_0)$ for any $t\in[n]\setminus[m]$ and $(i,j)\in\mathcal{S}_l$. To guarantee the identification of $\btheta_0$, we impose the following regularity conditions.

\begin{condition}\label{as:gambd}
{\rm(i)} There exists some universal constant $C_{1}>0$ such that
\[
\min_{t\in[n]\setminus[m]}\min_{i,j:\,1\leq i<j\leq p}\inf_{\btheta\in\bTheta}\gamma_{i,j}^{t-1}(\btheta)\{1-\gamma_{i,j}^{t-1}(\btheta)\}\geq C_1\,.
\]
{\rm (ii)} For any $1\leq i<j\leq p$ and $t\in[n]\setminus[m]$, $\gamma_{i,j}^{t-1}(\btheta)$ is thrice continuously
differentiable with respect to $\btheta\in\bTheta$. Furthermore, there exists some universal constant
$C_{2}>0$ such that
\[
\max_{t\in[n]\setminus[m]}\max_{i,j:\,1\leq i<j\leq p}\sup_{\btheta\in\bTheta}\bigg|\frac{\partial^k \gamma_{i,j}^{t-1}(\btheta)}{\partial\btheta^k}\bigg|_\infty\leq C_2
\]
for any $k\in[3]$.
\end{condition}
%\begin{remark}

Condition~{\rm\ref{as:gambd}} specifies conditions for  the parameter space
$\bTheta$. We assume $C_1$ is a universal constant for simplifying the presentation. Our proposed methods also work when $C_1$ converges slowly to zero as $p\rightarrow\infty$.  Recall that $
\gamma_{i,j}^{t-1}(\btheta) =\alpha_{i,j}^{t-1}(\btheta) +X_{i,j}^{t-1} \{1 - \alpha_{i,j}^{t-1}(\btheta) -\beta_{i,j}^{t-1}(\btheta)\}
$. Due to $X_{i,j}^{t-1}\in\{0,1\}$, Condition {\rm\ref{as:gambd}(i)} holds if there exist four universal constants $c_1,c_2,c_3,c_4\in(0,1)$ with $c_1<c_2$ and $c_3<c_4$ such that
\begin{align*}
c_1\leq\alpha_{i,j}^{t-1}(\btheta)\leq c_2~~~\textrm{and}~~~1-c_4\leq \beta_{i,j}^{t-1}(\btheta)\leq 1-c_3
\end{align*}
for any $\btheta\in\bTheta$, $t\in[n]\setminus[m]$ and $1\leq i<j\leq p$. Also, Condition {\rm\ref{as:gambd}(ii)} holds provided that
\[
\bigg|\frac{\partial^k \alpha_{i,j}^{t-1}(\btheta)}{\partial\btheta^k}\bigg|_\infty\leq C_2~~~\textrm{and}~~~\bigg|\frac{\partial^k \beta_{i,j}^{t-1}(\btheta)}{\partial\btheta^k}\bigg|_\infty\leq C_2
\]
for any $\btheta\in\bTheta$, $t\in[n]\setminus[m]$ and $1\leq i<j\leq p$. Based on the explicit forms of $\alpha_{i,j}^{t-1}(\btheta)$ and $\beta_{i,j}^{t-1}(\btheta)$ in the specific models, we can identify the associated restrictions for the parameter space $\bTheta$. See
 Section \ref{app:add.discussion} of the supplementary material for further discussion on Condition {\rm\ref{as:gambd}}.

%For the density-dependent model introduced in Section {\rm\ref{sec:density}}, Condition~{\rm\ref{as:gambd}} holds if $(\xi_i \xi_j -1) \exp(a_0 +a_1) <2 $ and $(\eta_i \eta_j -1) \exp(b_0 +b_1) <2$. For the persistence model introduced in Section {\rm\ref{sec:persistence}}, Condition~{\rm\ref{as:gambd}} holds  if   $\xi_i\xi_j < e$ and $\eta_i\eta_j < e$. For the transitivity model introduced in Section {\rm\ref{sec:transitivity}}, Condition~{\rm\ref{as:gambd}} holds  if $(\xi_i\xi_j-1) \exp(a) < 2$ and $(\eta_i\eta_j-1) \exp(b) < 2$. See Section {\rm\ref{sec:fr_con}} of the supplementary material for details.
%\end{remark}

For any $1\leq i<j\leq p$, we define 
\begin{align*}
\mathcal{I}_{i,j}=\big\{l\in [q]:\gamma_{i,j}^{t-1}(\btheta)~\textrm{involves}~\theta_l~\textrm{for any}~t\in[n]\setminus[m]\big\}\,.
%\mathcal{J}_{i,j}=&~\bigg\{(l_1,l_2,l_3)\subset[q]^3:\frac{\partial^3h_{i,j}^{t-1}(\btheta)}{\partial\theta_{l_1}\partial\theta_{l_2}\partial\theta_{l_3}}\equiv0~\textrm{for any}~t\in[n]\setminus[m]~\textrm{and}~\btheta\in\bTheta\bigg\}\,,
\end{align*}
%where
%\begin{align*}
%h_{i,j}^{t-1}(\btheta)
%=\log \{1- \gamma_{i,j}^{t-1}(\btheta)\}+\gamma_{i,j}^{t-1}(\btheta_0)\log \bigg\{\frac{\gamma_{i,j}^{t-1}(\btheta)}{1- \gamma_{i,j}^{t-1}(\btheta)}\bigg\}\,.
%\end{align*}

\begin{condition}\label{as:bdSij}
There exists a universal constant $s\geq1$ such that
$
\max_{1\leq i<j\leq p}|\mathcal{I}_{i,j}|\leq s
$.
\end{condition}

%\begin{remark}
Condition~{\rm\ref{as:bdSij}} requires that the dynamics of each edge process $\{X_{i,j}^t\}_{t\ge 1}$ be driven by a finite number of  parameters. Hence, the number of global parameters is finite while the total number of local parameters may diverge together with $p$. For the degree heterogeneity model introduced in Section {\rm\ref{sec:density}}, we have $\btheta_{\mathcal{I}_{i,j}}  = (a_0,a_1,b_0,b_1,\xi_i,\xi_j, \eta_i,\eta_j)^{\T}$ with $s = 8$. For both the persistence model and transitivity model introduced in Sections {\rm\ref{sec:persistence}} and {\rm\ref{sec:transitivity}}, we have $\btheta_{\mathcal{I}_{i,j}}  = (a,b,\xi_i,\xi_j,\eta_i,\eta_j)^{\T}$ with $s=6$.
%\end{remark}

\begin{condition}\label{as:eigen}
There exists a universal constant $C_3>0$ such that
\[
\min_{i,j:\,1\leq i<j\leq p}\lambda_{\min}\bigg\{\frac{1}{n-m}\sum_{t=m+1}^n\frac{\partial\gamma_{i,j}^{t-1}(\btheta_0)}{\partial\btheta_{\mathcal{I}_{i,j}}}\frac{\partial\gamma_{i,j}^{t-1}(\btheta_0)}{\partial\btheta_{\mathcal{I}_{i,j}}^{\T}}\bigg\}\geq C_3
\]
%and
%\[
%\max_{i,j:\,i\neq j}\lambda_{\max}\bigg\{\frac{1}{n-m}\sum_{t=m+1}^n\frac{\partial\gamma_{i,j}^{t-1}(\btheta_0)}{\partial\btheta_{\mathcal{I}_{i,j}}}\frac{\partial\gamma_{i,j}^{t-1}(\btheta_0)}{\partial\btheta_{\mathcal{I}_{i,j}}^{\T}}\bigg\}\leq C_4
%\]
with probability approaching one when $n\to \infty$. 
\end{condition}

\begin{proposition}\label{pn:ident}
Let Conditions {\rm\ref{as:gambd}--\ref{as:eigen}} hold, and $C_*=2(2C_1^{-2}+C_1^{-3})C_2^3+3(C_1^{-1}+C_1^{-2})C_2^2+C_1^{-1}C_2$ with $(C_1,C_2)$ specified in Condition {\rm\ref{as:gambd}}. Assume $\sup_{\btheta\in\bTheta}|\btheta-\btheta_0|_\infty<2C_3/(C_*s^3)$.
As $n\rightarrow\infty$, it holds 
with probability approaching one that 
\begin{align*}
\ell_{n,p}^{(l)}(\btheta_0)-\ell_{n,p}^{(l)}(\btheta)\geq \frac{\bar{C}}{|\mathcal{S}_l|}\sum_{(i,j)\in\mathcal{S}_l}|\btheta_{\mathcal{I}_{i,j}}-\btheta_{0,\mathcal{I}_{i,j}}|_2^2
\end{align*}
for any $\btheta\in\bTheta$ and $l\in[q]$, where $\bar{C}>0$ is a universal constant.  
\end{proposition}

The proof of Proposition~\ref{pn:ident} is given in Section~\ref{pr:pn:ident} of the supplementary material.
Notice that $|\mathcal{S}_{ l}\cap\mathcal{S}_{l'}|=|\mathcal{S}_{ l}|$ for any $l' \in\mathcal{G} \cup \{l\}$.
By Proposition~\ref{pn:ident}, it holds with probability approaching one that for any $\btheta\in\bTheta$ and $l \in [q]$,
\begin{equation}\label{eq:ident}
    \begin{split}
        \ell_{n,p}^{(l)}(\btheta_0)-\ell_{n,p}^{(l)}(\btheta)\geq&~\frac{\bar C}{|\mathcal{S}_l|}\sum_{(i,j) \in \mathcal{S}_l}\sum_{l'\in\mathcal{I}_{i,j}}|\theta_{l'}-\theta_{0,l'}|_2^2
=~\frac{{\bar C}}{|\mathcal{S}_l|}\sum_{l'=1}^q\sum_{(i,j)\in\mathcal{S}_{l}\cap \mathcal{S}_{l'}}|\theta_{l'}-\theta_{0,l'}|^2\\
=&~{\bar C}\sum_{l'\in\mathcal{G}\cup\{l\}}|\theta_{l'}-\theta_{0,l'}|^2+{\bar C}\sum_{l'\in\mathcal{G}^{\btc} \setminus \{l\}}\frac{|\mathcal{S}_{l}\cap\mathcal{S}_{l'}||\theta_{l'}-\theta_{0,l'}|^2}{|\mathcal{S}_l|}\,.
    \end{split}
\end{equation}
% \begin{align}\label{eq:ident}
% \end{align}
 Hence, for any $l\in[q]$, the function  $\ell_{n,p}^{(l)}(\cdot)$ defined as \eqref{eq:idelocf} is a good candidate for identifying $\theta_{0,l}$ and the global parameter vector $\btheta_{0,\mathcal{G}}$ but is powerless in identifying $\theta_{0,l'}$ with $l'\in\mathcal{G}^{\btc}\setminus\{l\}$ if $|\mathcal{S}_l\cap\mathcal{S}_{l'}|\ll|\mathcal{S}_l|$. In the degree heterogeneity model, persistence model, and transitivity model introduced in Section~\ref{sec:3examples}, given   $l\in[q]$ and $l'\in\mathcal{G}^{\btc}\setminus\{l\}$, we have (i) 
 $|\mathcal{S}_l\cap\mathcal{S}_{l'}| = 1$ and $|\mathcal{S}_l| = p-1$ if $l \in\mathcal{G}^{\btc}$, and (ii) $|\mathcal{S}_l\cap\mathcal{S}_{l'}| = p-1$ and $|\mathcal{S}_l| = p(p-1)/2$ if  $l\in\mathcal{G}$, reaffirming  that   $\ell_{n,p}^{(l)}(\cdot)$  is informative for the identification of $\theta_{0,l}$ and the global parameters, yet ineffective for identifying $\theta_{0,l'}$ with $l'\in\mathcal{G}^{\btc}\setminus\{l\}$. Therefore, for the three models introduced in Section \ref{sec:3examples}, 
we need to employ $\ell_{n,p}^{(l)}(\cdot)$ to identify each $\theta_{0,l}$ for $l\in[q]$.

\subsection{Initial estimation for $\btheta_0$} \label{sec:initialEstimation}

With available observations $\bX_1, \ldots, \bX_n$, since $\{\bX_t\}_{t\geq1}$ is a Markov chain with order $m$, the likelihood function for $\btheta$, conditionally on
$\bX_1, \ldots, \bX_m$, admits the form
\[
\mathcal{L}_{n,p}(\bX_n, \ldots, \bX_{m+1}\,|\,\bX_m, \ldots, \bX_1; \btheta) =
\prod_{t=m+1}^n L_{t,p}(\bX_t \,|\, \bX_{t-1}, \ldots, \bX_{t-m};\btheta)\,,
\]
where $L_{t,p}(\bX_t \,|\, \bX_{t-1}, \ldots, \bX_{t-m};\btheta)$ is the transition probability of $\bX_t$ given $\bX_{t-1},\ldots,\bX_{t-m}$. 
By \eqref{def1}, the (normalized) log-likelihood admits the form
\begin{align} \label{eq:trueLoglikeli}
&\frac{2}{(n-m)p(p-1)}\log \mathcal{L}_{n,p}(\bX_n, \ldots, \bX_{m+1}\,|\,\bX_m, \ldots, \bX_1; \btheta)\\
&~~~~~~~~~=\frac{2}{(n-m)p(p-1)}\sum_{t=m+1}^n\sum_{i,j:\,1\leq i<j\leq p}\log\big[\{\gamma_{i,j}^{t-1}(\btheta)\}^{X_{i,j}^t}\{1-\gamma_{i,j}^{t-1}(\btheta)\}^{1-X_{i,j}^t}\big]\,,\notag
\end{align}
which is the sample version of $\ell_{n,p}^{(l)}(\btheta)$ defined as \eqref{eq:idelocf} with $l\in\mathcal{G}$.
As pointed out below (\ref{eq:ident}), we should not estimate the local parameters based on this full log-likelihood. Therefore, for each $l\in[q]$, we define 
\begin{align}\label{logL:local}
    \hat{\ell}_{n,p}^{(l)}(\btheta) =\frac{1}{(n-m)|\mathcal{S}_l|}  \sum_{t=m+1}^n\sum_{ (i,j) \in \mathcal{S}_l}\log\big[\{\gamma_{i,j}^{t-1}(\btheta)\}^{X_{i,j}^t}\{1-\gamma_{i,j}^{t-1}(\btheta)\}^{1-X_{i,j}^t}\big]\,,
%=&~\frac{1}{(n-m)p(p-1)}\sum_{t = m+1}^n \sum_{i\neq j} \bigg[\log \{1- \gamma_{i,j}^{t-1}(\btheta)\} +  X_{i,j}^t\log \bigg\{\frac{\gamma_{i,j}^{t-1}(\btheta)}{1- \gamma_{i,j}^{t-1}(\btheta)}\bigg\}\bigg]\,,
\end{align}
which contains only the terms depending on $\theta_l$ on the right-hand side of (\ref{eq:trueLoglikeli}) (with a rescaled normalized constant).

For any $l\in[q]$, Lemma \ref{la:uniformcov} in the supplementary material shows that $\hat{\ell}_{n,p}^{(l)}(\btheta)$ converges in probability to $\ell_{n,p}^{(l)}(\btheta)$ defined as \eqref{eq:idelocf} uniformly over $\btheta\in\bTheta$. 
Together with Proposition \ref{pn:ident}, we can estimate the global parameter vector $\btheta_{0,\mathcal{G}}$ by maximizing the full log-likelihood $\hat{\ell}_{n,p}^{(l')}(\btheta)$ with some $l'\in\mathcal{G}$, and estimate the local parameter $\theta_{0,l}$ with $l\in\mathcal{G}^{\btc}$ by maximizing the corresponding (partial) log-likelihood $\hat{\ell}_{n,p}^{(l)}(\btheta)$. More specifically, letting $$(\hat \theta_{*,1}^{(l)}, \dots,\hat \theta_{*,q}^{(l)} )^{\T} =\arg\max_{\btheta\in\bTheta}\hat{\ell}_{n,p}^{(l)}(\btheta)$$ for each $l\in[q]$, we define the initial estimator $\tilde{\btheta}=(\tilde{\btheta}_{\mathcal{G}}^{\T},\tilde{\btheta}_{\mathcal{G}^{\btc}}^{\T})^{\T}$ for $\btheta_0$ as
\begin{align}\label{eq:estimate}
\tilde{\btheta}_{\mathcal{G}}=(\hat{\theta}_{*,l}^{(l')})_{l\in\mathcal{G}}~~~\textrm{and}~~~\tilde{\btheta}_{\mathcal{G}^{\btc}}=(\hat{\theta}_{*,l}^{(l)})_{l\in\mathcal{G}^{\btc}}
\end{align}
for some $l'\in\mathcal{G}$. Due to $\mathcal{S}_l=\{(i,j):1\le i < j \le p\}$
for any $l\in\mathcal{G}$, we know $\hat{\ell}_{n,p}^{(l_1)}(\btheta)=\hat{\ell}_{n,p}^{(l_2)}(\btheta)$ for any $l_1,l_2\in\mathcal{G}$, which implies that the estimator $\tilde{\btheta}_{\mathcal{G}}$ given in \eqref{eq:estimate} does not depend on the selection of $l'\in\mathcal{G}$.

%We summarize the estimation procedure for $\btheta_0$ in Algorithm \ref{alg:1}. 

To investigate the theoretical properties of the  estimator $\tilde{\btheta}=(\tilde{\btheta}_{\mathcal{G}}^{\T},\tilde{\btheta}_{\mathcal{G}^{\btc}}^{\T})^{\T}$, we define \begin{align}\label{eq:cnggc}
\left\{\begin{aligned}
   c_{n,\mathcal{G}}^2=&~\frac{q\log(np)}{\sqrt{n}p}+\frac{q^{3/2}\log^{3/2}(np)}{\sqrt{n}p^2}\,, \\
   c_{n,\mathcal{G}^{\btc}}^2=&~\frac{q\log(nS_{\mathcal{G}^{\btc},\min})}{\sqrt{nS_{\mathcal{G}^{\btc},\min}}}+\frac{q^{3/2}\log^{3/2}(nS_{\mathcal{G}^{\btc},\min})}{\sqrt{n}S_{\mathcal{G}^{\btc},\min}}\,,
\end{aligned}\right.
\end{align}
where $S_{\mathcal{G}^{\btc},\min}=\min_{l\in\mathcal{G}^{\btc}}|\mathcal{S}_l|$. 
Theorem \ref{tm:consistency} shows that the convergence rate of the
initial estimator for
the local parameters is slower than that of the global parameters if
$S_{\mathcal{G}^{\btc},\min} \ll p^2$. The proof of Theorem~\ref{tm:consistency} is given in Section~\ref{pr:tm:consistency} of the supplementary material.

%\begin{algorithm}\caption{Two-step estimation procedure for $\btheta_0$}\label{alg:1}
%\begin{itemize}
%\item (Initial estimation) Calculate
% $\widehat{\btheta}_*^{(l)} =\arg\max_{\btheta\in\bTheta}\hat{\ell}_{n,p}^{(l)}(\btheta)$ for each $l\in[q]$. Write $\widehat{\btheta}_*^{(l)} = (\hat \theta_{*,1}^{(l)}, \dots,\hat \theta_{*,q}^{(l)} )^{\T}$ for each $l \in [q]$.
 
% \item (Final estimation) Let $ \widetilde{\btheta}_{\mathcal{G}^{\btc}}= (\hat \theta_{*,l}^{(l)})_{l \in \mathcal{G}^c}$ and $\widetilde{\btheta}_{-l}^{(l)}=\widehat{\btheta}_{*,\mathcal{G}}^{(l_*)}\times (\hat{\theta}_{*,l'}^{(l')})_{l'\in\mathcal{G}^{\btc}\setminus\{l\}}$ for each $l\in\mathcal{G}^{\btc}$ and some $l_*\in\mathcal{G}$. Denote by $\bTheta_{\mathcal{G}}$ and $\bTheta_l$, respectively, the parameter spaces of the global parameter vector $\btheta_{\mathcal{G}}$ and the local parameter $\theta_l$. Then
%\begin{align*}
%\widehat{\btheta}_{\mathcal{G}}=\arg\max_{\btheta_{\mathcal{G}}\in\bTheta_{\mathcal{G}}}\hat{\ell}_{n,p}^{(l_*)}(\btheta_\mathcal{G}; \widetilde{\btheta}_{\mathcal{G}^{\btc}})
%\end{align*}
%provides the estimate of $\btheta_{0,\mathcal{G}}$, and 
%\begin{align*}
%\hat{\theta}_l=\arg\max_{\theta_l\in\bTheta_l}\hat{\ell}_{n,p}^{(l)}(\theta_l;\widetilde{\btheta}_{-l}^{(l)})
%\end{align*}
%provides the estimate of $\theta_{0,l}$ with $l\in\mathcal{G}^{\btc}$. 
%\end{itemize}
%\end{algorithm}

\begin{theorem}\label{tm:consistency}
Let the conditions of Proposition {\rm\ref{pn:ident}} hold. Then $|\tilde{\btheta}_{\mathcal{G}}-\btheta_{0,\mathcal{G}}|_2=O_{\rm p}(c_{n,\mathcal{G}})$ and $|\tilde{\btheta}_{\mathcal{G}^{\btc}}-\btheta_{0,\mathcal{G}^{\btc}}|_\infty=O_{\rm p}(c_{n,\mathcal{G}^{\btc}})$. 
\end{theorem}

\begin{remark} \label{rm:consistency}
By Theorem {\rm\ref{tm:consistency}},  the initial estimator $\tilde{\btheta}_{\mathcal{G}}$ for the global parameters is consistent provided that 
\[
q\ll \min\bigg\{\frac{\sqrt{n}p}{\log(np)},\,\frac{n^{1/3}p^{4/3}}{\log(np)}\bigg\}\,,
\]
and  the initial estimator $\tilde{\btheta}_{\mathcal{G}^{\btc}}$ for the local parameters is consistent provided that
\[
q\ll \min\bigg\{\frac{\sqrt{nS_{\mathcal{G}^{\btc},\min}}}{\log(nS_{\mathcal{G}^{\btc},\min})},\,\frac{n^{1/3}S_{\mathcal{G}^{\btc},\min}^{2/3}}{\log(nS_{\mathcal{G}^{\btc},\min})}\bigg\}\,.
\]
For the degree heterogeneity model introduced in Section {\rm\ref{sec:density}}, we have $q=2p+4$ and $S_{\mathcal{G}^{\btc},\min}=p-1$. For both the persistence model and transitivity model introduced in Sections {\rm\ref{sec:persistence}} and {\rm\ref{sec:transitivity}}, we have $q=2p+2$ and $S_{\mathcal{G}^{\btc},\min}=p-1$. Hence, for these three models, Theorem {\rm\ref{tm:consistency}} gives the convergence rates of $\tilde{\btheta}_{\mathcal{G}}$ and $\tilde{\btheta}_{\mathcal{G}^{\btc}}$ as follows: 
\begin{align*}
|\tilde{\btheta}_{\mathcal{G}}-\btheta_{0,\mathcal{G}}|_2=&~ O_{\rm p}\bigg\{\frac{\log^{1/2}(np)}{n^{1/4}}\vee\frac{\log^{3/4}(np)}{(np)^{1/4}}\bigg\}\,,  \\
|\tilde{\btheta}_{\mathcal{G}^{\btc}}-\btheta_{0,\mathcal{G}^{\btc}}|_\infty=&~ O_{\rm p}\bigg\{\frac{p^{1/4}\log^{1/2}(np)}{n^{1/4}}\vee\frac{p^{1/4}\log^{3/4}(np)}{n^{1/4}}\bigg\}\,,
\end{align*}
which implies the consistency of $\tilde{\btheta}_{\mathcal{G}}$ provided that $\log p \ll n^{1/2}$, and the consistency of $\tilde{\btheta}_{\mathcal{G}^{\btc}}$ provided that $p\ll n(\log n)^{-3}$. 
\end{remark}
\begin{remark} \label{rmk.complexity}
Motivated by \eqref{eq:ident}, we can approximate $\tilde{\btheta}_{\mathcal{G}}$ well by
solving an alternative optimization problem involving $|\mathcal{G}|$ parameters, and each $\tilde{\theta}_l$ for $l\in\mathcal{G}^{\btc}$ through a univariate optimization. See Section~{\rm\ref{sec:impl.details}} of the supplementary material for a detailed discussion. Since each optimization is of finite dimension, standard numerical methods such as the Newton–Raphson or a Quasi-Newton algorithm can be efficiently applied. For each $l \in [q]$,
evaluating $\hat{\ell}_{n,p}^{(l)}(\btheta)$ in \eqref{logL:local}, together with its gradient and Hessian, requires $O(n|\mathcal{S}_l|)$ operations per iteration. The subsequent step, i.e., solving a small linear system in Newton--Raphson or updating the inverse-Hessian approximation in Quasi-Newton, does not introduce any higher-order computational cost.
Therefore, for the three models  introduced in Section~{\rm\ref{sec:3examples}}, the per-iteration complexity is $O(np^2)$ for the global parameters and $O(np)$ for each local parameter.
\end{remark}

%\subsection{Inference for $\btheta_0$}
\subsection{Improved estimation for $\btheta_0$} \label{sec:improvedEstimation}

Recall $\btheta=(\theta_1,\ldots,\theta_q)^{\T}$. The initial estimator $\tilde{\btheta}$ specified in \eqref{eq:estimate} suffers from slow convergence rates due to the high dimensionality of $\btheta$.
In this section, we improve the estimation for each component $\theta_{0,l}$ by projecting the score function onto certain direction. See (\ref{eq:hatanl}) below for details. An improved estimator for $\theta_{0,l}$ is then obtained by solving the projected score function while  letting $\btheta_{-l} =
\tilde{\btheta}_{-l}$. The projection mitigates the impact 
of $\tilde{\btheta}_{-l}$ in the improved estimation for $\theta_{0,l}$. This strategy was initially proposed by \cite{ChangChenTangWu2021} and \cite{ChangShiZhang2023} for constructing the valid confidence regions of some low-dimensional subvector of the whole parameters in high-dimensional models with removing the impact of the high-dimensional nuisance parameters.

For $(c_{n,\mathcal{G}},c_{n,\mathcal{G}^{\btc}})$ defined as \eqref{eq:cnggc}, put
\begin{align}\label{eq:Deltan}
\Delta_n=\max\big\{|\mathcal{G}|c_{n,\mathcal{G}}^2,|\mathcal{G}^{\btc}|^2c_{n,\mathcal{G}^{\btc}}^2\big\}\,.
\end{align}
For any $t\in[n]\setminus[m]$, $l\in[q]$ and $\btheta\in\bTheta$, we define
\[
\bg_t^{(l)}(\btheta)=\frac{1}{|\mathcal{S}_l|}\sum_{(i,j)\in\mathcal{S}_l}\frac{X_{i,j}^t-\gamma_{i,j}^{t-1}(\btheta)}{\gamma_{i,j}^{t-1}(\btheta)\{1-\gamma_{i,j}^{t-1}(\btheta)\}}\frac{\partial\gamma_{i,j}^{t-1}(\btheta)}{\partial\btheta}\,.
\]
Then the score function can be written as
\[
\frac{\partial\hat{\ell}_{n,p}^{(l)}(\btheta)}{\partial\btheta}=\frac{1}{n-m}\sum_{t=m+1}^n\bg^{(l)}_t(\btheta)\,.
\]
To estimate $\theta_{0,l}$, $\btheta_{-l}$ %(components of $\btheta$ do not include $\theta_l$)
can be treated as a nuisance parameter vector. Following  \cite{ChangChenTangWu2021} and \cite{ChangShiZhang2023}, % for reducing the impact of nuisance parameter, 
we  project $\bg_t^{(l)}(\btheta)$ to form a new estimating function:
\[
\hat{f}_t^{(l)}(\btheta)=\hat{\bvarphi}_{l}^{\T}\bg_t^{(l)}(\btheta)\,,
\]
where $\hat{\bvarphi}_{l}$ is defined as
\begin{align}\label{eq:hatanl}
\hat{\bvarphi}_{l}=\arg\min_{\bu\in\mathbb{R}^q}|\bu|_1~~~~\textrm{s.t.}~~\bigg|\bigg\{\frac{1}{n-m}\sum_{t=m+1}^n\frac{\partial\bg_t^{(l)}(\tilde{\btheta})}{\partial\btheta}\bigg\}^{\T}\bu-\be_l\bigg|_\infty\leq\tau\,.
\end{align}
In the above expression, $\tau>0$ is a tuning parameter satisfying $\tau\lesssim \Delta_n^{1/2}$ with $\Delta_n$ defined as \eqref{eq:Deltan}, $\tilde{\btheta}=(\tilde{\theta}_1,\ldots,\tilde{\theta}_q)^{\T}$ is the initial estimator defined as \eqref{eq:estimate}, and $\be_l$ is a $q$-dimensional vector with the $l$-th component being $1$ and other components being $0$. Then we can re-estimate $\btheta_0$ by $\check{\btheta}=(\check{\theta}_1,\ldots,\check{\theta}_q)^{\T}$, where 
\begin{equation} \label{eq:theta_check}
    \check{\theta}_l=\arg\min_{\theta_l\in B(\tilde{\theta}_l, \tilde r)}\bigg|\frac{1}{n-m}\sum_{t=m+1}^n\hat{f}_t^{(l)}(\theta_l,\tilde{\btheta}_{-l})\bigg|^2
\end{equation}
for some $\tilde{r}>0$ satisfying $\max\{c_{n,\mathcal{G}},c_{n,\mathcal{G}^{\btc}}\}\ll \tilde{r}\ll 1$. 

To construct the convergence rate of $|\check{\btheta}-\btheta_{0}|_\infty$, we need the following regularity condition, which is analogous to Condition 1 of \cite{ChangChenTangWu2021} and Condition 7 of \cite{ChangShiZhang2023}. See the discussion there for the validity of such condition.
\begin{condition}\label{as:al}
  For each $l\in[q]$, there is a nonrandom vector $\bvarphi_l\in\mathbb{R}^q$ 
  such that $|\bvarphi_l|_1\leq C_4$ for some universal constant $C_4>0$, and $\max_{l\in[q]}|\hat{\bvarphi}_l-\bvarphi_l|_1=O_{\rm p}(\omega_n)$ for some $\omega_n\rightarrow0$ satisfying $\omega_n(\log q)^{1/2}\log(qn)=o(1)$.   
\end{condition}

 Proposition \ref{pn:reconv} shows that $\check{\btheta}$ has faster convergence rate than the initial estimator $\tilde{\btheta}$ given in \eqref{eq:estimate}.  The proof of Proposition~\ref{pn:reconv}
is given in Section~\ref{pr:pn:reconv} of the supplementary material.
\begin{proposition}\label{pn:reconv}
Let the conditions of Proposition {\rm\ref{pn:ident}} and Condition {\rm\ref{as:al}} hold. Then
$|\check{\btheta}-\btheta_{0}|_\infty=O_{\rm p}(\Delta_n)$, where $\Delta_n$ is defined as \eqref{eq:Deltan}.
\end{proposition}

Based on the obtained $\check{\btheta}$, we consider the final estimate $\hat{\btheta}=(\hat{\theta}_1,\ldots,\hat{\theta}_q)^{\T}$ for $\btheta_0$ defined as follows:
\begin{align}\label{eq:thetal}
\hat{\theta}_l=\arg\min_{\theta_l\in B(\check{\theta}_l,\check{r})}\bigg|\frac{1}{n-m}\sum_{t=m+1}^n\hat{f}_t^{(l)}(\theta_l,\check{\btheta}_{-l})\bigg|^2
\end{align}
for some $\check{r}>0$ satisfying $q\Delta_n\ll \check{r}\ll 1$ with $\Delta_n$ defined as \eqref{eq:Deltan}.

\begin{remark} 
Given the initial estimate $\tilde{\btheta}$, there are three tuning parameters $(\tau,\tilde{r},\check{r})$ for deriving our final estimate $\hat{\btheta}$. 
For the degree heterogeneity model introduced in Section {\rm\ref{sec:density}}, we have $|\mathcal{G}| = 4$ and $|\mathcal{G}^\btc| = 2p$. For both the persistence model and transitivity model introduced in Sections {\rm\ref{sec:persistence}} and {\rm\ref{sec:transitivity}}, we have $|\mathcal{G}| = 2$ and $|\mathcal{G}^\btc| = 2p$.  Together with Remark {\rm\ref{rm:consistency}}, we have $\Delta_n = n^{-1/2}p^{5/2}\log^{3/2}(np)$ for these three models. The improved estimation procedure thus requires $\tau\lesssim n^{-1/4}p^{5/4}\log^{3/4}(np)$, $n^{-1/4}p^{1/4}\log^{3/4}(np)\ll \tilde{r}\ll 1$ and $n^{-1/2}p^{7/2}\log^{3/2}(np)\ll \check{r}\ll 1$, which suggests $p\ll n^{1/7}(\log n)^{-3/7}$.
In practice, for the three models introduced in Section {\rm\ref{sec:3examples}}, we compute the final estimate $\hat{\btheta}$ with $\tau$ proportional to $n^{-1/4}p^{5/4}\log^{3/4}(np)$ and adopting reasonably large $\tilde{r}$ and $ \check{r}$. Numerical experiments in Section {\rm\ref{sec:real}}  and Section~{\rm\ref{sec:simulation}} of the supplementary material validate the robustness of our proposed estimation procedure regarding the selections of $\tilde{r}$ and $\check{r}$ as long as $\btheta_0$ falls within the defined search range.
% and adopting reasonably large $\tilde{r}$ and $ \check{r}$ so that $\btheta_0$ falls within the search range. This yields good estimation performance across all scenarios in Sections {\rm\ref{sec:simulation}} and {\rm\ref{sec:real}}.

%{\color{red} (Add remarks on how to choose three tuning parameters
%in practice: $\tau, \tilde r$ and $\check{r}$? -- QY)}
\end{remark}

\begin{remark}
 In terms of computational complexity, for each $l \in [q]$, $\hat{\bvarphi}_l$ in \eqref{eq:hatanl} is obtained by solving a well-studied and efficiently solvable linear programming (LP) problem that involves computing a Hessian matrix.
    Recall that  Condition~{\rm\ref{as:bdSij}} implies that the dynamics of each edge process $\{X_{i,j}^t\}_{t\ge1}$ depend on a finite number of parameters. Evaluating each Hessian matrix within the LP thus requires $O(n|\mathcal{S}_l|)$ operations.  Similar to Remark~{\rm\ref{rmk.complexity}},
we use the Newton--Raphson or a Quasi-Newton algorithm to obtain $\check \theta_l$ in \eqref{eq:theta_check} and $\hat \theta_l$ in \eqref{eq:thetal} with $O(nq|\mathcal{S}_l|)$ computations incurred per iteration for $l \in [q]$. Hence,
for the three models introduced in Section~{\rm\ref{sec:3examples}}, each iteration has a computational cost of $O(np^3)$ for each global parameter and $O(np^2)$ for each local parameter. 
Compared with the detailed computational cost of the initial estimation for these three models as reported in Remark~{\rm\ref{rmk.complexity}}, the additional factor of $p$ for both the global and local parameters arises from the extra projection operation required for $\hat{f}_t^{(l)}(\boldsymbol{\theta})$, which eventually leads to faster convergence rates of $\hat \btheta$ and substantially improved empirical performance. See Theorem~{\rm\ref{tm:asymnormal}}, Remark~{\rm\ref{rmk6}}, and Section~{\rm\ref{sec:sim.est}} of the supplementary material for further details.
Importantly, the sequence of computations in \eqref{eq:hatanl}, \eqref{eq:theta_check}, and \eqref{eq:thetal} can be performed in parallel across $l$ within each step to improve computational efficiency.   Further discussion on scalability is provided in Section~{\rm\ref{app:add.discussion}} of the supplementary material.
\end{remark}

For any $\btheta\in\bTheta$ and $l\in[q]$, define
\begin{align} \nonumber
\zeta_{n,l}(\btheta)=\frac{1}{(n-m)|\mathcal{S}_l|}\sum_{t=m+1}^n\sum_{(i,j)\in\mathcal{S}_l}\frac{1}{\gamma_{i,j}^{t-1}(\btheta)\{1-\gamma_{i,j}^{t-1}(\btheta)\}}\bigg\{\bvarphi_l^{\T}\frac{\partial\gamma_{i,j}^{t-1}(\btheta)}{\partial\btheta}\bigg\}^2\,,
\end{align}
where $\bvarphi_l$ is given in Condition {\rm\ref{as:al}}. Under Conditions \ref{as:gambd} and \ref{as:al}, we have $|\zeta_{n,l}(\btheta)|\leq C_1^{-1}C_2^2C_4^2$, which implies that, for any $\btheta\in\bTheta$, $\zeta_{n,l}(\btheta)$ is a bounded random variable. To construct the asymptotic distribution of each $\hat{\theta}_l$, we require the following condition. 

\begin{condition}\label{as:conv}
For each $l\in[q]$, there exists some random variable $\kappa_l\geq0$ such that
$
\zeta_{n,l}(\btheta_0)\rightarrow\kappa_l$ in probability as $n\rightarrow\infty$.
\end{condition}

\begin{remark}
For each $l\in [q]$ and $t\geq m+1$,
let \[
v_{l}^{t-1}=\frac{1}{|\mathcal{S}_l|}\sum_{(i,j)\in\mathcal{S}_l}\frac{1}{\gamma_{i,j}^{t-1}(\btheta_0)\{1-\gamma_{i,j}^{t-1}(\btheta_0)\}}\bigg\{\bvarphi_l^{\T}\frac{\partial\gamma_{i,j}^{t-1}(\btheta_0)}{\partial\btheta}\bigg\}^2\,.
\]
As $\{\zeta_{n,l}(\btheta_0)\}_{n\geq m+1}$ is a bounded sequence of random variables for each $l\in[q]$, Condition {\rm\ref{as:conv}} is mild and $\kappa_l$ is a random variable in general. Generally speaking, the asymptotic distribution
of $\hat \theta_l$ is a mixture of normal distributions. See Theorem {\rm\ref{tm:asymnormal}} below for details. However, if the long-run variance of $\{v_l^{t-1}\}_{t=m+1}^n$ satisfies the condition
\begin{align}\label{eq:cond1}
{\rm Var}\bigg(\frac{1}{\sqrt{n-m}}\sum_{t=m+1}^nv_l^{t-1}\bigg)=o(\sqrt{n})\,,
\end{align}
$\kappa_l$ is reduced to a constant 
$$
\kappa_l=\lim_{n\rightarrow\infty}\mathbb{E}\bigg(\frac{1}{n-m}\sum_{t=m+1}^nv_{l}^{t-1}\bigg)\,.
$$
Then Theorem {\rm\ref{tm:asymnormal}} implies that $\hat \theta_l$ is asymptotically normal distributed. When the sequence
$\{v_{l}^{t}\}_{t\ge m} $ is $\alpha$-mixing with the mixing coefficients 
attaining certain convergence rates, \eqref{eq:cond1} holds automatically. 
\end{remark}

\begin{theorem}\label{tm:asymnormal}
Let the conditions of Proposition {\rm\ref{pn:ident}} and Conditions {\rm\ref{as:al}} and {\rm\ref{as:conv}} hold. For each $l\in[q]$, if $\sqrt{n|\mathcal{S}_l|}\max\{q\Delta_n^{3/2},q^2\Delta_{n}^2\}=o(1)$ with $\Delta_n$ defined as \eqref{eq:Deltan}, it then holds that
\begin{align*}
\sqrt{n|\mathcal{S}_l|}(\hat{\theta}_l-\theta_{0,l})\rightarrow\sqrt{\kappa_l}\cdot Z
 \end{align*}
in distribution as $n\rightarrow\infty$, where $Z$ is a standard normally distributed random variable independent of $\kappa_l$ specified in Condition {\rm\ref{as:conv}}.
\end{theorem}

The proof of Theorem~{\rm\ref{tm:asymnormal}}
is given in Section~{\rm\ref{pr:tm:asymnormal}} of the supplementary material.

\begin{remark} \label{rmk6}
    %Recall $\Delta_n$ in \eqref{eq:Deltan}.  
    % $\sqrt{n|\mathcal{S}_l|}\max\big\{q^2|\mathcal{G}|^{3/2}c_{n,\mathcal{G}}^3,q^2|\mathcal{G}^{\btc}|^3c_{n,\mathcal{G}^{\btc}}^3\big\} \ll 1$ for each $l \in [q]$.
{\rm(i)} Theorem {\rm\ref{tm:asymnormal}} shows that, for the global parameter $\theta_{l}$ with $l \in \mathcal{G}$,
    \begin{align*}
        |\hat{\theta}_l-\theta_{0,l}| = O_{\rm p}\bigg(\frac{1}{\sqrt{n}p}\bigg)\,,
     \end{align*}   
provided that  
 \begin{align*}
     q\ll &\min\bigg\{\frac{n^{1/10}p^{1/5}}{|\mathcal{G}|^{3/5}\log^{3/5}(np)},\,\frac{n^{1/13}p^{8/13}}{|\mathcal{G}|^{6/13}\log^{9/13}(np)},
\\&~~~~~~~~~
\frac{n^{1/10}S_{\mathcal{G}^{\btc},\min}^{3/10}}{|\mathcal{G}^{\btc}|^{6/5}p^{2/5}\log^{3/5}(nS_{\mathcal{G}^{\btc},\min})},\,\frac{n^{1/13}S_{\mathcal{G}^{\btc},\min}^{6/13}}{|\mathcal{G}^{\btc}|^{12/13}p^{4/13}\log^{9/13}(nS_{\mathcal{G}^{\btc},\min})}\bigg\}\,,
 \end{align*}
 and for the local parameter $\theta_l$ with $l \in \mathcal{G}^{\btc},$
        \begin{align*}
         |\hat{\theta}_l-\theta_{0,l}| = O_{\rm p}\bigg(\frac{1}{\sqrt{n |\mathcal{S}_l|}}\bigg)\,, 
         \end{align*}
  provided that 
%  \begin{align*}
%      q\ll &\min\bigg\{\frac{n^{1/10}p^{3/5}}{\log^{3/5}(np)|\mathcal{G}|^{3/5}|\mathcal{S}_l|^{1/5}},\,\frac{n^{1/13}p^{12/13}}{\log^{9/13}(np)|\mathcal{G}|^{6/13}|\mathcal{S}_l|^{2/13}},
% \\&~~~~~~~~~
% \frac{n^{1/10}S_{\mathcal{G}^{\btc},\min}^{3/10}}{\log^{3/5}(nS_{\mathcal{G}^{\btc},\min})|\mathcal{G}^{\btc}|^{6/5}|\mathcal{S}_l|^{1/5}},\,\frac{n^{1/13}S_{\mathcal{G}^{\btc},\min}^{6/13}}{\log^{9/13}(nS_{\mathcal{G}^{\btc},\min})|\mathcal{G}^{\btc}|^{12/13}|\mathcal{S}_l|^{2/13}}\bigg\}\,.
%  \end{align*}
  \begin{align*}
     q\ll &\min\bigg\{\frac{n^{1/10}p^{3/5}}{|\mathcal{G}|^{3/5}|\mathcal{S}_l|^{1/5}\log^{3/5}(np)},\,\frac{n^{1/13}p^{12/13}}{|\mathcal{G}|^{6/13}|\mathcal{S}_l|^{2/13}\log^{9/13}(np)},\,
\frac{n^{1/8}p^{1/2}}{|\mathcal{G}|^{1/2}|\mathcal{S}_l|^{1/8}\log^{1/2}(np)},
\\&~~~~~~~~~
\frac{n^{1/10}S_{\mathcal{G}^{\btc},\min}^{3/10}}{|\mathcal{G}^{\btc}|^{6/5}|\mathcal{S}_l|^{1/5}\log^{3/5}(nS_{\mathcal{G}^{\btc},\min})},\,\frac{n^{1/13}S_{\mathcal{G}^{\btc},\min}^{6/13}}{|\mathcal{G}^{\btc}|^{12/13}|\mathcal{S}_l|^{2/13}\log^{9/13}(nS_{\mathcal{G}^{\btc},\min})}\bigg\}\,.
 \end{align*}
In particular, for the three models introduced in Section {\rm\ref{sec:3examples}},  the estimators satisfy $|\hat{\theta}_l-\theta_{0,l}| = O_{\rm p}(n^{-1/2}p^{-1})$ for $l \in \mathcal{G}$ if $p \ll n^{1/23} (\log n)^{-9/23}$, and $|\hat{\theta}_l-\theta_{0,l}| = O_{\rm p}\{(np)^{-1/2} \}$ for $l \in \mathcal{G}^{\btc}$ if $p \ll n^{1/21} (\log n)^{-3/7}$. Compared with the results in Theorem {\rm\ref{tm:consistency}}, the improved estimator $\hat \btheta$ achieves a faster convergence rate than the initial estimator $\tilde \btheta$. 

{\rm(ii)} For each $l\in[q]$, write 
\[
\hat{\zeta}_{n,l}(\hat{\btheta})=\frac{1}{(n-m)|\mathcal{S}_l|}\sum_{t=m+1}^n\sum_{(i,j)\in\mathcal{S}_l}\frac{1}{\gamma_{i,j}^{t-1}(\hat{\btheta})\{1-\gamma_{i,j}^{t-1}(\hat{\btheta})\}}\bigg\{\hat{\bvarphi}_l^{\T}\frac{\partial\gamma_{i,j}^{t-1}(\hat{\btheta})}{\partial\btheta}\bigg\}^2
\]
with $\hat{\bvarphi}_l$ defined as \eqref{eq:hatanl}. 
Since $\hat{\zeta}_{n,l}(\hat{\btheta})-\zeta_{n,l}(\btheta_0)\rightarrow0$ in probability as $n\rightarrow\infty$, by Corollary {\rm3.2} of {\rm\cite{HallHeyde1980}}, it holds that 
\begin{align}\label{eq:asy1}
\sqrt{\frac{n|\mathcal{S}_l|}{\hat{\zeta}_{n,l}(\hat{\btheta})}}(\hat{\theta}_l-\theta_{0,l})\rightarrow \mathcal{N}(0,1)
\end{align}
in distribution as $n\rightarrow\infty$, provided that $\mathbb{P}(\kappa_l>0)=1$. We can use \eqref{eq:asy1} to construct the confidence interval for each $\theta_l$.

% To be specific, for the three models introduced in Section {\rm\ref{sec:3examples}}, Theorem \ref{tm:asymnormal} provides the convergence rate as  $|\hat{\theta}_l-\theta_{0,l}| = o_p\{n^{-3/4}p^{23/4} \log^{9/4}(np)\}$,  which indicates $p\ll n^{3/23}(\log n)^{-9/23}$ for the consistency of the improved estimator.

\end{remark}

\section{Real data analysis: Email interactions} \label{sec:real}
%\subsection{Email interactions} \label{subsec:emails}

We studied estimation and coverage properties in the transitivity model (\ref{trans-model}) and an extended version of that model.  See Sections~\ref{sec:simulation} and~\ref{app:simulation} of the supplementary material for a summary of results. In this section, we apply the transitivity model (\ref{trans-model}) to a dynamic network dataset of email interactions in a medium-sized Polish manufacturing company, from January to September 2010 \citep{michalski2011}.
We analyze a subset of the data among $p=106$ of the most active participants out of an original 167 employees.
The organizational tree of direct reports in the company is also available for these employees. 
Each of the $n=39$ network snapshots corresponds to a non-overlapping time window, with $X_{i,j}^t=1$ if participants $i$ and $j$ exchanged at least one email in week $t$. Binarizing the communications at a weekly scale removes periodic effects and irregular behaviours which are present at higher frequencies.

In Section~\ref{subsec:emails_appendix} of the supplementary material, we inspect the stationarity of the network and the effective sample size. 
Our preliminary plots of edge density and dynamic activity suggest a change point in the network sequence, hence we fit our model separately to the first 13 and last 26 snapshots, referred to as ``period 1'' and ``period 2''.  
Overall, the proportion of non-edges that form in the next snapshot is about 5\%, while the proportion of existing edges which persist in the next snapshot is about 55\%, clear evidence of temporal edge dependence.

\begin{figure}[t] 
     \vspace{-0.2cm}
    \centering
     \begin{subfigure}[b]{0.45\textwidth}
         \centering
         \includegraphics[width=\textwidth]{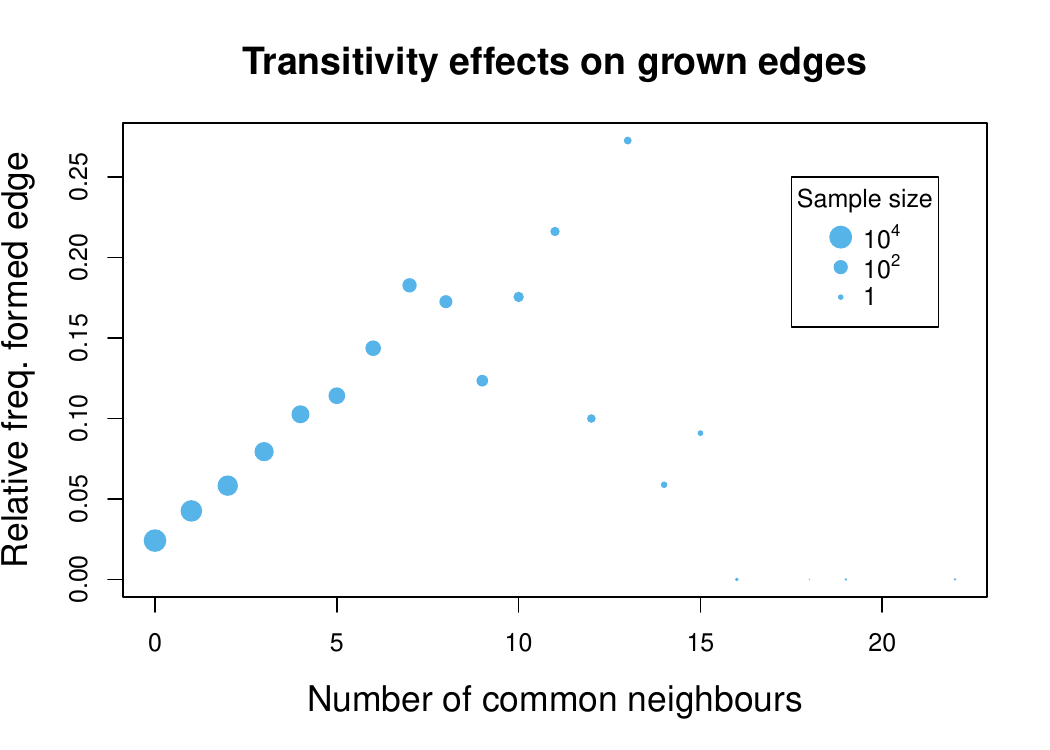}
     \end{subfigure}   
          \begin{subfigure}[b]{0.45\textwidth}
         \centering
         \includegraphics[width=\textwidth]{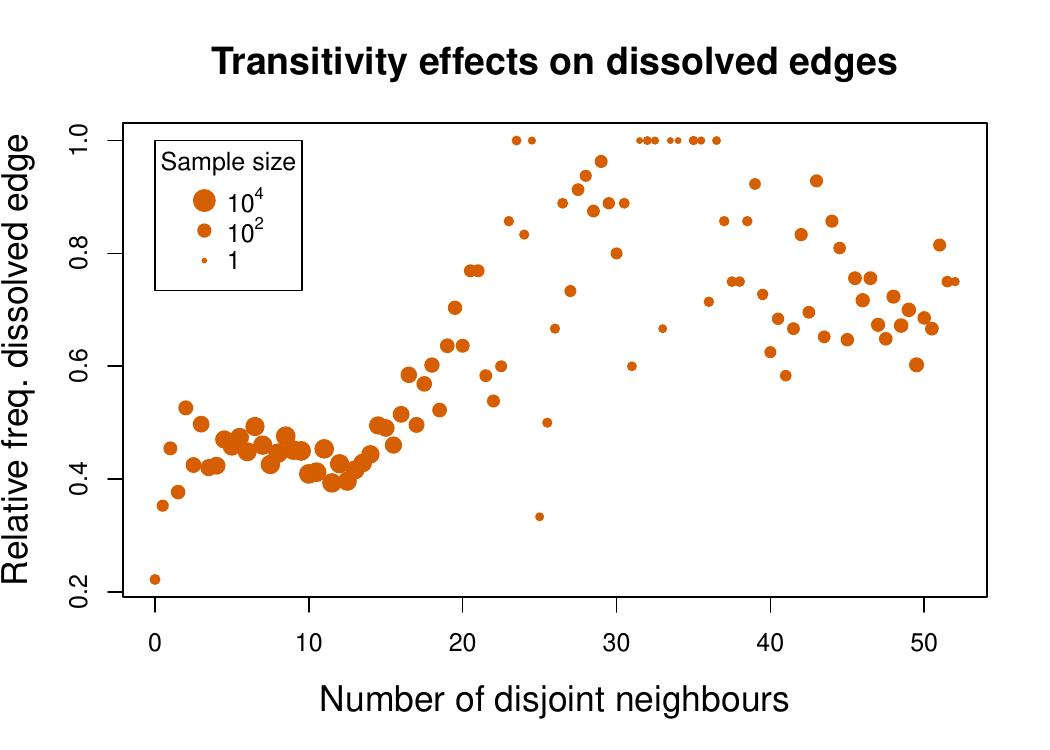}
     \end{subfigure}
          \vspace{-0.2cm}
     \caption{\label{fig:man_trans} Left panel: the plot of relative edge frequency $|\mathcal{U}_\ell^1|/|\mathcal{U}_\ell|$ against $\ell$, email interaction networks. Right panel: the plot of relative non-edge frequency $|\mathcal{V}_\ell^0|/|\mathcal{V}_\ell|$ against $\ell$, email interaction networks. In both panels, point size is proportional to the log sample sizes $\log|\mathcal{U}_\ell|$ and $\log|\mathcal{V}_\ell|$ respectively.}
          \vspace{-0.2cm}
\end{figure}

Basic summaries also identify empirical evidence for transitivity effects, demonstrated in Figure~\ref{fig:man_trans}.
To construct these plots, we partition the edge variables as follows: for each integer $\ell\geq0$, define
\begin{align*}
    &\mathcal{U}_\ell = \big\{(i,j,t) : 1 \leq i < j \leq p\,, ~t \in [n]\setminus\{1\}\,, ~X_{i,j}^{t-1}=0\,, ~U_{i,j}^{t-1} = \ell/(p-2) \big\}\,, \\
    &\mathcal{V}_\ell = \big\{(i,j,t) : 1 \leq i < j \leq p\,, ~t \in[n]\setminus\{1\}\,, ~X_{i,j}^{t-1}=1\,, ~V_{i,j}^{t-1} = \ell/(p-2) \big\}\,,\\
    &\mathcal{U}_\ell^1 = \big\{(i,j,t) \in \mathcal{U}_\ell, \; X_{i,j}^t=1 \big\}, \quad 
    \mathcal{V}_\ell^0 = \big\{(i,j,t) \in \mathcal{V}_\ell,\; X_{i,j}^t=0
    \big\}\,,
\end{align*}
where $U_{i,j}^{t-1}$ and $V_{i,j}^{t-1}$ are given in \eqref{UVij}. 
The left panel of Figure~\ref{fig:man_trans} plots the relative frequency $|\mathcal{U}_\ell^1|/|\mathcal{U}_\ell|$ against $\ell$ for $\ell =0,1, \ldots$, showing that this frequency of grown edges tends to be higher for node pairs with more common neighbours in the previous snapshot. 
The right panel of Figure~\ref{fig:man_trans} analogously plots the relative frequency  $|\mathcal{V}_\ell^0|/|\mathcal{V}_\ell|$ against $\ell$, and shows a similar increasing relationship between disjoint neighbours and frequency of dissolved edges.

% Moved the preliminary summaries to the appendix essentially as is

\begin{figure}[t] 
     \vspace{-0.2cm}
    \centering
     \begin{subfigure}[b]{0.45\textwidth}
         \centering
         \includegraphics[width=\textwidth]{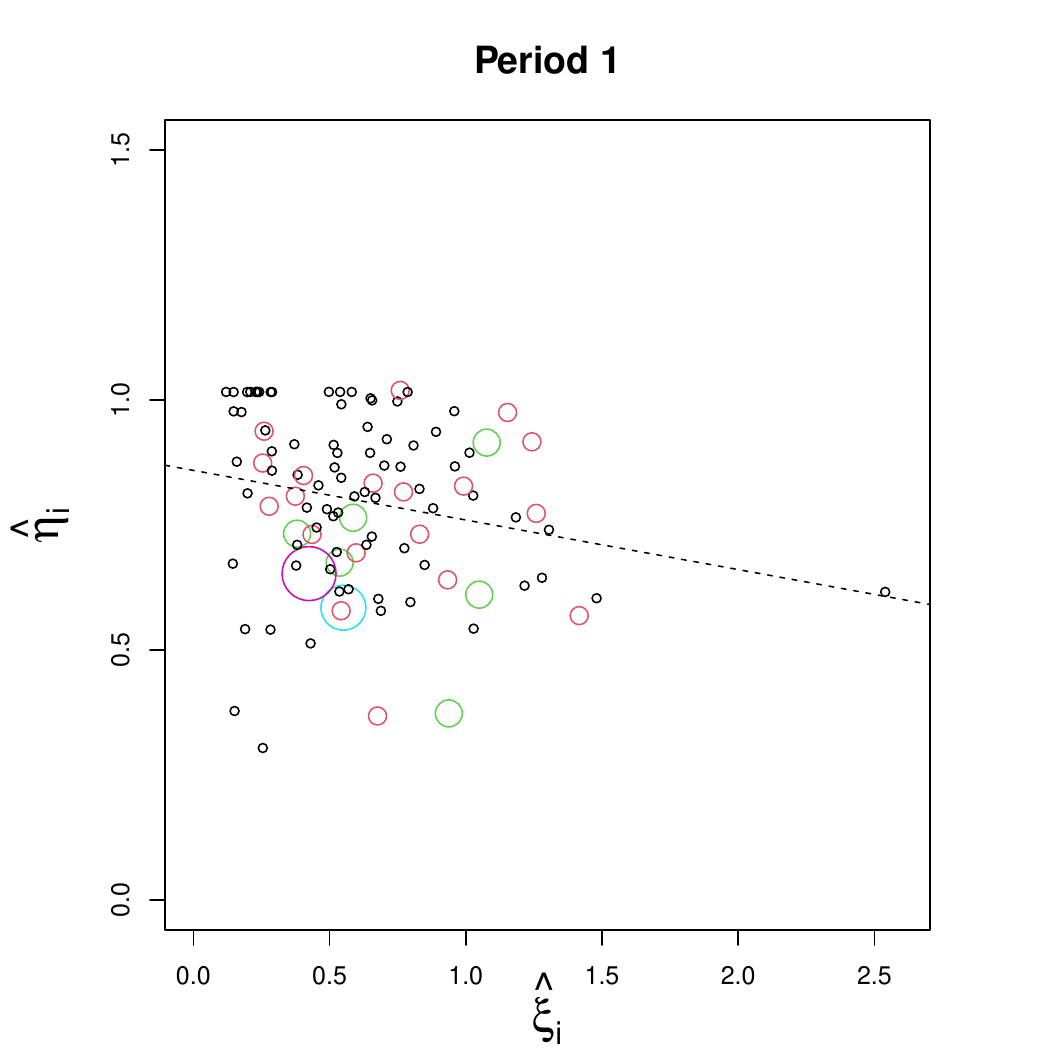}
     \end{subfigure}   
          \begin{subfigure}[b]{0.45\textwidth}
         \centering
         \includegraphics[width=\textwidth]{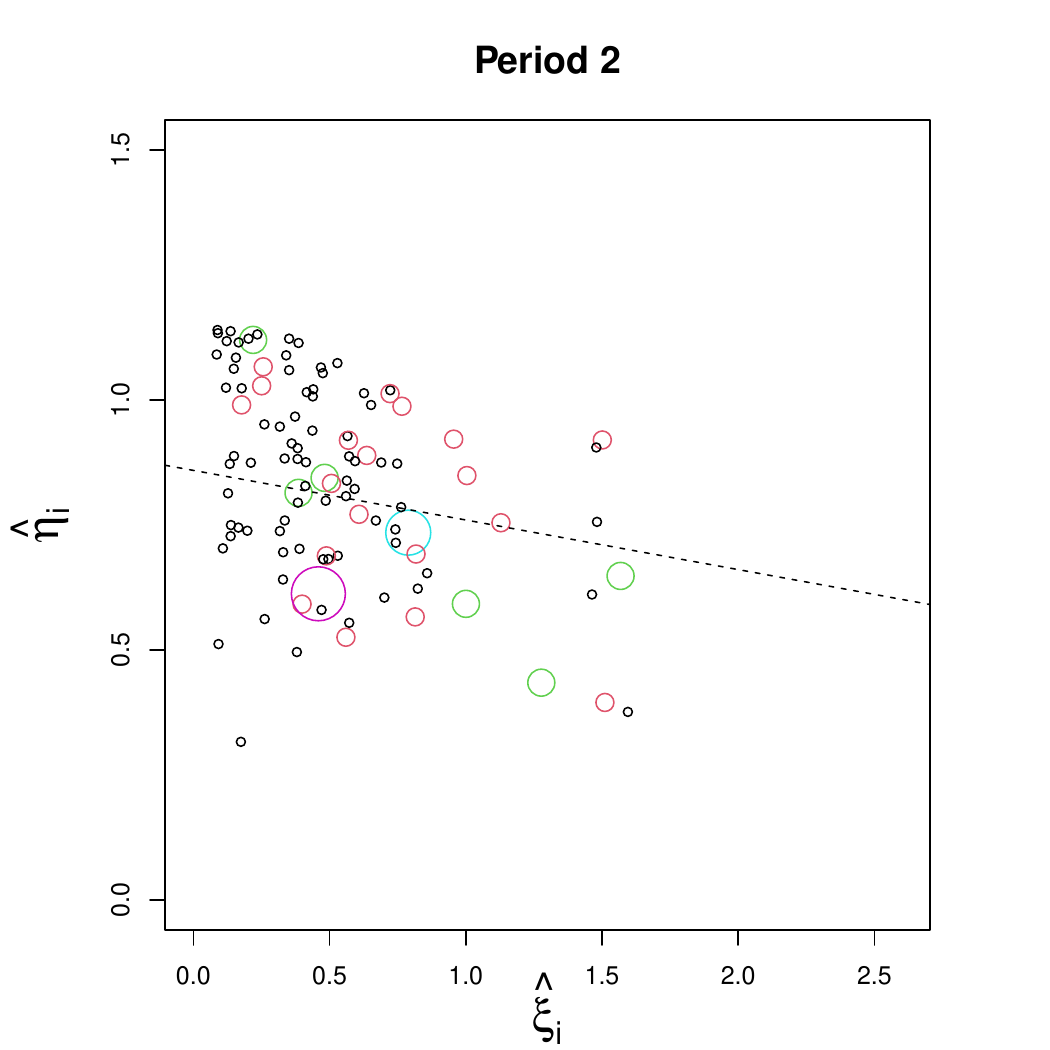}
     \end{subfigure}
          \vspace{-0.2cm}
     \caption{\label{fig:man_scatters}{Scatter plots of estimates $\{\hat{\xi}_i\}_{i=1}^{106}$ and $\{\hat{\eta}_i\}_{i=1}^{106}$ for periods 1 and 2, email interaction data. Circles are sized and coloured according to their level in the company organizational tree from 1 (no direct reports) to 5 (CEO). Level 1: black, smallest; level 2: red; level 3: green; level 4: cyan; level 5: purple, largest.}}
          \vspace{-0.2cm}
\end{figure}

The presence of transitivity effects is confirmed by the fit of our model parameters, using the estimation algorithm described in Section~\ref{sec:simulation} of the supplementary material. This algorithm applies the estimation procedure introduced in Section~\ref{sec:estimation}
and is implemented in our development R package \texttt{arnetworks}. For period 1, we estimate the global parameters $\hat{a}=13.64$ and $\hat{b}=9.51$, suggesting a tendency towards edge growth given more common neighbours, and edge dissolution given more distinct neighbours.
%which agrees with the empirical evidence in Figure~\ref{fig:man_trans}. 
We interpret the estimates of the local parameters $\{\xi_i\}_{i=1}^{106}$ and $\{\eta_i\}_{i=1}^{106}$ in the left panel of Figure~\ref{fig:man_scatters}.
The estimates $\{\hat{\xi}_i\}_{i=1}^{106}$ have mean $0.62$ and skew towards the right, implying node heterogeneity in the edge growth. Conversely, the estimates $\{\hat{\eta}_i\}_{i=1}^{106}$ have mean $0.80$ and skew towards the left. 
There is a decreasing relationship between the paired parameters: employees who tend to grow new edges also tend to maintain existing edges.
Finally, there is an observed relationship between email behavior and company hierarchy: managers (non-leaf nodes in the organizational tree) tend to have larger estimates $\{\hat{\xi}_i\}_{i=1}^{106}$ compared to non-managers (means $0.72$ and $0.59$ respectively), implying that managers are more likely to grow edges. 
However, this increasing pattern does not continue at higher levels of the organizational tree.

The model fit to period 2 shows many of the same patterns. We estimate $\hat{a}=22.01$ and $\hat{b}=11.31$ and summarize the estimates $\{\hat{\xi}_i\}_{i=1}^{106}$ and $\{\hat{\eta}_i\}_{i=1}^{106}$ in the right  panel of Figure~\ref{fig:man_scatters}.
Relative to period 1, the larger estimate of $a$ implies a stronger transitivity effect in this time period.
The estimates $\{\hat{\xi}_i\}_{i=1}^{106}$ now have mean $0.52$ and the estimates $\{\hat{\eta}_i\}_{i=1}^{106}$ have mean $0.84$, to model overall lower edge density than in period 1.
The decreasing relationship between the paired parameters is stronger, and the means of ${\hat{\xi}_i}$ for managers and
non-managers are, respectively, $0.74$ and $0.45$.
Along with the stronger transitivity effect, we interpret that the decreased edge density in period 2 has led to a concentration of email activity among a smaller group of employees, many of whom are managers.

We compare our model to some competing models from the literature in terms of Akaike and Bayesian information criteria (AIC, BIC). 
To briefly describe these competitors: the ``global AR model'' and ``edgewise AR model'' fit the model of \cite{Jiang2023}, with two global switching parameters or two parameters for each edge, respectively.  
The ``edgewise mean model'' assumes 
$
    X_{i,j}^t \overset{\mathrm{i.i.d.}}{\sim} \operatorname{Bernoulli}(P_{i,j})
$
with no temporal dependence, and estimates the edge probability $\{P_{i,j}\}_{i,j:\,i<j}$ for each node pair by its relative frequency of having an edge in the observed network snapshots; and the ``degree parameter mean model'' 
assumes 
$
    X_{i,j}^t \overset{\mathrm{i.i.d.}}{\sim} \operatorname{Bernoulli}(\nu_i\nu_j)
$
and estimates the degree parameters $\{\nu_i\}_{i=1}^{106}$ by fitting 1-dimensional adjacency spectral embedding \citep{Athreya2017} to the mean adjacency matrix over the observed network snapshots. 
The edgewise mean model has $O(p^2)$ parameters, while the degree parameter model has $O(p)$ parameters, like our AR network model with transitivity.
All of these models can be directly compared using the likelihoods under their respective AR network formulations, although only our AR network model with transitivity incorporates edge dependence, and the final two models do not incorporate any temporal dependence. 
Results for both periods are reported in Table~\ref{tab:aic_man}.

\begin{table}[]
     \vspace{-0.2cm}
\begin{center}
\resizebox{4in}{!} {
\begin{tabular}{|l|ll|ll|}
\hline
\textbf{}                   & \multicolumn{2}{l|}{\textbf{Period 1}}               & \multicolumn{2}{l|}{\textbf{Period 2}}               \\ \hline
\textbf{Model}              & \multicolumn{1}{l|}{\textbf{AIC}}   & \textbf{BIC}   & \multicolumn{1}{l|}{\textbf{AIC}}   & \textbf{BIC}   \\ \hline
Transitivity AR model       & \multicolumn{1}{l|}{33462} & \textbf{35412} & \multicolumn{1}{l|}{53221}          & \textbf{55327} \\ \hline
Global AR model             & \multicolumn{1}{l|}{36309}          & 36327          & \multicolumn{1}{l|}{58267}          & 58287          \\ \hline
Edgewise AR model           & \multicolumn{1}{l|}{42717}          & 144102         & \multicolumn{1}{l|}{55840}          & 165394         \\ \hline
Edgewise mean model         & \multicolumn{1}{l|}{\textbf{33248}}        & 83941          & \multicolumn{1}{l|}{\textbf{47133}} & 101910         \\ \hline
Degree parameter mean model & \multicolumn{1}{l|}{41730}          & 42695          & \multicolumn{1}{l|}{68969}          & 70013          \\ \hline
\end{tabular}
}
\end{center}
     \vspace{-0.2cm}
\caption{\label{tab:aic_man}{AIC and BIC performance for email interaction data, periods 1 and 2.}}
     \vspace{-0.2cm}
\end{table}

In both periods, our AR network model with transitivity is outperformed only by the edgewise mean model in terms of AIC, and achieves the lowest BIC, as it uses fewer parameters. 
This reduction of the parameter space is important for modelling sparse dynamic network data: although there is clear temporal edge dependence in this data, the edgewise mean model outperforms the edgewise AR model, as there is low effective sample size to estimate the local edge dissolution parameters. 

\begin{figure}[t] 
    \centering
     \begin{subfigure}[b]{0.3\textwidth}
         \centering
         \includegraphics[width=\textwidth]{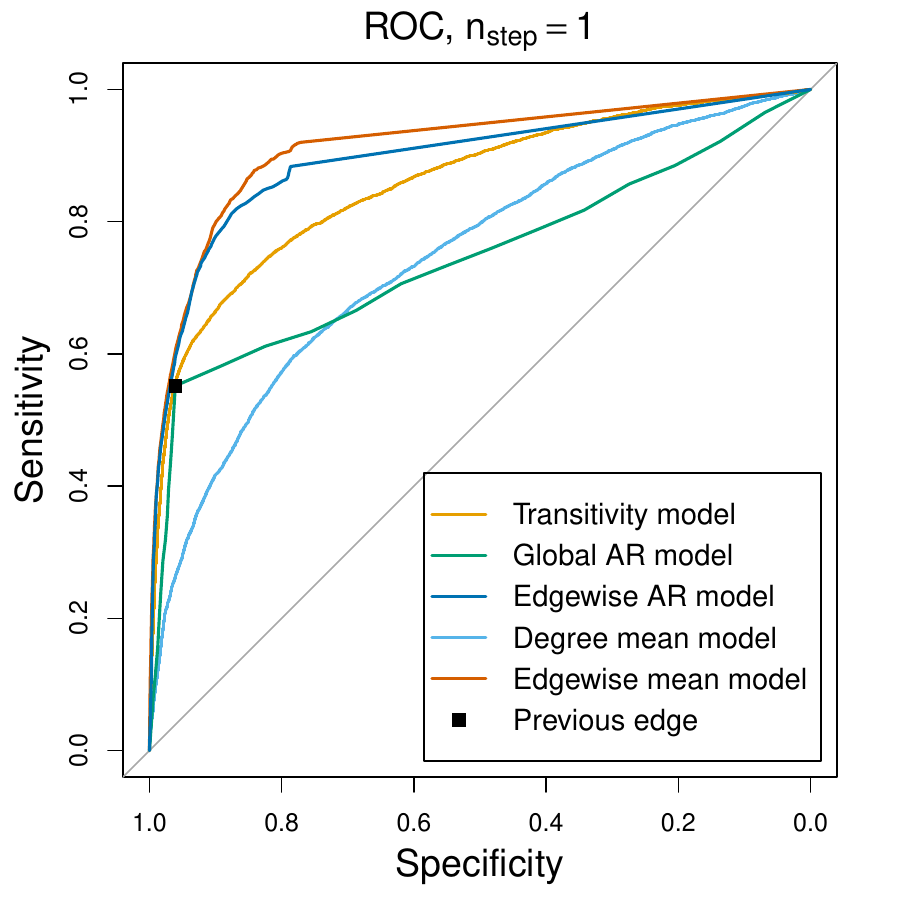}
     \end{subfigure}   
          \begin{subfigure}[b]{0.3\textwidth}
         \centering
         \includegraphics[width=\textwidth]{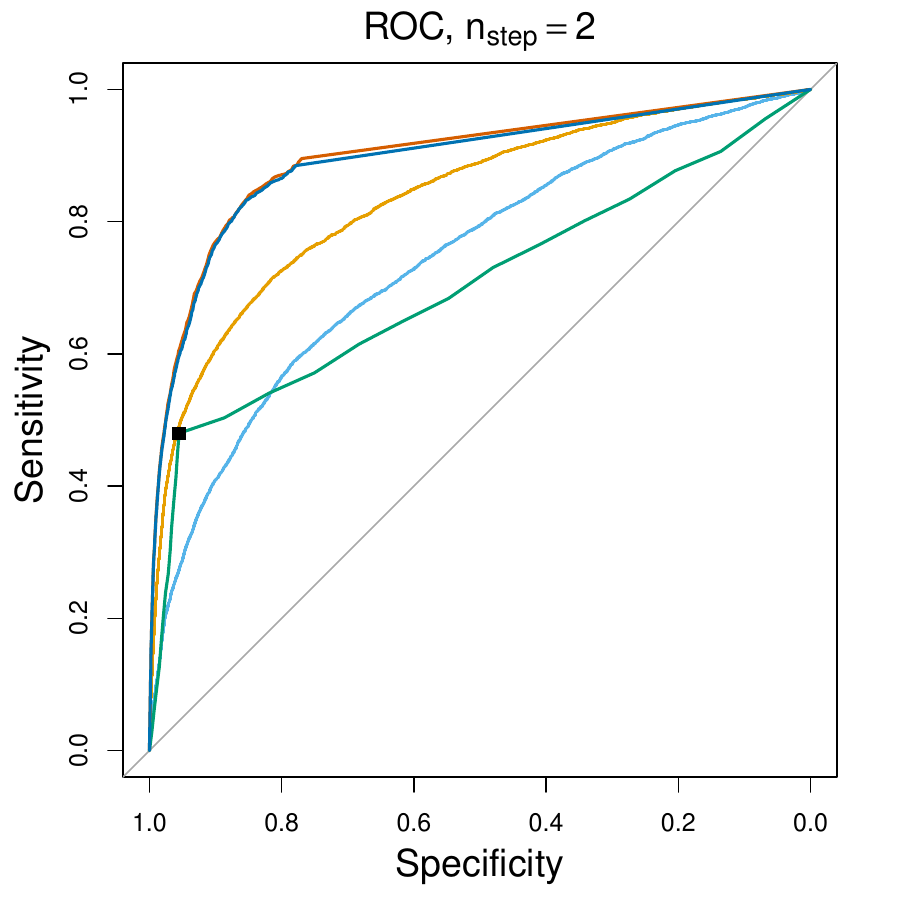}
     \end{subfigure}
              \begin{subfigure}[b]{0.3\textwidth}
         \centering
         \includegraphics[width=\textwidth]{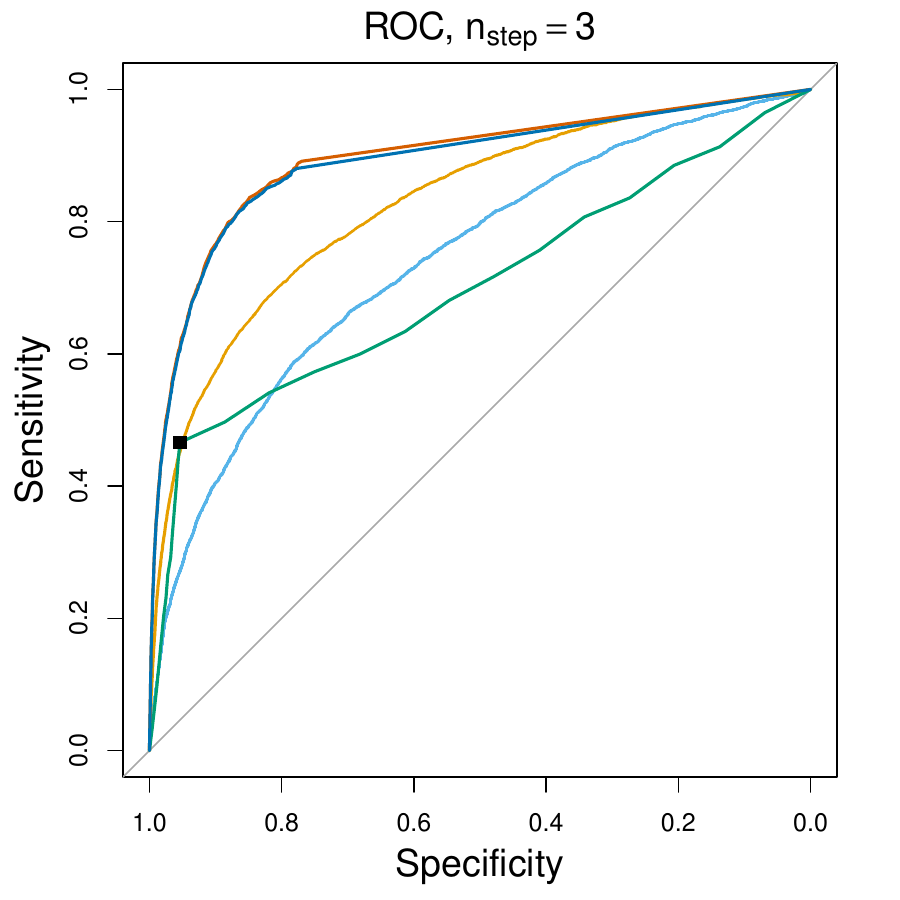}
     \end{subfigure}
      \vspace{-0.2cm}
     \caption{\label{fig:man_roc}{ROC curves for link prediction performance, email interaction data.}}
      \vspace{-0.2cm}
\end{figure}

Finally, we compare the performances of those  models in an edge forecasting task on the final 26 network snapshots (period 2). 
For $n_{\mathrm{train}}=10, \ldots, 23$, we train these models on the first $n_{\mathrm{train}}$ snapshots of period 2, then forecast the state of each edge $n_{\mathrm{step}}$ steps forward, for $n_{\mathrm{step}}=1, 2, 3$. 
The combined results are presented in Figure~\ref{fig:man_roc} as receiver operating characteristic (ROC) curves. We also include a single point summarizing the performance of naively forecasting the state of each edge by its state in the $n_{\mathrm{train}}$-th (last available observed) network snapshot (``previous edge'').

For all choices of $n_{\mathrm{step}}$, the ROC curve for our AR network model with transitivity dominates or is competitive with the global AR model, degree mean model, and the naive previous edge prediction. 
However, the two highly parameterized edgewise models achieve better performance in terms of area under the ROC, which suggests the presence of higher order structure in this network that cannot be modeled with only two parameters per node.
The edgewise mean and edgewise AR models give very similar, but not identical edge predictions; due to network sparsity the edgewise AR model has a low effective sample size to estimate the dissolution parameters, leading to slightly worse forecasting performance.

 Last but not least, the analysis above is based on the combined information for each pair $(i,j)$ in each week into a simple binary $X_{i,j}^t$, though the original data set contains more information including the precise time stamp on each email exchange.
While we did not make the full use of the available data, the information accumulation over one week makes data more stationary.
See Figure \ref{fig:man_dens_day} in comparison with Figure \ref{fig:man_dens} in Section~\ref{subsec:emails_appendix} of the supplementary material. 
Filtering out some irregular behaviours at higher frequencies
also reveals more clearly the transitivity patterns in this data set. 
See Figure \ref{fig:man_trans_day} in comparison with Figure \ref{fig:man_trans}.

\bigskip 

% \noindent
% {\bf Acknowledgements}\\
% The authors are grateful to the Co-Editor, an Associate Editor and two referees for their helpful suggestions.

% \bigskip 

% \noindent
% {\bf Conflict of interest}\\ 
% The authors declare that there are no conflicts of interest.

\bigskip

\noindent
{\bf Funding} \\
Jinyuan Chang was supported in part by the National Natural Science Foundation of China (Grant nos. 72495122 and 72125008). Eric D. Kolaczyk was supported in part by U.S. National Science Foundation grant SES-2120115, Canadian NSERC RGPIN-2023-03566 and NSERC DGDND-2023-03566. Peter W. MacDonald was partially supported by a Postdoctoral Fellowship from the Natural Sciences and Engineering Research Council of Canada. Qiwei Yao was partially supported by the U.K. Engineering and Physical Sciences Research Council (Grants nos. EP/V007556/1 and EP/X002195/1). 

\bigskip

\noindent
{\bf Data availability}\\
The R package \texttt{arnetworks}, which implements the proposed estimation and inference procedure, along with the data used in this paper, is publicly available at: \url{https://github.com/peterwmacd/arnetworks}.

\bigskip

\noindent
{\bf Supplementary material}\\
Supplementary material 
contains an extended literature review, detailed analysis on the relationship
between the proposed AR models and TERGMs, simulation studies with transitivity model, the method-of-moments-based estimation, additional real data results and all technical proofs.

\begin{appendix}
\end{appendix}
%%%%%%%%%%%%%%%%%%%%%%%%%%%%%%%%%%%%%%%%%%%%%%%%%%%%%%%%%%%%%%%%%%%%%
\linespread{1.1}\selectfont
\bibliography{paperbib}
\bibliographystyle{dcu}

% \clearpage
% \listoffigures

\newpage
\linespread{1.4}\selectfont
	\begin{center}
		\title{\large \bf Supplementary material to ``Autoregressive Networks with Dependent Edges"}	
\end{center}
\begin{center}
	{\noindent Jinyuan Chang, Qin Fang, Eric D. Kolaczyk, Peter W. MacDonald, and Qiwei Yao}
\end{center}
\bigskip

	\renewcommand{\theequation}{\thesection.\arabic{equation}}

 \renewcommand{\thepage}{S\arabic{page}}

%\bigskip

\setcounter{page}{1}
\setcounter{equation}{0}
\setcounter{section}{0}
\setcounter{table}{0}
\renewcommand{\thesection}{\Alph{section}}	
\renewcommand{\thefigure}{S\arabic{figure}}
\renewcommand{\thetable}{S\arabic{table}}
\counterwithin{equation}{section}
% Optionally reset the figure counter if starting a new numbering scheme
\setcounter{figure}{0}

This supplementary material contains an extended literature review (Section~\ref{app:extendedlit}), detailed analysis on the relationship
between the proposed AR models and TERGMs (Section \ref{app:relation2ERGM}), simulation studies with transitivity model (Section \ref{sec:simulation}), additional simulation results (Section \ref{app:simulation}), the method-of-moments-based estimation (Section \ref{sec:estimation2}), the analysis of an additional dynamic network dataset (Section~\ref{app:add.real}), additional discussion on  model specification and scalability (Section~\ref{app:add.discussion})
and all the technical proofs (Section \ref{app:proofs}).

\section{Extended literature review}
\label{app:extendedlit}

To place our contribution in the literature properly, and to clarify the novelty of our new theoretical results, we give a non-exhaustive overview of the growing research area of dynamic network modelling. We roughly taxonomize some competing classes of dynamic network models in terms of the patterns of dependence which they can capture among the edge variables.

In this work, we study a class of network models which explicitly model marginal edge dependence.
This contrasts with many well-studied network models which implictly model edge dependence through latent variables, including the stochastic block model \citep[SBM,][]{karrer2011stochastic_supp}, and many classes of latent space model \citep[LSM,][]{matias2014modeling_supp}.
% Now in main manucript:
%While latent variable models do induce marginal edge dependence, this dependence is implicit and cannot be configured easily to accommodate stylized features of real network data (e.g.~transitivity).

Our modelling framework is designed for dynamic network data collected as undirected binary edges on a common node set, with the status of edges changing over time.
We adopt a discrete-time view, and assume that we observe time-indexed network {\em snapshots} at regular intervals, or can coerce relational event data into network snapshots sampled at a constant time resolution.
When relational event data with continuous time stamps are available, dynamic network models can be specified which describe the instantaneous effects of network statistics on the rates, or relative rates of events between different node pairs.
Extensive work has been done in this setting, %on continuous-time models for network data, typically when data are observed as relational events between nodes with time stamps.
on SAOM for edge formation and dissolution events \citep{snijders2017stochastic_supp, koskinen2023multilevel_supp}, and relational event models (network-structured multivariate point process models) for instantaneous relational events \citep[]{PerryWolfe2013_supp,matias2018semiparametric_supp,butts2023relational_supp,bianchi2024relational_supp}. 
% In manuscript:
From our discrete-time viewpoint, with the exception of SAOMs, which have been used to analyze discretely observed network snapshots, we do not compare further to these approaches here. 
When both discrete and continuous-time modelling is possible, continuous-time model classes may be preferred, depending on the nature of the data and the goals of the analysis.
We show in Section~\ref{sec:real} that relational event data coerced to discrete-time snapshots at a suitable resolution can be analyzed effectively within our modelling framework. In this case, accumulation of the data over regular intervals can be used to smooth out very non-stationary or periodic behaviour at higher frequencies which may be present in the underlying continuous-time process (see Appendix~\ref{app:add.real} of the supplementary material).
%SAOM model fitting has been implemented for discrete-time network snapshots by treating underlying continuous-time events as missing data, but requires computationally expensive imputation of intermediate edge events \citep{snijders2017stochastic}.

In the discrete-time setting, existing classes of models can be loosely characterized in terms of two types of edge dependence: {\em within-snapshot} and {\em between-snapshot}.
Within-snapshot edge dependence describes the dependence patterns possible among the edge variables in a fixed snapshot, thus we refer to previous work on static networks.
The simplest pattern is {\em full independence}, for instance in the classical Erdos-Renyi (ER) model.
Exponential random graph models \citep[ERGM,][]{robins2007_supp} in principle allow for arbitrary {\em full dependence} patterns among the edge variables; in practice these patterns must be specified through a parametric model.
Note that ERGMs can also specify conditional models for edge variables, which depend on exogenous node or edge attributes; we do not pursue this direction in our modelling framework.
In recent years, advances have been made to overcome some practical issues with ERGM inference \citep{blackburn2023practical_supp}. 
%: lack of computational scalability, unstable estimation algorithms, and concentration on extreme subspaces of graph space \citep{blackburn2023practical}. 
% Clark/Handcock (dropped) is not strictly an ERGM paper, but it proposes a similar 'LOLOG' model class which is easier to work with
However, asymptotic analysis of ERGM parameter estimates remains limited. % cf the recent Fischer et al paper which looks at LD-ERGM
In recent work, \citet{schweinberger2020_supp} have considered a restriction of ERGMs to an intermediate {\em local dependence} regime (LD-ERGM): under a known partition of the network nodes into communities, they model within-community edges and allow for arbitrary dependence among edge variables in the same community, facilitating theoretical advances for parameter estimation and inference \citep{stewart2024rates_supp}. %,fischer2025stein}.

Between-snapshot edge dependence describes the dependence patterns possible between the edge variables in different snapshots.
As in the within-snapshot case, the simplest pattern is {\em full independence}, which holds if edges in different snapshots are all mutually independent.
Previous works on AR network models (\citeauthor{Jiang2023_supp}, \citeyear{Jiang2023_supp}; \citeauthor{Jiang2025_supp}, \citeyear{Jiang2025_supp}) have extended this to {\em edgewise dependence}: in this case, the only dependence between snapshots is between edge variables with the same indices $\{i,j\}$.
To our knowledge, more general forms of local dependence based on communities like in \cite{schweinberger2020_supp} have not yet been studied in the discrete dynamic setting. % and present an interesting avenue  for future research.
Seminal work on (S)TERGM \citep[]{Hanneke2010_supp, Krivitsky2014_supp, statnet2025_supp}, LNR, \citep{almquist2014logistic_supp}, as well as our new class of dependent-edge AR models allow for the specification of {\em full dependence} patterns among the edge variables between different snapshots.
All three frameworks can capture arbitrary patterns of between-snapshot dependence through parametric models which specify the transition distributions of network snapshots.
%Finally, SAOMs can be used to model dynamic network snapshots in discrete time: the underlying continuous time events are treated as missing data and must be imputed as part of the model fitting procedure \citep{snijders2017stochastic}.
%This missing data problem means that discrete time modelling does not mitigate the computational challenges of SAOM fitting described above, instead requiring computationally expensive imputation of intermediate edge events.
%SAOMs in discrete time can also allow for full dependence between snapshots.

Between-snapshot dependence can also be induced though latent variable modelling.
Both the SBM \citep{matias2017_supp,pensky2019spectral_supp} and LSM \citep{sewell2015latent_supp,Durante2016_supp,gallagher2021spectral_supp,zhang2024change_supp} have been extended to the the dynamic regime, typically by assuming full independence of edge variables (both within and between snapshots), conditional on an evolving sequence of latent variables, as in a hidden Markov model. % Not 100% certain this is the Zhang et al 2024 mentioned explicitly in the response letter
We note that latent variable modelling does not preclude additional explicit modelling of edge dependence, see for instance the {\em co-evolving} network model of \citet{zhu2023_supp} which specifies both edgewise persistence effects and evolving latent variables which depend on previous network snapshots; or the {\em block logistic autoregressive} model of \citet{suveges2023_supp}, which specifies edgewise autoregressive terms, as well as evolving sequences of latent community memberships for each node.

We summarize this discussion in Table~\ref{tab:model_classes}, which organizes existing marginal dependence models in terms of patterns of within and between snapshot edge dependence. Our new dependent-edge AR models, (S)TERGMs, and SAOMs are all able to specify arbitrary between-snapshot edge dependence. 
%${\color{red} (BUT WE DO NEED CONDITIONAL INDEPENDENCE -- QY.))}
General (S)TERGMs and SAOMs are more flexible, and can also specify arbitrary within-snapshot dependence.
In the case of SAOMs, edge dependence occurs as a consequence of temporal dependence among the edge events which are intermediate to the discrete observation times.
Our framework assumes independence of the within-snapshot edges conditional on past snapshots. %(see Section~\ref{sec:sec2}).
The same conditional independence assumption is made for LNRs by \cite{almquist2014logistic_supp} and TERGMs by \cite{Hanneke2010_supp}.
%\cite{hanneke2010} present within-snapshot edge independence as an additional structural assumption which can prevent degeneracy and make estimation more stable and computationally efficient.
%The subset of (S)TERGMs, including those studied by \cite{hanneke2010}, which make this within-snapshot independence assumption is in fact a (strict) subset of our dependent-edge AR models.
%We make this relationship to (S)TERGMs precise in Section~\ref{sec: relation2ERGM} and Section~\ref{app:relation2ERGM} of the supplementary material.

Table~\ref{tab:model_classes} also clearly identifies avenues for future work.
Expanding the allowable patterns of within-snapshot edge dependence would be valuable from a modelling perspective, allowing us to capture behavior currently only possible with (S)TERGMs or SAOMs, for instance, transitivity effects on edge formation, even when common neighbours are not yet present in the previous snapshot.
For modelling network snapshots conditional on the past, the LD-ERGM specification and methodology extends naturally by treating past snapshots as exogenous node or edge attributes. 
This form of {\em local within-snapshot} dependence, combined with our high-dimensional asymptotic techniques, may provide a path to formal theoretical analysis when the number of networks and number of unknown parameters diverges together with the network size.

%There is also a related body of literature on (discrete-time) network autoregression, which uses the network structure to regulate vector autoregressive models. The focus there is to model the dynamic changes of the features associated with each node of a network instead of the interactions across different nodes. See, among others,  \cite{Knight2016}, \cite{Zhu2017}, \cite{Yin2023}, and \cite{Butts2023}. 

\section{Relationship to TERGMs}
\label{app:relation2ERGM}

A dynamic network sequence follows a TERGM of order $m$ if it satisfies
\begin{equation}\label{eq:tergmdef}
  \mathbb{P}({\bf X}_t \,|\, {\bf X}_{t-1}, \ldots, {\bf X}_{t-m} ; \btheta)\, \propto \exp\{ \bvarsigma(\btheta)^{\T} \bvarrho({\bf X}_t, {\bf X}_{t-1}, \ldots, {\bf X}_{t-m})\}\,,
\end{equation}
where $\bvarsigma: \mathbb{R}^q \rightarrow \mathbb{R}^p$ maps the parameter vector $\btheta$ to the vector of natural parameters, and $\bvarrho$ maps the data, including the past network snapshots, to the corresponding sufficient statistics.

As in Equation (2) of \cite{Hanneke2010_supp}, 
suppose $\bvarrho$ factors over the edges of the present snapshot,
\begin{equation*}
  \bvarrho({\bf X}_t, {\bf X}_{t-1}, \ldots, {\bf X}_{t-m}) = \sum_{i,j:\,i < j} \bvarrho_{i,j}(X^t_{i,j} ; {\bf X}_{t-1}, \ldots, {\bf X}_{t-m})\,.
\end{equation*}
Then 
\begin{align}
  &\mathbb{P}({\bf X}_t \,|\, {\bf X}_{t-1}, \ldots, {\bf X}_{t-m} ; \btheta) \notag\\
  &~~~~~~~~~~~~~~~~\propto \prod_{i,j:\,i < j} \exp \big\{ \bvarsigma(\btheta)^{\T} \bvarrho_{i,j}(X^t_{i,j} ; X_{i,j}^{t-1} , {\bf X}_{t-1} \setminus X_{i,j}^{t-1}, {\bf X}_{t-2}, \ldots, {\bf X}_{t-m} )\big\}\,,\label{eq:TERGM}
\end{align}
which implies ${\bf X}_t$ will have mutually independent edges conditional on the past snapshots. We refer to this as the edge conditional independence assumption, which is a property of AR network models defined in Definition 1 of the main document.
%Hanneke et al. (2010) also show that under this setting and a condition similar to our Condition 1 (i), the resulting TERGM is non-degenerate.
We will show that any edge conditionally independent TERGM can be rewritten as an AR network model.

%Write

%where we also highlight the dependence of the sufficient statistic $\varrho_{i,j}$ on $X_{i,j}^{t-1}$, the previous edge value for the same node pair.
%Edge independence makes the normalizing constant of this conditional distribution tractable.

Denote the logit function by $\sigma(x) = \log\{x/(1-x)\}$, and specify an AR network model defined in Definition 1 by setting
\begin{align}
  \alpha_{i,j}^{t-1} %=&~ f_{i,j}({\bf X}_{t-1} \setminus X_{i,j}^{t-1}, {\bf X}_{t-2}, \ldots, {\bf X}_{t-m} ; \btheta) \\
  =&~ \sigma^{-1} \big[ \bvarsigma(\btheta)^{\T} \big\{ \bvarrho_{i,j}(1 ; 0 , {\bf X}_{t-1} \setminus X_{i,j}^{t-1}, {\bf X}_{t-2}, \ldots, {\bf X}_{t-m}) \notag\\
  &~~~~~~~~~~~~~~~~~~~~~~~~~- \bvarrho_{i,j}(0 ; 0 , {\bf X}_{t-1} \setminus X_{i,j}^{t-1}, {\bf X}_{t-2}, \ldots, {\bf X}_{t-m} )\big\} \big] \notag\\
  :=&~ \sigma^{-1}[ \bvarsigma(\btheta)^{\T}(\bvarrho_{i,j,10}^{t-1} - \bvarrho_{i,j,00}^{t-1})]\,,\label{eq:aijsupp}\\
  \beta_{i,j}^{t-1} %=&~ g_{i,j}({\bf X}_{t-1} \setminus X_{i,j}^{t-1}, {\bf X}_{t-2}, \ldots, {\bf X}_{t-m} ; \btheta) \\
  =&~ \sigma^{-1} \big[ \bvarsigma(\btheta)^{\T} \big\{ \bvarrho_{i,j}(0 ; 1 , {\bf X}_{t-1} \setminus X_{i,j}^{t-1}, {\bf X}_{t-2}, \ldots, {\bf X}_{t-m})\notag \\
  &~~~~~~~~~~~~~~~~~~~~~~~~~- \bvarrho_{i,j}(1 ; 1 , {\bf X}_{t-1} \setminus X_{i,j}^{t-1}, {\bf X}_{t-2}, \ldots, {\bf X}_{t-m} )\big\} \big]\notag \\
  :=&~ \sigma^{-1}[ \bvarsigma(\btheta)^{\T}(\bvarrho_{i,j,01}^{t-1} - \bvarrho_{i,j,11}^{t-1})]\,.\label{eq:bijsupp}
\end{align}
With renormalizion, we have
\begin{align*}
  \alpha_{i,j}^{t-1} &= \frac{\exp \{ \bvarsigma(\btheta)^{\T} \bvarrho_{i,j,10}^{t-1} \}}{\exp \{  \bvarsigma(\btheta)^{\T} \bvarrho_{i,j,10}^{t-1}\} + \exp \{  \bvarsigma(\btheta)^{\T} \bvarrho_{i,j,00}^{t-1}\}} \,,\\
  1 - \alpha_{i,j}^{t-1} &= \frac{\exp \{ \bvarsigma(\btheta)^{\T} \bvarrho_{i,j,00}^{t-1} \}}{\exp \{  \bvarsigma(\btheta)^{\T} \bvarrho_{i,j,10}^{t-1}\} + \exp \{  \bvarsigma(\btheta)^{\T} \bvarrho_{i,j,00}^{t-1}\}}\,, \\
  \beta_{i,j}^{t-1} &= \frac{\exp \{ \bvarsigma(\btheta)^{\T} \bvarrho_{i,j,01}^{t-1} \}}{\exp \{  \bvarsigma(\btheta)^{\T} \bvarrho_{i,j,01}^{t-1}\} + \exp \{  \bvarsigma(\btheta)^{\T} \bvarrho_{i,j,11}^{t-1}\}}\,, \\
  1 - \beta_{i,j}^{t-1} &= \frac{\exp \{ \bvarsigma(\btheta)^{\T} \bvarrho_{i,j,11}^{t-1} \}}{\exp \{  \bvarsigma(\btheta)^{\T} \bvarrho_{i,j,01}^{t-1}\} + \exp \{  \bvarsigma(\btheta)^{\T} \bvarrho_{i,j,11}^{t-1}\}}\,. 
\end{align*}
For the AR network model with $(\alpha_{i,j}^{t-1},\beta_{i,j}^{t-1})$ specified in \eqref{eq:aijsupp} and \eqref{eq:bijsupp}, it holds that
\begin{align*}
   &\mathbb{P}({\bf X}_t\,|\, {\bf X}_{t-1}, \ldots, {\bf X}_{t-m} ; \btheta)
  = \prod_{i,j:\,i < j} \mathbb{P}(X_{i,j}^t \,|\, {\bf X}_{t-1}, \ldots, {\bf X}_{t-m} ; \btheta) \\
  %&= \prod_{i < j : X_{i,j}^{t-1}=0} (\alpha_{i,j}^{t-1})^{X_{i,j}^t} (1 - \alpha_{i,j}^{t-1})^{1 - X_{i,j}^t} \cdot \prod_{i < j : X_{i,j}^{t-1}=1} (1- \beta_{i,j}^{t-1})^{X_{i,j}^t} (\beta_{i,j}^{t-1})^{1 - X_{i,j}^t} \\
  &~~~~~= \prod_{i,j:\,i < j,\, X_{i,j}^{t}=0, \, X_{i,j}^{t-1}=0} (1 - \alpha_{i,j}^{t-1})\cdot \prod_{i,j:\,i < j,\, X_{i,j}^{t}=1,\, X_{i,j}^{t-1}=0} \alpha_{i,j}^{t-1} \\
  &~~~~~~~~~~~~~~~\cdot \prod_{i,j:\,i < j,\,X_{i,j}^{t}=0, \, X_{i,j}^{t-1}=1} \beta_{i,j}^{t-1} \cdot \prod_{i,j:\,i < j,\, X_{i,j}^{t}=1,\, X_{i,j}^{t-1}=1} (1- \beta_{i,j}^{t-1}) \\
  &~~~~~~\propto \prod_{i,j:\,i < j,\, X_{i,j}^{t}=0, \, X_{i,j}^{t-1}=0} \exp \{ \bvarsigma(\btheta)^{\T} \bvarrho_{i,j,00}^{t-1} \} \cdot \prod_{i,j:\,i < j,\, X_{i,j}^{t}=1,\, X_{i,j}^{t-1}=0} \exp \{ \bvarsigma(\btheta)^{\T} \bvarrho_{i,j,10}^{t-1} \} \\
  &~~~~~~~~~~~~~~~\cdot \prod_{i,j:\,i < j,\, X_{i,j}^{t}=0,\, X_{i,j}^{t-1}=1} \exp \{ \bvarsigma(\btheta)^{\T} \bvarrho_{i,j,01}^{t-1} \} \cdot \prod_{i,j:\,i < j,\, X_{i,j}^{t}=1,\, X_{i,j}^{t-1}=1} \exp \{ \bvarsigma(\btheta)^{\T} \bvarrho_{i,j,11}^{t-1} \} \\
  &~~~~~= \prod_{i,j:\,i < j} \exp \big\{ \bvarsigma(\btheta)^{\T} \bvarrho_{i,j}(X^t_{i,j} ; X_{i,j}^{t-1} , {\bf X}_{t-1} \setminus X_{i,j}^{t-1}, {\bf X}_{t-2}, \ldots, {\bf X}_{t-m} )\big\}\,,
\end{align*}
which shares the same form as \eqref{eq:TERGM}.
Thus, for any edge conditionally independent TERGM, we can specify an AR network model with the same distribution.

Conversely, suppose we have specified an AR network model such that
\begin{align*}
  \sigma(\alpha_{i,j}^{t-1}) &= \bphi(\btheta)^{\T} \bu_{i,j}({\bf X}_{t-1} \setminus X_{i,j}^{t-1}, {\bf X}_{t-2}, \ldots, {\bf X}_{t-m} ) \,,\\
  \sigma(\beta_{i,j}^{t-1}) &= \bpsi(\btheta)^{\T} \bv_{i,j}({\bf X}_{t-1} \setminus X_{i,j}^{t-1}, {\bf X}_{t-2}, \ldots, {\bf X}_{t-m} )\,,
\end{align*}
for some functions of the  parameter vector $\btheta$, and the past network behavior. We claim that this AR network model can be written as an edge conditionally independent TERGM \eqref{eq:tergmdef} with $\bvarsigma(\btheta) = (\bphi(\btheta)^{\T}, \bpsi(\btheta)^{\T})^{\T}$ and $
  \bvarrho({\bf X}_t , {\bf X}_{t-1}, \ldots, {\bf X}_{t-m}) = (\bvarrho_{\alpha}^{\T} , \bvarrho_{\beta}^{\T})^{\T}$, where
$\bvarrho_{\alpha} = \sum_{i,j:\,i < j} \bvarrho_{\alpha,i,j}$ and $\bvarrho_{\beta} = \sum_{i,j:\,i < j} \bvarrho_{\beta,i,j}$ with
\begin{align*}
  &\bvarrho_{\alpha,i,j}(X_{i,j}^t ; X_{i,j}^{t-1}, {\bf X}_{t-1} \setminus X_{i,j}^{t-1}, {\bf X}_{t-2}, \ldots, {\bf X}_{t-m})\\
  &~~~~~~~~~= \bu_{i,j}({\bf X}_{t-1} \setminus X_{i,j}^{t-1}, {\bf X}_{t-2}, \ldots, {\bf X}_{t-m} )I(X_{i,j}^{t}=1, X_{i,j}^{t-1}=0) \,,\\
  &\bvarrho_{\beta,i,j}(X_{i,j}^t ; X_{i,j}^{t-1}, {\bf X}_{t-1} \setminus X_{i,j}^{t-1}, {\bf X}_{t-2}, \ldots, {\bf X}_{t-m}) \\&~~~~~~~~~=\bv_{i,j}({\bf X}_{t-1} \setminus X_{i,j}^{t-1}, {\bf X}_{t-2}, \ldots, {\bf X}_{t-m} )I(X_{i,j}^{t}=0, X_{i,j}^{t-1}=1)\,.
\end{align*}

%{\bf Summary:} any conditionally independent TERGM can be written as an AR network model. Any AR network model which is linear on the logistic scale can be written as a TERGM with conditionally independent edges.

\cite{Krivitsky2014_supp}
define the concept of separability of a dynamic model for binary networks.
%A dmodel is said to be separable if the edge formation and dissolution processes are independent conditional on the past, and the parameters also split over these two processes.
%Define two subnetworks ${\bf X}^+_t = \{X_{i,j}^t : X_{i,j}^{t-1}=0\}$ or \textcolor{red}{ $(X_{i,j}^t)_{i,j: X_{i,j}^{t-1} = 0}$} and ${\bf X}^-_t = \{X_{i,j}^t : X_{i,j}^{t-1}=1\}$ \textcolor{red}{or $(X_{i,j}^t)_{i,j:  X_{i,j}^{t-1} = 1}$ to be confirmed}.
Define two subnetworks ${\bf X}^+_t$ and ${\bf X}^-_t$, where
\begin{align} \label{eq_X+-}
    X_{i,j}^{t,+} = 1 - I(X_{i,j}^t=0,X_{i,j}^{t-1}=0) ~~\text{and}~~X_{i,j}^{t,-}= I(X_{i,j}^t=1,X_{i,j}^{t-1}=1)\,,
\end{align}
for all $i,j$.
A dynamic network model is said to be separable if ${\bf X}^+_t$ and ${\bf X}^-_t$ are independent conditional on the past, and do not share any parameters.
In particular, a separable TERGM (STERGM) can be specified by a product of a formation model and a dissolution model:
\begin{align}
  &\mathbb{P}({\bf X}_t \,|\, {\bf X}_{t-1}, \ldots, {\bf X}_{t-m} ; \btheta) \notag \\& ~~~~~~\propto \exp \big\{ \bvarsigma^+(\btheta^+)^{\T} \bvarrho^+({\bf X}^+_t, {\bf X}_{t-1}, \ldots, {\bf X}_{t-m}) + \bvarsigma^-(\btheta^-)^{\T} \bvarrho^-({\bf X}_t^-, {\bf X}_{t-1}, \ldots, {\bf X}_{t-m})\big\}\,, \label{stergm_def}
\end{align}
where $\bvarsigma^+$ and $\bvarsigma^-$ map parameter vectors $\btheta^+$ and $\btheta^-$ to the vectors of natural parameters, and $\bvarrho^+$ and $\bvarrho^-$ map the data, including the past network snapshots, to the corresponding sufficient statistics.

%Under the edge conditional independence assumption, we can write
%\begin{equation*}
%  \varrho^+({\bf X}^+_t, {\bf X}_{t-1}, \ldots, {\bf X}_{t-m}) = \sum_{i < j : X_{i,j}^{t-1}=0} \varrho_{i,j}(X_{i,j}^t, {\bf X}_{t-1}, \ldots, {\bf X}_{t-m}),
%\end{equation*}
%or equivalently
%\begin{equation*}
%  \varrho^+({\bf X}^+_t, {\bf X}_{t-1}, \ldots, {\bf X}_{t-m}) = \sum_{i < j} \tilde{\varrho}_{i,j}^+(X_{i,j}^t, {\bf X}_{t-1}, \ldots, {\bf X}_{t-m}),
%\end{equation*}
%%where for all $i < j$, $\tilde{\varrho}^+_{i,j} = 0$ if $X_{i,j}^{t-1} = 1$.
%Define $\tilde{\varrho}^-_{i,j}$ analogously for the dissolution model.

%NEW:
Under the edge conditional independence assumption, we can write
\begin{equation*} %\label{stergm_stat}
  \bvarrho^+({\bf X}^+_t, {\bf X}_{t-1}, \ldots, {\bf X}_{t-m}) = \sum_{i < j} \bvarrho_{i,j}^+(X_{i,j}^{t,+}; X_{i,j}^{t-1}, {\bf X}_{t-1} \setminus X_{i,j}^{t-1}, \ldots, {\bf X}_{t-m})\,,
\end{equation*}
or equivalently
\begin{equation*}
  \bvarrho^+({\bf X}^+_t, {\bf X}_{t-1}, \ldots, {\bf X}_{t-m}) = \sum_{i < j} \tilde{\bvarrho}_{i,j}^+(X_{i,j}^t; X_{i,j}^{t-1}, {\bf X}_{t-1} \setminus X_{i,j}^{t-1}, \ldots, {\bf X}_{t-m})\,,
\end{equation*}
since for all $i < j$, $X_{i,j}^{t,+}$ can be recovered from $X_{i,j}^t$ and $X_{i,j}^{t-1}$.
Define $\tilde{\bvarrho}^-_{i,j}$ analogously for the dissolution model.

In this way $X_t$ is an edge conditionally independent TERGM with parameter vector $(\btheta^+,\btheta^-)$, natural parameter
\begin{equation*}
  (\varsigma^+(\btheta^+),\varsigma^-(\btheta^-))\,,
\end{equation*}
and sufficient statistic
\begin{equation*}
  \bigg( \sum_{i < j} \tilde{\bvarrho}^+_{i,j}, \sum_{i < j} \tilde{\bvarrho}^-_{i,j} \bigg)\,.
\end{equation*}

%Following the above construction to rewrite this as an AR network model, we see that $\alpha_{i,j}^{t-1}$ is free of $\btheta^-$ and $\beta_{i,j}^{t-1}$ is free of $\btheta^+$.

%NEW:
Following the above construction to rewrite this as an AR network model, we can write
\begin{align*}
    \alpha_{i,j}^{t-1} %=&~ f_{i,j}({\bf X}_{t-1} \setminus X_{i,j}^{t-1}, {\bf X}_{t-2}, \ldots, {\bf X}_{t-m} ; \btheta) \\
  =&~ \sigma^{-1} \big[ \bvarsigma^+(\btheta^+)^{\T} \big\{ \tilde{\bvarrho}^+_{i,j}(1 ; 0 , {\bf X}_{t-1} \setminus X_{i,j}^{t-1}, {\bf X}_{t-2}, \ldots, {\bf X}_{t-m}) \\
  &~~~~~~~~~~~~~~~~- \tilde{\bvarrho}^+_{i,j}(0 ; 0 , {\bf X}_{t-1} \setminus X_{i,j}^{t-1}, {\bf X}_{t-2}, \ldots, {\bf X}_{t-m} )\big\} \\
  &~~~~~~~~+\bvarsigma^-(\btheta^-)^{\T} \big\{ \tilde{\bvarrho}^-_{i,j}(1 ; 0 , {\bf X}_{t-1} \setminus X_{i,j}^{t-1}, {\bf X}_{t-2}, \ldots, {\bf X}_{t-m}) \\
  &~~~~~~~~~~~~~~~~~~~- \tilde{\bvarrho}^-_{i,j}(0 ; 0 , {\bf X}_{t-1} \setminus X_{i,j}^{t-1}, {\bf X}_{t-2}, \ldots, {\bf X}_{t-m} )\big\} \big]\,,
\end{align*} 
and note that 
\begin{align*}
    \tilde{\bvarrho}^-_{i,j}(1 ; 0 , {\bf X}_{t-1} \setminus X_{i,j}^{t-1}, {\bf X}_{t-2}, \ldots, {\bf X}_{t-m}) &= \tilde{\bvarrho}^-_{i,j}(0 ; 0 , {\bf X}_{t-1} \setminus X_{i,j}^{t-1}, {\bf X}_{t-2}, \ldots, {\bf X}_{t-m} ) \\
    &= \bvarrho^-_{i,j}(0 ; 0 , {\bf X}_{t-1} \setminus X_{i,j}^{t-1}, {\bf X}_{t-2}, \ldots, {\bf X}_{t-m} )
\end{align*}
for all $i < j$ and $t$, since $X_{i,j}^{t,-}=0$ whenever $X_{i,j}^{t-1}=0$.
Thus $\alpha_{i,j}^{t-1}$ is free of $\btheta^-$. Similarly,
we can write 
\begin{align*}
    \beta_{i,j}^{t-1} %=&~ f_{i,j}({\bf X}_{t-1} \setminus X_{i,j}^{t-1}, {\bf X}_{t-2}, \ldots, {\bf X}_{t-m} ; \btheta) \\
  =&~ \sigma^{-1} \big[ \bvarsigma^+(\btheta^+)^{\T} \big\{ \tilde{\bvarrho}^+_{i,j}(0 ; 1 , {\bf X}_{t-1} \setminus X_{i,j}^{t-1}, {\bf X}_{t-2}, \ldots, {\bf X}_{t-m}) \\
  &~~~~~~~~~~~~~~~~- \tilde{\bvarrho}^+_{i,j}(1 ; 1 , {\bf X}_{t-1} \setminus X_{i,j}^{t-1}, {\bf X}_{t-2}, \ldots, {\bf X}_{t-m} )\big\} \\
  &~~~+\bvarsigma^-(\btheta^-)^{\T} \big\{ \tilde{\bvarrho}^-_{i,j}(0 ; 1 , {\bf X}_{t-1} \setminus X_{i,j}^{t-1}, {\bf X}_{t-2}, \ldots, {\bf X}_{t-m}) \\
  &~~~~~~~~~~~~~~~~- \tilde{\bvarrho}^-_{i,j}(1 ; 1 , {\bf X}_{t-1} \setminus X_{i,j}^{t-1}, {\bf X}_{t-2}, \ldots, {\bf X}_{t-m} )\big\} \big]\,,
\end{align*} 
and
\begin{align*}
    \tilde{\bvarrho}^+_{i,j}(0 ; 1 , {\bf X}_{t-1} \setminus X_{i,j}^{t-1}, {\bf X}_{t-2}, \ldots, {\bf X}_{t-m}) &= \tilde{\bvarrho}^+_{i,j}(1 ; 1 , {\bf X}_{t-1} \setminus X_{i,j}^{t-1}, {\bf X}_{t-2}, \ldots, {\bf X}_{t-m} ) \\
    &= \bvarrho^+_{i,j}(1 ; 1 , {\bf X}_{t-1} \setminus X_{i,j}^{t-1}, {\bf X}_{t-2}, \ldots, {\bf X}_{t-m} )
\end{align*}
for all $i < j$ and $t$, since $X_{i,j}^{t,+}=1$ whenever $X_{i,j}^{t-1}=1$.
Then $\beta_{i,j}^{t-1}$ is free of $\btheta^+$.
Hence, any edge conditionally independent STERGM can be written as an AR network model with separable parameters.

Conversely, suppose we have specified an AR network model such that
\begin{align*}
  \sigma(\alpha_{i,j}^{t-1}) &= \bphi(\btheta_{\alpha})^{\T} \bu_{i,j}({\bf X}_{t-1} \setminus X_{i,j}^{t-1}, {\bf X}_{t-2}, \ldots, {\bf X}_{t-m} ) \\
  \sigma(\beta_{i,j}^{t-1}) &= \bpsi(\btheta_{\beta})^{\T} \bv_{i,j}({\bf X}_{t-1} \setminus X_{i,j}^{t-1}, {\bf X}_{t-2}, \ldots, {\bf X}_{t-m} )
\end{align*}
with separable parameters $\btheta_{\alpha}$ and $\btheta_{\beta}$.

We follow the same construction as above to rewrite this model as an edge conditionally independent TERGM, with parameter vector $\bvarsigma(\btheta) = (\bphi(\btheta_{\alpha})^{\T}, \bpsi(\btheta_{\beta})^{\T})^{\T}$ and sufficient statistics $
  \bvarrho({\bf X}_t , {\bf X}_{t-1}, \ldots, {\bf X}_{t-m}) = (\bvarrho_{\alpha}^{\T} , \bvarrho_{\beta}^{\T})^{\T}$, where
$\bvarrho_{\alpha} = \sum_{i,j:\,i < j} \bvarrho_{\alpha,i,j}$ and $\bvarrho_{\beta} = \sum_{i,j:\,i < j} \bvarrho_{\beta,i,j}$ with
\begin{align}
  &\bvarrho_{\alpha,i,j}(X_{i,j}^t ; X_{i,j}^{t-1}, {\bf X}_{t-1} \setminus X_{i,j}^{t-1}, {\bf X}_{t-2}, \ldots, {\bf X}_{t-m}) \notag \\
  &~~~~~~~~~= \bu_{i,j}({\bf X}_{t-1} \setminus X_{i,j}^{t-1}, {\bf X}_{t-2}, \ldots, {\bf X}_{t-m} )I(X_{i,j}^{t}=1, X_{i,j}^{t-1}=0) \notag \\
  &~~~~~~~~~= \bu_{i,j}({\bf X}_{t-1} \setminus X_{i,j}^{t-1}, {\bf X}_{t-2}, \ldots, {\bf X}_{t-m} )I(X_{i,j}^{t,+}=1, X_{i,j}^{t-1}=0) \,, \label{stergm_statA} \\
  &\bvarrho_{\beta,i,j}(X_{i,j}^t ; X_{i,j}^{t-1}, {\bf X}_{t-1} \setminus X_{i,j}^{t-1}, {\bf X}_{t-2}, \ldots, {\bf X}_{t-m}) \notag \\&~~~~~~~~~=\bv_{i,j}({\bf X}_{t-1} \setminus X_{i,j}^{t-1}, {\bf X}_{t-2}, \ldots, {\bf X}_{t-m} )I(X_{i,j}^{t}=0, X_{i,j}^{t-1}=1)\notag \\
  &~~~~~~~~~=\bv_{i,j}({\bf X}_{t-1} \setminus X_{i,j}^{t-1}, {\bf X}_{t-2}, \ldots, {\bf X}_{t-m} )I(X_{i,j}^{t,-}=0, X_{i,j}^{t-1}=1)\,. \label{stergm_statB}
\end{align}

Note that in \eqref{stergm_statA} and \eqref{stergm_statB}, the sufficient statistics depend on $X_{i,j}^t$ only through $X_{i,j}^{t,+}$ and $X_{i,j}^{t,-}$ respectively. 
It follows that when written as a conditionally independent TERGM, the distribution factors into a product of a formation model and dissolution model, as in \eqref{stergm_def}. 
Thus, this AR network model is an edge conditionally independent STERGM.

%{\color{red}note that the sufficient statistics satisfy $\bvarrho_{\alpha,i,j} = 0$ if $X_{i,j}^{t-1}=1$ and $\bvarrho_{\beta,i,j} = 0$ if $X_{i,j}^{t-1}=1$}, so they can be reparameterized to depend on $X_{i,j}^t$ only through $X_{i,j}^{t,+}$ and $X_{i,j}^{t,-}$ respectively, as in \eqref{stergm_stat}.
%They also separate across ${\bf X}_t^+$ and ${\bf X}_t^-$, and by construction the natural parameter $\bvarsigma(\btheta) = (\bphi(\btheta_{\alpha}), \bpsi(\btheta_{\beta}))$ is separable.
%Thus, this AR network model is also an edge conditionally independent STERGM.

\section{Simulations with transitivity model} \label{sec:simulation}

In this section, we use the transitivity model introduced in Section \ref{sec:transitivity} as an example to illustrate numerical behaviour of the initial estimation proposed in Section \ref{sec:initialEstimation} and the improved estimation suggested in Section \ref{sec:improvedEstimation}.
% {\color{blue}and an alternative, computationally efficient method-of-moments-based estimation to be introduced in Section~\ref{sec:estimation2}}. 

\subsection{Implementation details}
\label{sec:impl.details}

Network data $\{\bX_1, \ldots, \bX_n\}$ used in the experiments described below are generated according to (\ref{alphaij}),\,(\ref{betaij}) and (\ref{trans-model}). %, except in Section~\ref{sec:more.gen.model}, where a slight generalization is used. 
For each sample, we generate a sequence of  length $n+200$, and discard the first 200 observations.

Regarding the implementation of our estimation procedures in Sections \ref{sec:initialEstimation} and \ref{sec:improvedEstimation}, recall  that $\mathcal{G}$ and $\mathcal{G}^{\btc}$ are, respectively, the index sets of the global parameters and the local parameters. For the transitivity model \eqref{trans-model}, we have the parameter vector $\btheta=(a,b,\xi_1,\ldots,\xi_p,\eta_1,\ldots,\eta_p)^{\T}$, where $a$ and $b$ are the global parameters and $\{\xi_i\}_{i=1}^p$ and $\{\eta_i\}_{i=1}^p$ are the local parameters.  Hence, for this model we have $|\mathcal{S}_l| = p(p-1)/2$  and $|\mathcal{S}_{l}\cap\mathcal{S}_{l'}| = p-1$ when  $l\in \mathcal{G}$ and $l'\in\mathcal{G}^{\btc}$. By \eqref{eq:ident}, for each given $l\in\mathcal{G}$ and $\btheta\in\bTheta$, it holds with probability approaching one that $\ell_{n,p}^{(l)}(\btheta_0)-\ell_{n,p}^{(l)}(\btheta)\geq \bar{C}|\btheta_{\mathcal{G}}-\btheta_{0,\mathcal{G}}|_2^2+2\bar{C}p^{-1}|\btheta_{\mathcal{G}^{\btc}}-\btheta_{0,\mathcal{G}^{\btc}}|_2^2$ for some universal constant $\bar{C}>0$ independent of $\btheta$, which means that the function  $\ell_{n,p}^{(l)}(\cdot)$ defined as \eqref{eq:idelocf} exhibits robustness against fluctuations in the values of local parameters when $l\in\mathcal{G}$ and $p$ is large. 

Motivated by this fact, when we compute the initial estimator $\tilde{\btheta}_{\mathcal{G}}$  for the global parameter vector $\btheta_{0,\mathcal{G}}$, we can just approximate $\tilde{\btheta}_{\mathcal{G}}$ by $\tilde{\btheta}_{\mathcal{G}}^{({\rm app})}=\arg\max_{\btheta_{\mathcal{G}}}\hat{\ell}_{n,p}^{(l)}(\btheta_{\mathcal{G}},\bar{\btheta}_{\mathcal{G}^{\btc}})$, for some given $\bar{\btheta}_{\mathcal{G}^{\btc}}$ and $\hat{\ell}_{n,p}^{(l)}(\btheta)$ defined as \eqref{logL:local} for some $l\in\mathcal{G}$. This simple idea can significantly improve the computational efficiency. Specifically, note that computing the original $\tilde{\btheta}_{\mathcal{G}}$ requires solving an optimization problem with $2p+2$ variables while this alternative approach only requires solving an optimization problem with two variables. Our above discussion guarantees $\tilde{\btheta}_{\mathcal{G}}^{({\rm app})}$ can approximate $\tilde{\btheta}_{\mathcal{G}}$ well. 
Similarly, when we compute the initial estimator $\tilde{\theta}_l$ for the local parameter $\theta_{0,l}$ with $l\in\mathcal{G}^{\btc}$, 
we can  approximate it by $\tilde{\theta}_{l}^{({\rm app})}=\arg\max_{\theta_l}\hat{\ell}_{n,p}^{(l)}(\tilde{\btheta}_{\mathcal{G}}^{({\rm app})},  \theta_l,  \bar{\btheta}_{\mathcal{G}^{\btc}\setminus\{l\}})$ with some given $ \bar{\btheta}_{\mathcal{G}^{\btc}\setminus\{l\}}$.

In practice, we first estimate the global parameters $a$ and $b$ via the Quasi-Newton method, 
 given certain initial values for the local parameters $\{\xi_i\}_{i=1}^p$ and $\{\eta_i\}_{i=1}^p$. To be specific, we consider 9 different sets of the initial values between 0.5 and 0.9 for $\{\xi_i\}_{i=1}^p$ and $\{\eta_i\}_{i=1}^p$, and compute $\tilde a^{(\nu)}$ and $\tilde b^{(\nu)}$  for the $\nu$-th initial setting.
% compute $\tilde a^{(\nu)}$ and $\tilde b^{(\nu)}$  by maximizing the global log-likelihood function via the Quasi-Newton method 
%  using the $\nu$-th of the
% 9 different sets of the initial values between 0.5 and 0.9 for local parameters $(\xi_i, \eta_i)$. 
%
%For the $\nu$-th initial setting,  We obtain the first-step global parameter estimates, denoted as $\tilde a^{(\nu)}$ and $\tilde b^{(\nu)}$.
%These estimates then serve as inputs for computing the first-step local parameter estimates: 
With $ a=\tilde a^{(\nu)}$ and $b=\tilde b^{(\nu)}$, we then compute $\tilde \xi_1^{(\nu)}, \dots, \tilde \xi_p^{(\nu)},$ $ \tilde \eta_1^{(\nu)}, \dots, \tilde \eta_p^{(\nu)}$ through maximizing each of the associated $\hat{\ell}_{n,p}^{(l)}(\btheta)$ with $l\in\mathcal{G}^{\btc}$. 
To reduce sensitivity to initial values and enhance numerical stability, we subsequently apply a computationally efficient method-of-moments-based estimation, which leverages the separable structure between the global and local components of the transition probabilities \eqref{trans-model} in the transitivity model. See Section~\ref{sec:estimation2} below for implementation details.
The improved estimates $\hat a^{(\nu)}, \hat b^{(\nu)},$ $ \hat \xi_1^{(\nu)}, \dots, \hat \xi_p^{(\nu)}, $ $\hat \eta_1^{(\nu)}, \dots, \hat \eta_p^{(\nu)}$ are then obtained according to \eqref{eq:thetal}
with $(\tilde r, \check r) = (0.2,0.05)$ for the local parameters  and $(\tilde r, \check r) =  (10, 2)$ for the global parameters. The tuning parameter $\tau$ required in \eqref{eq:hatanl} is selected individually for each parameter $\theta_l$ based on the principle of minimizing $\var(\check \theta_l)$, where $\check \theta_l$ is defined in \eqref{eq:theta_check}. Specifically, for each $l \in [2+2p]$, we consider a sequence of $\tilde \tau_l$ values. For each candidate value, we compute the corresponding projection vector $\hat \bvarphi_{l}$ via \eqref{eq:hatanl} with $\tau = \tilde \tau_{l} \Delta_n^{1/2}$ and subsequently obtain the updated estimator $\check \theta_l$ via \eqref{eq:theta_check}. The optimal $\tau$ is then chosen to minimize the  empirical variance approximation, i.e.,
\[
\frac{
\hat{\bvarphi}_l^\top \, \sum_{t=2}^n \big[ \bg_t^{(l)}(\check \theta_l, \tilde{\boldsymbol{\theta}}_{-l}) \, \bg_t^{(l)}(\check \theta_l, \tilde{\boldsymbol{\theta}}_{-l})^\top \big] \, \hat{\bvarphi}_l
}{
\big[\hat{\bvarphi}_l^\top \, \sum_{t=2}^n \nabla_{\theta_l} \bg_t^{(l)}(\check \theta_l, \tilde{\boldsymbol{\theta}}_{-l})  \big]^2
}\,.
\]
Finally, we apply \eqref{eq:asy1} to construct the $95\%$ confidence interval for each $\theta_l$.

\subsection{Estimation errors and coverage probabilities}
\label{sec:sim.est}
Here we report results on experiments exploring the behavior of the
initial estimator $\tilde{\btheta}$ given in \eqref{eq:estimate}, the improved estimator $\hat{\btheta}$ given in \eqref{eq:thetal}  and the coverage probabilities of the associated confidence intervals constructed using \eqref{eq:asy1}.  For simplicity, we  set all the true values for  $\{\xi_i\}_{i=1}^p$ to be the same, and those for $\{\eta_i\}_{i=1}^p$ also to be the same. The same three sets of parameter values were used as in Section~\ref{sec:simul1}.
We set $n\in\{100, 200\}$ and $p\in\{ 50, 100\}$.
For each setting, we replicate the estimation 100 times. The simulations in this section utilise our development R package \texttt{arnetworks}, which provides a user-friendly implementation of the estimation and inference procedure described above. 

Table~\ref{table.dynamic.err1} presents the means  and the standard errors, over the 100 replications, of the relative mean absolute errors (rMAE):
%over the nine initial levels. For illustrative purposes, the rMAEs for $\hat\theta_i$ and $\hat a$ are calculated as follows:
$$
 \text{rMAE}(\hat \xi_i)= \frac{1}{9}\sum_{\nu = 1}^9\frac{1}{p}\sum_{i =1}^{p} \bigg|\frac{\hat \xi_i^{(\nu)} - \xi_i}{\xi_i}\bigg| ~~\text{and}~~ \text{rMAE}(\hat a)= \frac{1}{9}\sum_{\nu = 1}^9\bigg|\frac{\hat a^{(\nu)} - a}{a}\bigg|\,,
$$
where the sum over $\nu$ corresponds to taking the average over the 9 initial values discussed in Section~\ref{sec:impl.details}.  Also reported are the coverage probabilities of the $95\%$ confidence intervals constructed based on the limiting distribution in \eqref{eq:asy1}, for local parameters $\xi_i$ and $\eta_i$ with $i=p$, and global parameters $a$ and $b$.  Since each replication includes results from 9 different initializations, the coverage probabilities are based on a total of 900 experiments per setting.
%When sample size $n$ increases, the errors and standard deviations for both the estimators decrease.
%Furthermore, 
 Several conclusions can be drawn here. First, the improved estimator $\hat \btheta$ significantly outperforms the initial estimator $\tilde{\btheta}$ in terms of estimation accuracy across all settings. For example, under the setting $(0.8,0.9,10,10)$, the rMAE for estimating 
$a$ decreases from 1.248 to 0.049 when $n = 100$ and $p = 100$, reflecting an approximate 95\% improvement. 
This demonstrates the effectiveness of the projection strategy in \eqref{eq:hatanl} for computing $\hat{\btheta}$, as it successfully reduces the impact of nuisance parameter vector $\btheta_{-l}$  and substantially improves the estimation of the target parameter $\theta_{0,l}$.
Second, the improved estimator $\hat \btheta$ achieves higher accuracy as $n$ increases. Notably, when $n = 200$, we observe comparable or even better estimation performance as the network size $p$ increases, which highlights the robustness and scalability of our method in high-dimensional settings.
Lastly,  the coverage probabilities remain consistently close to the nominal level $95\%$, indicating the validity of the proposed asymptotic normal approximation in \eqref{eq:asy1}.

% Lastly, the setting $(0.7,0.8,30,15)$ attains the lowest overall estimation errors, with IMoM  $\mathring{\btheta}$ outperforming the improved estimator $\widehat{\btheta}$.
% This is well-expected, as this is the most dynamic setting among the four settings considered.
 %See Figs.\ref{trans_plot_50} and \ref{trans_plot_100}--\ref{trans_plot_150} of the Supplementary Material. 

%For the transitivity model \eqref{trans-model}, the local parameters $\{\xi_i\}_{i=1}^p$ and $\{\eta_i\}_{i=1}^p$ are separable in the sense that $\{\xi_i\}_{i=1}^p$ only appear in $\alpha_{i,j}^{t-1}(\btheta)$'s while $\{\eta_i\}_{i=1}^p$ only appear in $\beta_{i,j}^{t-1}(\btheta)$'s.
%Due to $\gamma_{i,j}^{t-1}(\btheta) =\alpha_{i,j}^{t-1}(\btheta) +X_{i,j}^{t-1}\{1 - \alpha_{i,j}^{t-1}(\btheta) -\beta_{i,j}^{t-1}(\btheta)\}$, then $\gamma_{i,j}^{t-1}(\btheta)$ does not include the information of $(\xi_i,\xi_j)$ if $X_{i,j}^{t-1}=1$. For given observations $\bX_1,\ldots,\bX_n$, to estimate $\xi_i$, the effective samples for calculating the associated $\hat{\ell}_{n,p}^{(l)}(\btheta)$ defined as \eqref{logL:local} consist of those $X_{i,j}^{t-1}$, $t\in[n]\setminus[m]$ and $j\neq i$, such that $X_{i,j}^{t-1}=0$. Analogously, to estimate $\eta_i$, the effective samples for calculating the associated  $\hat{\ell}_{n,p}^{(l)}(\btheta)$ defined as \eqref{logL:local} consist of those $X_{i,j}^{t-1}$, $t\in[n]\setminus[m]$ and $j\neq i$, such that $X_{i,j}^{t-1}=1$.

\begin{table}[tbp]
\vspace{-1cm}
 %Means (standard errors) of 
	\caption{The means and  STDs (in parenthesis) of rMAEs for estimating parameters in the transitivity model (\ref{trans-model}) over 100 replications,  and the corresponding coverage probabilities of the 95\% confidence intervals.\label{table.dynamic.err1}}
	\begin{center}
		\vspace{-0.5cm}
		\resizebox{6.7in}{!}{
\begin{tabular}{ccccccccccccccc}
\hline
\multirow{2}{*}{$(\xi_i,   \eta_i, a, b)$} &  & \multirow{2}{*}{$p$} &  & \multirow{2}{*}{} &  & \multicolumn{4}{c}{$n=100$}                                   &  & \multicolumn{4}{c}{$n=200$}                                   \\ \cline{7-10} \cline{12-15} 
                                           &  &                      &  &                   &  & $\xi_i$       & $\eta_i$      & $a$           & $b$           &  & $\xi_i$       & $\eta_i$      & $a$           & $b$           \\ \hline
(0.8, 0.9, 10, 10)                         &  & 50                   &  & Initial           &  & 0.131 (0.001) & 0.133 (0.010)  & 0.743 (0.045) & 0.157 (0.002) &  & 0.128 (0.001) & 0.134 (0.008) & 0.737 (0.031) & 0.157 (0.001) \\
                                           &  &                      &  & Improved          &  & 0.083 (0.006) & 0.033 (0.003) & 0.063 (0.044) & 0.020 (0.012)  &  & 0.068 (0.005) & 0.023 (0.002) & 0.042 (0.033) & 0.019 (0.011) \\
                                           &  &                      &  & Coverage          &  & 0.951         & 0.969         & 0.933         & 0.948         &  & 0.927         & 0.958         & 0.958         & 0.949         \\
                                           &  & 100                  &  & Initial           &  & 0.140 (0.001)  & 0.122 (0.005) & 1.248 (0.034) & 0.154 (0.001) &  & 0.139 (0.001) & 0.123 (0.004) & 1.246 (0.025) & 0.154 (0.001) \\
                                           &  &                      &  & Improved          &  & 0.088 (0.011) & 0.049 (0.011) & 0.049 (0.027) & 0.058 (0.035) &  & 0.055 (0.006) & 0.024 (0.002) & 0.058 (0.019) & 0.026 (0.003) \\
                                           &  &                      &  & Coverage          &  & 0.953         & 0.922         & 0.955         & 0.931         &  & 0.953         & 0.934         & 0.941         & 0.950          \\ \hline
(0.8, 0.9, 25, 15)                         &  & 50                   &  & Initial           &  & 0.129 (0.002) & 0.144 (0.012) & 0.305 (0.015) & 0.140 (0.001)  &  & 0.126 (0.001) & 0.145 (0.008) & 0.304 (0.011) & 0.140 (0.001)  \\
                                           &  &                      &  & Improved          &  & 0.077 (0.008) & 0.031 (0.003) & 0.031 (0.022) & 0.017 (0.011) &  & 0.058 (0.006) & 0.022 (0.002) & 0.022 (0.017) & 0.015 (0.008) \\
                                           &  &                      &  & Coverage          &  & 0.956         & 0.942         & 0.950          & 0.954         &  & 0.940          & 0.952         & 0.963         & 0.944         \\
                                           &  & 100                  &  & Initial           &  & 0.142 (0.001) & 0.133 (0.008) & 0.455 (0.020)  & 0.141 (0.001) &  & 0.140 (0.001)  & 0.133 (0.005) & 0.457 (0.015) & 0.141 (0.001) \\
                                           &  &                      &  & Improved          &  & 0.064 (0.003) & 0.026 (0.002) & 0.025 (0.015) & 0.014 (0.005) &  & 0.055 (0.002) & 0.018 (0.001) & 0.018 (0.011) & 0.012 (0.005) \\
                                           &  &                      &  & Coverage          &  & 0.937         & 0.946         & 0.944         & 0.952         &  & 0.950          & 0.940          & 0.947         & 0.959         \\ \hline
(0.7, 0.9, 15, 10)                         &  & 50                   &  & Initial           &  & 0.159 (0.002) & 0.122 (0.010)  & 0.433 (0.015) & 0.196 (0.004) &  & 0.155 (0.001) & 0.122 (0.007) & 0.432 (0.011) & 0.196 (0.003) \\
                                           &  &                      &  & Improved          &  & 0.082 (0.007) & 0.039 (0.003) & 0.046 (0.035) & 0.019 (0.011) &  & 0.066 (0.005) & 0.026 (0.002) & 0.031 (0.024) & 0.019 (0.012) \\
                                           &  &                      &  & Coverage          &  & 0.938         & 0.954         & 0.946         & 0.946         &  & 0.943         & 0.955         & 0.967         & 0.943         \\
                                           &  & 100                  &  & Initial           &  & 0.170 (0.001)  & 0.111 (0.004) & 0.653 (0.013) & 0.196 (0.002) &  & 0.168 (0.001) & 0.110 (0.003)  & 0.651 (0.009) & 0.196 (0.001) \\
                                           &  &                      &  & Improved          &  & 0.131 (0.011) & 0.086 (0.007) & 0.121 (0.026) & 0.146 (0.020)  &  & 0.048 (0.003) & 0.029 (0.003) & 0.027 (0.015) & 0.025 (0.006) \\
                                           &  &                      &  & Coverage          &  & 0.964         & 0.931         & 0.958         & 0.953         &  & 0.959         & 0.932         & 0.934         & 0.947  \\ \hline     
\end{tabular}
}

	\end{center}
	 \vspace{-0.5cm}
\end{table}

\section{Additional simulation results}  \label{app:simulation}

Further to Section \ref{sec:simulation}, we report more simulation results on the transitivity model.

\begin{skip}
\subsection{More plots for the simulated transitivity models}
 \label{app:transitivity}
We provide two additional figures for the simulated transitivity models
in Section~\ref{sec:simul1}.
Figure~\ref{trans_plot_100} and \ref{trans_plot_150} display the time series plots of simulated  $\{D_t\}_{t=2}^{200}$, $\{D_{1,t}\}_{t=2}^{200}$, $\{D_{0,t}\}_{t=2}^{200}$, $\{\bar D_t\}_{t=2}^{200}$, $\{\bar D_{1,t}\}_{t=2}^{200}$ and $\{\bar D_{0,t}\}_{t=2}^{200}$ when $p=100$ and $150$, respectively. All settings are stationary, among which $(\xi_i, \eta_i, a, b) = (0.7,0.8,30,15)$ shows the most frequent edge-changing behaviour.

\begin{figure}[!htbp]
\centering
\begin{subfigure}{.3\linewidth}
    \centering
    \centering
    \includegraphics[width=5cm]{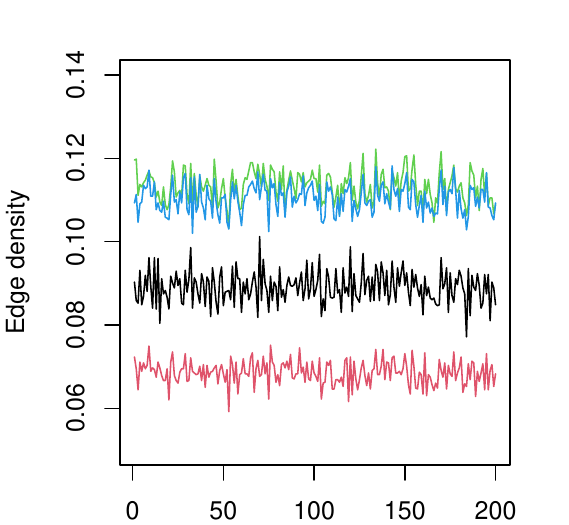}
    \caption{$D_t$}
\end{subfigure}
\begin{subfigure}{.3\linewidth}
    \centering
    \includegraphics[width=5cm]{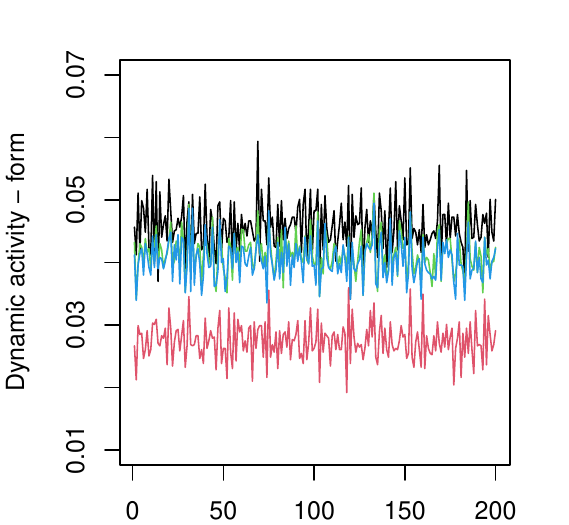}
    \caption{$D_{1,t}$}
\end{subfigure}
\begin{subfigure}{.3\linewidth}
    \centering
    \includegraphics[width=5cm]{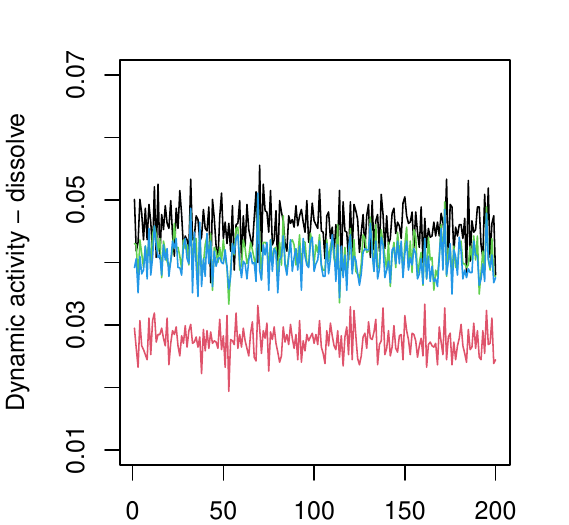}
    \caption{$D_{0,t}$}
\end{subfigure}

\begin{subfigure}{.3\linewidth}
    \centering
    \centering
    \includegraphics[width=5cm]{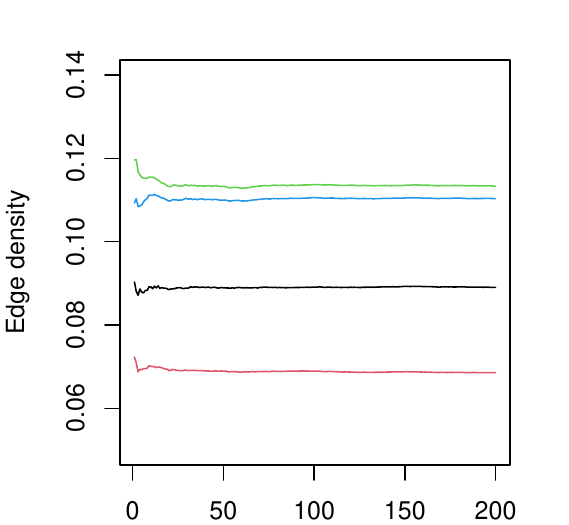}
    \caption{$\bar D_t$}
\end{subfigure}
\begin{subfigure}{.3\linewidth}
    \centering
    \includegraphics[width=5cm]{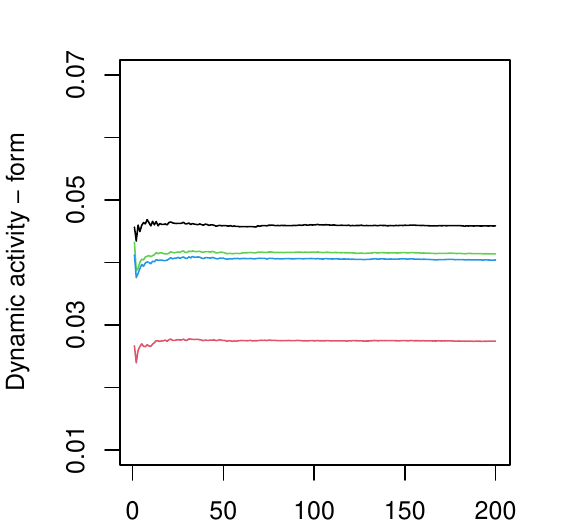}
    \caption{$\bar D_{1,t}$}
\end{subfigure}
\begin{subfigure}{.3\linewidth}
    \centering
    \includegraphics[width=5cm]{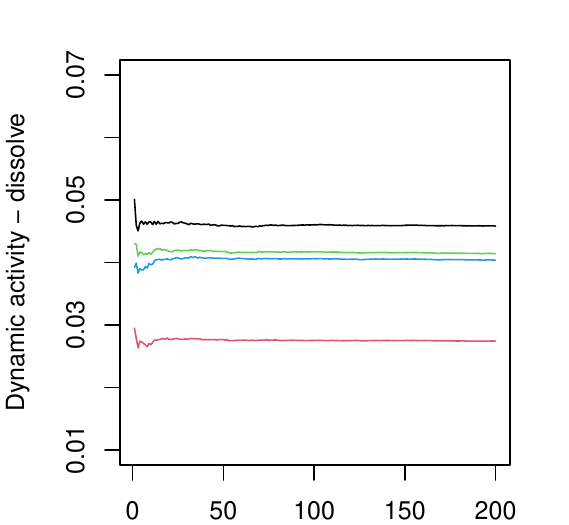}
    \caption{$\bar D_{0,t}$}
\end{subfigure}

\centering
\caption{\label{trans_plot_100}{Time series plots of $\{D_t\}_{t=2}^{200}$, $\{D_{1,t}\}_{t=2}^{200}$, $\{D_{0,t}\}_{t=2}^{200}$, $\{\bar D_t\}_{t=2}^{200}$, $\{\bar D_{1,t}\}_{t=2}^{200}$ and $\{\bar D_{0,t}\}_{t=2}^{200}$ for the four simulated settings  for $p = 100$. The black, red, green and blue curves correspond to the setting $(\xi_i,  \eta_i, a, b) = (0.7,0.8,30,15),$ $(0.6,0.7,20,20)$, $(0.6,0.7,15,10)$ and $(0.6,0.7,10,10)$, respectively. }}
\end{figure}

\begin{figure}[!htbp]
\centering
\begin{subfigure}{.3\linewidth}
    \centering
    \centering
    \includegraphics[width=5cm]{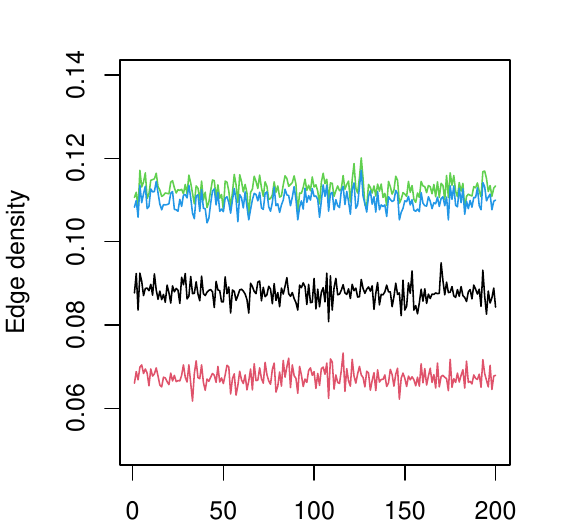}
    \caption{$D_t$}
\end{subfigure}
\begin{subfigure}{.3\linewidth}
    \centering
    \includegraphics[width=5cm]{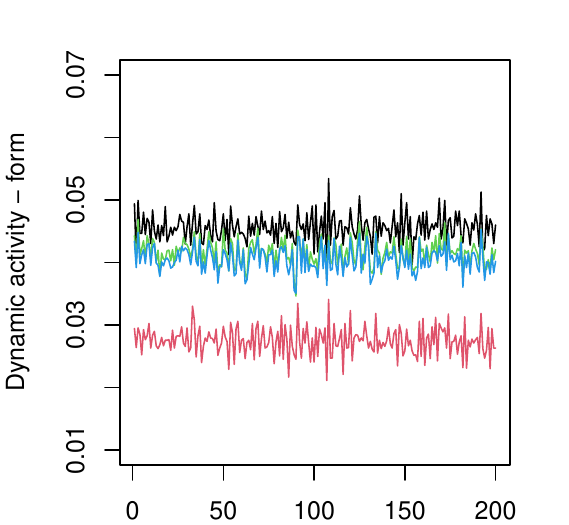}
    \caption{$D_{1,t}$}
\end{subfigure}
\begin{subfigure}{.3\linewidth}
    \centering
    \includegraphics[width=5cm]{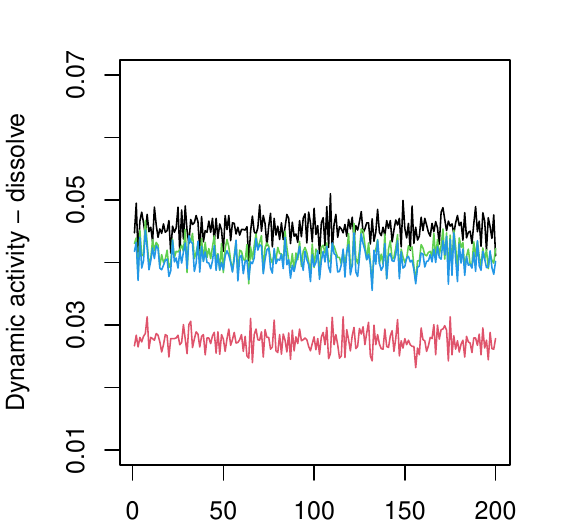}
    \caption{$D_{0,t}$}
\end{subfigure}

\begin{subfigure}{.3\linewidth}
    \centering
    \centering
    \includegraphics[width=5cm]{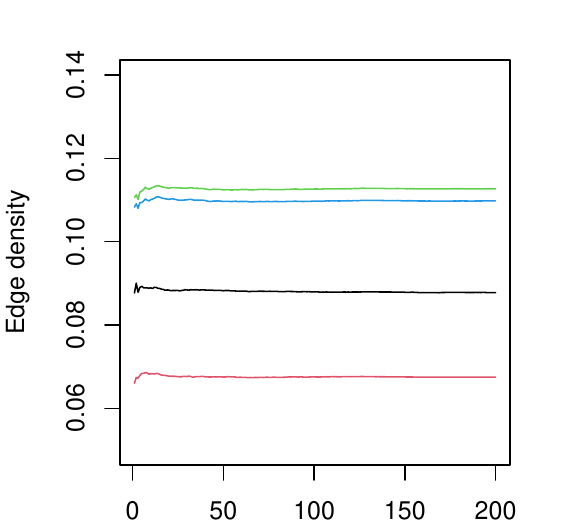}
    \caption{$\bar D_t$}
\end{subfigure}
\begin{subfigure}{.3\linewidth}
    \centering
    \includegraphics[width=5cm]{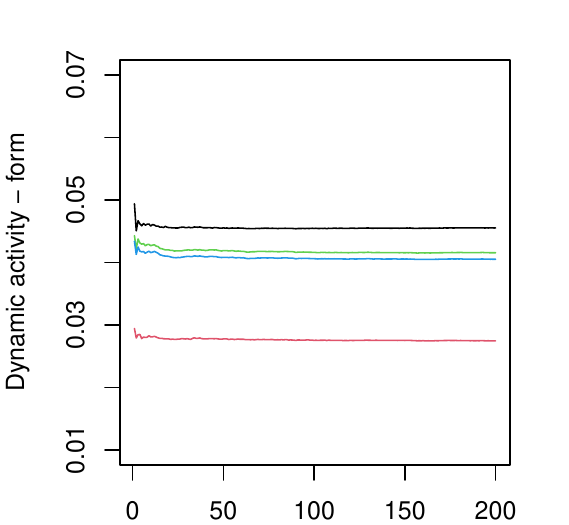}
    \caption{$\bar D_{1,t}$}
\end{subfigure}
\begin{subfigure}{.3\linewidth}
    \centering
    \includegraphics[width=5cm]{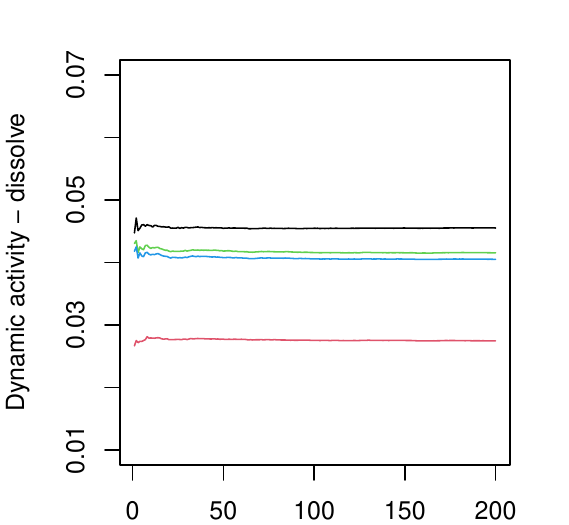}
    \caption{$\bar D_{0,t}$}
\end{subfigure}

\centering
\caption{\label{trans_plot_150}{Time series plots of $\{D_t\}_{t=2}^{200}$, $\{D_{1,t}\}_{t=2}^{200}$, $\{D_{0,t}\}_{t=2}^{200}$, $\{\bar D_t\}_{t=2}^{200}$, $\{\bar D_{1,t}\}_{t=2}^{200}$ and $\{\bar D_{0,t}\}_{t=2}^{200}$ for the four simulated settings  for $p = 150$. The black, red, green and blue curves correspond to the setting $(\xi_i,  \eta_i, a, b) = (0.7,0.8,30,15),$ $(0.6,0.7,20,20)$, $(0.6,0.7,15,10)$ and $(0.6,0.7,10,10)$, respectively. }}
\end{figure}
\subsection{Stationarity and ergodicity}% of degree heterogeneity models}
\label{app:density-dependent}
In this section, we investigate the stationarity and ergodicity of our proposed dynamic models with more synthetic data. We simulate $\{\bX_1, \dots, \bX_n\}$ from two dependent AR(1) networks, i.e., degree heterogeneity model as in \eqref{densityalphaij}, and transitivity model in \eqref{trans-model}.
For each sample, we generate a sequence of length $n +200$ and discard the first $200$ observations. For the degree heterogeneity model, we consider 625 settings with $a_0=a_1\in\{0.2,0.4,0.6,0.8,1\}$, $ b_0=b_1\in\{0.2,0.4,0.6,0.8,1\}$, $\xi_1 = \dots= \xi_p \in \{0.1,0.3,0.5,0.7,0.9\}$ and $ \eta_1 = \dots =\eta_p \in \{0.1,0.3,0.5,0.7,0.9\}$. For transitivity model, we focus on another 625 settings with $a \in\{5,15,25,35,45\}, b \in\{5,15,25,35,45\}$, $\xi_1 = \dots =  \xi_p \in \{0.1,0.3,0.5,0.7,0.9\}$ and $ \eta_1 = \dots =\eta_p \in \{0.1,0.3,0.5,0.7,0.9\}$. 
% It is worth mentioning that we adopt a wider range of $a$ and $b$ in the transitivity model aiming to record the edge-changing behavior of dependent dynamic networks over a broader spectrum of the parameter space.

For each setting, we present the  time series plots of simulated  $\{D_t\}_{t=2}^{500}$, $\{D_{1,t}\}_{t=2}^{500}$, $\{D_{0,t}\}_{t=2}^{500}$, $\{\bar D_t\}_{t=2}^{500}$, $\{\bar D_{1,t}\}_{t=2}^{500}$ and $\{\bar D_{0,t}\}_{t=2}^{500}$, defined as in Section~\ref{sec:simulation}, when $p = 50$. Figures \ref{checkdd1}--\ref{checkdd_end} and Figures \ref{checkts1}--\ref{checkts_end} report the dynamic patterns of the degree heterogeneity model and transitivity models, respectively. Several conclusions can be drawn. First,  all simulated series $\{D_t\}_{t\geq2}$, $\{D_{1,t}\}_{t\geq2}$ and $\{D_{0,t}\}_{t\geq2}$ exhibit stationary behavior. Similar patterns are discussed in Section~\ref{sec:simulation}. Second,  in contrast to many exponential random graph models (ERGMs), the proposed autoregressive (AR) model design is free from degeneracy problem when $\alpha_{i,j}^{t-1}$ and $\beta_{i,j}^{t-1}$ are strictly between 0 and 1.
Note that we do observe instances of nearly complete or empty graphs when $\xi_i = 0.1$ or $ \eta_i = 0.1$. Such degenerate graphs fall outside the scope of our discussion as the node-specific parameters are set to values as low as $0.01$. Examples of this can be seen in Figures \ref{checkdd1}--\ref{checkdd_er2} and \ref{checkts1}--\ref{checkts_er2}.

\begin{figure}[t]
	\centering
\includegraphics[width=\textwidth]{plot/check/dd/check_setting_dd_01.pdf}
    \\ 
   \includegraphics[width=\textwidth]{plot/check/dd/check_setting_dd_02.pdf}
    \\[\smallskipamount]
    \includegraphics[width=\textwidth]{plot/check/dd/check_setting_dd_03.pdf}
    \caption{\label{checkdd1}{Time series plots of $\{D_t\}_{t=2}^{500}$, $\{D_{1,t}\}_{t=2}^{500}$, $\{D_{0,t}\}_{t=2}^{500}$ (left to right) for degree heterogeneity model when $\xi_i = 0.1$ and $\eta_i \in (0.1, 0.3, 0.5).$ The 25 colored curves in each subplot correspond to 25 settings where $a_0=a_1\in\{0.2,0.4,0.6,0.8,1\}$ and $b_0=b_1\in\{0.2,0.4,0.6,0.8,1\}$.}}
\end{figure}

\begin{figure}[t]
	\centering
\includegraphics[width=\textwidth]{plot/check/dd/check_setting_dd_er_01.pdf}
    \\[\smallskipamount]
   \includegraphics[width=\textwidth]{plot/check/dd/check_setting_dd_er_02.pdf}
    \\[\smallskipamount]
    \includegraphics[width=\textwidth]{plot/check/dd/check_setting_dd_er_03.pdf}
    \caption{\label{checkdd_er1}{Time series plots of  $\{\bar D_t\}_{t=2}^{500}$, $\{\bar D_{1,t}\}_{t=2}^{500}$ and $\{\bar D_{0,t}\}_{t=2}^{500}$ (left to right) for degree heterogeneity model when $\xi_i = 0.1$ and $\eta_i \in (0.1, 0.3, 0.5).$ The 25 colored curves in each subplot correspond to 25 settings where $a_0=a_1\in\{0.2,0.4,0.6,0.8,1\}$ and $b_0=b_1\in\{0.2,0.4,0.6,0.8,1\}$.}}
\end{figure}

\begin{figure}[t]
	\centering
    \includegraphics[width=\textwidth]{plot/check/dd/check_setting_dd_04.pdf}
    \\[\smallskipamount]
    \includegraphics[width=\textwidth]{plot/check/dd/check_setting_dd_05.pdf}
    \caption{{Time series plots of $\{D_t\}_{t=2}^{500}$, $\{D_{1,t}\}_{t=2}^{500}$, $\{D_{0,t}\}_{t=2}^{500}$ (left to right) for degree heterogeneity model when $\xi_i = 0.1$ and $\eta_i \in (0.7, 0.9).$ The 25 colored curves in each subplot correspond to 25 settings where $a_0=a_1\in\{0.2,0.4,0.6,0.8,1\}$ and $b_0=b_1\in\{0.2,0.4,0.6,0.8,1\}$.}}
\end{figure}

\begin{figure}[t]
	\centering
    \includegraphics[width=\textwidth]{plot/check/dd/check_setting_dd_er_04.pdf}
    \\[\smallskipamount]
    \includegraphics[width=\textwidth]{plot/check/dd/check_setting_dd_er_05.pdf}
    \caption{\label{checkdd_er2}{Time series plots of  $\{\bar D_t\}_{t=2}^{500}$, $\{\bar D_{1,t}\}_{t=2}^{500}$ and $\{\bar D_{0,t}\}_{t=2}^{500}$ (left to right) for degree heterogeneity model when $\xi_i = 0.1$ and $\eta_i \in (0.7, 0.9).$ The 25 colored curves in each subplot correspond to 25 settings where $a_0=a_1\in\{0.2,0.4,0.6,0.8,1\}$ and $b_0=b_1\in\{0.2,0.4,0.6,0.8,1\}$.}}
\end{figure}

\begin{figure}[t]
	\centering
\includegraphics[width=\textwidth]{plot/check/dd/check_setting_dd_06.pdf}
    \\[\smallskipamount]
   \includegraphics[width=\textwidth]{plot/check/dd/check_setting_dd_07.pdf}
    \\[\smallskipamount]
    \includegraphics[width=\textwidth]{plot/check/dd/check_setting_dd_08.pdf}
    \caption{{Time series plots of $\{D_t\}_{t=2}^{500}$, $\{D_{1,t}\}_{t=2}^{500}$, $\{D_{0,t}\}_{t=2}^{500}$ (left to right)  for degree heterogeneity model when $\xi_i = 0.3$ and $\eta_i \in (0.1, 0.3, 0.5).$ The 25 colored curves in each subplot correspond to 25 settings where $a_0=a_1\in\{0.2,0.4,0.6,0.8,1\}$ and $b_0=b_1\in\{0.2,0.4,0.6,0.8,1\}$.}}
\end{figure}

\begin{figure}[t]
	\centering
\includegraphics[width=\textwidth]{plot/check/dd/check_setting_dd_er_06.pdf}
    \\[\smallskipamount]
   \includegraphics[width=\textwidth]{plot/check/dd/check_setting_dd_er_07.pdf}
    \\[\smallskipamount]
    \includegraphics[width=\textwidth]{plot/check/dd/check_setting_dd_er_08.pdf}
    \caption{{Time series plots of  $\{\bar D_t\}_{t=2}^{500}$, $\{\bar D_{1,t}\}_{t=2}^{500}$ and $\{\bar D_{0,t}\}_{t=2}^{500}$ (left to right)  for degree heterogeneity model when $\xi_i = 0.3$ and $\eta_i \in (0.1, 0.3, 0.5).$ The 25 colored curves in each subplot correspond to 25 settings where $a_0=a_1\in\{0.2,0.4,0.6,0.8,1\}$ and $b_0=b_1\in\{0.2,0.4,0.6,0.8,1\}$.}}
\end{figure}

\begin{figure}[t]
	\centering
    \includegraphics[width=\textwidth]{plot/check/dd/check_setting_dd_09.pdf}
    \\[\smallskipamount]
    \includegraphics[width=\textwidth]{plot/check/dd/check_setting_dd_10.pdf}
    \caption{{Time series plots of $\{D_t\}_{t=2}^{500}$, $\{D_{1,t}\}_{t=2}^{500}$, $\{D_{0,t}\}_{t=2}^{500}$ (left to right)  for degree heterogeneity model when $\xi_i = 0.3$ and $\eta_i \in (0.7, 0.9).$ The 25 colored curves in each subplot correspond to 25 settings where $a_0=a_1\in\{0.2,0.4,0.6,0.8,1\}$ and $b_0=b_1\in\{0.2,0.4,0.6,0.8,1\}$.}}
\end{figure}

\begin{figure}[t]
	\centering
    \includegraphics[width=\textwidth]{plot/check/dd/check_setting_dd_er_09.pdf}
    \\[\smallskipamount]
    \includegraphics[width=\textwidth]{plot/check/dd/check_setting_dd_er_10.pdf}
    \caption{{Time series plots of  $\{\bar D_t\}_{t=2}^{500}$, $\{\bar D_{1,t}\}_{t=2}^{500}$ and $\{\bar D_{0,t}\}_{t=2}^{500}$ (left to right)  for degree heterogeneity model when $\xi_i = 0.3$ and $\eta_i \in (0.7, 0.9).$ The 25 colored curves in each subplot correspond to 25 settings where $a_0=a_1\in\{0.2,0.4,0.6,0.8,1\}$ and $b_0=b_1\in\{0.2,0.4,0.6,0.8,1\}$.}}
\end{figure}

\begin{figure}[t]
	\centering
\includegraphics[width=\textwidth]{plot/check/dd/check_setting_dd_11.pdf}
    \\[\smallskipamount]
   \includegraphics[width=\textwidth]{plot/check/dd/check_setting_dd_12.pdf}
    \\[\smallskipamount]
    \includegraphics[width=\textwidth]{plot/check/dd/check_setting_dd_13.pdf}
    \caption{{Time series plots of $\{D_t\}_{t=2}^{500}$, $\{D_{1,t}\}_{t=2}^{500}$, $\{D_{0,t}\}_{t=2}^{500}$ (left to right) for degree heterogeneity model when $\xi_i = 0.5$ and $\eta_i \in (0.1, 0.3, 0.5).$ The 25 colored curves in each subplot correspond to 25 settings where $a_0=a_1\in\{0.2,0.4,0.6,0.8,1\}$ and $b_0=b_1\in\{0.2,0.4,0.6,0.8,1\}$.}}
\end{figure}

\begin{figure}[t]
	\centering
\includegraphics[width=\textwidth]{plot/check/dd/check_setting_dd_er_11.pdf}
    \\[\smallskipamount]
   \includegraphics[width=\textwidth]{plot/check/dd/check_setting_dd_er_12.pdf}
    \\[\smallskipamount]
    \includegraphics[width=\textwidth]{plot/check/dd/check_setting_dd_er_13.pdf}
    \caption{{Time series plots of  $\{\bar D_t\}_{t=2}^{500}$, $\{\bar D_{1,t}\}_{t=2}^{500}$ and $\{\bar D_{0,t}\}_{t=2}^{500}$ (left to right) for degree heterogeneity model when $\xi_i = 0.5$ and $\eta_i \in (0.1, 0.3, 0.5).$ The 25 colored curves in each subplot correspond to 25 settings where $a_0=a_1\in\{0.2,0.4,0.6,0.8,1\}$ and $b_0=b_1\in\{0.2,0.4,0.6,0.8,1\}$.}}
\end{figure}

\begin{figure}[t]
	\centering
    \includegraphics[width=\textwidth]{plot/check/dd/check_setting_dd_14.pdf}
    \\[\smallskipamount]
    \includegraphics[width=\textwidth]{plot/check/dd/check_setting_dd_15.pdf}
    \caption{{Time series plots of $\{D_t\}_{t=2}^{500}$, $\{D_{1,t}\}_{t=2}^{500}$, $\{D_{0,t}\}_{t=2}^{500}$ (left to right) for degree heterogeneity model when $\xi_i = 0.5$ and $\eta_i \in (0.7, 0.9).$ The 25 colored curves in each subplot correspond to 25 settings where $a_0=a_1\in\{0.2,0.4,0.6,0.8,1\}$ and $b_0=b_1\in\{0.2,0.4,0.6,0.8,1\}$.}}
\end{figure}

\begin{figure}[t]
	\centering
    \includegraphics[width=\textwidth]{plot/check/dd/check_setting_dd_er_14.pdf}
    \\[\smallskipamount]
    \includegraphics[width=\textwidth]{plot/check/dd/check_setting_dd_er_15.pdf}
    \caption{{Time series plots of  $\{\bar D_t\}_{t=2}^{500}$, $\{\bar D_{1,t}\}_{t=2}^{500}$ and $\{\bar D_{0,t}\}_{t=2}^{500}$ (left to right) for degree heterogeneity model when $\xi_i = 0.5$ and $\eta_i \in (0.7, 0.9).$ The 25 colored curves in each subplot correspond to 25 settings where $a_0=a_1\in\{0.2,0.4,0.6,0.8,1\}$ and $b_0=b_1\in\{0.2,0.4,0.6,0.8,1\}$.}}
\end{figure}

\begin{figure}[t]
	\centering
\includegraphics[width=\textwidth]{plot/check/dd/check_setting_dd_16.pdf}
    \\[\smallskipamount]
   \includegraphics[width=\textwidth]{plot/check/dd/check_setting_dd_17.pdf}
    \\[\smallskipamount]
    \includegraphics[width=\textwidth]{plot/check/dd/check_setting_dd_18.pdf}
    \caption{{Time series plots of $\{D_t\}_{t=2}^{500}$, $\{D_{1,t}\}_{t=2}^{500}$, $\{D_{0,t}\}_{t=2}^{500}$ (left to right) for degree heterogeneity model when $\xi_i = 0.7$ and $\eta_i \in (0.1, 0.3, 0.5).$ The 25 colored curves in each subplot correspond to 25 settings where $a_0=a_1\in\{0.2,0.4,0.6,0.8,1\}$ and $b_0=b_1\in\{0.2,0.4,0.6,0.8,1\}$.}}
\end{figure}

\begin{figure}[t]
	\centering
\includegraphics[width=\textwidth]{plot/check/dd/check_setting_dd_er_16.pdf}
    \\[\smallskipamount]
   \includegraphics[width=\textwidth]{plot/check/dd/check_setting_dd_er_17.pdf}
    \\[\smallskipamount]
    \includegraphics[width=\textwidth]{plot/check/dd/check_setting_dd_er_18.pdf}
    \caption{{Time series plots of  $\{\bar D_t\}_{t=2}^{500}$, $\{\bar D_{1,t}\}_{t=2}^{500}$ and $\{\bar D_{0,t}\}_{t=2}^{500}$ (left to right) for degree heterogeneity model when $\xi_i = 0.7$ and $\eta_i \in (0.1, 0.3, 0.5).$ The 25 colored curves in each subplot correspond to 25 settings where $a_0=a_1\in\{0.2,0.4,0.6,0.8,1\}$ and $b_0=b_1\in\{0.2,0.4,0.6,0.8,1\}$.}}
\end{figure}

\begin{figure}[t]
	\centering
    \includegraphics[width=\textwidth]{plot/check/dd/check_setting_dd_19.pdf}
    \\[\smallskipamount]
    \includegraphics[width=\textwidth]{plot/check/dd/check_setting_dd_20.pdf}
    \caption{{Time series plots of $\{D_t\}_{t=2}^{500}$, $\{D_{1,t}\}_{t=2}^{500}$, $\{D_{0,t}\}_{t=2}^{500}$ (left to right) for degree heterogeneity model when $\xi_i = 0.7$ and $\eta_i \in (0.7, 0.9).$ The 25 colored curves in each subplot correspond to 25 settings where $a_0=a_1\in\{0.2,0.4,0.6,0.8,1\}$ and $b_0=b_1\in\{0.2,0.4,0.6,0.8,1\}$.}}
\end{figure}

\begin{figure}[t]
	\centering
    \includegraphics[width=\textwidth]{plot/check/dd/check_setting_dd_er_19.pdf}
    \\[\smallskipamount]
    \includegraphics[width=\textwidth]{plot/check/dd/check_setting_dd_er_20.pdf}
    \caption{{Time series plots of  $\{\bar D_t\}_{t=2}^{500}$, $\{\bar D_{1,t}\}_{t=2}^{500}$ and $\{\bar D_{0,t}\}_{t=2}^{500}$ (left to right) for degree heterogeneity model when $\xi_i = 0.7$ and $\eta_i \in (0.7, 0.9).$ The 25 colored curves in each subplot correspond to 25 settings where $a_0=a_1\in\{0.2,0.4,0.6,0.8,1\}$ and $b_0=b_1\in\{0.2,0.4,0.6,0.8,1\}$.}}
\end{figure}

\begin{figure}[t]
	\centering
\includegraphics[width=\textwidth]{plot/check/dd/check_setting_dd_21.pdf}
    \\[\smallskipamount]
   \includegraphics[width=\textwidth]{plot/check/dd/check_setting_dd_22.pdf}
    \\[\smallskipamount]
    \includegraphics[width=\textwidth]{plot/check/dd/check_setting_dd_23.pdf}
    \caption{{Time series plots of $\{D_t\}_{t=2}^{500}$, $\{D_{1,t}\}_{t=2}^{500}$, $\{D_{0,t}\}_{t=2}^{500}$ (left to right) for degree heterogeneity model when $\xi_i = 0.9$ and $\eta_i \in (0.1, 0.3, 0.5).$ The 25 colored curves in each subplot correspond to 25 settings where $a_0=a_1\in\{0.2,0.4,0.6,0.8,1\}$ and $b_0=b_1\in\{0.2,0.4,0.6,0.8,1\}$.}}
\end{figure}

\begin{figure}[t]
	\centering
\includegraphics[width=\textwidth]{plot/check/dd/check_setting_dd_er_21.pdf}
    \\[\smallskipamount]
   \includegraphics[width=\textwidth]{plot/check/dd/check_setting_dd_er_22.pdf}
    \\[\smallskipamount]
    \includegraphics[width=\textwidth]{plot/check/dd/check_setting_dd_er_23.pdf}
    \caption{{Time series plots of  $\{\bar D_t\}_{t=2}^{500}$, $\{\bar D_{1,t}\}_{t=2}^{500}$ and $\{\bar D_{0,t}\}_{t=2}^{500}$ (left to right) for degree heterogeneity model when $\xi_i = 0.9$ and $\eta_i \in (0.1, 0.3, 0.5).$ The 25 colored curves in each subplot correspond to 25 settings where $a_0=a_1\in\{0.2,0.4,0.6,0.8,1\}$ and $b_0=b_1\in\{0.2,0.4,0.6,0.8,1\}$.}}
\end{figure}

\begin{figure}[t]
	\centering
    \includegraphics[width=\textwidth]{plot/check/dd/check_setting_dd_24.pdf}
    \\[\smallskipamount]
    \includegraphics[width=\textwidth]{plot/check/dd/check_setting_dd_25.pdf}
    \caption{{Time series plots of $\{D_t\}_{t=2}^{500}$, $\{D_{1,t}\}_{t=2}^{500}$, $\{D_{0,t}\}_{t=2}^{500}$ (left to right) for degree heterogeneity model when $\xi_i = 0.9$ and $\eta_i \in (0.7, 0.9).$ The 25 colored curves in each subplot correspond to 25 settings where $a_0=a_1\in\{0.2,0.4,0.6,0.8,1\}$ and $b_0=b_1\in\{0.2,0.4,0.6,0.8,1\}$.}}
\end{figure}

\begin{figure}[t]
	\centering
    \includegraphics[width=\textwidth]{plot/check/dd/check_setting_dd_er_24.pdf}
    \\[\smallskipamount]
    \includegraphics[width=\textwidth]{plot/check/dd/check_setting_dd_er_25.pdf}
    \caption{\label{checkdd_end}{Time series plots of  $\{\bar D_t\}_{t=2}^{500}$, $\{\bar D_{1,t}\}_{t=2}^{500}$ and $\{\bar D_{0,t}\}_{t=2}^{500}$ (left to right) for degree heterogeneity model when $\xi_i = 0.9$ and $\eta_i \in (0.7, 0.9).$ The 25 colored curves in each subplot correspond to 25 settings where $a_0=a_1\in\{0.2,0.4,0.6,0.8,1\}$ and $b_0=b_1\in\{0.2,0.4,0.6,0.8,1\}$.}}
\end{figure}

% transitivity model

\begin{figure}[t]
	\centering
\includegraphics[width=\textwidth]{plot/check/ts/check_setting_ts_01.pdf}
    \\[\smallskipamount]
   \includegraphics[width=\textwidth]{plot/check/ts/check_setting_ts_02.pdf}
    \\[\smallskipamount]
    \includegraphics[width=\textwidth]{plot/check/ts/check_setting_ts_03.pdf}
    \caption{\label{checkts1}{Time series plots of $\{D_t\}_{t=2}^{500}$, $\{D_{1,t}\}_{t=2}^{500}$, $\{D_{0,t}\}_{t=2}^{500}$ (left to right)  for transitivity model  when $\xi_i = 0.1$ and $\eta_i \in (0.1, 0.3, 0.5).$ The 25 colored curves in each subplot correspond to 25 settings where $a, b \in \{5, 15, 25, 35, 45\}$.}}
\end{figure}

\begin{figure}[t]
	\centering
\includegraphics[width=\textwidth]{plot/check/ts/check_setting_ts_er_01.pdf}
    \\[\smallskipamount]
   \includegraphics[width=\textwidth]{plot/check/ts/check_setting_ts_er_02.pdf}
    \\[\smallskipamount]
    \includegraphics[width=\textwidth]{plot/check/ts/check_setting_ts_er_03.pdf}
    \caption{\label{checkts_er1}{Time series plots of  $\{\bar D_t\}_{t=2}^{500}$, $\{\bar D_{1,t}\}_{t=2}^{500}$ and $\{\bar D_{0,t}\}_{t=2}^{500}$ (left to right)  for transitivity model  when $\xi_i = 0.1$ and $\eta_i \in (0.1, 0.3, 0.5).$ The 25 colored curves in each subplot correspond to 25 settings where $a, b \in \{5, 15, 25, 35, 45\}$.}}
\end{figure}

\begin{figure}[t]
	\centering
    \includegraphics[width=\textwidth]{plot/check/ts/check_setting_ts_04.pdf}
    \\[\smallskipamount]
    \includegraphics[width=\textwidth]{plot/check/ts/check_setting_ts_05.pdf}
    \caption{{Time series plots of $\{D_t\}_{t=2}^{500}$, $\{D_{1,t}\}_{t=2}^{500}$, $\{D_{0,t}\}_{t=2}^{500}$ (left to right) for transitivity model when $\xi_i = 0.1$ and $\eta_i \in (0.7, 0.9).$ The 25 colored curves in each subplot correspond to 25 settings where $a, b \in \{5, 15, 25, 35, 45\}$. }}
\end{figure}

\begin{figure}[t]
	\centering
    \includegraphics[width=\textwidth]{plot/check/ts/check_setting_ts_er_04.pdf}
    \\[\smallskipamount]
    \includegraphics[width=\textwidth]{plot/check/ts/check_setting_ts_er_05.pdf}
    \caption{\label{checkts_er2}{Time series plots of  $\{\bar D_t\}_{t=2}^{500}$, $\{\bar D_{1,t}\}_{t=2}^{500}$ and $\{\bar D_{0,t}\}_{t=2}^{500}$ (left to right) for transitivity model when $\xi_i = 0.1$ and $\eta_i \in (0.7, 0.9).$ The 25 colored curves in each subplot correspond to 25 settings where $a, b \in \{5, 15, 25, 35, 45\}$. }}
\end{figure}

\begin{figure}[t]
	\centering
\includegraphics[width=\textwidth]{plot/check/ts/check_setting_ts_06.pdf}
    \\[\smallskipamount]
   \includegraphics[width=\textwidth]{plot/check/ts/check_setting_ts_07.pdf}
    \\[\smallskipamount]
    \includegraphics[width=\textwidth]{plot/check/ts/check_setting_ts_08.pdf}
    \caption{{Time series plots of $\{D_t\}_{t=2}^{500}$, $\{D_{1,t}\}_{t=2}^{500}$, $\{D_{0,t}\}_{t=2}^{500}$ (left to right) for transitivity model when $\xi_i = 0.3$ and $\eta_i \in (0.1, 0.3, 0.5).$ The 25 colored curves in each subplot correspond to 25 settings where $a, b \in \{5, 15, 25, 35, 45\}$.}}
\end{figure}

\begin{figure}[t]
	\centering
\includegraphics[width=\textwidth]{plot/check/ts/check_setting_ts_er_06.pdf}
    \\[\smallskipamount]
   \includegraphics[width=\textwidth]{plot/check/ts/check_setting_ts_er_07.pdf}
    \\[\smallskipamount]
    \includegraphics[width=\textwidth]{plot/check/ts/check_setting_ts_er_08.pdf}
    \caption{{Time series plots of  $\{\bar D_t\}_{t=2}^{500}$, $\{\bar D_{1,t}\}_{t=2}^{500}$ and $\{\bar D_{0,t}\}_{t=2}^{500}$ (left to right) for transitivity model when $\xi_i = 0.3$ and $\eta_i \in (0.1, 0.3, 0.5).$ The 25 colored curves in each subplot correspond to 25 settings where $a, b \in \{5, 15, 25, 35, 45\}$.}}
\end{figure}

\begin{figure}[t]
	\centering
    \includegraphics[width=\textwidth]{plot/check/ts/check_setting_ts_09.pdf}
    \\[\smallskipamount]
    \includegraphics[width=\textwidth]{plot/check/ts/check_setting_ts_10.pdf}
    \caption{{Time series plots of $\{D_t\}_{t=2}^{500}$, $\{D_{1,t}\}_{t=2}^{500}$, $\{D_{0,t}\}_{t=2}^{500}$ (left to right) for transitivity model when $\xi_i = 0.3$ and $\eta_i \in (0.7, 0.9).$ The 25 colored curves in each subplot correspond to 25 settings where $a, b \in \{5, 15, 25, 35, 45\}$.}}
\end{figure}

\begin{figure}[t]
	\centering
    \includegraphics[width=\textwidth]{plot/check/ts/check_setting_ts_er_09.pdf}
    \\[\smallskipamount]
    \includegraphics[width=\textwidth]{plot/check/ts/check_setting_ts_er_10.pdf}
    \caption{{Time series plots of  $\{\bar D_t\}_{t=2}^{500}$, $\{\bar D_{1,t}\}_{t=2}^{500}$ and $\{\bar D_{0,t}\}_{t=2}^{500}$ (left to right) for transitivity model when $\xi_i = 0.3$ and $\eta_i \in (0.7, 0.9).$ The 25 colored curves in each subplot correspond to 25 settings where $a, b \in \{5, 15, 25, 35, 45\}$.}}
\end{figure}

\begin{figure}[t]
	\centering
\includegraphics[width=\textwidth]{plot/check/ts/check_setting_ts_11.pdf}
    \\[\smallskipamount]
   \includegraphics[width=\textwidth]{plot/check/ts/check_setting_ts_12.pdf}
    \\[\smallskipamount]
    \includegraphics[width=\textwidth]{plot/check/ts/check_setting_ts_13.pdf}
    \caption{{Time series plots of $\{D_t\}_{t=2}^{500}$, $\{D_{1,t}\}_{t=2}^{500}$, $\{D_{0,t}\}_{t=2}^{500}$ (left to right) for transitivity model when $\xi_i = 0.5$ and $\eta_i \in (0.1, 0.3, 0.5).$ The 25 colored curves in each subplot correspond to 25 settings where $a, b \in \{5, 15, 25, 35, 45\}$.}}
\end{figure}

\begin{figure}[t]
	\centering
\includegraphics[width=\textwidth]{plot/check/ts/check_setting_ts_er_11.pdf}
    \\[\smallskipamount]
   \includegraphics[width=\textwidth]{plot/check/ts/check_setting_ts_er_12.pdf}
    \\[\smallskipamount]
    \includegraphics[width=\textwidth]{plot/check/ts/check_setting_ts_er_13.pdf}
    \caption{{Time series plots of  $\{\bar D_t\}_{t=2}^{500}$, $\{\bar D_{1,t}\}_{t=2}^{500}$ and $\{\bar D_{0,t}\}_{t=2}^{500}$ (left to right) for transitivity model when $\xi_i = 0.5$ and $\eta_i \in (0.1, 0.3, 0.5).$ The 25 colored curves in each subplot correspond to 25 settings where $a, b \in \{5, 15, 25, 35, 45\}$.}}
\end{figure}

\begin{figure}[t]
	\centering
    \includegraphics[width=\textwidth]{plot/check/ts/check_setting_ts_14.pdf}
    \\[\smallskipamount]
    \includegraphics[width=\textwidth]{plot/check/ts/check_setting_ts_15.pdf}
    \caption{{Time series plots of $\{D_t\}_{t=2}^{500}$, $\{D_{1,t}\}_{t=2}^{500}$, $\{D_{0,t}\}_{t=2}^{500}$ (left to right) for transitivity model when $\xi_i = 0.5$ and $\eta_i \in (0.7, 0.9).$ The 25 colored curves in each subplot correspond to 25 settings where $a, b \in \{5, 15, 25, 35, 45\}$.}}
\end{figure}

\begin{figure}[t]
	\centering
    \includegraphics[width=\textwidth]{plot/check/ts/check_setting_ts_er_14.pdf}
    \\[\smallskipamount]
    \includegraphics[width=\textwidth]{plot/check/ts/check_setting_ts_er_15.pdf}
    \caption{{Time series plots of  $\{\bar D_t\}_{t=2}^{500}$, $\{\bar D_{1,t}\}_{t=2}^{500}$ and $\{\bar D_{0,t}\}_{t=2}^{500}$ (left to right) for transitivity model when $\xi_i = 0.5$ and $\eta_i \in (0.7, 0.9).$ The 25 colored curves in each subplot correspond to 25 settings where $a, b \in \{5, 15, 25, 35, 45\}$.}}
\end{figure}

\begin{figure}[t]
	\centering
\includegraphics[width=\textwidth]{plot/check/ts/check_setting_ts_16.pdf}
    \\[\smallskipamount]
   \includegraphics[width=\textwidth]{plot/check/ts/check_setting_ts_17.pdf}
    \\[\smallskipamount]
    \includegraphics[width=\textwidth]{plot/check/ts/check_setting_ts_18.pdf}
    \caption{{Time series plots of $\{D_t\}_{t=2}^{500}$, $\{D_{1,t}\}_{t=2}^{500}$, $\{D_{0,t}\}_{t=2}^{500}$ (left to right) for transitivity model when $\xi_i = 0.7$ and $\eta_i \in (0.1, 0.3, 0.5).$ The 25 colored curves in each subplot correspond to 25 settings where $a, b \in \{5, 15, 25, 35, 45\}$.}}
\end{figure}

\begin{figure}[t]
	\centering
\includegraphics[width=\textwidth]{plot/check/ts/check_setting_ts_er_16.pdf}
    \\[\smallskipamount]
   \includegraphics[width=\textwidth]{plot/check/ts/check_setting_ts_er_17.pdf}
    \\[\smallskipamount]
    \includegraphics[width=\textwidth]{plot/check/ts/check_setting_ts_er_18.pdf}
    \caption{{Time series plots of  $\{\bar D_t\}_{t=2}^{500}$, $\{\bar D_{1,t}\}_{t=2}^{500}$ and $\{\bar D_{0,t}\}_{t=2}^{500}$ (left to right) for transitivity model when $\xi_i = 0.7$ and $\eta_i \in (0.1, 0.3, 0.5).$ The 25 colored curves in each subplot correspond to 25 settings where $a, b \in \{5, 15, 25, 35, 45\}$.}}
\end{figure}

\begin{figure}[t]
	\centering
    \includegraphics[width=\textwidth]{plot/check/ts/check_setting_ts_19.pdf}
    \\[\smallskipamount]
    \includegraphics[width=\textwidth]{plot/check/ts/check_setting_ts_20.pdf}
    \caption{{Time series plots of $\{D_t\}_{t=2}^{500}$, $\{D_{1,t}\}_{t=2}^{500}$, $\{D_{0,t}\}_{t=2}^{500}$ (left to right) for transitivity model  when $\xi_i = 0.7$ and $\eta_i \in (0.7, 0.9).$ The 25 colored curves in each subplot correspond to 25 settings where $a, b \in \{5, 15, 25, 35, 45\}$.}}
\end{figure}

\begin{figure}[t]
	\centering
    \includegraphics[width=\textwidth]{plot/check/ts/check_setting_ts_er_19.pdf}
    \\[\smallskipamount]
    \includegraphics[width=\textwidth]{plot/check/ts/check_setting_ts_er_20.pdf}
    \caption{{Time series plots of  $\{\bar D_t\}_{t=2}^{500}$, $\{\bar D_{1,t}\}_{t=2}^{500}$ and $\{\bar D_{0,t}\}_{t=2}^{500}$ (left to right) for transitivity model  when $\xi_i = 0.7$ and $\eta_i \in (0.7, 0.9).$ The 25 colored curves in each subplot correspond to 25 settings where $a, b \in \{5, 15, 25, 35, 45\}$.}}
\end{figure}

\begin{figure}[t]
	\centering
\includegraphics[width=\textwidth]{plot/check/ts/check_setting_ts_21.pdf}
    \\[\smallskipamount]
   \includegraphics[width=\textwidth]{plot/check/ts/check_setting_ts_22.pdf}
    \\[\smallskipamount]
    \includegraphics[width=\textwidth]{plot/check/ts/check_setting_ts_23.pdf}
    \caption{{Time series plots of $\{D_t\}_{t=2}^{500}$, $\{D_{1,t}\}_{t=2}^{500}$, $\{D_{0,t}\}_{t=2}^{500}$ (left to right) for transitivity model when $\xi_i = 0.9$ and $\eta_i \in (0.1, 0.3, 0.5).$ The 25 colored curves in each subplot correspond to 25 settings where $a, b \in \{5, 15, 25, 35, 45\}$.}}
\end{figure}

\begin{figure}[t]
	\centering
\includegraphics[width=\textwidth]{plot/check/ts/check_setting_ts_er_21.pdf}
    \\[\smallskipamount]
   \includegraphics[width=\textwidth]{plot/check/ts/check_setting_ts_er_22.pdf}
    \\[\smallskipamount]
    \includegraphics[width=\textwidth]{plot/check/ts/check_setting_ts_er_23.pdf}
    \caption{{Time series plots of  $\{\bar D_t\}_{t=2}^{500}$, $\{\bar D_{1,t}\}_{t=2}^{500}$ and $\{\bar D_{0,t}\}_{t=2}^{500}$ (left to right) for transitivity model when $\xi_i = 0.9$ and $\eta_i \in (0.1, 0.3, 0.5).$ The 25 colored curves in each subplot correspond to 25 settings where $a, b \in \{5, 15, 25, 35, 45\}$.}}
\end{figure}

\begin{figure}[t]
	\centering
    \includegraphics[width=\textwidth]{plot/check/ts/check_setting_ts_24.pdf}
    \\[\smallskipamount]
    \includegraphics[width=\textwidth]{plot/check/ts/check_setting_ts_25.pdf}
    \caption{{Time series plots of $\{D_t\}_{t=2}^{500}$, $\{D_{1,t}\}_{t=2}^{500}$, $\{D_{0,t}\}_{t=2}^{500}$ (left to right) for transitivity model when $\xi_i = 0.9$ and $\eta_i \in (0.7, 0.9).$ The 25 colored curves in each subplot correspond to 25 settings where $a, b \in \{5, 15, 25, 35, 45\}$.}}
\end{figure}

\begin{figure}[t]
	\centering
    \includegraphics[width=\textwidth]{plot/check/ts/check_setting_ts_er_24.pdf}
    \\[\smallskipamount]
    \includegraphics[width=\textwidth]{plot/check/ts/check_setting_ts_er_25.pdf}
    \caption{\label{checkts_end}{Time series plots of  $\{\bar D_t\}_{t=2}^{500}$, $\{\bar D_{1,t}\}_{t=2}^{500}$ and $\{\bar D_{0,t}\}_{t=2}^{500}$ (left to right) for transitivity model when $\xi_i = 0.9$ and $\eta_i \in (0.7, 0.9).$ The 25 colored curves in each subplot correspond to 25 settings where $a, b \in \{5, 15, 25, 35, 45\}$.}}
\end{figure}
\end{skip}

\subsection{Stationarity and ergodicity} \label{sec:simul1}

As stated in the last  paragraph of Section \ref{sec:models}, for each fixed constant $p$, $\{\bX_t\}_{t\geq1}$ defined by the transitivity model (\ref{trans-model}) is stationary and ergodic as long as $\alpha_{i,j}^{t-1}$ and $\beta_{i,j}^{t-1}$ are strictly between 0 and 1. Nevertheless, it is a Markov chain with $2^{p(p-1)/2}$ states. When $p$ is a fixed constant, the ergodicity of $\bX_t$ 
(i.e., the average in time converges to the average over the state space) may take a long time to be observed; see Theorem 4.1 in Chapter 3 of \cite{Bremaud1998_supp}. 
However, the ergodicity of some scalar summary statistics of $\bX_t$ can be observed in much short time spans, as indicated in the simulation reported below.

We  consider
the following three network statistics measures at each time $t$:
\begin{equation}
    \label{eq:3density}
D_t= {\sum_{i,j:\,i < j} X_{i,j}^t \over p(p-1)/2}\,, ~~
D_{1,t} = \frac{\sum_{i,j:\,i < j} (1-X_{i,j}^{t-1}) X_{i,j}^t}{p(p-1)/2}\,,
~~ 
D_{0,t} = \frac{\sum_{i,j:\,i< j} X_{i,j}^{t-1} (1-X_{i,j}^t)}{p(p-1)/2}\,,
\end{equation}
where $D_t$ is the edge density at time $t$, and $D_{1,t}$ and $D_{0,t}$ are, 
respectively, the (normalized) number of newly formed edges and newly dissolved edges at time $t$. If $\{\bX_t\}_{t\geq1}$ is stationary, all three density sequences $\{D_t\}_{t\geq1}$, $\{D_{1,t}\}_{t\geq2}$ and $\{D_{0,t}\}_{t\geq2}$ are also  stationary. We also plot 
\[
\bar D_t = {1 \over t} \sum_{u=1}^t D_u\,, \qquad 
\bar{D}_{1,t} = {1 \over t} \sum_{u=1}^t D_{1,u}\,, \qquad
\bar{D}_{0,t} = {1 \over t} \sum_{u=1}^t D_{0,u}\,,
\]
against $t$, for $t\geq2$, to see how quickly the ergodicity can be observed. These are sample means of one-dimensional network summaries. We expect that their convergences are much faster than that for the sample mean of $p\times p$ network $\bX_t$ itself.

\begin{figure}[t]
\centering
\begin{subfigure}{.3\linewidth}
    \centering
    \centering
    \includegraphics[width=5cm]{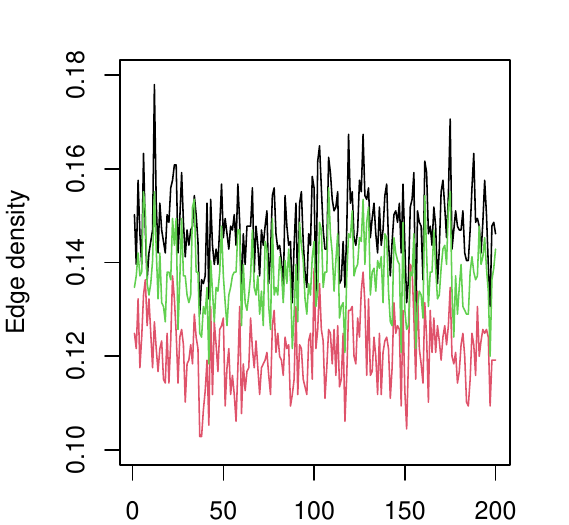}
    \caption{$D_t$}
\end{subfigure}
\begin{subfigure}{.3\linewidth}
    \centering
    \includegraphics[width=5cm]{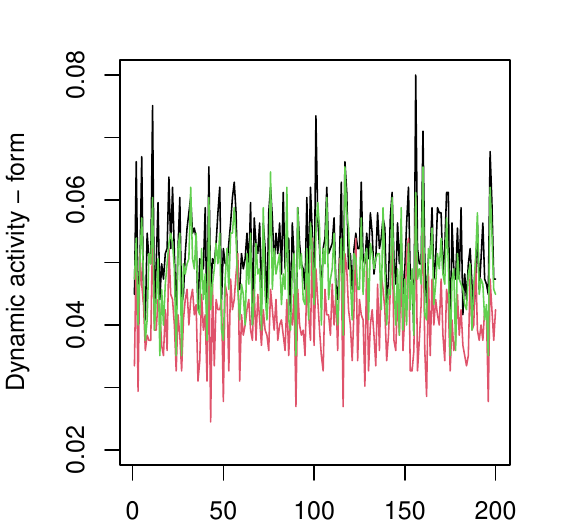}
    \caption{$D_{1,t}$}
\end{subfigure}
\begin{subfigure}{.3\linewidth}
    \centering
    \includegraphics[width=5cm]{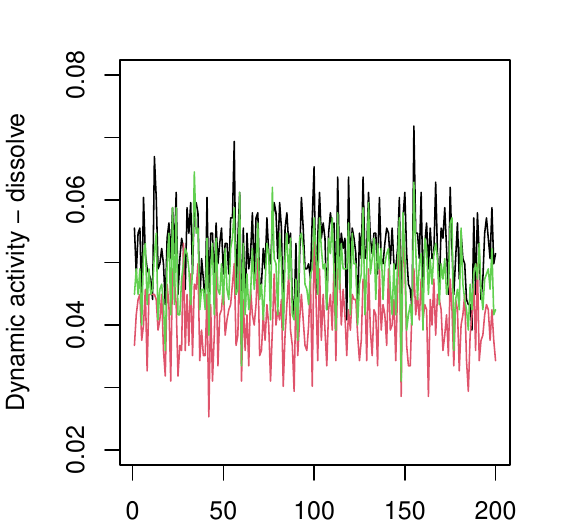}
    \caption{$D_{0,t}$}
\end{subfigure}

\begin{subfigure}{.3\linewidth}
    \centering
    \centering
    \includegraphics[width=5cm]{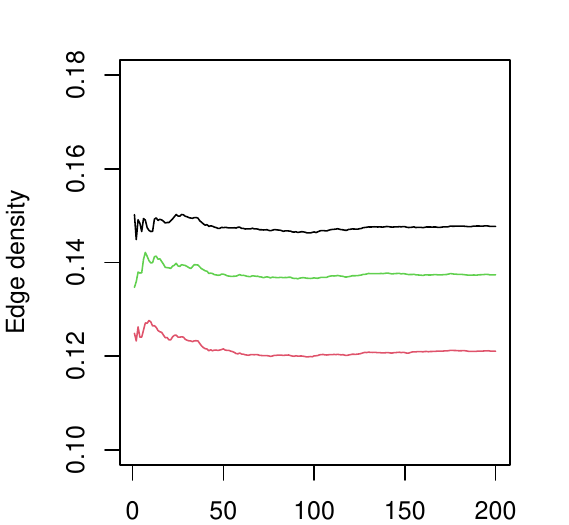}
    \caption{$\bar D_t$}
\end{subfigure}
\begin{subfigure}{.3\linewidth}
    \centering
    \includegraphics[width=5cm]{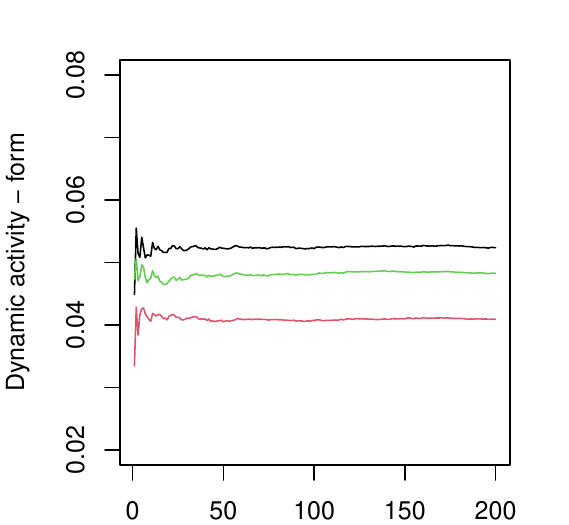}
    \caption{$\bar D_{1,t}$}
\end{subfigure}
\begin{subfigure}{.3\linewidth}
    \centering
    \includegraphics[width=5cm]{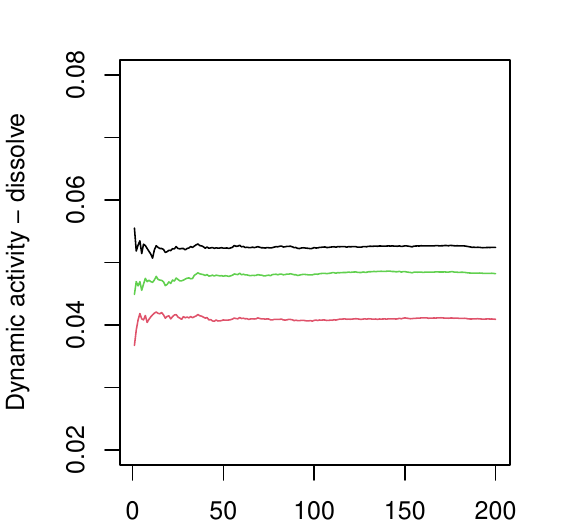}
    \caption{$\bar D_{0,t}$}
\end{subfigure}

\centering

\caption{\label{trans_plot_50}{Time series plots of $\{D_t\}_{t=2}^{200}$, $\{D_{1,t}\}_{t=2}^{200}$, $\{D_{0,t}\}_{t=2}^{200}$, $\{\bar D_t\}_{t=2}^{200}$, $\{\bar D_{1,t}\}_{t=2}^{200}$ and $\{\bar D_{0,t}\}_{t=2}^{200}$ for the three simulated settings with $p = 50$. The black, red and green curves correspond to the settings $(\xi_i,  \eta_i, a, b) = ( 0.8,0.9,10,10),$ $(0.8,0.9,25,15)$ and $(0.7,0.9,15,10)$, respectively. }}
\vspace{-0.2cm}
\end{figure}

Setting $\xi_1=\cdots=\xi_p$ and $\eta_1=\cdots=\eta_p$, we let
$(\xi_i, \eta_i, a, b)$ take three different sets of values: 
$(0.8,0.9,10,10),$ $(0.8,0.9,25,15)$ and $(0.7,0.9,15,10)$.
%To evaluate the effectiveness of the estimation procedure, we focus our analysis on settings with 
Figure~\ref{trans_plot_50} displays the time series plots of simulated  $\{D_t\}_{t=2}^{200}$, $\{D_{1,t}\}_{t=2}^{200}$, $\{D_{0,t}\}_{t=2}^{200}$, $\{\bar D_t\}_{t=2}^{200}$, $\{\bar D_{1,t}\}_{t=2}^{200}$ and $\{\bar D_{0,t}\}_{t=2}^{200}$ when $p=50$. %The same plots with $p = 100$ and $150$ are, respectively, shown in Figures~\ref{trans_plot_100} and \ref{trans_plot_150} in Section \ref{app:transitivity} of the supplementary material. 
As expected, all simulated series $\{D_t\}_{t\geq2}$, $\{D_{1,t}\}_{t\geq2}$ and $\{D_{0,t}\}_{t\geq2}$ exhibit patterns in line with stationarity. The convergence of their sample means is observed with the sample sizes greater than 50. %As presented in Figures \ref{trans_plot_50}, %\ref{trans_plot_100} and \ref{trans_plot_150} in Section \ref{app:transitivity} of the supplementary material,
%$D_t \in (0.02,0.14)$, $D_{1,t} \in (0.015, 0.07)$, and $ D_{0,t} \in (0.02,0.07)$ in the current settings with $p\in\{50, 100, 150\}$. We also notice the similar ranges of $D_t, D_{1,t}$ and $D_{0,t}$
%for  the manufacturing email network data  analysed later in Section~\ref{sec:real}. Hence, the proposed transitivity model in \eqref{trans-model} can be used to capture the sparsity and the temporal dynamics  observed in this real data.
% In particular, $(\xi_i, \eta_i, a, b) = (0.8,0.9,10,10)$ displays the most dynamic edge-changing behavior, while $(\xi_i, \eta_i, a, b) =(0.8,0.9,25,15)$ is the least dynamic among the three settings.  

%Additional simulation results, on the degree heterogeneity model (\ref{densityalphaij}), are reported in Section \ref{app:density-dependent} of the supplementary material.

\subsection{The transitivity model with $p = 500$}
 We further examine larger networks with $p = 500$ and  $n =100$ under the
transitivity model, and
report in Table~\ref{table:1} the means and the standard errors of rMAEs over the 100 replications. The initial values of $\{\xi_i\}_{i=1}^p$ and $\{\eta_i\}_{i=1}^p$ are sampled independently from $\text{Uniform}[0.5, 0.9]$. It can be seen that the proposed estimation procedure remains effective when $p = 500$, with the improved estimator $\hat{\btheta}$ significantly
outperforming the initial estimator $\tilde{\btheta}$.

\begin{minipage}{1\linewidth}
\begin{table}[H]
 %Means (standard errors) of 
	\caption{The means and  STDs (in parenthesis) of rMAEs for estimating parameters in the transitivity model (\ref{trans-model}) with $n = 100$  and $p = 500$ over 100 replications. \label{table:1}}
    \vspace{-0.5cm}
	\begin{center}
		\resizebox{5in}{!}{
\begin{tabular}{cccccccc}
\hline
$(\xi_i,   \eta_i, a, b)$ &  & Estimation &  & $\xi_i$       & $\eta_i$      & $a$           & $b$           \\ \hline
(0.8, 0.9, 10, 10)        &  & Initial    &  & 0.105 (0.049) & 0.048 (0.001) & 3.506 (0.077) & 0.187 (0.006) \\
                          &  & Improved   &  & 0.069 (0.037) & 0.035 (0.012) & 0.678 (0.138) & 0.126 (0.027) \\
(0.8, 0.9, 25, 15)        &  & Initial    &  & 0.112 (0.051) & 0.047 (0.001) & 1.330 (0.034)  & 0.097 (0.007) \\
                          &  & Improved   &  & 0.082 (0.050)  & 0.036 (0.009) & 0.234 (0.069) & 0.094 (0.016) \\
(0.7, 0.9, 15, 10)        &  & Initial    &  & 0.138 (0.056) & 0.047 (0.001) & 1.890 (0.055)  & 0.285 (0.007) \\
                          &  & Improved   &  & 0.113 (0.027) & 0.014 (0.004) & 0.685 (0.133) & 0.032 (0.026)\\ \hline
\end{tabular}
}
\end{center}
\end{table}
  \end{minipage}

\subsection{A more general model}
\label{sec:more.gen.model}

As stated towards the end of Section \ref{sec:transitivity},
 a more general transitivity model admits the form:
\begin{align} \label{trans-model-ab}
    \alpha_{i,j}^{t-1}(\btheta) = \frac{\xi_i \xi_je^{a_1 U_{i,j}^{t-1}} }{1+ e^{a_1 U_{i,j}^{t-1}}+ e^{b_1 V_{i,j}^{t-1}}}\,,
    \qquad \beta_{i,j}^{t-1}(\btheta) =\frac{ \eta_i \eta_j e^{b_2 V_{i,j}^{t-1}}  }{1+ e^{a_2 U_{i,j}^{t-1}}+ e^{b_2 V_{i,j}^{t-1}}}\,,
\end{align}
with $\btheta=(a_1, b_1, a_2, b_2,\xi_1,\ldots,\xi_p, \eta_1,\ldots,\eta_p)^{\T}\in\mathbb{R}_+^{2p+4}$. Different from the transitivity model \eqref{trans-model} introduced in Section \ref{sec:transitivity}, we allow $a_1\neq a_2$ and $b_1\neq b_2$ in \eqref{trans-model-ab}. We adopt the same simulation settings as above, i.e., 
$(\xi_i,  a_1, b_1, \eta_i, a_2, b_2) \in \{(0.8,10,10, 0.9,10,10),$ $(0.8,25,15, 0.9,25,15),$ $(0.7,15,10,0.9,15,10)\}$.

%We implement the two-step estimation in a similar way as that in Section~\ref{sec:sim.est}. Using 9 different sets of the initial values, set  $(\tau, \tilde r, \check r) = (0.5 \Delta_n^{1/2},0.5,0.1)$ for local parameters $\{\xi_i,\eta_i\}$ and $\tau =  0.01 \Delta_n^{1/2}$ for global parameters $\{a_1,b_1\}$. As the ideal choice of $(\tilde r, \check r)$ shall cover the true value of parameters, we keep $(\tilde r, \check r) =  (10, 2)$ for $\{a_1,b_1\}$ and employee $(\tilde r, \check r) =  (100, 20)$ for $\{a_2,b_2\}$. 

We implement the estimation procedure in the same manner as in Section~\ref{sec:sim.est}.
Table~\ref{table.dynamic.err2} reports the resulting rMAEs over $100$ replications with $p = 50$ and $n \in \{100,200\}$. 
The estimation for the parameters in  $\alpha_{i,j}^{t-1}(\btheta)$, namely $\xi_i, a_1$ and $b_1$, exhibits the similar patterns as in Table \ref{table.dynamic.err1}. In contrast, the estimation for parameters in $\beta_{i,j}^{t-1}(\btheta)$ deteriorates significantly,  and especially for $a_2$ and $b_2$.
%As we mentioned in the last paragraph of Section \ref{sec:sim.est}, 
Note that only some components of $\bX_{t-1}$ with $X_{i,j}^{t-1}=1$, $t\in[n]\setminus[m]$ and $j\neq i$,  were used in estimating parameters $a_2$ and $b_2$. For sparse networks, the total number of those data points is small. This is the intrinsic difficulty in estimating the parameters in $\beta_{i,j}^{t-1}(\btheta)$.
See also the relevant discussion at the end of Section \ref{sec:transitivity}.

% Similar to the model degeneracy of static network analysis (the random graph model places almost all its mass on the empty or the complete graphs)\citep{snijders2002,handcock2003}, 

\begin{table}[tbp]
	\caption{\label{table.dynamic.err2} %Means (standard errors) of 
 The means and  STDs (in parenthesis) of rMAEs for estimating parameters in transitivity model (\ref{trans-model-ab}) with 100 replications.}% The lowest values are in bold font.}
 \vspace{-0.5cm}
	\begin{center}
		\resizebox{6.6in}{!}{
\begin{tabular}{ccccccccccccc}
\hline
\multirow{2}{*}{$(\xi_i,   a_1, b_1,  \eta_i, a_2, b_2)$} &  & \multirow{2}{*}{$n$} &  & \multirow{2}{*}{Estimation} &  & \multicolumn{3}{c}{$\alpha_{i,j}^{t-1}(\btheta)$} &  & \multicolumn{3}{c}{$\beta_{i,j}^{t-1}(\btheta)$}  \\ \cline{7-9} \cline{11-13} 
                                                          &  &                      &  &                             &  & $\xi_i$         & $a_1$          & $b_1$          &  & $\eta_i$      & $a_2$           & $b_2$           \\ \hline
(0.8, 10, 10, 0.9, 10, 10)                                &  & 100                  &  & Initial                     &  & 0.183 (0.001)   & 0.611 (0.042)  & 0.285 (0.005)  &  & 0.131 (0.021) & 11.953 (21.723) & 52.369 (26.797) \\
                                                          &  &                      &  & Improved                    &  & 0.059 (0.006)   & 0.069 (0.052)  & 0.020 (0.015)   &  & 0.103 (0.015) & 11.843 (22.698) & 52.629 (29.730)  \\
                                                          &  & 200                  &  & Initial                     &  & 0.182 (0.001)   & 0.609 (0.026)  & 0.285 (0.003)  &  & 0.128 (0.022) & 10.082 (18.822) & 53.861 (29.054) \\
                                                          &  &                      &  & Improved                    &  & 0.045 (0.005)   & 0.043 (0.032)  & 0.013 (0.010)   &  & 0.104 (0.032) & 9.900 (19.676)    & 54.005 (32.168) \\ \hline
(0.8, 25, 15, 0.9, 25, 15)                                &  & 100                  &  & Initial                     &  & 0.173 (0.001)   & 0.331 (0.015)  & 0.216 (0.002)  &  & 0.117 (0.021) & 6.265 (10.585)  & 34.119 (23.249) \\ 
                                                          &  &                      &  & Improved                    &  & 0.062 (0.007)   & 0.025 (0.021)  & 0.015 (0.010)   &  & 0.060 (0.020)   & 6.238 (11.305)  & 33.400 (25.999)   \\
                                                          &  & 200                  &  & Initial                     &  & 0.172 (0.001)   & 0.333 (0.013)  & 0.216 (0.002)  &  & 0.113 (0.022) & 4.389 (7.712)   & 32.834 (25.108) \\
                                                          &  &                      &  & Improved                    &  & 0.048 (0.005)   & 0.018 (0.014)  & 0.081 (0.009)  &  & 0.065 (0.033) & 4.031 (8.166)   & 31.552 (28.022) \\ \hline
(0.7, 15, 10, 0.9, 15, 10)                                &  & 100                  &  & Initial                     &  & 0.177 (0.002)   & 0.379 (0.012)  & 0.248 (0.004)  &  & 0.142 (0.025) & 12.506 (18.162) & 54.621 (29.036) \\ 
                                                          &  &                      &  & Improved                    &  & 0.060 (0.006)    & 0.051 (0.036)  & 0.022 (0.017)  &  & 0.112 (0.014) & 12.747 (18.630)  & 55.040 (31.396)  \\
                                                          &  & 200                  &  & Initial                     &  & 0.172 (0.001)   & 0.375 (0.010)   & 0.247 (0.002)  &  & 0.135 (0.029) & 2.748 (5.513)   & 49.116 (30.803) \\
\multicolumn{1}{l}{}                                      &  & \multicolumn{1}{l}{} &  & Improved                    &  & 0.046 (0.004)   & 0.035 (0.025)  & 0.075 (0.005)  &  & 0.112 (0.041) & 2.415 (5.850)    & 48.528 (34.177) \\ \hline
\end{tabular}

		}	
	\end{center}
	\vspace{-0.5cm}
\end{table}

{

 \section{Method of moments estimation} \label{sec:estimation2}
% The maximum-likelihood-based estimation procedure described in Section~\ref{sec:impl.details} is general and does not require any specific functional forms for the transition probabilities $\alpha_{i,j}^{t-1}(\btheta)$ and $\beta_{i,j}^{t-1}(\btheta)$ in \eqref{eq:alphaijbetaij:theta}. 

 In this section, we leverage the
 structural similarities across the three AR network models introduced in Section~\ref{sec:3examples} and propose a computationally efficient estimation procedure, termed Iterative Method of Moments (IMoM). 

 \subsection{IMoM algorithm}
 Consider a simplified formulation of the transition probabilities, given by 
\begin{equation}\label{eq:mom:alphaijbetaij}
\begin{split}
 \alpha_{i,j}^{t-1}(\btheta) = \xi_i \xi_j \tilde f_{i,j}(\bX_{t-1}\setminus X_{i,j}^{t-1}, \bX_{t-2}, \ldots, \bX_{t-m}; \btheta_{\mathcal{G}}) \,,\\
    \beta_{i,j}^{t-1}(\btheta) =  \eta_i \eta_j  \tilde g_{i,j}(\bX_{t-1}\setminus X_{i,j}^{t-1}, \bX_{t-2}, \ldots, \bX_{t-m}; \btheta_{\mathcal{G}})\,,
\end{split}
\end{equation}
where $\{\xi_i\}_{i=1}^p$ and $\{\eta_i\}_{i=1}^p$ are local parameters capturing node-specific heterogeneity in edge formation and dissolution, respectively,  and $\tilde f_{i,j}$'s and $\tilde g_{i,j}$'s are known functions that depend only on the global parameters $\btheta_{\mathcal{G}}$. This separable structure between the local and global components allows us to develop the IMoM estimation procedure as follows.

For illustrative purposes, we focus on the AR($m$) network with $m = 1$. We start with some initial values for the local parameters $\{\xi_i\}_{i=1}^p$ and $\{\eta_i\}_{i=1}^p$, denoted as $\{\mathring \xi_i^{(0)}\}_{i=1}^p$ and $\{\mathring \eta_i^{(0)}\}_{i=1}^p$, respectively. At the $(\iota+1)$-th iteration, we first update the global parameter vector 
by solving an optimization problem similar to the computation of $\tilde{\btheta}_{\mathcal{G}}^{({\rm app})}$ in Section~\ref{sec:impl.details}.
Specifically, we write  $\mathring \btheta_{\mathcal{G}^{\btc}}^{(\iota)}=(\mathring\xi_1^{(\iota)},\ldots,\mathring\xi_p^{(\iota)},\mathring\eta_1^{(\iota)},\ldots,\mathring\eta_p^{(\iota)})^{\T}$ 
with $\mathring\xi_1^{(\iota)},\ldots,\mathring\xi_p^{(\iota)},\mathring\eta_1^{(\iota)},\ldots,\mathring\eta_p^{(\iota)}$ denoting the $\iota$-th iterate of local parameters, and obtain $\mathring{\btheta}_{\mathcal{G}}^{(\iota+1)}=\arg\max_{\btheta_{\mathcal{G}}}\hat{\ell}_{n,p}^{(l)}(\btheta_{\mathcal{G}},\mathring \btheta_{\mathcal{G}^{\btc}}^{(\iota)})$, for $\hat{\ell}_{n,p}^{(l)}(\btheta)$ defined in \eqref{logL:local} for some $l\in\mathcal{G}$. 
As discussed in Section~\ref{sec:impl.details}, 
this low-dimensional optimization step is computationally more efficient and offers a reliable approximation of the global parameters.

Let $\Xi_{i,j} = \xi_i \xi_j$ and $\Gamma_{i,j} = \eta_i \eta_j$ denote the edgewise local parameters  for $1 \leq i \neq j \leq p$. 
We next turn to updating these parameters via a method-of-moments approach.
Denote by $I(\cdot)$ the indicator function.
Inspired by the fact that 
\begin{align*} 
     \mathbb{P}(X_{i,j}^t=1\,|\,X_{i,j}^{t-1} =0) &=\mathbb{E}\big\{ \mathbb{P}(X_{i,j}^t=1\,|\,X_{i,j}^{t-1} =0,  \bX_{t-1}\setminus X_{i,j}^{t-1}) | \, X_{i,j}^{t-1} = 0 \big\}
     \\ & =  \mathbb{E}(\alpha_{i,j}^{t-1} |\, X_{i,j}^{t-1} = 0) = \frac{\mathbb{E}\big\{\alpha_{i,j}^{t-1} I( X_{i,j}^{t-1} = 0)\big\}}{\mathbb{P}( X_{i,j}^{t-1} = 0)}\,, 
\end{align*}
 % the expectation is over $\bX_{t-1} \setminus X_{i,j}^{t-1}$   given $X_{i,j}^{t-1} = 0$
we compute the $(\iota+1)$-th iterate for $\Xi_{i,j}$ by 
\begin{equation} \label{eq:mom:xi}
    \mathring \Xi_{i,j}^{(\iota+1)}  = \, \frac{ \sum_{t =2}^n X_{i,j}^t (1- X_{i,j}^{t-1})}{\sum_{t =2}^n \tilde f_{i,j}\big\{\bX_{t-1}\setminus X_{i,j}^{t-1}; \mathring \btheta_{\mathcal{G}}^{(\iota+1)}\big\} ( 1-X_{i,j}^{t-1} )} \equiv \mathring \xi_{i}^{(\iota+1)}\mathring \xi_{j}^{(\iota+1)}\,.
\end{equation}
Similarly, we update 
\begin{equation} \label{eq:mom:eta}
    \mathring \Gamma_{i,j}^{(\iota+1)}  =\, \frac{ \sum_{t =2}^n (1-X_{i,j}^t)  X_{i,j}^{t-1}}{\sum_{t =2}^n \tilde g_{i,j}\big\{\bX_{t-1}\setminus X_{i,j}^{t-1}; \mathring \btheta_{\mathcal{G}}^{(\iota+1)}\big\} {X_{i,j}^{t-1} }}\equiv \mathring \eta_{i}^{(\iota+1)}\mathring \eta_{j}^{(\iota+1)} \,.
\end{equation}

Finally, we recover the individual local parameters $\{\xi_i\}_{i=1}^p$ and $\{\eta_i\}_{i=1}^p$ by solving a strongly convex optimization problem. Set $\mathring{\Xi}_{i,i}^{(\iota+1)} = \mathring{\Gamma}_{i,i}^{(\iota+1)} = 0$ for all $i \in [p]$ and define $\boldsymbol{\Xi}^{(\iota+1)} = \{\mathring{\Xi}_{i,j}^{(\iota+1)}\}_{i,j \in [p]}$ and $\boldsymbol{\Gamma}^{(\iota+1)} = \{\mathring{\Gamma}_{i,j}^{(\iota+1)}\}_{i,j \in [p]}$. Let $\mathbf{1}_p$ denote the $p$-dimensional column vector of ones. We then compute the log-transformed local parameters $\{\log(\xi_i)\}_{i \in [p]}$ and $\{\log(\eta_i)\}_{i \in [p]}$, respectively, as
$
\hat{\bx}_{\iota+1} = \arg\min_{\bx \in \mathbb{R}^p} \mathring \ell(\bx; \boldsymbol{\Xi}^{(\iota+1)} \mathbf{1}_p)$ and $ 
\hat{\by}_{\iota+1} = \arg\min_{\by \in \mathbb{R}^p} \mathring \ell(\by; \boldsymbol{\Gamma}^{(\iota+1)} \mathbf{1}_p)$,
where the function $\mathring \ell(\cdot\,;\cdot)$ is specified in Lemma~\ref{convex} below.  The strong convexity of $\mathring \ell(\cdot; \cdot)$ ensures that these minimizers are unique and can be  computed efficiently. 
See Section~\ref{supp:imom} below for further technical details. 
Write $\hat{\bx}_{\iota+1} = (\hat{x}_{1, \iota+1}, \dots, \hat{x}_{p, \iota+1})^\T$ and $\hat{\by}_{\iota+1} = (\hat{y}_{1, \iota+1}, \dots, \hat{y}_{p, \iota+1})^\T$. The updated local parameters are obtained via $\mathring{\xi}_i^{(\iota+1)} = \exp(\hat{x}_{i,\iota+1})$ and $\mathring{\eta}_i^{(\iota+1)} = \exp(\hat{y}_{i,\iota+1})$ for all $i \in [p]$.

We iterate the above three steps until convergence to obtain the final IMoM estimator $\mathring \btheta$. The complete estimation procedure is summarized in Algorithm~\ref{alg1}.
One possible alternative is to iterate only between Steps~\ref{alg1:2a} and \ref{alg1:2b}, deferring the recovery of the individual local parameters, i.e., Step~\ref{alg1:2c}, to the final iteration.  While such an approach introduces a slight reduction in computational effort per iteration, our numerical experiments indicate that it leads to lower estimation accuracy and slower convergence in practice.
We therefore incorporate Step~\ref{alg1:2c} into each iteration, as described in Algorithm~\ref{alg1}.
 
\begin{center}
	\begin{algorithm}[ht] \caption{\label{alg1}{IMoM for AR networks with transition probabilities in \eqref{eq:mom:alphaijbetaij}.}}
		\begin{enumerate}
			\item[1.] Input: Initial values for local parameters, i.e., $\{\mathring \xi_i^{(0)}\}_{i=1}^p$ and $\{\mathring \eta_i^{(0)}\}_{i=1}^p$.
			\item[2.] For $\iota = 0, 1, \ldots$ until convergence: 	
			\begin{enumerate}[label=(\roman*)]		
				\item  \label{alg1:2a} Update the global parameter vector via $$
    \mathring{\btheta}_{\mathcal{G}}^{(\iota+1)} = \arg\max_{\btheta_{\mathcal{G}}} \hat{\ell}_{n,p}^{(l)}(\btheta_{\mathcal{G}}, \mathring{\btheta}_{\mathcal{G}^{\btc}}^{(\iota)})\,,$$
        where $\mathring{\btheta}_{\mathcal{G}^{\btc}}^{(\iota)} = (\mathring{\xi}_1^{(\iota)}, \ldots, \mathring{\xi}_p^{(\iota)}, \mathring{\eta}_1^{(\iota)}, \ldots, \mathring{\eta}_p^{(\iota)})^{\T}$ and $\hat{\ell}_{n,p}^{(l)}(\btheta)$ is defined in \eqref{logL:local} for some $l \in \mathcal{G}$.
				\item  \label{alg1:2b} Compute the edgewise local parameters 
        $\mathring \Xi_{i,j}^{(\iota+1)} $ and $\mathring \Gamma_{i,j}^{(\iota+1)} $ using method-of-moments-based approach in \eqref{eq:mom:xi} and \eqref{eq:mom:eta}, respectively.
				\item \label{alg1:2c} Recover the log-transformed local parameters by solving 
\[
\hat{\bx}_{\iota+1} = \arg\min_{\bx\in \mathbb{R}^p} \mathring \ell(\bx; \boldsymbol{\Xi}^{(\iota+1)} \mathbf{1}_p)\,, \quad 
\hat{\by}_{\iota+1} = \arg\min_{\by\in \mathbb{R}^p} \mathring \ell(\by; \boldsymbol{\Gamma}^{(\iota+1)} \mathbf{1}_p)\,,
\]
where $\hat{\bx}_{\iota+1} = (\hat{x}_{1,\iota+1}, \ldots, \hat{x}_{p,\iota+1})^\T$, $\hat{\by}_{\iota+1} = (\hat{y}_{1,\iota+1}, \ldots, \hat{y}_{p,\iota+1})^\T$ and $\mathring \ell(\cdot\,;\cdot)$ is specified in Lemma~\ref{convex}.  Then obtain $\mathring{\xi}_i^{(\iota+1)} = \exp(\hat{x}_{i,\iota+1})$, $\mathring{\eta}_i^{(\iota+1)} = \exp(\hat{y}_{i,\iota+1})$ for $i \in [p]$.

			\end{enumerate}
			%\end{itemize}
			%\item[] end do.	
			\item[3.] Output: Final estimator $$\mathring\btheta = \big(\{\mathring \btheta_{\mathcal{G}}^{(\iota+1)} \}^{\T}, \{\mathring \btheta_{\mathcal{G}^{\btc}}^{(\iota+1)}\}^{\T}\big)^{\T}\,,$$ where $\mathring{\btheta}_{\mathcal{G}^{\btc}}^{(\iota+1)} = (\mathring{\xi}_1^{(\iota+1)}, \ldots, \mathring{\xi}_p^{(\iota+1)}, \mathring{\eta}_1^{(\iota+1)}, \ldots, \mathring{\eta}_p^{(\iota+1)})^{\T}$.
		\end{enumerate}
	\end{algorithm}		
\end{center}

}

\subsection{Technical support}
\label{supp:imom}

We present Lemma~\ref{convex} to establish the strong convexity of $\mathring \ell(\bx; \bv)$ and provide relevant discussion to justify Step~\ref{alg1:2c} in Algorithm~\ref{alg1}.

\begin{lemma} \label{convex}
    For any $\bx = (x_1, \dots, x_p)^\T \in \mathbb{R}^p$ and $\bv = (v_1, \dots, v_p)^\T \in \mathbb{R}^p$, define 
\[
\mathring \ell(\bx; \bv) = \sum_{i,j:\,1 \le i < j \le p} \exp(x_i + x_j) - \sum_{i=1}^p x_i v_i\,.
\]
Then $\mathring \ell(\bx; \bv)$ is strongly convex in $\bx \in \mathbb{R}^p$.
\end{lemma}
\noindent
\textit{{Proof}}. 
Let $\bH(\bx) = \{H_{i,j}(\bx)\}_{i,j \in [p]}$ denote the Hessian matrix of $\mathring \ell(\bx; \bv)$. By basic calculations, we obtain that
\[
H_{i,i}(\bx) = \sum_{j : \, j \ne i} \exp(x_i + x_j)\,, \quad \text{for~} i \in [p]\,,
\]
\[
H_{i,j}(\bx) = \exp(x_i + x_j)\,, \quad \text{for~} 1 \le i \ne j \le p\,.
\]
For any $\bx \in \mathbb{R}^p$ and $\bz = (z_1, \dots, z_p)^\T \in \mathbb{R}^p$, it follows that
\[
\mathbf{z}^\top \mathbf{H}(\bx) \mathbf{z} = \sum_{i,j:\,1 \le i < j \le p} (z_i + z_j)^2 \exp(x_i + x_j) \ge 0\,,
\]
with equality if and only if $z_1 = \cdots = z_p = 0$. We thus complete the
proof of Lemma~\ref{convex}. $\hfill\Box$

Lemma~\ref{convex} indicates that $\mathring \ell(\bx; \bv)$ admits a unique minimizer. Observe that 
\begin{equation}\label{eq:lm5:dr}
    \frac{\partial \mathring \ell(\bx; \bv)}{\partial x_i} = \sum_{j :\, j \ne i} \exp(x_i + x_j) - v_i\,.
\end{equation}
Recall $\bXi^{(\iota+1)} =  \{\mathring \Xi_{i,j}^{(\iota+1)}\}_{i,j \in [p]}$ and note that 
\begin{equation} \label{eq:lm5:model}
    \sum_{j :\, j \neq i} \exp\big\{\log(\mathring \xi_{i}^{(\iota+1)}) +\log(\mathring \xi_{j}^{(\iota+1)})\big\} = \sum_{j=1}^p \mathring \Xi_{i,j}^{(\iota+1)}\,.
\end{equation}
Comparing \eqref{eq:lm5:dr} and \eqref{eq:lm5:model}, the  log-transformed local parameters $\{\log(\xi_i)\}_{i \in [p]}$ can be identified by solving 
 $$\hat{\bx}_{\iota+1} = (\hat{x}_{1,\iota+1}, \dots, \hat{x}_{p,\iota+1})^\top = \arg\min_{\bx \in \mathbb{R}^p} \mathring \ell(\bx; \boldsymbol{\Xi}^{(\iota+1)} \mathbf{1}_p)\,,$$ and then we update $\{\xi_i\}_{i \in [p]}$ via $\mathring{\xi}_i^{(\iota+1)} = \exp(\hat{x}_{i,\iota+1})$ for all $i \in [p]$. Analogously, we compute $$\hat{\by}_{\iota+1} = (\hat{y}_{1,\iota+1}, \dots, \hat{y}_{p,\iota+1})^\top = \arg\min_{\by \in \mathbb{R}^p} \mathring \ell(\by; \boldsymbol{\Gamma}^{(\iota+1)} \mathbf{1}_p)\,,$$ leading to the update $\mathring{\eta}_i^{(\iota+1)} = \exp(\hat{y}_{i,\iota+1})$ for all $i \in [p]$.

\section{Additional real data analysis}
\label{app:add.real}

\subsection{Preliminary summaries of email interaction} \label{subsec:emails_appendix}

\begin{figure}
     \vspace{-0.2cm}
    \centering
     \begin{subfigure}[b]{0.45\textwidth}
         \centering
         \includegraphics[width=\textwidth]{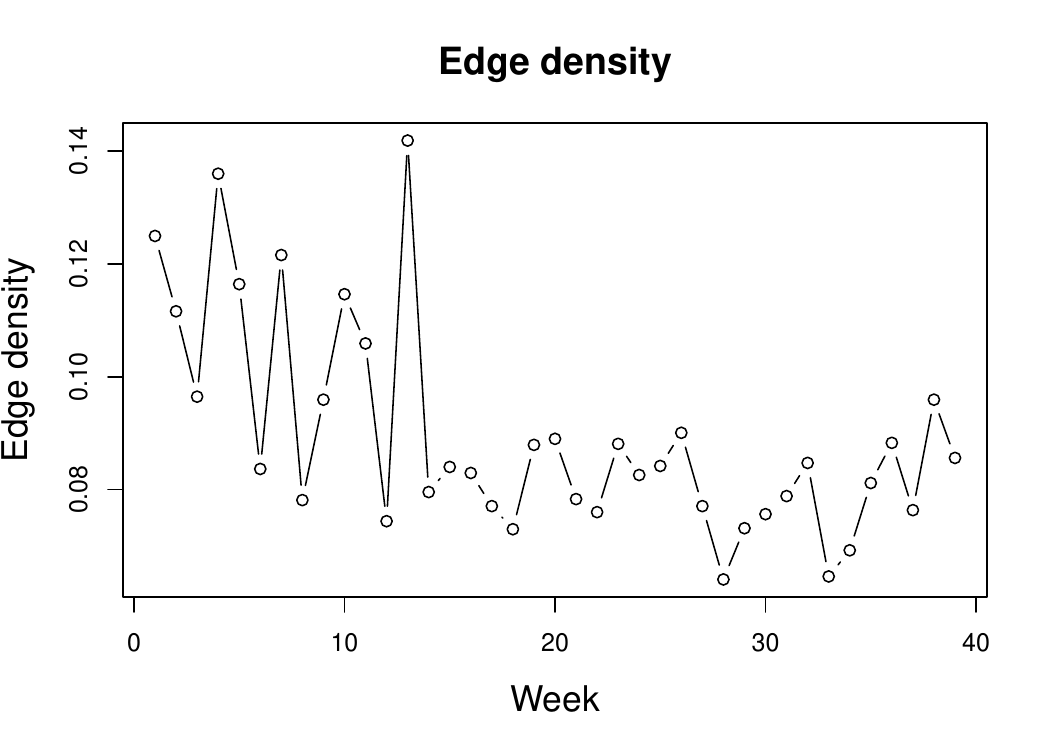}
     \end{subfigure}   
     %      \begin{subfigure}[b]{0.45\textwidth}
     %     \centering
     %     \includegraphics[width=\textwidth]{plot/real_data/man_grds_norm.pdf}
     % \end{subfigure}
        \begin{subfigure}[b]{0.45\textwidth}
         \centering
         \includegraphics[width=\textwidth]{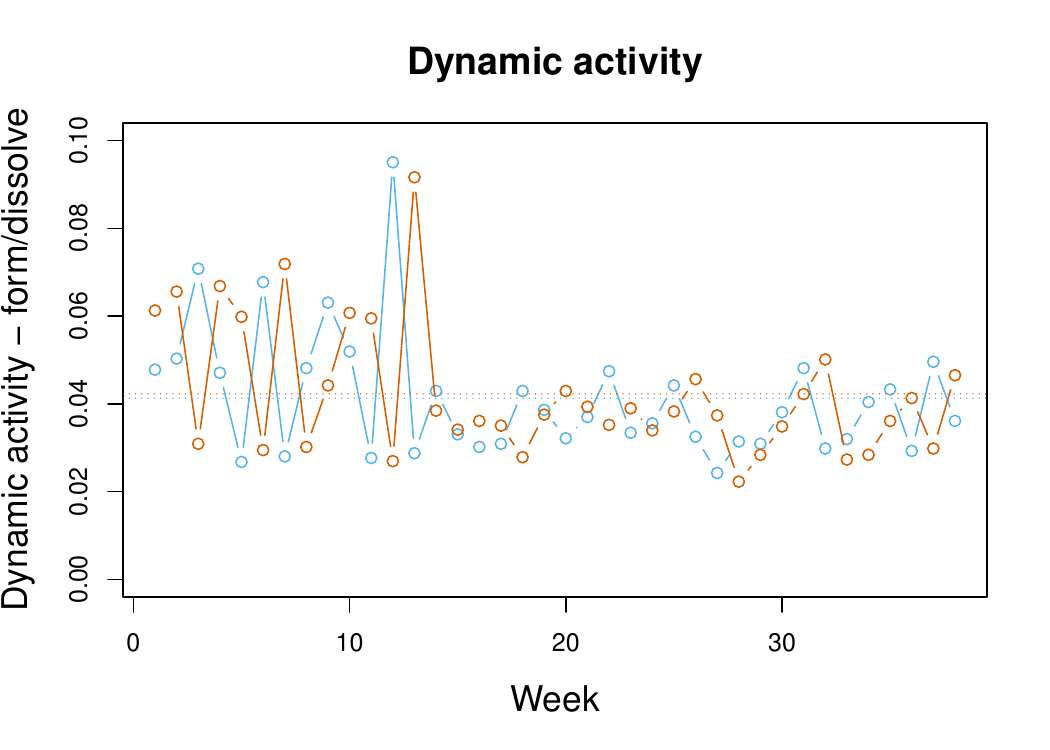}
     \end{subfigure}
          \vspace{-0.2cm}
     \caption{\label{fig:man_dens}{Evolution of edge density $D_t$ (left panel), percentage of grown $D_{1,t}$ (blue) and dissolved $D_{0,t}$ (orange) edges (right panel), email interaction networks.}}
     \vspace{-0.2cm}
\end{figure}

In this section, we present some preliminary summaries of the dynamic email network dataset analyzed in Section~\ref{sec:real}.
%We first present some preliminary summaries of the data to inspect the stationarity of the network and the effective sample size. 
The behavior shown in Figure~\ref{fig:man_dens} suggests a change point in the network behavior, in terms of both edge density $D_t$ and two dynamic (normalized) formation and dissolution measures $D_{1,t}$ and $D_{0,t}$ defined by
\begin{equation*}
    %\label{eq:3density}
D_t= {\sum_{i,j:\,i < j} X_{i,j}^t \over p(p-1)/2}\,, ~~
D_{1,t} = \frac{\sum_{i,j:\,i < j} (1-X_{i,j}^{t-1}) X_{i,j}^t}{p(p-1)/2}\,,
~~ 
D_{0,t} = \frac{\sum_{i,j:\,i< j} X_{i,j}^{t-1} (1-X_{i,j}^t)}{p(p-1)/2}\,.
\end{equation*}
%Hence in the following analysis, we fit the model separately to the first 13 and last 26 snapshots, referred to as ``period 1'' and ``period 2''. 
See also \eqref{eq:3density}. Based on this analysis, we split the dataset into the first 13 and last 26 snapshots, referred to as ``period 1'' and ``period 2''. 
In the right panel, about 4\% of node pairs see a grown edge or a dissolved edge between consecutive snapshots. 
However, dissolution can only occur for an active edge, and as we see in the left panel of Figure~\ref{fig:man_dens}, this network has relatively low edge density. 
On average (over time), the proportion of non-edge that form in the next snapshot is about 5\%, while the proportion of existing edges which persist in the next snapshot is about 55\%.
%clear evidence of temporal edge dependence.

\begin{figure} 
     \vspace{-0.2cm}
    \centering
     \begin{subfigure}[b]{0.85\textwidth}
         \centering
         \includegraphics[width=\textwidth]{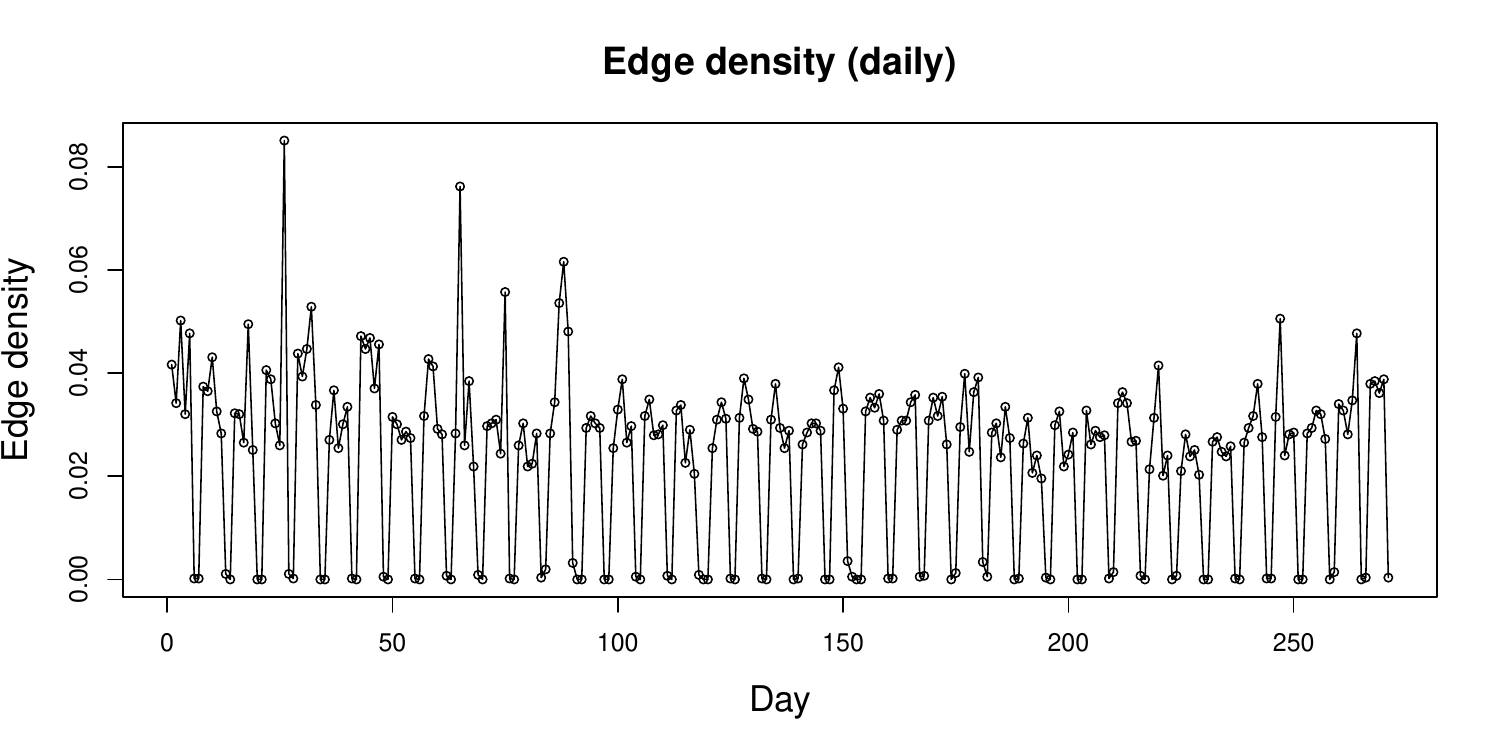}
     \end{subfigure}   
     
        \begin{subfigure}[b]{0.85\textwidth}
         \centering
         \includegraphics[width=\textwidth]{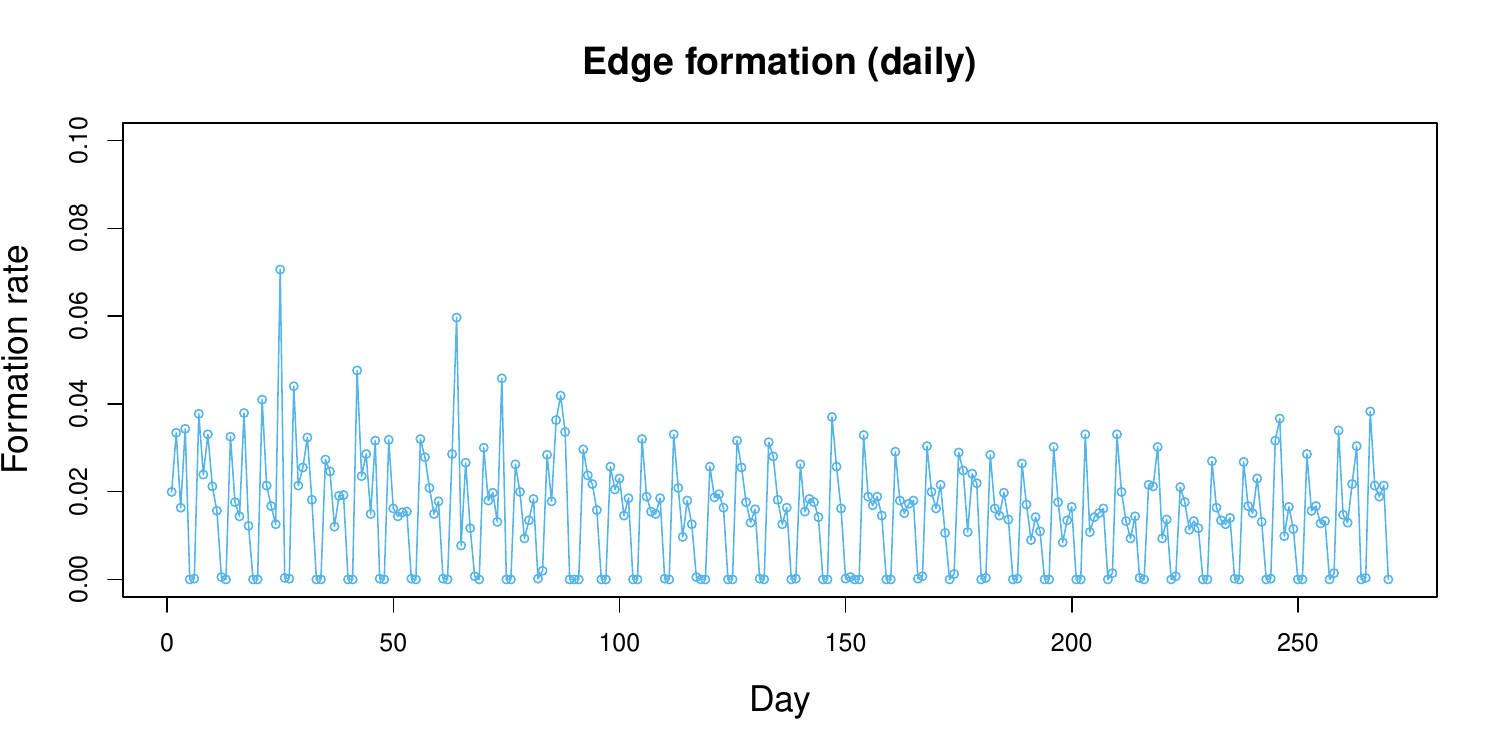}
     \end{subfigure}

     \begin{subfigure}[b]{0.85\textwidth}
         \centering
         \includegraphics[width=\textwidth]{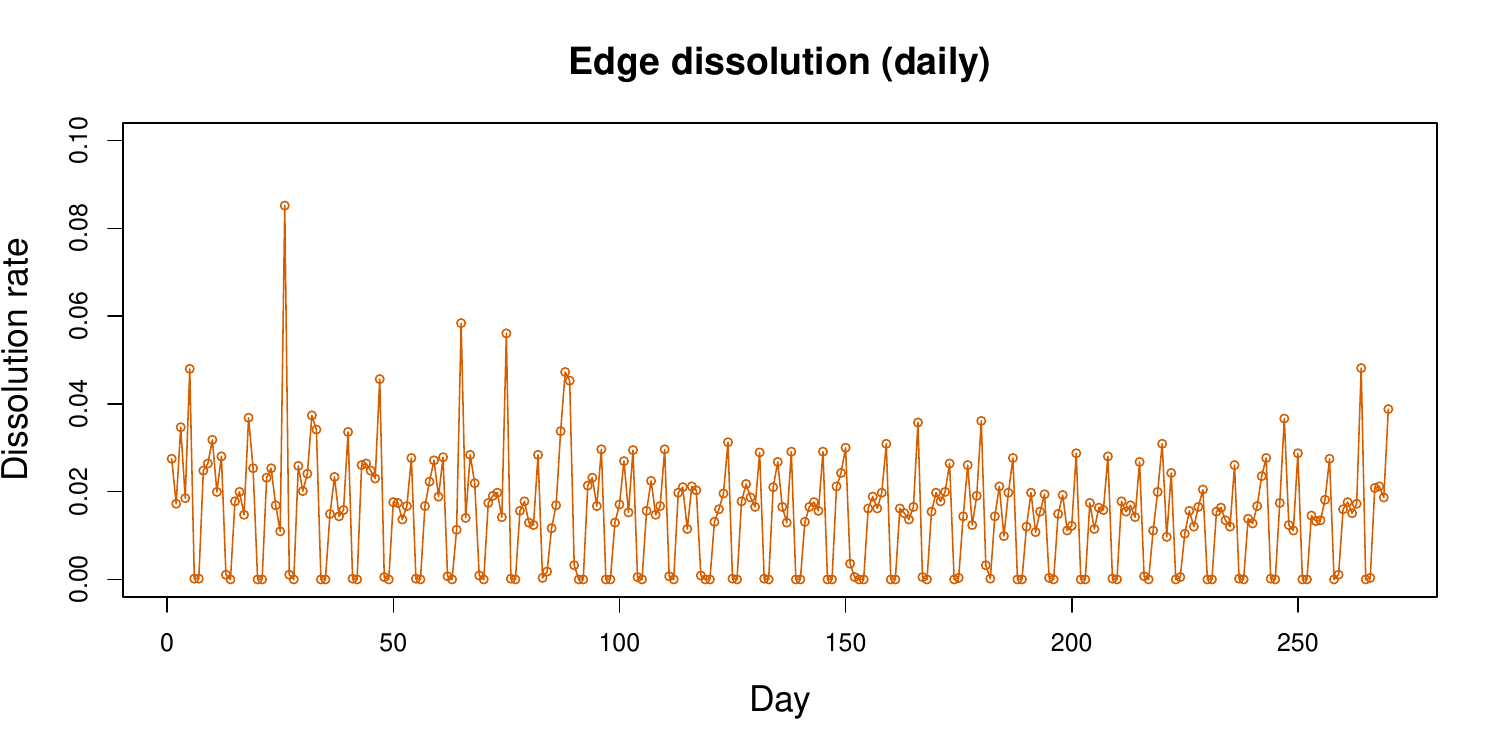}
     \end{subfigure}
          \vspace{-0.2cm}
     \caption{\label{fig:man_dens_day}{Evolution of edge density $D_t$ (top panel), percentage of grown $D_{1,t}$ (blue, center panel) and dissolved $D_{0,t}$ (orange, bottom panel) edges, manufacturing email networks constructed on a daily level.}} \label{dayEmail}
     \vspace{-0.2cm}
\end{figure}

\begin{figure}[t] 
     \vspace{-0.2cm}
    \centering
     \begin{subfigure}[b]{0.45\textwidth}
         \centering
         \includegraphics[width=\textwidth]{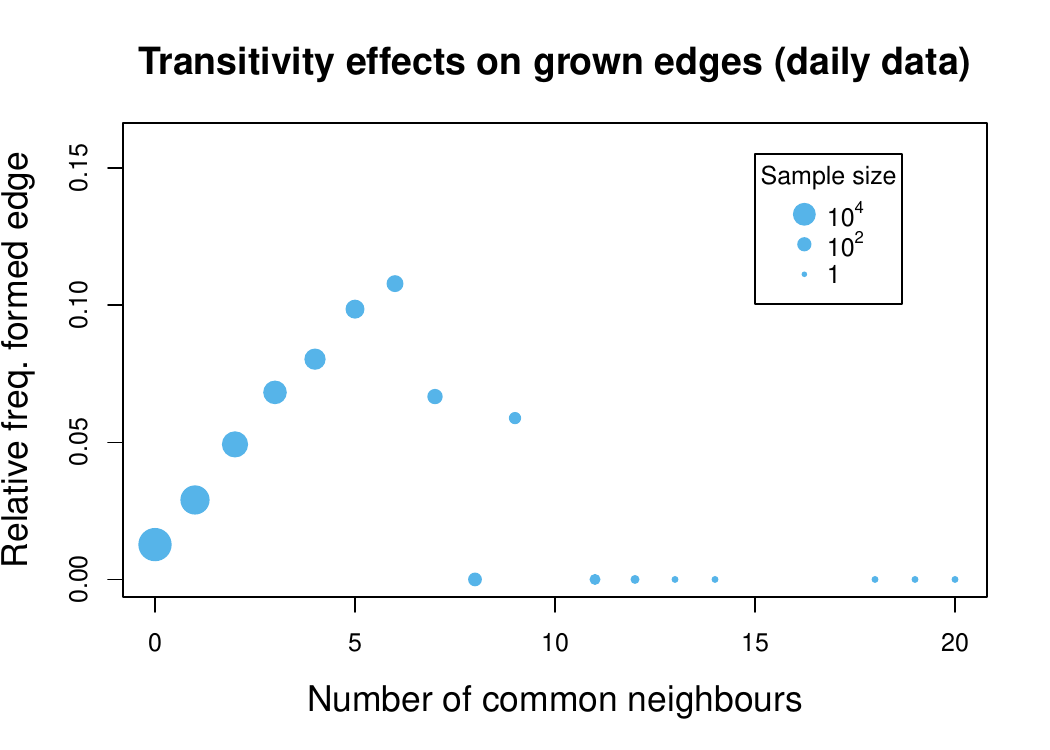}
     \end{subfigure}   
          \begin{subfigure}[b]{0.45\textwidth}
         \centering
         \includegraphics[width=\textwidth]{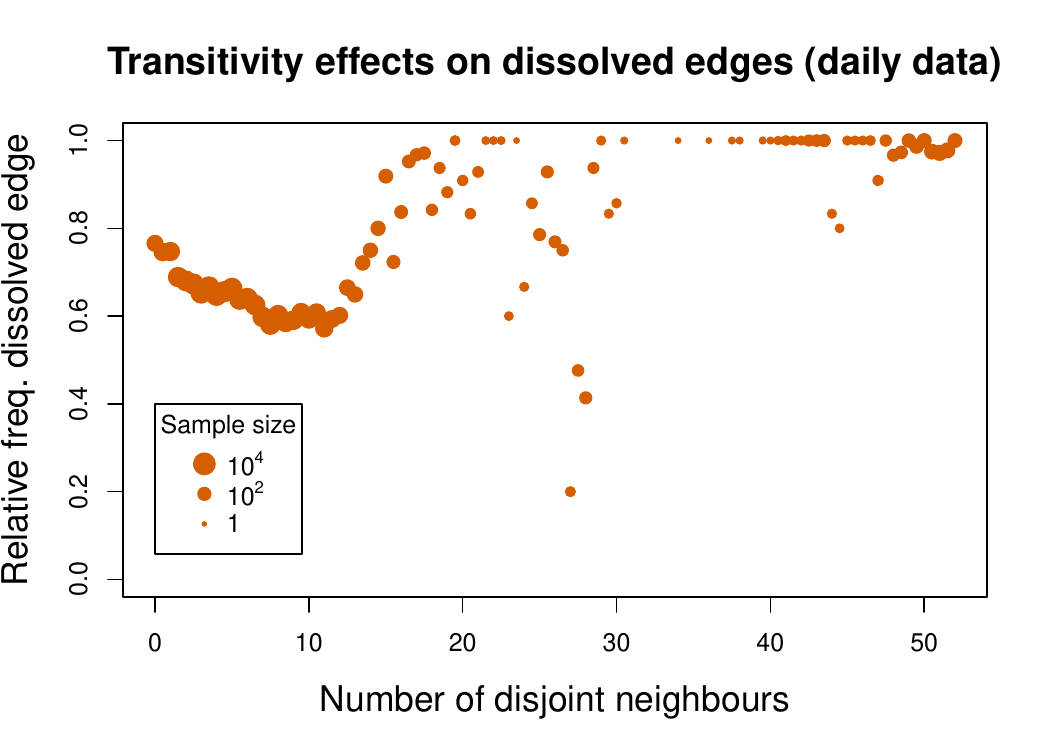}
     \end{subfigure}
          \vspace{-0.2cm}
     \caption{\label{fig:man_trans_day} Left panel: the plot of relative edge frequency $|\mathcal{U}_\ell^1|/|\mathcal{U}_\ell|$ against $\ell$, email interaction networks. Right panel: the plot of relative non-edge frequency $|\mathcal{V}_\ell^0|/|\mathcal{V}_\ell|$ against $\ell$, email interaction networks. In both panels, point size is proportional to the log sample sizes $\log|\mathcal{U}_\ell|$ and $\log|\mathcal{V}_\ell|$ respectively. The plots are constructed based on daily data.}
          \vspace{-0.2cm}
\end{figure}

The analysis above is based on the weekly email network data. As the original data contain the time stamp of each email exchange, we can also construct daily networks based on if there is any email exchange or not between two individuals on each day. Figure \ref{fig:man_dens_day} plots the daily edge density, the daily density of forming new edges and that of dissolving existing edges. By comparing it with Figure \ref{fig:man_dens}, the daily data exhibit more nonstationary features (such as effects of weekends, holidays, and periodicity through the week). The difference between the first 13 weeks and the last 26 weeks is more clearly exhibited in the weekly data. Figure \ref{fig:man_trans_day} reproduces Figure \ref{fig:man_trans} for the daily network data. Though the transitivity is still present with the daily data, it is mixed with more irregular signals. 
In the right panel of Figure \ref{fig:man_trans_day}, we do not see a clear effect of number of disjoint neighbours on edge dissolution; the daily networks are much sparser than the weekly ones, implying a lack of effective sample size to see meaningful patterns in edge dissolution.
By accumulating information on a weekly level, the data is more stationary and the empirical transitivity patterns are clearer.

\subsection{Comparative STERGM analysis of email interactions} \label{sec:email.tergm}

\begin{figure}[t]
    \centering
    \includegraphics[width=0.8\linewidth]{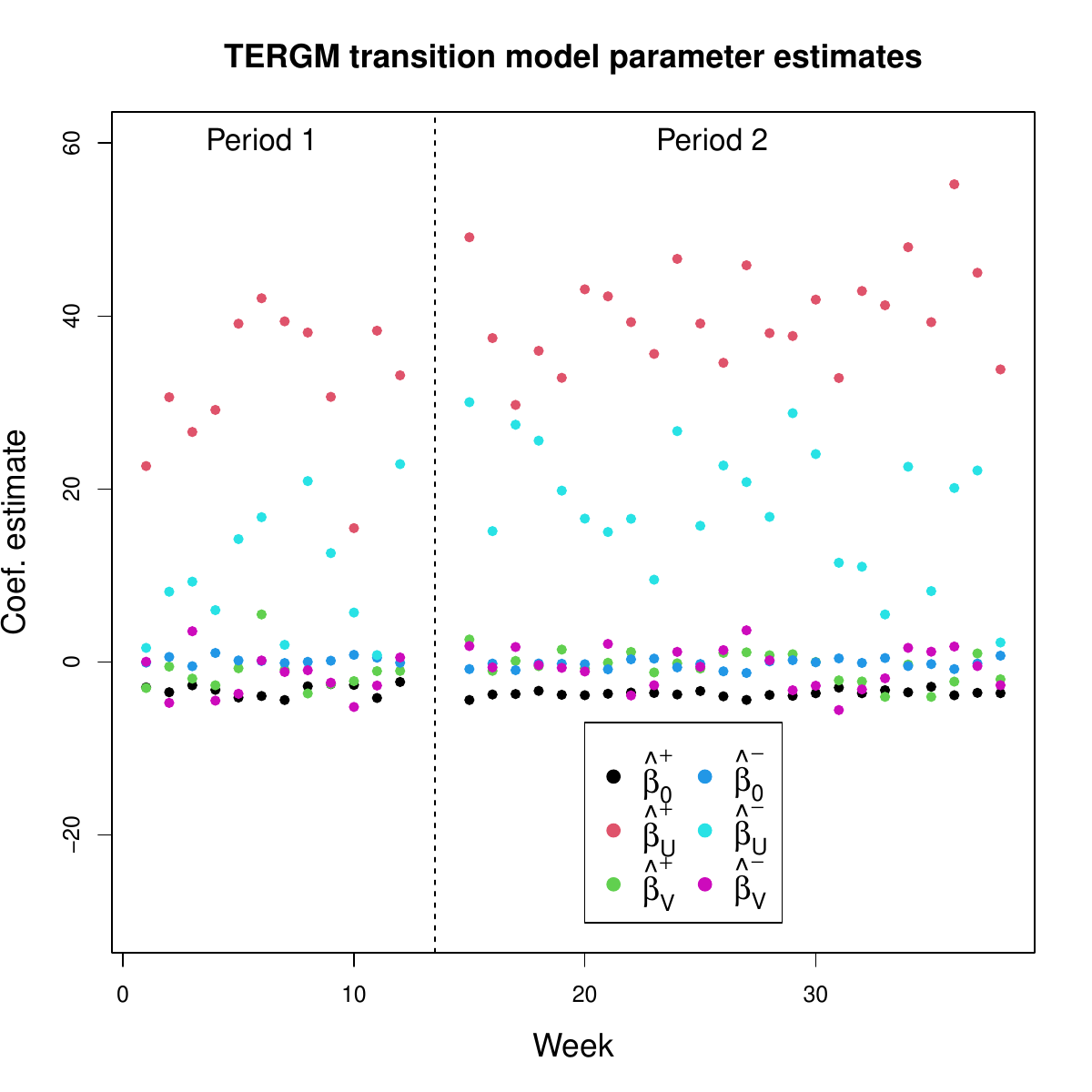}
    \caption{Sequence of TERGM parameter estimates for each network transition.}
    \label{fig:tergm_sequence}
\end{figure}

In order to compare the effects of the transitivity statistics $\{U^{t-1}_{i,j}\}$ and $\{V^{t-1}_{i,j}\}$ on the manufacturing email networks, we fit STERGM models using an implementation in the R package \texttt{tergm} from the \texttt{statnet} suite of software \citep{statnet2025_supp}.
First we fit a sequence of models, one for each transition from ${\bf X}_{t-1}$ to ${\bf X}_t$ for each $t$ within periods 1 and 2. 
Following \cite{Krivitsky2014_supp} and the notation introduced in our Appendix~\ref{app:relation2ERGM}, for each $t$, we specify a formation model by
$$
    \mathbb{P}({\bf X}^+_t ~\vert~ {\bf X}_{t-1})\, \propto \exp\bigg( \beta^{+}_0 \sum_{i,j:\, i \neq j} X_{i,j}^{t,+} + \beta^+_U \sum_{i,j:\, i \neq j} U_{i,j}^{t-1} X_{i,j}^{t,+} + \beta^+_V \sum_{i,j:\, i \neq j} V_{i,j}^{t-1} X_{i,j}^{t,+} \bigg)
$$
and ${\bf X}^+_t \in \mathcal{Y}^{+}({\bf X}_{t-1})$, where $\mathcal{Y}^{+}({\bf X}_{t-1})$ denotes the sample space of formation networks, i.e., all networks obtained by adding edges to ${\bf X}_{t-1}$, and $\bX_t^{+}= (X_{i,j}^{t,+})_{p\times p}$ with $X_{i,j}^{t,+}$ defined in \eqref{eq_X+-}.
% The edges present in ${\bf X}^+_t$ are added to the new network $X^t$.
We similarly specify a persistence model by
$$
    \mathbb{P}({\bf X}_t^- ~\vert~ {\bf X}_{t-1})\, \propto \exp\bigg( \beta^-_0 \sum_{i,j:\, i \neq j} X_{i,j}^{t,-} + \beta^-_U \sum_{i,j:\, i \neq j} U_{i,j}^{t-1} X_{i,j}^{t,-} + \beta^-_V \sum_{i,j:\, i \neq j} V_{i,j}^{t-1} X_{i,j}^{t,-} \bigg)
$$
and ${\bf X}_t^- \in \mathcal{Y}^{-}({\bf X}_{t-1})$, where $\mathcal{Y}^{-}({\bf X}_{t-1})$ denotes the sample space of dissolution networks, i.e., all networks obtained by removing edges from $\bX_{t-1}$, and $\bX_t^{-}= (X_{i,j}^{t,-})_{p\times p}$ with $X_{i,j}^{t,-}$ defined in \eqref{eq_X+-}.
Note that the edges present in $\bX_t^{-}$ are those that persist (i.e., are not dissolved) from $\bX_{t-1}$ to $\bX_t$, whereas our AR framework models edge dissolution.
To fit this separable model, we call \texttt{tergm} with the formula \texttt{$\sim$ Form($\sim$ edges + edgecov(U) + edgecov(V)) + Persist($\sim$ edges + edgecov(U) + edgecov(V))}. The resulting coefficient estimates are plotted in Figure~\ref{fig:tergm_sequence}. 
We expect the coefficients corresponding to \texttt{edgecov(U)} to influence edge formation and dissolution probabilities similarly to our AR transitivity model parameter $a$: positive values of $\beta^+_U$ and $\beta^-_U$ imply that more common neighbours in the previous snapshot increases probability of formation, and increases probability of persistence (decreases probability of dissolution). 
Analogously, we expect the coefficients corresponding to \texttt{edgecov(V)} to influence edge formation and dissolution probabilities inversely to our AR transitivity model parameter $b$ (see Section~\ref{sec:transitivity}).

Consistent with the analysis with the AR network model, we find that there is a strong effect of number of common neighbours in the previous time step on edge formation. \texttt{tergm} reports that all 38 of these parameter estimates are significantly different from 0 at the Bonferroni-corrected level $0.05/(6\times38)$. The signs of the estimated coefficients agree with our analysis (see Section~\ref{sec:real}), showing that number of common neighbours has an increasing relationship with edge formation probability. Also similar to our analysis, \texttt{tergm} finds that the magnitude of these effects are higher on average in period 2 than period 1. \texttt{tergm} reports (Bonferroni) significant coefficients for 19 out of 38 effects of $\{U^{t-1}_{i,j}\}$ on edge persistence (dissolution); the signs of these estimated coefficients are all positive, which agrees with our analysis. 
For the estimated \texttt{edgecov(V)} coefficients, \texttt{tergm} reports (Bonferroni) significant coefficients for only 2 out of 38 effects of $\{V^{t-1}_{i,j}\}$ on edge formation, and 7 out of 38 effects of $\{V^{t-1}_{i,j}\}$ on edge persistence (dissolution). Relative to our analysis, which pools information across time, and between formation and dissolution, this sequential analysis with \texttt{tergm} finds a less consistent effect of $\{V^{t-1}_{i,j}\}$ (normalized number of distinct neighbours) on edge formation and dissolution in $\bX_t$.

In order to inspect the node heterogeneity of the manufacturing email networks, we fit STERGM models to the sequences of snapshots from periods 1 and 2 with sociality effects, which fit one parameter per node to describe the degree heterogeneity in the formation and dissolution models.
The precise formation and persistence models are given by
\begin{align*}
    \mathbb{P}({\bf X}^+_t ~\vert~ {\bf X}_{t-1})\, &\propto \exp\bigg( \alpha^{+}_0 \sum_{i,j:\, i \neq j} X_{i,j}^{t,+} + \sum_{i=1}^p \alpha^+_i \sum_{j:\,j \neq i} X_{i,j}^{t,+} \bigg)\,, \quad {\bf X}^+_t \in \mathcal{Y}^{+}({\bf X}_{t-1})\,, \\
    \mathbb{P}({\bf X}_t^- ~\vert~ {\bf X}_{t-1})\, &\propto \exp\bigg(\alpha^{-}_0 \sum_{i,j:\, i \neq j} X_{i,j}^{t,-} + \sum_{i=1}^p \alpha^-_i \sum_{j:\,j \neq i} X_{i,j}^{t,-} \bigg)\,, \quad {\bf X}_t^- \in \mathcal{Y}^{-}({\bf X}_{t-1})\,.
\end{align*}
We fit this separable model by calling \texttt{tergm} with the formula \texttt{$\sim$ Form($\sim$ edges + sociality) + Persist($\sim$ edges + sociality)}. 
We expect the coefficients corresponding to the formation model sociality term to influence edge formation probabilities similarly to our local parameters $\{\xi_i\}$. We expect the coefficients corresponding to the persistence model sociality term to influence edge dissolution probabilities inversely to our local parameters $\{\eta_i\}$ (see Section~\ref{sec:transitivity}).

\begin{figure}[t]
    \centering
    \includegraphics[width=0.8\linewidth]{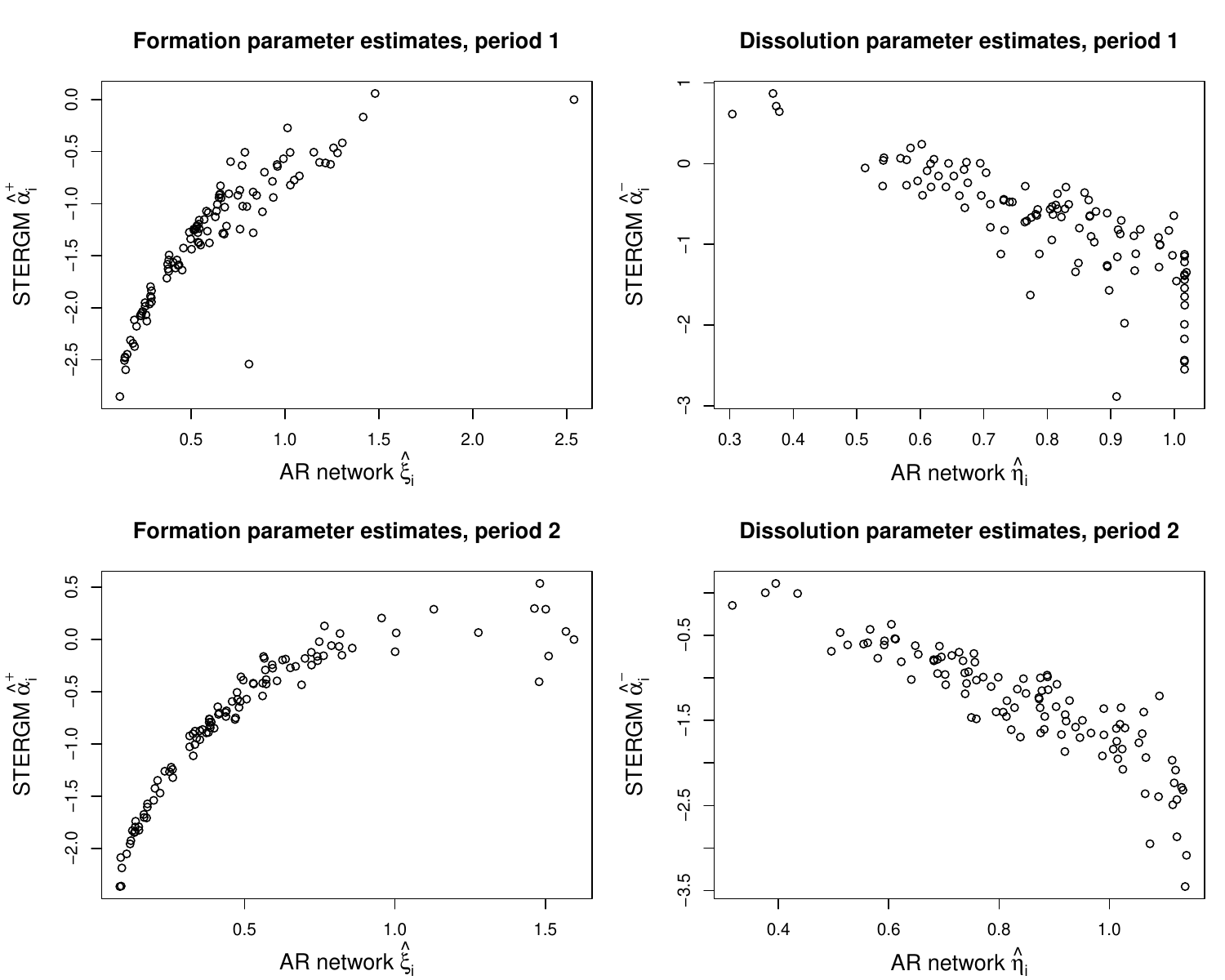}
    \caption{TERGM sociality parameter estimates plotted against AR network local parameter estimates.}
    \label{fig:tergm_degree}
\end{figure}

In Figure~\ref{fig:tergm_degree}, we summarize the fitted sociality effects by comparing them to the corresponding local parameters fit by our AR network model in Section~\ref{sec:real}.
Broadly, we see the expected increasing relationship between estimates of $\{\xi_i\}$ parameters and formation model sociality effects, and a decreasing relationship between estimates of $\{\eta_i\}$ parameters and persistence model sociality effects.

However, there are also subtle differences in interpretation of these parameters afforded by our new model formulation. The multiplicative structure of our model means that while each $\xi_i$ and $\eta_i$ are related to observed degree, their direct impact on the edge formation and dissolution is to modulate the strength of transitivity for each individual node.
In Figure~\ref{fig:tergm_degree}, we see a wide range of $\{\xi_i\}$ and $\{\eta_i\}$ estimates even for nodes with very similar estimates of their TERGM sociality parameters $\alpha^+_i$ and $\alpha^-_i$, implying that there is not a perfect correspondence between these two phenomena: nodes with similar socialities may still vary considerably in the strength of transitivity on their formed and dissolved connections. 

\subsection{Analysis of conference interactions}
\label{sec:conf.int.appl}

We also apply our dependent AR model with transitivity to an additional dynamic network dataset, in this case of face-to-face interactions among attendees of an academic conference.
The conference in question was the 2009 congress of the Soci\'{e}t\'{e} Fran\c{c}aise d'Hygi\`{e}ne Hospitali\`{e}re (SFHH) \citep{Cattuto2010_supp,Genois2018_supp}.
% \citep{cattuto2010,genois2018}.
The original data was collected automatically by RFID badges worn by the conference participants. 
We analyze a subset of the data corresponding to an active portion of the first day of the congress (June 4, 2009) from about 11:00AM to 6:00PM among the $p=200$  most active participants out of the total 403.
Each of the $n=22$ network snapshots corresponds to a non-overlapping time window, with $X_{i,j}^t=1$ if participants were in close proximity at any time during the prior 20 minutes.

%Summaries 

\begin{figure}
    \centering
     \begin{subfigure}[b]{0.45\textwidth}
         \centering
         \includegraphics[width=\textwidth]{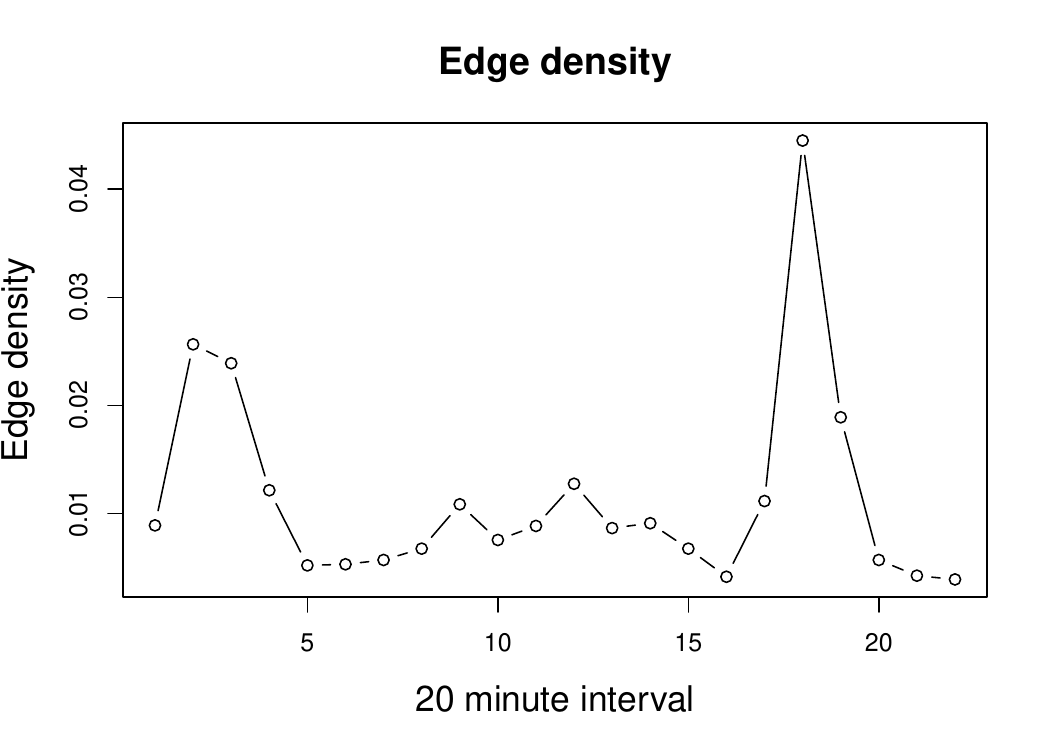}
     \end{subfigure}   
     %      \begin{subfigure}[b]{0.45\textwidth}
     %     \centering
     %     \includegraphics[width=\textwidth]{plot/real_data/sfhh_grds_norm.pdf}
     % \end{subfigure}
               \begin{subfigure}[b]{0.45\textwidth}
         \centering
         \includegraphics[width=\textwidth]{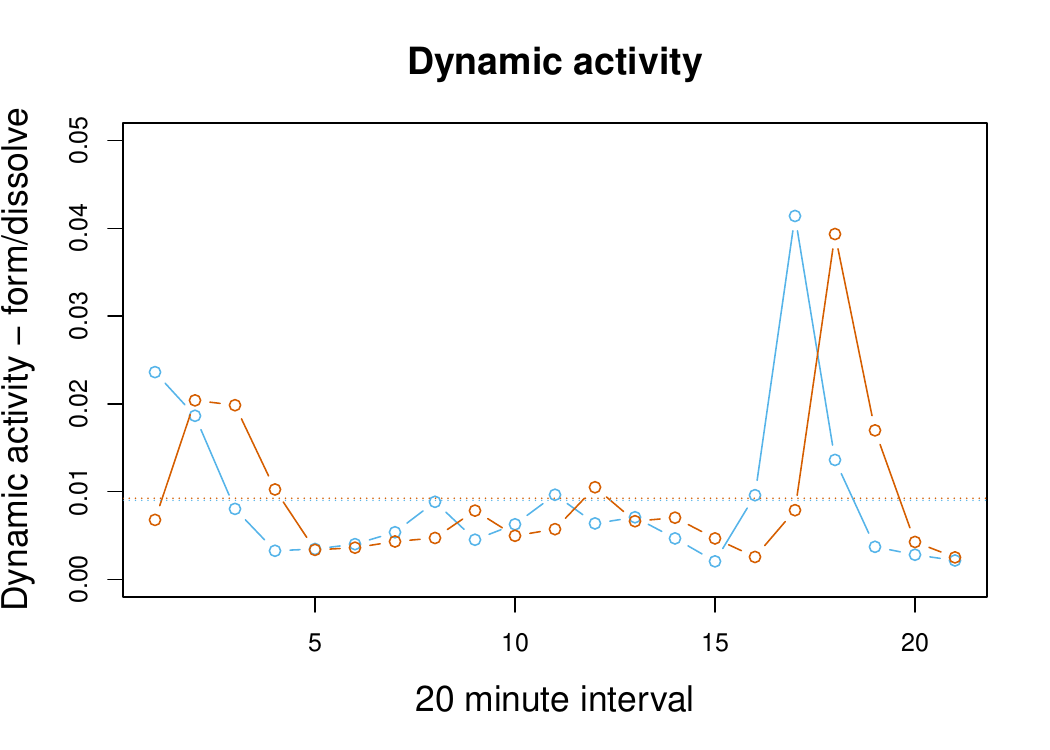}
     \end{subfigure}
     \caption{\label{fig:sfhh_dens}{Evolution of edge density $D_t$ (left panel), percentage of grown $D_{1,t}$ (blue) and dissolved $D_{0,t}$ (orange) edges (right panel), SFHH participant networks. The quantities $D_t$, $D_{0,t}$, and $D_{1,t}$ are defined as in Section~\ref{subsec:emails_appendix}}.}
\end{figure}

\begin{figure}[t] 
    \centering
     \begin{subfigure}[b]{0.45\textwidth}
         \centering
         \includegraphics[width=\textwidth]{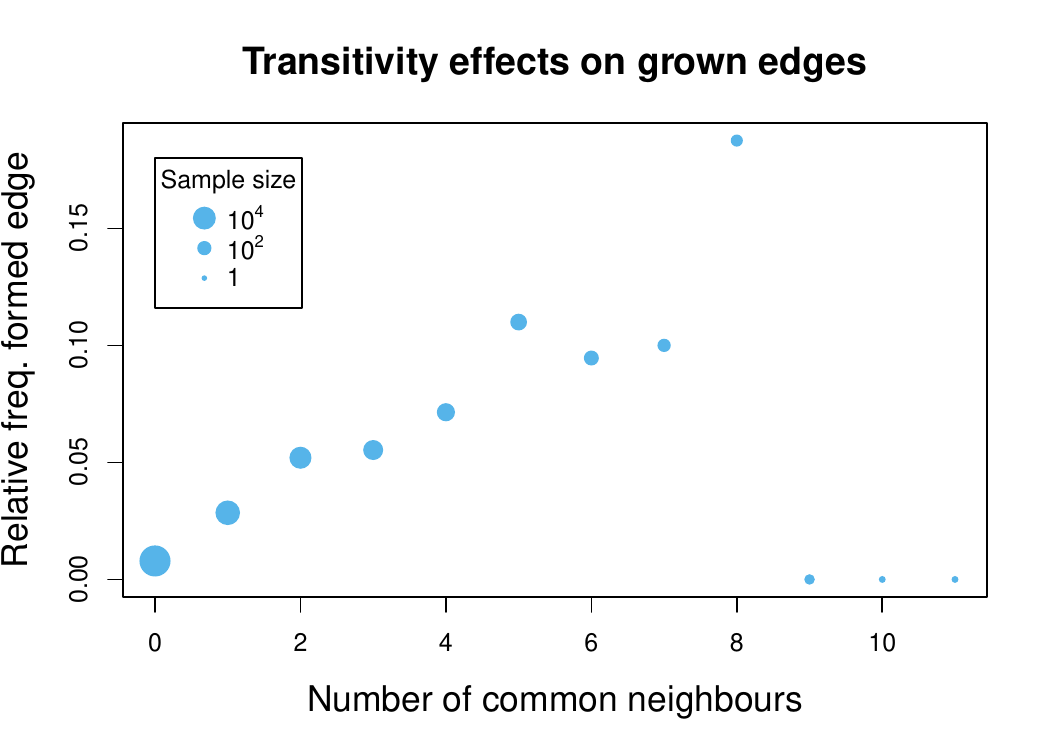}
     \end{subfigure}   
          \begin{subfigure}[b]{0.45\textwidth}
         \centering
         \includegraphics[width=\textwidth]{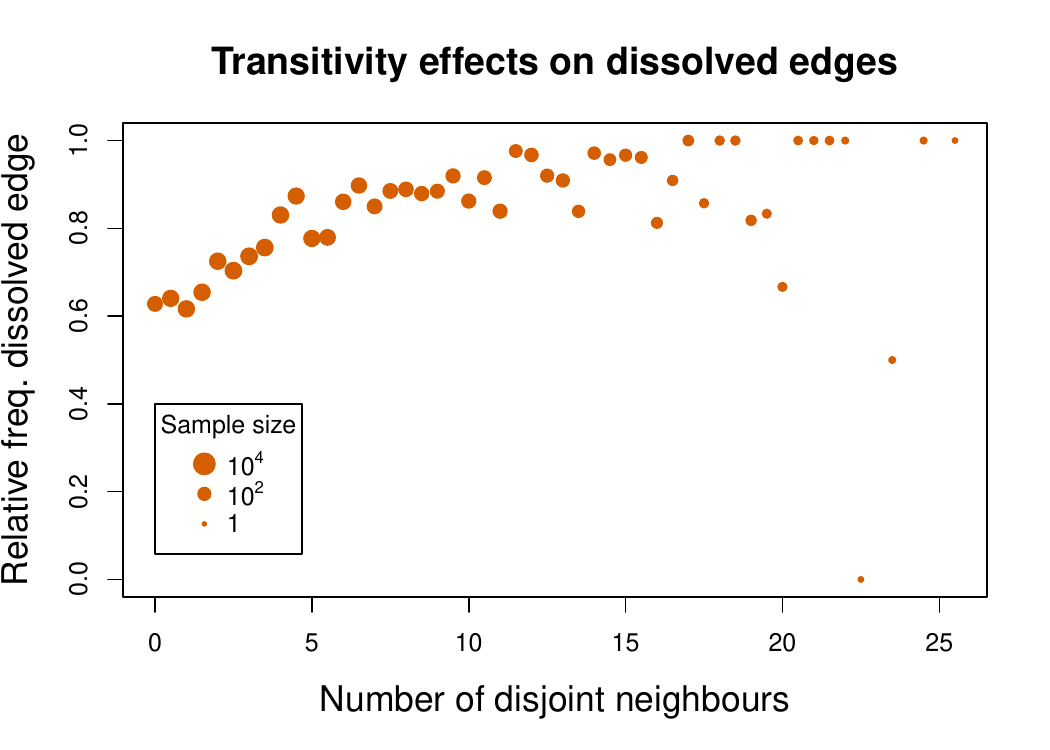}
     \end{subfigure}
     \caption{\label{fig:sfhh_trans} Left panel: the plot of relative edge frequency $|\mathcal{U}_\ell^1|/|\mathcal{U}_\ell|$ against $\ell$, SFHH participant networks. Right panel: the plot of relative non-edge frequency $|\mathcal{V}_\ell^0|/|\mathcal{V}_\ell|$ against $\ell$, SFHH participant networks. In both panels, point size is proportional to the log sample sizes $\log|\mathcal{U}_\ell|$ and $\log|\mathcal{V}_\ell|$ respectively. The quantities $\mathcal{U}_{\ell}$, $\mathcal{U}_{\ell}^1$, $\mathcal{V}_{\ell}$, and $\mathcal{V}_{\ell}^0$ are defined as in Section~\ref{subsec:emails_appendix}.}
          \vspace{-0.2cm}
\end{figure}

Similar to Appendix~\ref{subsec:emails_appendix}, we summarize some key features of this dataset.
Figure~\ref{fig:sfhh_dens} (left panel) shows that while there are some spikes in edge density, there is no clear increasing or decreasing pattern, so we choose to model this dataset with a single AR network model.
Figures~\ref{fig:sfhh_dens} and \ref{fig:sfhh_trans} show empirical evidence of temporal edge dependence, as well as transitivity effects: after accounting for edge density, edges persist at a higher rate than they grow, they more often grow for node pairs which had more common neighbours, and they more often dissolve for node pairs which had more disjoint neighbours.

%Model estimates 

\begin{figure}[t] 
    \centering
     \begin{subfigure}[b]{0.3\textwidth}
         \centering
         \includegraphics[width=\textwidth]{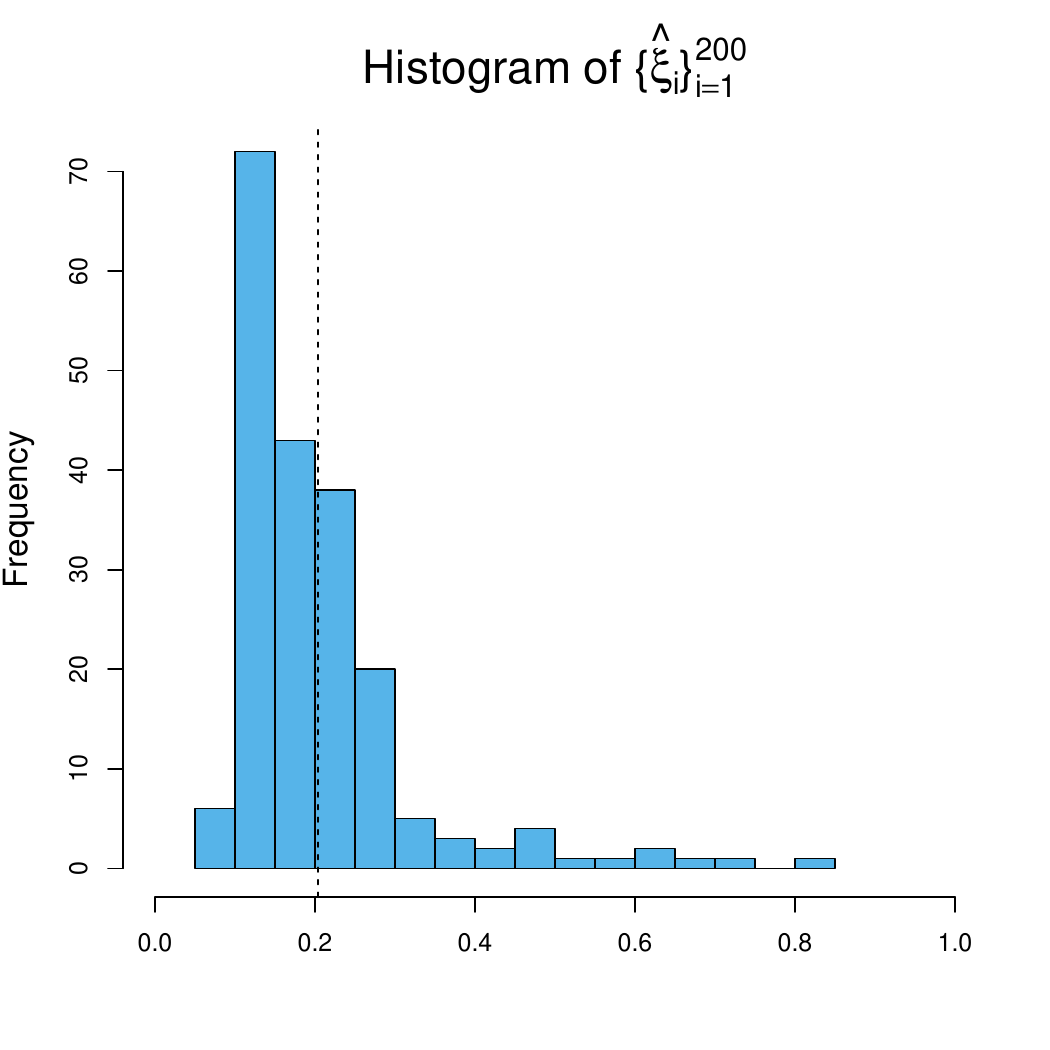}
     \end{subfigure}   
          \begin{subfigure}[b]{0.3\textwidth}
         \centering
         \includegraphics[width=\textwidth]{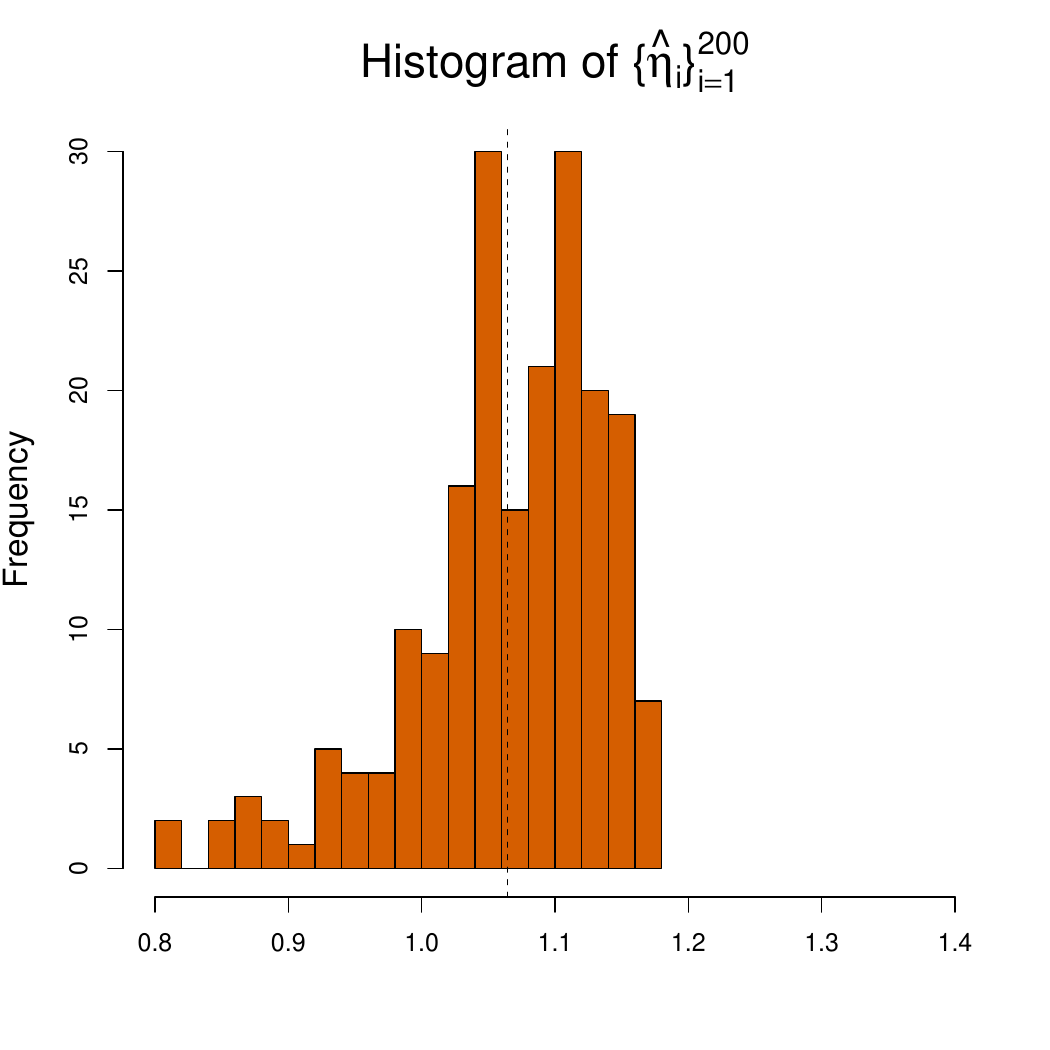}
     \end{subfigure}   
          \begin{subfigure}[b]{0.3\textwidth}
         \centering
         \includegraphics[width=\textwidth]{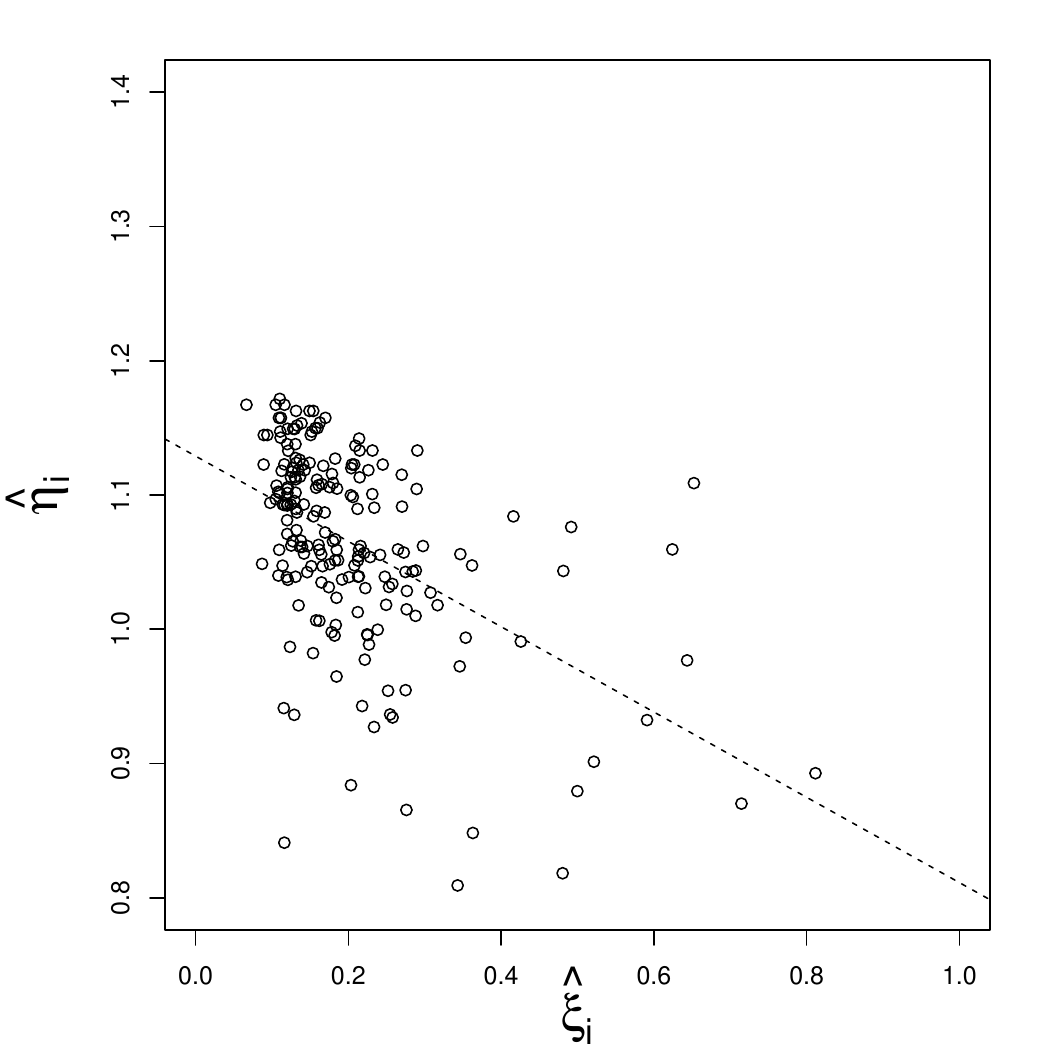}
     \end{subfigure}
     \caption{\label{fig:sfhh_scatters}{Histograms and scatter plot of estimates $\{\hat{\xi}_i\}_{i=1}^{200}$ and $\{\hat{\eta}_i\}_{i=1}^{200}$, SFHH participant networks.}}
\end{figure}

Fitting our AR network model with transitivity, we estimate $\hat{a}=23.20$ and $\hat{b}=19.08$, confirming these empirical dynamic effects of common and disjoint neighbours. 
We summarize the estimates of the local parameters $\{\xi_i\}_{i=1}^{200}$ and $\{\eta_i\}_{i=1}^{200}$ in Figure~\ref{fig:sfhh_scatters}.

The estimates $\{\hat{\xi}_i\}_{i=1}^{200}$ have mean 0.20 and a longer right tail, while the estimates $\{\hat{\eta}_i\}_{i=1}^{200}$ have mean 1.06 and a longer left tail. 
Moreover, their scatter plot shows that there is a negative relationship between these estimates for a given node. All of these observations are consistent with overall degree heterogeneity of the network.
The estimates $\{\hat{\xi}_i\}_{i=1}^{200}$ have a wide range from 0.07 to 0.81, while the estimates $\{\hat{\eta}_i\}_{i=1}^{200}$ are all between 0.81 and 1.17.
This implies that conference attendees are more heterogeneous in their propensity to form new connections than in their propensity to extend the length of existing ones.

%AIC/BIC comparison

\begin{table}[]
\begin{center}
\begin{tabular}{|l|l|l|}
\hline
\textbf{Model}              & \textbf{AIC}   & \textbf{BIC}   \\ \hline
Transitivity AR model       & 48997 & 53397          \\ \hline
Global AR model             & \textbf{48037}          & \textbf{48059} \\ \hline
Edgewise AR model           & 109284         & 544815         \\ \hline
Edgewise mean model         & 71205          & 288970         \\ \hline
Degree parameter mean model & 54061          & 56250          \\ \hline
\end{tabular}
\end{center}
\caption{\label{tab:aic_sfhh}{AIC and BIC performance for SFHH participant networks}}
\end{table}

Finally, we compare our model to the same competing models described in Section~\ref{sec:real} in terms of AIC and BIC. These results are reported in Table~\ref{tab:aic_sfhh}. The global AR model, which only has 2 parameters, achieves the lowest AIC and BIC, followed closely by our AR network model with transitivity. 
%In both cases the best model incorporates temporal edge dependence.
Under either criterion, the two edgewise models require $O(p^2)$ parameters and thus perform poorly, as there are relatively many nodes ($p=200$) compared to network samples ($n=22$).
%Although our model and the degree parameter mean model still have $O(p)$ parameters, their criteria are the best or competitive, suggesting that degree heterogeneity parameters are an efficient parameterization to summarize the structure in this dataset.
Although the transitivity model has $O(p)$ parameters, it achieves the second lowest AIC and BIC,  suggesting that degree heterogeneity parameters and the imposed transitivity form are an effective parameterization to summarize the structure in this dataset.

\section{Discussion on model specification and scalability}
\label{app:add.discussion}

\subsection{Remarks on Condition~\ref{as:gambd}}

Recall that Condition~\ref{as:gambd} requires that the transition probabilities $\alpha_{i,j}^{t-1}(\btheta)$  and $\beta_{i,j}^{t-1}(\btheta)$ for $t\in[n]\setminus[m]$ and $1\leq i<j\leq p$ to lie strictly between 0 and 1. This condition implicitly defines the feasible parameter space for all parameters. To illustrate, consider the transitivity model introduced in Section \ref{sec:transitivity}, with transition probabilities given in \eqref{trans-model}.  If no explicit bounds are imposed on the global parameters $a$ and $b$, the logistic components 
\begin{align*}
    \frac{\exp(a U_{i,j}^{t-1}) }{1+ \exp(a U_{i,j}^{t-1})+ \exp(b V_{i,j}^{t-1})}~~\text{and}~~\frac{ \exp(b V_{i,j}^{t-1})}  {1+ \exp(a U_{i,j}^{t-1})+ \exp(b V_{i,j}^{t-1})}
\end{align*}
always lie between 0 and 1. One could then impose constraints
like $\max_i|\xi_i| \le 1$ and $\max_i|\eta_i| \le 1$ as a sufficient condition to guarantee that transition probabilities lie between 0 and 1. 
However, such restrictions are not necessary for ensuring valid transition probabilities in the transitivity model.
To maintain generality and flexibility in the theoretical framework, we thus do not impose explicit parameter constraints tailored to specific models. Instead, Condition~\ref{as:gambd}  ensures that all transition probabilities lie strictly between 0 and 1 while keeping the parameter space fully general.

\subsection{Comparison of scalability with continuous-time models}

Practically speaking, the proposed estimation procedures in Sections~\ref{sec:initialEstimation} and \ref{sec:improvedEstimation} are not as scalable as those developed for continuous-time network models in \cite{Vu2011a_supp, Vu2011b_supp,Salathe2013_supp}.
The approaches in these works achieved notable scalability by exploiting event-level sparsity in continuous-time network data, in the sense that only a small subset of nodes or edges are updated between two consecutive event times. 
Our framework adopts a discrete-time formulation that specifies transition probabilities for the formation and dissolution of edges between every pair of nodes.
To fully leverage the dynamics across time, for each parameter $l \in [p]$, all related $|\mathcal{S}_l|$ edges need to be evaluated, which leads to potentially heavier computation in practice.
Moreover, the models in \cite{Vu2011a_supp, Vu2011b_supp} estimate a  coefficient vector $\bbeta$, which captures the effects of network summary statistics such as out-degree, in-degree, reciprocity, and transitivity. This $\bbeta$ plays a similar role to the global parameters in our framework, though the two formulations differ fundamentally in how time is represented and how parameters are estimated.  In addition, our model includes local parameters to represent node-specific heterogeneity, whose number may diverge together with $p$. For the three models introduced in Section~\ref{sec:3examples}, the number of local parameters is of order $O(p)$, which generally requires greater computational effort during estimation.

\section{Technical proofs} \label{app:proofs}

\subsection{Proof of Proposition \ref{pn:ident}} \label{pr:pn:ident}

Recall 
\[
\ell_{n,p}^{(l)}(\btheta)=\frac{1}{(n-m)|\mathcal{S}_l|}\sum_{t=m+1}^n\sum_{(i,j)\in\mathcal{S}_l}h_{i,j}^{t-1}(\btheta)
\]
with
\begin{align*}
h_{i,j}^{t-1}(\btheta)
=\log \{1- \gamma_{i,j}^{t-1}(\btheta)\}+\gamma_{i,j}^{t-1}(\btheta_0)\log \bigg\{\frac{\gamma_{i,j}^{t-1}(\btheta)}{1- \gamma_{i,j}^{t-1}(\btheta)}\bigg\}\,.
\end{align*}
Then
\begin{align}
\frac{\partial h_{i,j}^{t-1}(\btheta)}{\partial\theta_{l_1}}=&~\frac{\gamma_{i,j}^{t-1}(\btheta_0)-\gamma_{i,j}^{t-1}(\btheta)}{\gamma_{i,j}^{t-1}(\btheta)\{1-\gamma_{i,j}^{t-1}(\btheta)\}}\frac{\partial \gamma_{i,j}^{t-1}(\btheta)}{\partial\theta_{l_1}}\,,\label{eq:h1der}\\
\frac{\partial^2 h_{i,j}^{t-1}(\btheta)}{\partial\theta_{l_1}\partial\theta_{l_2}}=&-\bigg[\frac{1}{\gamma_{i,j}^{t-1}(\btheta)\{1-\gamma_{i,j}^{t-1}(\btheta)\}}+\frac{\{\gamma_{i,j}^{t-1}(\btheta_0)-\gamma_{i,j}^{t-1}(\btheta)\}\{1-2\gamma_{i,j}^{t-1}(\btheta)\}}{\{\gamma_{i,j}^{t-1}(\btheta)\}^2\{1-\gamma_{i,j}^{t-1}(\btheta)\}^2}\bigg]\notag\\
&~~~~~\times\frac{\partial\gamma_{i,j}^{t-1}(\btheta)}{\partial \theta_{l_1}}\frac{\partial\gamma_{i,j}^{t-1}(\btheta)}{\partial \theta_{l_2}}\label{eq:h2der}\\
&+\frac{\gamma_{i,j}^{t-1}(\btheta_0)-\gamma_{i,j}^{t-1}(\btheta)}{\gamma_{i,j}^{t-1}(\btheta)\{1-\gamma_{i,j}^{t-1}(\btheta)\}}\frac{\partial^2\gamma_{i,j}^{t-1}(\btheta)}{\partial\theta_{l_1}\partial\theta_{l_2}}\,,\notag\\
\frac{\partial^3h_{i,j}^{t-1}(\btheta)}{\partial\theta_{l_1}\partial\theta_{l_2}\partial\theta_{l_3}}=&~2\bigg[\frac{1-2\gamma_{i,j}^{t-1}(\btheta)}{\{\gamma_{i,j}^{t-1}(\btheta)\}^2\{1-\gamma_{i,j}^{t-1}(\btheta)\}^2}+\frac{\gamma_{i,j}^{t-1}(\btheta_0)-\gamma_{i,j}^{t-1}(\btheta)}{\{\gamma_{i,j}^{t-1}(\btheta)\}^2\{1-\gamma_{i,j}^{t-1}(\btheta)\}^2}\notag\\
&~~~~~+\frac{\{\gamma_{i,j}^{t-1}(\btheta_0)-\gamma_{i,j}^{t-1}(\btheta)\}\{1-2\gamma_{i,j}^{t-1}(\btheta)\}^2}{\{\gamma_{i,j}^{t-1}(\btheta)\}^3\{1-\gamma_{i,j}^{t-1}(\btheta)\}^3}\bigg]\frac{\partial\gamma_{i,j}^{t-1}(\btheta)}{\partial\theta_{l_1}}\frac{\partial\gamma_{i,j}^{t-1}(\btheta)}{\partial\theta_{l_2}}\frac{\partial\gamma_{i,j}^{t-1}(\btheta)}{\partial\theta_{l_3}}\notag\\
&-\bigg[\frac{1}{\gamma_{i,j}^{t-1}(\btheta)\{1-\gamma_{i,j}^{t-1}(\btheta)\}}+\frac{\{\gamma_{i,j}^{t-1}(\btheta_0)-\gamma_{i,j}^{t-1}(\btheta)\}\{1-2\gamma_{i,j}^{t-1}(\btheta)\}}{\{\gamma_{i,j}^{t-1}(\btheta)\}^2\{1-\gamma_{i,j}^{t-1}(\btheta)\}^2}\bigg]\label{eq:h3der}\\
&~~~~~\times\bigg\{\frac{\partial^2\gamma_{i,j}^{t-1}(\btheta)}{\partial\theta_{l_1}\partial\theta_{l_2}}\frac{\partial\gamma_{i,j}^{t-1}(\btheta)}{\partial\theta_{l_3}}+\frac{\partial^2\gamma_{i,j}^{t-1}(\btheta)}{\partial\theta_{l_2}\partial\theta_{l_3}}\frac{\partial\gamma_{i,j}^{t-1}(\btheta)}{\partial\theta_{l_1}}+\frac{\partial^2\gamma_{i,j}^{t-1}(\btheta)}{\partial\theta_{l_3}\partial\theta_{l_1}}\frac{\partial\gamma_{i,j}^{t-1}(\btheta)}{\partial\theta_{l_2}}\bigg\}\notag\\
&+\frac{\gamma_{i,j}^{t-1}(\btheta_0)-\gamma_{i,j}^{t-1}(\btheta)}{\gamma_{i,j}^{t-1}(\btheta)\{1-\gamma_{i,j}^{t-1}(\btheta)\}}\frac{\partial^3\gamma_{i,j}^{t-1}(\btheta)}{\partial\theta_{l_1}\partial\theta_{l_2}\partial\theta_{l_3}}\,.\notag
\end{align}
By the triangle inequality and Condition \ref{as:gambd},
\begin{align}\label{eq:h3bd}
\bigg|\frac{\partial^3h_{i,j}^{t-1}(\btheta)}{\partial\theta_{l_1}\partial\theta_{l_2}\partial\theta_{l_3}}\bigg|\leq&~2(2C_1^{-2}+C_1^{-3})\bigg|\frac{\partial\gamma_{i,j}^{t-1}(\btheta)}{\partial\btheta}\bigg|_\infty^3+C_1^{-1}\bigg|\frac{\partial^3\gamma_{i,j}^{t-1}(\btheta)}{\partial\btheta^3}\bigg|_\infty\notag\\
&+3(C_1^{-1}+C_1^{-2})\bigg|\frac{\partial^2\gamma_{i,j}^{t-1}(\btheta)}{\partial\btheta^2}\bigg|_\infty\bigg|\frac{\partial\gamma_{i,j}^{t-1}(\btheta)}{\partial\btheta}\bigg|_\infty\\
\leq&~2(2C_1^{-2}+C_1^{-3})C_2^3+3(C_1^{-1}+C_1^{-2})C_2^2+C_1^{-1}C_2=:C_*\,.\notag
\end{align}
Write $\btheta=(\theta_1,\ldots,\theta_q)^{\T}$ and $\btheta_0=(\theta_{0,1},\ldots,\theta_{0,q})^{\T}$. By Taylor's theorem, \eqref{eq:h1der} and \eqref{eq:h2der}, %it holds that
\begin{align*}
h_{i,j}^{t-1}(\btheta)-h_{i,j}^{t-1}(\btheta_0)=&~\frac{\partial h_{i,j}^{t-1}(\btheta_0)}{\partial\btheta^{\T}}(\btheta-\btheta_0)+\frac{1}{2}(\btheta-\btheta_0)^{\T}\frac{\partial^2h_{i,j}^{t-1}(\btheta_0)}{\partial\btheta\partial\btheta^{\T}}(\btheta-\btheta_0)+R_{i,j}^{t-1}(\btheta)\\
=&-\frac{1}{2\gamma_{i,j}^{t-1}(\btheta_0)\{1-\gamma_{i,j}^{t-1}(\btheta_0)\}}\bigg|\frac{\partial\gamma_{i,j}^{t-1}(\btheta_0)}{\partial\btheta^{\T}}(\btheta_0-\btheta)\bigg|^2+R_{i,j}^{t-1}(\btheta)\,,
\end{align*}
where 
\begin{align*}
R_{i,j}^{t-1}(\btheta)=&~\frac{1}{2}\sum_{l_1=1}^q(\theta_{l_1}-\theta_{0,l_1})^3\int_0^1(1-v)^2\frac{\partial^3h_{i,j}^{t-1}(\btheta_v)}{\partial\theta_{l_1}^3}\,{\rm d}v\\
&+\frac{3}{2}\sum_{l_1\neq l_2}(\theta_{l_1}-\theta_{0,l_1})^2(\theta_{l_2}-\theta_{0,l_2})\int_0^1(1-v)^2\frac{\partial^3h_{i,j}^{t-1}(\btheta_v)}{\partial\theta_{l_1}^2\partial\theta_{l_2}}\,{\rm d}v\\
&+3\sum_{l_1\neq l_2\neq l_3}(\theta_{l_1}-\theta_{0,l_1})(\theta_{l_2}-\theta_{0,l_2})(\theta_{l_3}-\theta_{0,l_3})\int_0^1(1-v)^2\frac{\partial^3h_{i,j}^{t-1}(\btheta_v)}{\partial\theta_{l_1}\partial\theta_{l_2}\partial\theta_{l_3}}\,{\rm d}v
\end{align*}
with $\btheta_v=\btheta_0+v(\btheta-\btheta_0)$. Recall $\mathcal{I}_{i,j}=\{l\subset[q]:\gamma_{i,j}^{t-1}(\btheta)~\textrm{involves}~\theta_l~\textrm{for any}~t\in[n]\setminus[m]\}$. We have 
$\partial\gamma_{i,j}^{t-1}(\btheta)/\partial\theta_{l_1}\equiv0$ if $l_1\notin\mathcal{I}_{i,j}$, 
$\partial^2\gamma_{i,j}^{t-1}(\btheta)/\partial\theta_{l_1}\partial\theta_{l_2}\equiv0$ if $l_1\notin\mathcal{I}_{i,j}$ or $l_2\notin\mathcal{I}_{i,j}$, and $\partial^3\gamma_{i,j}^{t-1}(\btheta)/\partial\theta_{l_1}\partial\theta_{l_2}\partial\theta_{l_3}\equiv0$ if $l_1\notin\mathcal{I}_{i,j}$ or $l_2\notin\mathcal{I}_{i,j}$ or $l_3\notin\mathcal{I}_{i,j}$. By \eqref{eq:h3der}, it holds that
\[
\frac{\partial^3h_{i,j}^{t-1}(\btheta)}{\partial\theta_{l_1}\partial\theta_{l_2}\partial\theta_{l_3}}\equiv0~\textrm{if}~l_1\notin\mathcal{I}_{i,j}~\textrm{or}~ l_2\notin\mathcal{I}_{i,j}~\textrm{or}~l_3\notin\mathcal{I}_{i,j}\,.
\]
By Condition \ref{as:bdSij}, it holds that 
$
|R_{i,j}^{t-1}(\btheta)|\leq C_*s^3|\btheta_{\mathcal{I}_{i,j}}-\btheta_{0,\mathcal{I}_{i,j}}|_\infty^3$, 
which implies
\begin{align}
&\ell_{n,p}^{(l)}(\btheta_0)-\ell_{n,p}^{(l)}(\btheta)=\frac{1}{(n-m)|\mathcal{S}_l|}\sum_{t=m+1}^n\sum_{(i,j)\in\mathcal{S}_l}\big\{h_{i,j}^{t-1}(\btheta_0)-h_{i,j}^{t-1}(\btheta)\big\}\notag\\
&~~~~~\geq\frac{1}{(n-m)|\mathcal{S}_l|}\sum_{t=m+1}^n\sum_{(i,j)\in\mathcal{S}_l}\frac{1}{2\gamma_{i,j}^{t-1}(\btheta_0)\{1-\gamma_{i,j}^{t-1}(\btheta_0)\}}\bigg|\frac{\partial\gamma_{i,j}^{t-1}(\btheta_0)}{\partial\btheta^{\T}}(\btheta_0-\btheta)\bigg|^2\label{eq:diffbarell}\\
&~~~~~~~~~-\frac{C_*s^3}{|\mathcal{S}_l|}\sum_{(i,j)\in\mathcal{S}_l}|\btheta_{\mathcal{I}_{i,j}}-\btheta_{0,\mathcal{I}_{i,j}}|_\infty^3\notag\\
&~~~~~\geq\frac{2}{(n-m)|\mathcal{S}_l|}\sum_{t=m+1}^n\sum_{(i,j)\in\mathcal{S}_l}\bigg|\frac{\partial\gamma_{i,j}^{t-1}(\btheta_0)}{\partial\btheta^{\T}_{\mathcal{I}_{i,j}}}(\btheta_{0,\mathcal{I}_{i,j}}-\btheta_{\mathcal{I}_{i,j}})\bigg|^2-\frac{C_*s^3}{|\mathcal{S}_l|}\sum_{(i,j)\in\mathcal{S}_l}|\btheta_{\mathcal{I}_{i,j}}-\btheta_{0,\mathcal{I}_{i,j}}|_\infty^3\notag
\end{align}
for any $\btheta\in\bTheta$. By  Condition \ref{as:eigen}, it holds with probability approaching one that 
\[
\frac{1}{n-m}\sum_{t=m+1}^n\sum_{(i,j)\in\mathcal{S}_l}\bigg|\frac{\partial\gamma_{i,j}^{t-1}(\btheta_0)}{\partial\btheta^{\T}_{\mathcal{I}_{i,j}}}(\btheta_{0,\mathcal{I}_{i,j}}-\btheta_{\mathcal{I}_{i,j}})\bigg|^2\geq C_3\sum_{(i,j)\in\mathcal{S}_l}|\btheta_{0,\mathcal{I}_{i,j}}-\btheta_{\mathcal{I}_{i,j}}|_2^2
\]
for any $l\in[q]$, which implies
\begin{align*}
\ell_{n,p}^{(l)}(\btheta_0)-\ell_{n,p}^{(l)}(\btheta)\geq&~ \frac{2C_3}{|\mathcal{S}_l|}\sum_{(i,j)\in\mathcal{S}_l}|\btheta_{0,\mathcal{I}_{i,j}}-\btheta_{\mathcal{I}_{i,j}}|_2^2-\frac{C_*s^3}{|\mathcal{S}_l|}\sum_{(i,j)\in\mathcal{S}_l}|\btheta_{\mathcal{I}_{i,j}}-\btheta_{0,\mathcal{I}_{i,j}}|_\infty^3
\end{align*}
for any $\btheta\in\bTheta$ and $l\in[q]$ with probability approaching one. Since $d:=\sup_{\btheta\in\bTheta}|\btheta-\btheta_0|_\infty<2C_3/(C_*s^3)$, there exists a universal constant $\tilde{C}>0$ such that $d<\tilde{C}<2C_3/(C_*s^3)$. Hence, it holds with probability approaching one that
\begin{align*}
\ell_{n,p}^{(l)}(\btheta_0)-\ell_{n,p}^{(l)}(\btheta)\geq&~ \frac{2C_3-\tilde{C}C_*s^3}{|\mathcal{S}_l|}\sum_{(i,j)\in\mathcal{S}_l}|\btheta_{0,\mathcal{I}_{i,j}}-\btheta_{\mathcal{I}_{i,j}}|_2^2\\
&+\frac{\tilde{C}C_*s^3}{|\mathcal{S}_l|}\sum_{(i,j)\in\mathcal{S}_l}|\btheta_{0,\mathcal{I}_{i,j}}-\btheta_{\mathcal{I}_{i,j}}|_2^2-\frac{C_*s^3}{|\mathcal{S}_l|}\sum_{(i,j)\in\mathcal{S}_l}|\btheta_{\mathcal{I}_{i,j}}-\btheta_{0,\mathcal{I}_{i,j}}|_\infty^3\\
\geq&~\frac{2C_3-\tilde{C}C_*s^3}{|\mathcal{S}_l|}\sum_{(i,j)\in\mathcal{S}_l}|\btheta_{0,\mathcal{I}_{i,j}}-\btheta_{\mathcal{I}_{i,j}}|_2^2\\
&+\frac{\tilde{C}C_*s^3}{|\mathcal{S}_l|}\sum_{(i,j)\in\mathcal{S}_l}|\btheta_{0,\mathcal{I}_{i,j}}-\btheta_{\mathcal{I}_{i,j}}|_2^2-\frac{C_*s^3d}{|\mathcal{S}_l|}\sum_{(i,j)\in\mathcal{S}_l}|\btheta_{\mathcal{I}_{i,j}}-\btheta_{0,\mathcal{I}_{i,j}}|_\infty^2\\
\geq&~\frac{2C_3-\tilde{C}C_*s^3}{|\mathcal{S}_l|}\sum_{(i,j)\in\mathcal{S}_l}|\btheta_{0,\mathcal{I}_{i,j}}-\btheta_{\mathcal{I}_{i,j}}|_2^2
\end{align*}
for any $\btheta\in\bTheta$ and $l\in[q]$. 
We then have Proposition \ref{pn:ident} by selecting $\bar{C}=2C_3-\tilde{C}C_*s^3$. 
$\hfill\Box$

\subsection{Proof of Theorem \ref{tm:consistency}} \label{pr:tm:consistency}
Recall $S_{\mathcal{H},\min}=\min_{l\in\mathcal{H}}|\mathcal{S}_l|$ for any $\mathcal{H}\subset[q]$, and \begin{align*}
\left\{\begin{aligned}
   c_{n,\mathcal{G}}^2=&~\frac{q\log(n|\mathcal{S}_{l'}|)}{\sqrt{n|\mathcal{S}_{l'}|}}+\frac{q^{3/2}\log^{3/2}(n|\mathcal{S}_{l'}|)}{\sqrt{n}|\mathcal{S}_{l'}|}\,, \\
   c_{n,\mathcal{G}^{\btc}}^2=&~\frac{q\log(nS_{\mathcal{G}^{\btc},\min})}{\sqrt{nS_{\mathcal{G}^{\btc},\min}}}+\frac{q^{3/2}\log^{3/2}(nS_{\mathcal{G}^{\btc},\min})}{\sqrt{n}S_{\mathcal{G}^{\btc},\min}}\,,
\end{aligned}\right.
\end{align*}
where $l'\in\mathcal{G}$. 
To show Theorem \ref{tm:consistency}, we first present the following lemma whose proof is given in Section \ref{sec:pflauniformcov}. 

\begin{lemma}\label{la:uniformcov}
 Assume Conditions {\rm\ref{as:gambd}} and {\rm\ref{as:bdSij}} hold. Then  
\begin{align*}
\max_{l\in\mathcal{H}}\sup_{\btheta \in \bTheta}|\hat{\ell}_{n,p}^{(l)}(\btheta) - \ell_{n,p}^{(l)}(\btheta)|=O_{\rm p}\bigg\{\frac{q\log(nS_{\mathcal{H},\min})}{\sqrt{nS_{\mathcal{H},\min}}}\bigg\}+O_{\rm p}\bigg\{\frac{q^{3/2}\log^{3/2}(nS_{\mathcal{H},\min})}{\sqrt{n}S_{\mathcal{H},\min}}\bigg\}
\end{align*}
for any $\mathcal{H}\subset[q]$.
\end{lemma}

%\begin{lemma}\label{la:inirate}
%Assume Conditions {\rm\ref{as:gambd}}--{\rm\ref{as:eigen}} hold. Then $|\widetilde{\btheta}_{\mathcal{G}}-\btheta_{0,\mathcal{G}}|_2=O_{\rm p}(c_{n,\mathcal{G}})$ and $|\widetilde{\btheta}_{\mathcal{G}^{\btc}}-\btheta_{0,\mathcal{G}^{\btc}}|_\infty=O_{\rm p}(c_{n,\mathcal{G}^{\btc}})$. 
%\end{lemma}

Notice that $
\hat{\btheta}_*^{(l)}=(\hat \theta_{*,1}^{(l)}, \dots,\hat \theta_{*,q}^{(l)} )^{\T}=\arg\max_{\btheta\in\bTheta}\hat{\ell}_{n,p}^{(l)}(\btheta)$ for any $l\in[q]$. Then
\begin{align*}
{\ell}_{n,p}^{(l')}(\btheta_0)-\sup_{\btheta\in\bTheta}|\hat{\ell}_{n,p}^{(l')}(\btheta)-{\ell}_{n,p}^{(l')}(\btheta)|\leq&~\hat{\ell}_{n,p}^{(l')}(\btheta_0)\\
\leq&~\hat{\ell}_{n,p}^{(l')}(\hat{\btheta}_*^{(l')})\leq {\ell}_{n,p}^{(l')}(\hat{\btheta}_*^{(l')})+\sup_{\btheta\in\bTheta}|\hat{\ell}_{n,p}^{(l')}(\btheta)-{\ell}_{n,p}^{(l')}(\btheta)|\,,
\end{align*}
which implies
\begin{align*}%\label{eq:diffbd1}
0\leq {\ell}_{n,p}^{(l')}(\btheta_0)-{\ell}_{n,p}^{(l')}(\hat{\btheta}_*^{(l')})\leq 2\sup_{\btheta\in\bTheta}|\hat{\ell}_{n,p}^{(l')}(\btheta)-{\ell}_{n,p}^{(l')}(\btheta)|\,.
\end{align*}
Recall $\tilde{\btheta}_{\mathcal{G}}=\hat{\btheta}_{*,\mathcal{G}}^{(l')}$. Selecting $\mathcal{H}=\{l'\}$ in Lemma \ref{la:uniformcov}, we have $\sup_{\btheta\in\bTheta}|\hat{\ell}_{n,p}^{(l')}(\btheta)-{\ell}_{n,p}^{(l')}(\btheta)|=O_{\rm p}(c_{n,\mathcal{G}}^2)$, which implies 
\begin{align}\label{eq:diffbd11}
{\ell}_{n,p}^{(l')}(\btheta_0)-{\ell}_{n,p}^{(l')}(\hat{\btheta}_*^{(l')})=O_{\rm p}(c_{n,\mathcal{G}}^2)\,.
\end{align}
For any diverging $\epsilon_{n,p}>0$, if $|\tilde{\btheta}_{\mathcal{G}}-\btheta_{0,\mathcal{G}}|_2\geq \epsilon_{n,p}c_{n,\mathcal{G}}$, Proposition \ref{pn:ident} yields that
\[
{\ell}_{n,p}^{(l')}(\btheta_0)-{\ell}_{n,p}^{(l')}(\hat{\btheta}_*^{(l)})\geq\bar{C}\epsilon_{n,p}^2c_{n,\mathcal{G}}^2 
\]
with probability approaching one, which contradicts with \eqref{eq:diffbd11} and then implies $|\tilde{\btheta}_{\mathcal{G}}-\btheta_{0,\mathcal{G}}|_2=O_{\rm p}( \epsilon_{n,p}c_{n,\mathcal{G}})$. Notice that we can select arbitrary slowly diverging $\epsilon_{n,p}$. Following a standard result from
probability theory, we have $|\tilde{\btheta}_{\mathcal{G}}-\btheta_{0,\mathcal{G}}|_2=O_{\rm p}(c_{n,\mathcal{G}})$. 

For any $l\in\mathcal{G}^{\btc}$, %since $
%\widehat{\btheta}_*^{(l)}=(\hat \theta_{*,1}^{(l)}, \dots,\hat \theta_{*,q}^{(l)} )^{\T}=\arg\max_{\btheta\in\bTheta}\hat{\ell}_{n,p}^{(l)}(\btheta)$, then
%\begin{align*}
%{\ell}_{n,p}^{(l)}(\btheta_0)-\sup_{\btheta\in\bTheta}|\hat{\ell}_{n,p}^{(l)}(\btheta)-{\ell}_{n,p}^{(l)}(\btheta)|\leq&~\hat{\ell}_{n,p}^{(l)}(\btheta_0)\\
%\leq&~\hat{\ell}_{n,p}^{(l)}(\widehat{\btheta}_*^{(l)})\leq {\ell}_{n,p}^{(l)}(\widehat{\btheta}_*^{(l)})+\sup_{\btheta\in\bTheta}|\hat{\ell}_{n,p}^{(l)}(\btheta)-{\ell}_{n,p}^{(l)}(\btheta)|\,,
%\end{align*}
%which implies
we also have 
\begin{align*}%\label{eq:diffbd1}
0\leq {\ell}_{n,p}^{(l)}(\btheta_0)-{\ell}_{n,p}^{(l)}(\hat{\btheta}_*^{(l)})\leq 2\sup_{\btheta\in\bTheta}|\hat{\ell}_{n,p}^{(l)}(\btheta)-{\ell}_{n,p}^{(l)}(\btheta)|\,.
\end{align*}
Recall $\tilde{\btheta}_{\mathcal{G}^{\btc}}=(\hat{\theta}_{*,l}^{(l)})_{l\in\mathcal{G}^{\btc}}$. Selecting $\mathcal{H}=\mathcal{G}^{\btc}$ in Lemma \ref{la:uniformcov}, we have $\max_{l\in\mathcal{G}^{\btc}}\sup_{\btheta\in\bTheta}|\hat{\ell}_{n,p}^{(l)}(\btheta)-{\ell}_{n,p}^{(l)}(\btheta)|=O_{\rm p}(c_{n,\mathcal{G}^{\btc}}^2)$, which implies 
\begin{align}\label{eq:diffbd1}
\max_{l\in\mathcal{G}^{\btc}}\big\{{\ell}_{n,p}^{(l)}(\btheta_0)-{\ell}_{n,p}^{(l)}(\hat{\btheta}_*^{(l)})\big\}=O_{\rm p}(c_{n,\mathcal{G}^{\btc}}^2)\,.
\end{align}
For any diverging $\epsilon_{n,p}>0$, if $|\tilde{\btheta}_{\mathcal{G}^{\btc}}-\btheta_{0,\mathcal{G}^{\btc}}|_\infty\geq \epsilon_{n,p}c_{n,\mathcal{G}^{\btc}}$, Proposition \ref{pn:ident} yields that
\[
\max_{l\in\mathcal{G}^{\btc}}\big\{{\ell}_{n,p}^{(l)}(\btheta_0)-{\ell}_{n,p}^{(l)}(\hat{\btheta}_*^{(l)})\big\}\geq\bar{C}\epsilon_{n,p}^2c_{n,\mathcal{G}^{\btc}}^2 
\]
with probability approaching one, which contradicts with \eqref{eq:diffbd1} and then implies $|\tilde{\btheta}_{\mathcal{G}^{\btc}}-\btheta_{0,\mathcal{G}^{\btc}}|_\infty=O_{\rm p}( \epsilon_{n,p}c_{n,\mathcal{G}^{\btc}})$. Notice that we can select arbitrary slowly diverging $\epsilon_{n,p}$. Following a standard result from
probability theory, we have $|\tilde{\btheta}_{\mathcal{G}^{\btc}}-\btheta_{0,\mathcal{G}^{\btc}}|_\infty=O_{\rm p}(c_{n,\mathcal{G}^{\btc}})$. We complete the proof of Theorem \ref{tm:consistency}.  $\hfill\Box$

%{\rm (ii)} For any $\mathcal{H}\subset[q]$, if $n$ is fixed but $S_{\mathcal{H},\min}\rightarrow\infty$, then 
%\begin{align*}
%\max_{l\in\mathcal{H}}\sup_{\btheta \in \bTheta}|\hat{\ell}_{n,p}^{(l)}(\btheta) - \ell_{n,p}^{(l)}(\btheta)|=O_{\rm p}\bigg(\sqrt{\frac{q\log S_{\mathcal{H},\min}}{S_{\mathcal{H},\min}}}\bigg)+O_{\rm p}\bigg(\frac{q\log S_{\mathcal{H},\min}}{S_{\mathcal{H},\min}}\bigg)\,.
%\end{align*}

%\begin{lemma}\label{la:uniformcov:local}
%Assume Conditions {\rm\ref{as:gambd}} and {\rm\ref{as:bdSij}} hold. 

%{\rm (i)} If $n\rightarrow\infty$, then 
%\[
%\sup_{\btheta \in \bTheta}|\ell_{n,p}^{(l)}(\btheta) - \bar{\ell}_{n,p}^{(l)}(\btheta)|=O_{\rm p}\bigg\{\frac{q\log(|\mathcal{S}_l|n)}{\sqrt{n|\mathcal{S}_l|}}\bigg\}+O_{\rm p}\bigg\{\frac{q^{3/2}\log^{3/2}(|\mathcal{S}_l|n)}{|\mathcal{S}_l|\sqrt{n}}\bigg\}\,.
%\]

%{\rm (ii)} If $n$ is fixed but $p\rightarrow\infty$, then 
%\[
%\sup_{\btheta \in \bTheta}|\ell_{n,p}^{(l)}(\btheta) - \bar{\ell}_{n,p}^{(l)}(\btheta)|=O_{\rm p}\bigg(\sqrt{\frac{q\log |\mathcal{S}_l|}{|\mathcal{S}_l|}}\bigg)+O_{\rm p}\bigg(\frac{q\log |\mathcal{S}_l|}{|\mathcal{S}_l|}\bigg)\,.
%\]
%\end{lemma}

\subsection{Proof of Proposition \ref{pn:reconv}} \label{pr:pn:reconv}

Given $\bvarphi_l$ specified in Condition \ref{as:al}, define
\[
f_t^{(l)}(\btheta)=\bvarphi_l^{\T}\bg_t^{(l)}(\btheta)
\]
for any $t\in[n]\setminus[m]$ and $\btheta\in\bTheta$. 
To show Proposition \ref{tm:asymnormal}, we need the following lemmas whose proofs are given in Sections \ref{pf:laasymp1}--\ref{pf:laasymp3}.

\begin{lemma}\label{la:asym1}
Assume Condition {\rm\ref{as:gambd}} holds. Then
\begin{align*}
\bigg|\frac{1}{n-m}\sum_{t=m+1}^n\big\{\hat{f}_t^{(l)}(\theta_{0,l},\tilde{\btheta}_{-l})-\hat{f}_t^{(l)}(\btheta_0)\big\}\bigg|\leq \tau |\tilde{\btheta}-\btheta_{0}|_1+C_*|\hat{\bvarphi}_l|_1 |\tilde{\btheta}-\btheta_{0}|_1^2
\end{align*}
for any $l\in[q]$, where $C_*$ is a universal constant specified in \eqref{eq:h3bd}. 
\end{lemma}

\begin{lemma}\label{la:asym2}
Assume Conditions {\rm\ref{as:gambd}} and {\rm\ref{as:al}} hold. Then
\begin{align*}
&\max_{l\in\mathcal{H}}\bigg|\frac{1}{n-m}\sum_{t=m+1}^n\hat{f}_t^{(l)}(\btheta_0)\bigg|\\
&~~~~~~~=O_{\rm p}\bigg\{\frac{\sqrt{\log(1+ |\mathcal{H}|)\log(n|\mathcal{H}|)}}{\sqrt{nS_{\mathcal{H},\min}}}\bigg\}+O_{\rm p}\bigg\{\frac{\sqrt{\log(1+ |\mathcal{H}|)}\log(n|\mathcal{H}|)}{\sqrt{n}S_{\mathcal{H},\min}}\bigg\}
\end{align*}
for any $\mathcal{H}\subset[q]$. 
\end{lemma}

\begin{lemma}\label{la:asym3}
Assume Condition {\rm\ref{as:gambd}} holds. Then
\[
\sup_{\theta_l\in B(\tilde{\theta}_l,\tilde{r})}\bigg|\frac{1}{n-m}\sum_{t=m+1}^n\frac{\hat{f}_t^{(l)}(\theta_l,\tilde{\btheta}_{-l})}{\partial\theta_l}-1\bigg|\leq \tau+C_*\tilde{r}|\hat{\bvarphi}_l|_1
\]
for any $l\in[q]$, where $C_*$ is a universal constant specified in \eqref{eq:h3bd}. 
\end{lemma}

Recall $|\tilde{\btheta}_{\mathcal{G}}-\btheta_{0,\mathcal{G}}|_2=O_{\rm p}(c_{n,\mathcal{G}})$ and $|\tilde{\btheta}_{\mathcal{G}^{\btc}}-\btheta_{0,\mathcal{G}^{\btc}}|_\infty=O_{\rm p}(c_{n,\mathcal{G}^{\btc}})$ with \begin{align*}
\left\{\begin{aligned}
   c_{n,\mathcal{G}}^2=&~\frac{q\log(n|\mathcal{S}_{l'}|)}{\sqrt{n|\mathcal{S}_{l'}|}}+\frac{q^{3/2}\log^{3/2}(n|\mathcal{S}_{l'}|)}{\sqrt{n}|\mathcal{S}_{l'}|}\,, \\
   c_{n,\mathcal{G}^{\btc}}^2=&~\frac{q\log(nS_{\mathcal{G}^{\btc},\min})}{\sqrt{nS_{\mathcal{G}^{\btc},\min}}}+\frac{q^{3/2}\log^{3/2}(nS_{\mathcal{G}^{\btc},\min})}{\sqrt{n}S_{\mathcal{G}^{\btc},\min}}\,,
\end{aligned}\right.
\end{align*}
where $l'\in\mathcal{G}$. Then $|\tilde{\btheta}-\btheta_0|_1=|\tilde{\btheta}_{\mathcal{G}}-\btheta_{0,\mathcal{G}}|_1+|\tilde{\btheta}_{\mathcal{G}^{\btc}}-\btheta_{0,\mathcal{G}^{\btc}}|_1=O_{\rm p}(|\mathcal{G}|^{1/2}c_{n,\mathcal{G}})+O_{\rm p}(|\mathcal{G}^{\btc}|c_{n,\mathcal{G}^{\btc}})=O_{\rm p}(\Delta_n^{1/2})$. Based on Lemmas \ref{la:asym1} and \ref{la:asym2}, due to $\tau\lesssim \Delta_n^{1/2}$ and $|\mathcal{G}|+|\mathcal{G}^{\btc}|=q$, we have
\begin{align*}
&\max_{l\in\mathcal{G}}\bigg|\frac{1}{n-m}\sum_{t=m+1}^n\hat{f}_t^{(l)}(\theta_{0,l},\tilde{\btheta}_{-l})\bigg|\\
&~~~~~~~~~=O_{\rm p}(\Delta_n)+O_{\rm p}\bigg\{\frac{\sqrt{\log(1+|\mathcal{G}|)\log(n|\mathcal{G}|)}}{\sqrt{n|\mathcal{S}_{l'}|}}\bigg\}+O_{\rm p}\bigg\{\frac{\sqrt{\log(1+|\mathcal{G}|)}\log(n|\mathcal{G}|)}{\sqrt{n}|\mathcal{S}_{l'}|}\bigg\}\\
&~~~~~~~~~=O_{\rm p}(\Delta_n)
\end{align*}
for some $l'\in\mathcal{G}$, and 
\begin{align*}
&\max_{l\in\mathcal{G}^{\btc}}\bigg|\frac{1}{n-m}\sum_{t=m+1}^n\hat{f}_t^{(l)}(\theta_{0,l},\tilde{\btheta}_{-l})\bigg|\\
&~~~~~~~~~=O_{\rm p}(\Delta_n)+O_{\rm p}\bigg\{\frac{\sqrt{\log(1+|\mathcal{G}^{\btc}|)\log(n|\mathcal{G}^{\btc}|)}}{\sqrt{nS_{\mathcal{G}^{\btc},\min}}}\bigg\}+O_{\rm p}\bigg\{\frac{\sqrt{\log(1+|\mathcal{G}^{\btc}|)}\log(n|\mathcal{G}^{\btc}|)}{\sqrt{n}S_{\mathcal{G}^{\btc},\min}}\bigg\}\\
&~~~~~~~~~=O_{\rm p}(\Delta_n)\,.
\end{align*}
Due to $\check{\theta}_l=\arg\min_{\theta_l\in B(\tilde{\theta}_l,\tilde{r})}|(n-m)^{-1}\sum_{t=m+1}^n\hat{f}_t^{(l)}(\theta_l,\tilde{\btheta}_{-l})|^2$ and $\theta_{0,l}\in B(\tilde{\theta}_l,\tilde{r})$ with probability approaching one, then
\[
\bigg|\frac{1}{n-m}\sum_{t=m+1}^n\hat{f}_t^{(l)}(\check{\theta}_{l},\tilde{\btheta}_{-l})\bigg|\leq \bigg|\frac{1}{n-m}\sum_{t=m+1}^n\hat{f}_t^{(l)}(\theta_{0,l},\tilde{\btheta}_{-l})\bigg|
\]
with probability approaching one, 
which implies
\begin{align}
\max_{l\in[q]}\bigg|\frac{1}{n-m}\sum_{t=m+1}^n\big\{\hat{f}_t^{(l)}(\check{\theta}_{l},\tilde{\btheta}_{-l})-\hat{f}_t^{(l)}(\theta_{0,l},\tilde{\btheta}_{-l})\big\}\bigg|=O_{\rm p}(\Delta_n)\,.\label{eq:asymbd5}
\end{align}
Notice that
\begin{align*}%\label{eq:asymbd6}
\frac{1}{n-m}\sum_{t=m+1}^n\big\{\hat{f}_t^{(l)}(\check{\theta}_{l},\tilde{\btheta}_{-l})-\hat{f}_t^{(l)}(\theta_{0,l},\tilde{\btheta}_{-l})\big\}=\bigg\{\frac{1}{n-m}\sum_{t=m+1}^n\frac{\partial\hat{f}_t^{(l)}(\bar{\theta}_l,\tilde{\btheta}_{-l})}{\partial\theta_l}\bigg\}(\check{\theta}_l-\theta_{0,l})
\end{align*}
for some $\bar{\theta}_l$ on the joint line between $\check{\theta}_l$ and $\theta_{0,l}$. 
By Lemma \ref{la:asym3} and Condition \ref{as:al}, due to $\tau=o(1)$ and $\tilde{r}=o(1)$, we know
\[
\min_{l\in[q]}\bigg\{\frac{1}{n-m}\sum_{t=m+1}^n\frac{\partial\hat{f}_t^{(l)}(\bar{\theta}_l,\tilde{\btheta}_{-l})}{\partial\theta_l}\bigg\}>\frac{1}{2}
\]
with probability approaching one. Hence, by \eqref{eq:asymbd5}, we have $|\check{\btheta}-\btheta_{0}|_\infty=O_{\rm p}(\Delta_n)$. We complete the proof of Proposition \ref{pn:reconv}. $\hfill\Box$

\subsection{Proof of Theorem \ref{tm:asymnormal}} \label{pr:tm:asymnormal}
Recall $\hat{f}_t^{(l)}(\btheta)=\hat{\bvarphi}_{l}^{\T}\bg_t^{(l)}(\btheta)$. By the definition of $\hat{\bvarphi}_{l}$ given in \eqref{eq:hatanl}, we have 
\[
\max_{l\in[q]}\bigg|\frac{1}{n-m}\sum_{t=m+1}^n\frac{\partial\hat{f}_t^{(l)}(\tilde{\btheta})}{\partial\btheta_{-l}}\bigg|_\infty\leq \tau\,.
\]
It follows from the Taylor expansion that
\begin{align*}
&\bigg|\frac{1}{n-m}\sum_{t=m+1}^n\bigg\{\frac{\partial \hat{f}_t^{(l)}(\check{\btheta})}{\partial\btheta_{-l}}-\frac{\partial\hat{f}_t^{(l)}(\tilde{\btheta})}{\partial\btheta_{-l}}\bigg\}\bigg|_\infty\leq \bigg|\frac{1}{n-m}\sum_{t=m+1}^n\frac{\partial^2\hat{f}_t^{(l)}(\dot{\btheta})}{\partial\btheta_{-l}\partial\btheta^{\T}}\bigg|_\infty|\check{\btheta}-\tilde{\btheta}|_1\,,
\end{align*}
where $\dot{\btheta}$ is on the joint line between $\check{\btheta}$ and $\tilde{\btheta}$. %By Theorem \ref{tm:consistency} and Proposition \ref{pn:reconv}, we have 
%\begin{align*}
%|\widecheck{\btheta}-\widetilde{\btheta}|_1\leq&~ |\widecheck{\btheta}-\btheta_0|_1+|\widetilde{\btheta}-\btheta_0|_1\\
%=&~O_{\rm p}(q|\mathcal{G}|c_{n,\mathcal{G}}^2)+O_{\rm p}(q|\mathcal{G}^{\btc}|^2c_{n,\mathcal{G}^{\btc}}^2)+O_{\rm p}(|\mathcal{G}|^{1/2}c_{n,\mathcal{G}})+O_{\rm p}(|\mathcal{G}^{\btc}|c_{n,\mathcal{G}^{\btc}})\\
%=&~o_{\rm p}(1)\,.
%\end{align*}
Notice that $\dot{\btheta}\in B(\tilde{\theta}_1,\tilde{r})\times\cdots\times B(\tilde{\theta}_q,\tilde{r})$. 
Following the same arguments for deriving \eqref{eq:asybd2} in Section \ref{pf:laasymp1} for the proof of Lemma \ref{la:asym1}, we know 
\[
\sup_{\btheta\in B(\tilde{\theta}_1,\tilde{r})\times\cdots\times B(\tilde{\theta}_q,\tilde{r})}\max_{l\in[q]}\bigg|\frac{1}{n-m}\sum_{t=m+1}^n\frac{\partial^2\hat{f}_t^{(l)}({\btheta})}{\partial\btheta_{-l}\partial\btheta^{\T}}\bigg|_\infty\leq C_*\max_{l\in[q]}|\hat{\bvarphi}_l|_1
\]
for some universal constant $C_*$ specified in \eqref{eq:h3bd}, which implies
\begin{align*}
\max_{l\in[q]}\bigg|\frac{1}{n-m}\sum_{t=m+1}^n\bigg\{\frac{\partial \hat{f}_t^{(l)}(\check{\btheta})}{\partial\btheta_{-l}}-\frac{\partial\hat{f}_t^{(l)}(\tilde{\btheta})}{\partial\btheta_{-l}}\bigg\}\bigg|_\infty\leq C_*|\check{\btheta}-\tilde{\btheta}|_1\max_{l\in[q]}|\hat{\bvarphi}_l|_1\,.
\end{align*}
By Condition \ref{as:al}, we have $\max_{l\in[q]}|\hat{\bvarphi}_l|_1=O_{\rm p}(1)$, which implies 
\begin{align*}
\max_{l\in[q]}\bigg|\frac{1}{n-m}\sum_{t=m+1}^n\bigg\{\frac{\partial \hat{f}_t^{(l)}(\check{\btheta})}{\partial\btheta_{-l}}-\frac{\partial\hat{f}_t^{(l)}(\tilde{\btheta})}{\partial\btheta_{-l}}\bigg\}\bigg|_\infty\leq|\check{\btheta}-\tilde{\btheta}|_1\cdot O_{\rm p}(1)\,.
\end{align*}
Hence, it holds that 
\begin{align}\label{eq:asymbd10}
\max_{l\in[q]}\bigg|\frac{1}{n-m}\sum_{t=m+1}^n\frac{\partial \hat{f}_t^{(l)}(\check{\btheta})}{\partial\btheta_{-l}}\bigg|_\infty\leq \tau+ |\check{\btheta}-\tilde{\btheta}|_1\cdot O_{\rm p}(1)\,.
\end{align}
Repeating the arguments for deriving Lemma \ref{la:asym1} in Section \ref{pf:laasymp1}, we can also show
\begin{align}
&\max_{l\in[q]}\bigg|\frac{1}{n-m}\sum_{t=m+1}^n\big\{\hat{f}_t^{(l)}(\theta_{0,l},\check{\btheta}_{-l})-\hat{f}_t^{(l)}(\btheta_0)\big\}\bigg|\notag\\
&~~~~~~~~~~~~~~~~~\leq \big\{\tau+|\check{\btheta}-\tilde{\btheta}|_1\cdot O_{\rm p}(1)\big\} |\check{\btheta}-\btheta_{0}|_1+ |\check{\btheta}-\btheta_{0}|_1^2\cdot O_{\rm p}(1)\label{eq:asym15}\\
&~~~~~~~~~~~~~~~~~\leq \big\{\tau+|\tilde{\btheta}-\btheta_0|_1\cdot O_{\rm p}(1)\big\} |\check{\btheta}-\btheta_{0}|_1+ |\check{\btheta}-\btheta_{0}|_1^2\cdot O_{\rm p}(1)\notag\,.
\end{align}
Together with Lemma \ref{la:asym2}, it holds that
\[
\max_{l\in[q]}\bigg|\frac{1}{n-m}\sum_{t=m+1}^n\hat{f}_t^{(l)}(\theta_{0,l},\check{\btheta}_{-l})\bigg|\leq O_{\rm p}(q\Delta_n)\,.
\]
 Since $\hat{\theta}_l=\arg\min_{\theta_l\in B(\check{\theta}_l,\check{r})}|(n-m)^{-1}\sum_{t=m+1}^n\hat{f}_t^{(l)}(\theta_l,\check{\btheta}_{-l})|^2$ and $\theta_{0,l}\in B(\check{\theta}_l,\check{r})$ with probability approaching one, then
 \[
\bigg|\frac{1}{n-m}\sum_{t=m+1}^n\hat{f}_t^{(l)}(\hat{\theta}_{l},\check{\btheta}_{-l})\bigg|\leq \bigg|\frac{1}{n-m}\sum_{t=m+1}^n\hat{f}_t^{(l)}(\theta_{0,l},\check{\btheta}_{-l})\bigg|
\]
with probability approaching one, which implies 
\begin{align}\label{eq:asym11}
\max_{l\in[q]}\bigg|\frac{1}{n-m}\sum_{t=m+1}^n\big\{\hat{f}_t^{(l)}(\hat{\theta}_{l},\check{\btheta}_{-l})-\hat{f}_t^{(l)}(\theta_{0,l},\check{\btheta}_{-l})\big\}\bigg|=O_{\rm p}(q\Delta_n)\,.
\end{align}
Following the same arguments for deriving Lemma \ref{la:asym3} in Section \ref{pf:laasymp3} and noting \eqref{eq:asymbd10}, we can also have
\begin{align}\label{eq:asym12}
\sup_{\theta_l\in B(\check{\theta}_l,\check{r})}\bigg|\frac{1}{n-m}\sum_{t=m+1}^n\frac{\hat{f}_t^{(l)}(\theta_l,\check{\btheta}_{-l})}{\partial\theta_l}-1\bigg|\leq \tau+o_{\rm p}(1)+C_*\check{r}|\hat{\bvarphi}_l|_1
\end{align}
for any $l\in[q]$. Due to $\tau=o(1)$, $\check{r}=o(1)$ and $\max_{l\in[q]}|\hat{\bvarphi}_l|_1=O_{\rm p}(1)$, we have
\begin{align}\label{eq:asym13}
\max_{l\in[q]}\sup_{\theta_l\in B(\check{\theta}_l,\check{r})}\bigg|\frac{1}{n-m}\sum_{t=m+1}^n\frac{\hat{f}_t^{(l)}(\theta_l,\check{\btheta}_{-l})}{\partial\theta_l}-1\bigg|=o_{\rm p}(1)\,.
\end{align}
Due to
\begin{align}\label{eq:asymbd14}
\frac{1}{n-m}\sum_{t=m+1}^n\big\{\hat{f}_t^{(l)}(\hat{\theta}_{l},\check{\btheta}_{-l})-\hat{f}_t^{(l)}(\theta_{0,l},\check{\btheta}_{-l})\big\}=\bigg\{\frac{1}{n-m}\sum_{t=m+1}^n\frac{\partial\hat{f}_t^{(l)}(\bar{\theta}_l,\check{\btheta}_{-l})}{\partial\theta_l}\bigg\}(\hat{\theta}_l-\theta_{0,l})
\end{align}
for some $\bar{\theta}_l$ on the joint line between $\hat{\theta}_l$ and $\theta_{0,l}$, by \eqref{eq:asym11} and \eqref{eq:asym12}, we have $|\hat{\btheta}-\btheta_0|_\infty=O_{\rm p}(q\Delta_n)$, which implies $|\hat{\btheta}-\check{\btheta}|_\infty\leq |\hat{\btheta}-\btheta_0|_\infty+|\check{\btheta}-\btheta_0|_\infty=O_{\rm p}(q\Delta_n)$. Hence, $|\hat{\btheta}-\check{\btheta}|_\infty\ll\check{r}$ with probability approaching one. Therefore, 
\[
0=\bigg\{\frac{1}{n-m}\sum_{t=m+1}^n\hat{f}_t^{(l)}(\hat{\theta}_l,\check{\btheta}_{-l})\bigg\}\bigg\{\frac{1}{n-m}\sum_{t=m+1}^n\frac{\partial\hat{f}_t^{(l)}(\hat{\theta}_l,\check{\btheta}_{-l})}{\partial\theta_l}\bigg\}
\]
with probability approaching one. Together with \eqref{eq:asym13}, we have
\[
\frac{1}{n-m}\sum_{t=m+1}^n\hat{f}_t^{(l)}(\hat{\theta}_l,\check{\btheta}_{-l})=0
\]
with probability approaching one.

By \eqref{eq:asymbd14} and \eqref{eq:asym15}, due to $|\check{\btheta}-\btheta_{0}|_1=O_{\rm p}(q\Delta_n)$, $|\tilde{\btheta}-\btheta_0|_1=O_{\rm p}(\Delta_n^{1/2})$ and $\tau\lesssim \Delta_n^{1/2}$, then
\begin{align*}
\hat{\theta}_l-\theta_{0,l}=&-\bigg\{\frac{1}{n-m}\sum_{t=m+1}^n\frac{\partial\hat{f}_t^{(l)}(\bar{\theta}_l,\check{\btheta}_{-l})}{\partial\theta_l}\bigg\}^{-1}\bigg\{\frac{1}{n-m}\sum_{t=m+1}^n\hat{f}_t^{(l)}(\theta_{0,l},\check{\btheta}_{-l})\bigg\}\\
=&-\bigg\{\frac{1}{n-m}\sum_{t=m+1}^n\frac{\partial\hat{f}_t^{(l)}(\bar{\theta}_l,\check{\btheta}_{-l})}{\partial\theta_l}\bigg\}^{-1}\bigg\{\frac{1}{n-m}\sum_{t=m+1}^n\hat{f}_t^{(l)}(\btheta_{0})\bigg\}\\
&+O_{\rm p}(q\Delta_n^{3/2})+O_{\rm p}(q^2\Delta_n^2)
\end{align*}
with probability approaching one. As shown in \eqref{eq:asymbd4} in Section \ref{pf:laasymp2} for the proof of Lemma \ref{la:asym2}, 
\begin{align*}
&\bigg|\frac{1}{n-m}\sum_{t=m+1}^n\big\{\hat{f}_t^{(l)}(\btheta_0)-f_t^{(l)}(\btheta_0)\big\}\bigg|\\
&~~~~~~~~~=O_{\rm p}\bigg\{\frac{\omega_n\sqrt{(\log q)\log(qn)}}{\sqrt{n|\mathcal{S}_{l}|}}\bigg\}+O_{\rm p}\bigg\{\frac{\omega_n(\log q)^{1/2}\log(qn)}{\sqrt{n}|\mathcal{S}_{l}|}\bigg\}=o_{\rm p}\bigg(\frac{1}{\sqrt{n|\mathcal{S}_l|}}\bigg)\,,
\end{align*}
where the last step is based on Condition \ref{as:al}. 
Hence, \begin{align}
\hat{\theta}_l-\theta_{0,l}
=&-\bigg\{\frac{1}{n-m}\sum_{t=m+1}^n\frac{\partial\hat{f}_t^{(l)}(\bar{\theta}_l,\check{\btheta}_{-l})}{\partial\theta_l}\bigg\}^{-1}\bigg\{\frac{1}{n-m}\sum_{t=m+1}^nf_t^{(l)}(\btheta_{0})\bigg\}\notag\\
&+O_{\rm p}(q\Delta_n^{3/2})+O_{\rm p}(q^2\Delta_n^2)+o_{\rm p}\bigg(\frac{1}{\sqrt{n|\mathcal{S}_l|}}\bigg)\label{eq:asymbd8}
\end{align}
with probability approaching one. 

Write 
\[
\mathring{Q}_{l,t}=\frac{1}{|\mathcal{S}_l|^{1/2}}\sum_{(i,j)\in\mathcal{S}_l}\frac{X_{i,j}^t-\gamma_{i,j}^{t-1}(\btheta_0)}{\gamma_{i,j}^{t-1}(\btheta_0)\{1-\gamma_{i,j}^{t-1}(\btheta_0)\}}\bigg\{\bvarphi_l^{\T}\frac{\partial\gamma_{i,j}^{t-1}(\btheta_0)}{\partial\btheta}\bigg\}\,.
\]
Then
\begin{align}
\frac{1}{n-m}\sum_{t=m+1}^nf_t^{(l)}(\btheta_0)
=\frac{1}{(n-m)|\mathcal{S}_l|^{1/2}}\sum_{t=m+1}^n\mathring{Q}_{l,t}\,.\label{eq:asymbd9}
\end{align}
In the sequel, we will use the martingale central limit theorem to establish the asymptotic distribution of $(n-m)^{-1/2}\sum_{t=m+1}^n\mathring{Q}_{l,t}$. 
Denote by $\mathbb{P}_{\mathcal{F}_{t-1}}(\cdot)$ and $\mathbb{E}_{\mathcal{F}_{t-1}}(\cdot)$, respectively, the conditional probability measure and the conditional expectation given $\mathcal{F}_{t-1}$ with the unknown true parameter vector $\btheta_0$. By Conditions \ref{as:gambd} and \ref{as:al}, it holds that
\[
\max_{(i,j)\in\mathcal{S}_l}\bigg|\frac{X_{i,j}^t-\gamma_{i,j}^{t-1}(\btheta_0)}{\gamma_{i,j}^{t-1}(\btheta_0)\{1-\gamma_{i,j}^{t-1}(\btheta_0)\}}\bigg\{\bvarphi_l^{\T}\frac{\partial\gamma_{i,j}^{t-1}(\btheta_0)}{\partial\btheta}\bigg\}\bigg|\leq C_1^{-1}C_2C_4=:C_{**}\,.
\] 
%Notice that 
%\[
%\mathring{Q}_{l,t}=\frac{2}{|\mathcal{S}_l|^{1/2}}\sum_{(i,j)\in\mathcal{S}_l:\,i<j}\frac{X_{i,j}^t-\gamma_{i,j}^{t-1}(\btheta_0)}{\gamma_{i,j}^{t-1}(\btheta_0)\{1-\gamma_{i,j}^{t-1}(\btheta_0)\}}\bigg\{\bvarphi_l^{\T}\frac{\partial\gamma_{i,j}^{t-1}(\btheta_0)}{\partial\btheta}\bigg\}\,.
%\]
It follows from the Bernstein inequality that
\begin{align*}
\mathbb{P}_{\mathcal{F}_{t-1}}\big(|\mathring{Q}_{l,t}|\geq x\big)\leq 2\exp\bigg(-\frac{3|\mathcal{S}_l|^{1/2}x^2}{6|\mathcal{S}_l|^{1/2}C_{**}^2+2xC_{**}}\bigg)
\end{align*}
for any $x>0$, %Furthermore, due to $\mathbb{P}(|\mathring{Q}_{l,t}|\geq x)=\mathbb{E}\{\mathbb{P}_{\mathcal{F}_{t-1}}(|\mathring{Q}_{l,t}|\geq x)\}$, we have
%\begin{align*}
%\mathbb{P}\big(|\mathring{Q}_{l,t}|\geq x\big)\leq 2\exp\bigg(-\frac{3|\mathcal{S}_l|^{1/2}x^2}{6|\mathcal{S}_l|^{1/2}C_{**}^2+2xC_{**}}\bigg)
%\end{align*}
%for any $x>0$. 
which implies, for any $\delta>0$,
\begin{align*}
&\sum_{t=m+1}^n\mathbb{E}_{\mathcal{F}_{t-1}}\bigg\{\bigg(\frac{\mathring{Q}_{l,t}}{\sqrt{n-m}}\bigg)^2I\bigg(\frac{|\mathring{Q}_{l,t}|}{\sqrt{n-m}}\geq\delta\bigg)\bigg\}\notag\\
&~~~~~~~=\frac{1}{n-m}\sum_{t=m+1}^n\mathbb{E}_{\mathcal{F}_{t-1}}\big\{\mathring{Q}_{l,t}^2I(|\mathring{Q}_{l,t}|\geq\delta\sqrt{n-m})\big\}\rightarrow0\label{eq:asymbd7}
\end{align*}
in probability as $n\rightarrow\infty$. Meanwhile, by Condition \ref{as:conv}, we also have
\begin{align*}
&\frac{1}{n-m}\sum_{t=m+1}^n\mathbb{E}_{\mathcal{F}_{t-1}}(\mathring{Q}_{l,t}^2)\\
&~~~~~~~~~~~=\frac{1}{(n-m)|\mathcal{S}_l|}\sum_{t=m+1}^n\sum_{(i,j)\in\mathcal{S}_l}\frac{1}{\gamma_{i,j}^{t-1}(\btheta_0)\{1-\gamma_{i,j}^{t-1}(\btheta_0)\}}\bigg\{\bvarphi_l^{\T}\frac{\partial\gamma_{i,j}^{t-1}(\btheta_0)}{\partial\btheta}\bigg\}^2\rightarrow\kappa_l
\end{align*}
% we know
%\[
%\frac{1}{n-m}\sum_{t=m+1}^n\mathbb{E}_{\mathcal{F}_{t-1}}(\mathring{Q}_{l,t}^2)\rightarrow \kappa_l
%\]
in probability as $n\rightarrow\infty$. By Conditions \ref{as:gambd} and \ref{as:al}, we know $\kappa_l$ is a almost surely bounded random variable. Corollary 3.1 of \cite{Hall1980_supp} 
implies
\[
\frac{1}{\sqrt{n-m}}\sum_{t=m+1}^n\mathring{Q}_{l,t}\rightarrow\sqrt{\kappa_l}\cdot Z
\]
in distribution as $n\rightarrow\infty$, where $Z$ is a standard normally distributed random variable independent of $\kappa_l$. By \eqref{eq:asymbd8} and \eqref{eq:asymbd9}, due to  
\[
\frac{1}{n-m}\sum_{t=m+1}^n\frac{\partial\hat{f}_t^{(l)}(\bar{\theta}_l,\check{\btheta}_{-l})}{\partial\theta_l}\rightarrow1
\]
in probability which is obtained in \eqref{eq:asym13}, it holds that
 \begin{align*}
\sqrt{n|\mathcal{S}_l|}(\hat{\theta}_l-\theta_{0,l})\rightarrow \sqrt{\kappa_l}\cdot Z
 \end{align*}
in distribution as $n\rightarrow\infty$, provided that $\sqrt{n|\mathcal{S}_l|}\max\{q\Delta_n^{3/2},q^2\Delta_n^2\}=o(1)$. We complete the proof of Theorem \ref{tm:asymnormal}. $\hfill\Box$

\subsection{Proofs of auxiliary lemmas}

\subsubsection{Proof of Lemma \ref{la:uniformcov}}\label{sec:pflauniformcov}
Without loss of generality, we assume $\bTheta=[-C,C]^q$ for some constant $C>0$. For given $\epsilon>0$ which will be specified later, we partition $[-C,C]$ into $K=\lceil 2C/\epsilon\rceil$ sub-intervals $B_1,\ldots,B_K$ with equal length, where the length of each $B_k$ does not exceed $\epsilon$. Based on such defined $B_1,\ldots,B_K$, we can partition $\bTheta$ as follows:
\[
\bTheta=\bigcup_{k_1=1}^K\cdots\bigcup_{k_q=1}^KB_{k_1}\times\cdots \times B_{k_q}\,,
\]
which includes $K^q$ hyper-rectangles $\mathcal{B}_1,\ldots,\mathcal{B}_{K^q}$. For each given $u\in[K^q]$, there exists $(k_{1,u},\ldots,k_{q,u})\in[K]^q$ such that $\mathcal{B}_u=B_{k_{1,u}}\times\cdots\times B_{k_{q,u}}$. Let $\btheta_u$ be the center of $\mathcal{B}_u$. 

For each $\btheta \in \mathcal{B}_u$, since $\gamma_{i,j}^{t-1}(\btheta)$ only depends on $\theta_l$ with $l\in\mathcal{I}_{i,j}$, it follows from the Taylor expansion that  
\begin{align*}
&\hat{\ell}_{n,p}^{(l)}(\btheta)-\hat{\ell}_{n,p}^{(l)}(\btheta_u)\\
&~~~~~~~=\bigg[\frac{1}{(n-m)|\mathcal{S}_l|}\sum_{t=m+1}^n\sum_{(i,j)\in\mathcal{S}_l}\frac{X_{i,j}^{t}-\gamma_{i,j}^{t-1}(\tilde{\btheta}_u)}{\gamma_{i,j}^{t-1}(\tilde{\btheta}_u)\{1-\gamma_{i,j}^{t-1}(\tilde{\btheta}_u)\}}\frac{\partial\gamma_{i,j}^{t-1}(\tilde{\btheta}_u)}{\partial\btheta^{\T}}\bigg](\btheta-\btheta_u)\\
&~~~~~~~=\frac{1}{(n-m)|\mathcal{S}_l|}\sum_{t=m+1}^n\sum_{(i,j)\in\mathcal{S}_l}\frac{X_{i,j}^{t}-\gamma_{i,j}^{t-1}(\tilde{\btheta}_u)}{\gamma_{i,j}^{t-1}(\tilde{\btheta}_u)\{1-\gamma_{i,j}^{t-1}(\tilde{\btheta}_u)\}}\frac{\partial\gamma_{i,j}^{t-1}(\tilde{\btheta}_u)}{\partial\btheta_{\mathcal{I}_{i,j}}^{\T}}(\btheta_{\mathcal{I}_{i,j}}-\btheta_{u,\mathcal{I}_{i,j}})\,,
\end{align*}
where $\tilde{\btheta}_u \in \mathcal{B}_u$ is on the joint line between $\btheta$ and $\btheta_u$. Write $\btheta=(\theta_1,\ldots,\theta_q)^{\T}$ and $\btheta_u=(\theta_{u,1},\ldots,\theta_{u,q})^{\T}$. 
By Conditions \ref{as:gambd} and \ref{as:bdSij}, it holds that 
\begin{align*}
    &|\hat{\ell}_{n,p}^{(l)}(\btheta) - \hat{\ell}_{n,p}^{(l)}(\btheta_u)|
       \\&~~~~~~ \leq \frac{1}{(n-m)|\mathcal{S}_l|}\sum_{t = m+1}^n \sum_{(i,j)\in\mathcal{S}_l}\sum_{l\in\mathcal{I}_{i,j}} \bigg| \frac{X_{i,j}^{t}- \gamma_{i,j}^{t-1}(\tilde{ \btheta}_u)}{\gamma_{i,j}^{t-1}(\tilde{\btheta}_u)\{1- \gamma_{i,j}^{t-1}(\tilde{ \btheta}_u)\} }\bigg|\bigg|\frac{\partial \gamma_{i,j}^{t-1}(\tilde{\btheta}_u)}{\partial \theta_l}\bigg| |\theta_l - \theta_{u,l}| 
       \\&~~~~~~ \leq \frac{C_1^{-1}C_2\epsilon}{(n-m)|\mathcal{S}_l|}\sum_{t = m+1}^n \sum_{(i,j)\in\mathcal{S}_l}\sum_{l\in\mathcal{I}_{i,j}}  |X_{i,j}^{t}- \gamma_{i,j}^{t-1}(\tilde{\btheta}_u)|\leq sC_1^{-1}C_2\epsilon
\end{align*}
for any $\btheta\in\mathcal{B}_u$ and $l\in[q]$. Analogously, we also have 
\[
\sup_{\btheta\in\mathcal{B}_u}|\ell_{n,p}^{(l)}(\btheta)-\ell_{n,p}^{(l)}(\btheta_u)|\leq sC_1^{-1}C_2\epsilon\,.
\]
Therefore, by the triangle inequality, it holds that
\begin{align}\label{eq:unicov1}
     &\max_{l\in\mathcal{H}}\sup_{\btheta \in \bTheta}|\hat{\ell}_{n,p}^{(l)}(\btheta) - \ell_{n,p}^{(l)}(\btheta)| =\max_{l\in\mathcal{H}}\max_{u\in[K^q]} \sup_{\btheta \in \mathcal{B}_u} |\hat{\ell}_{n,p}^{(l)}(\btheta) - \ell_{n,p}^{(l)}(\btheta)| \notag
     \\
     &~~~~~~\leq\max_{l\in\mathcal{H}}\max_{u\in[K^q]} |\hat{\ell}_{n,p}^{(l)}(\btheta_u) - \ell_{n,p}^{(l)}(\btheta_u)| +\max_{l\in\mathcal{H}}\max_{u\in[K^q]}\sup_{\btheta\in\mathcal{B}_u}|\hat{\ell}_{n,p}^{(l)}(\btheta) - \hat{\ell}_{n,p}^{(l)}(\btheta_u)|\notag\\
     &~~~~~~~~~~+\max_{l\in\mathcal{H}}\max_{u\in[K^q]}\sup_{\btheta\in\mathcal{B}_u}|{\ell}_{n,p}^{(l)}(\btheta)-{\ell}_{n,p}^{(l)}(\btheta_u)|\\
     &~~~~~~\leq  \max_{l\in\mathcal{H}}\max_{u\in[K^q]} |\hat{\ell}_{n,p}^{(l)}(\btheta_u) - \ell_{n,p}^{(l)}(\btheta_u)| +  2sC_1^{-1}C_2\epsilon\,.\notag
\end{align}
 For each $l\in[q]$, $u\in[K^q]$ and $t\in[n]\setminus[m]$, define
\begin{align*}
Q_{u,t}^{(l)}=&~\frac{1}{|\mathcal{S}_l|^{1/2}}\sum_{(i,j)\in\mathcal{S}_l} \bigg[ \{X_{i,j}^t - \gamma_{i,j}^{t-1}(\btheta_0)\}\log \bigg\{\frac{\gamma_{i,j}^{t-1}(\btheta_u)}{1- \gamma_{i,j}^{t-1}(\btheta_u)}\bigg\}  \bigg]\,.
%=&~\frac{2}{|\mathcal{S}_l|^{1/2}}\sum_{(i,j)\in\mathcal{S}_l:\,i<j} \bigg[ \{X_{i,j}^t - \gamma_{i,j}^{t-1}(\btheta_0)\}\log \bigg\{\frac{\gamma_{i,j}^{t-1}(\btheta_u)}{1- \gamma_{i,j}^{t-1}(\btheta_u)}\bigg\}  \bigg]\,.
\end{align*}
Then 
\begin{align}\label{eq:ellexpa}
    \hat{\ell}_{n,p}^{(l)}(\btheta_u) - \ell_{n,p}^{(l)}(\btheta_u) =\frac{1}{(n-m)|\mathcal{S}_l|^{1/2}}\sum_{t = m+1}^n Q_{u,t}^{(l)}\,.
\end{align}
%In the sequel, we consider the convergence rate of $\max_{l\in\mathcal{H}}\max_{u\in[K^q]} |\hat{\ell}_{n,p}^{(l)}(\btheta_u) - \ell_{n,p}^{(l)}(\btheta_u)|$.

%\subsubsection{Case 1: $n\rightarrow\infty$}
Notice that   
\begin{align*}
|X_{i,j}^t-\gamma_{i,j}^{t-1}(\btheta_0)|\bigg|\log\bigg\{\frac{\gamma_{i,j}^{t-1}(\btheta_u)}{1-\gamma_{i,j}^{t-1}(\btheta_u)}\bigg\}\bigg|\leq \max_{i\neq j}\bigg|\log\bigg\{\frac{\gamma_{i,j}^{t-1}(\btheta_u)}{1-\gamma_{i,j}^{t-1}(\btheta_u)}\bigg\}\bigg|\,.
\end{align*}
Due to 
\begin{align*}
\log\big[\gamma_{i,j}^{t-1}(\btheta_u)\{1-\gamma_{i,j}^{t-1}(\btheta_u)\}\big]<\log\bigg\{\frac{\gamma_{i,j}^{t-1}(\btheta_u)}{1-\gamma_{i,j}^{t-1}(\btheta_u)}\bigg\}<-\log\big[\gamma_{i,j}^{t-1}(\btheta_u)\{1-\gamma_{i,j}^{t-1}(\btheta_u)\}\big]\,,
\end{align*}
by Condition \ref{as:gambd}, we have
\begin{align}
&\max_{u\in[K^q]}\max_{t\in[n]\setminus[m]}|X_{i,j}^t-\gamma_{i,j}^{t-1}(\btheta_0)|\bigg|\log\bigg\{\frac{\gamma_{i,j}^{t-1}(\btheta_u)}{1-\gamma_{i,j}^{t-1}(\btheta_u)}\bigg\}\bigg| \notag\\
&~~~~~~~~~~~~~~~\leq\max_{u\in[K^q]}\max_{t\in[n]\setminus[m]}\max_{i\neq j}\big|\log\big[\gamma_{i,j}^{t-1}(\btheta_u)\{1-\gamma_{i,j}^{t-1}(\btheta_u)\}\big]\big|\leq \log C_1^{-1}\,. \label{eq:boundQ}
\end{align}
Denote by $\mathbb{P}_{\mathcal{F}_{t-1}}(\cdot)$ the conditional probability measure given $\mathcal{F}_{t-1}$ with the unknown true parameter vector $\btheta_0$. It follows from the Bernstein inequality that
\begin{align}\label{eq:condtailQ}
&\max_{u\in[K^q]}\max_{t\in[n]\setminus[m]}\mathbb{P}_{\mathcal{F}_{t-1}}\big\{|Q_{u,t}^{(l)}|\geq x\big\}\notag\\
&~~~~~~~~~~~~~~\leq 2\exp\bigg\{-\frac{3|\mathcal{S}_l|^{1/2}x^2}{6|\mathcal{S}_l|^{1/2}(\log C_1^{-1})^2+2x\log C_1^{-1}}\bigg\}
\end{align}
for any $x>0$. Furthermore, due to $\mathbb{P}\{|Q_{u,t}^{(l)}|\geq x\}=\mathbb{E}[\mathbb{P}_{\mathcal{F}_{t-1}}\{|Q_{u,t}^{(l)}|\geq x\}]$, we also have 
\begin{align}\label{eq:tailQ}
\max_{u\in[K^q]}\max_{t\in[n]\setminus[m]}\mathbb{P}\big\{|Q_{u,t}^{(l)}|\geq x\big\}\leq 2\exp\bigg\{-\frac{3|\mathcal{S}_l|^{1/2}x^2}{6|\mathcal{S}_l|^{1/2}(\log C_1^{-1})^2+2x\log C_1^{-1}}\bigg\}
\end{align}
for any $x>0$. Let
\begin{align*}
\tilde{Q}_{u,t}^{(l)}=&~Q_{u,t}^{(l)}I\big\{|Q_{u,t}^{(l)}|\leq M\big\}-\mathbb{E}_{\mathcal{F}_{t-1}}\big[Q_{u,t}^{(l)}I\big\{|Q_{u,t}^{(l)}|\leq M\big\}\big]%\,,\\
%\check{Q}_{u,t}^{(l)}=&~Q_{u,t}^{(l)}I\big\{|Q_{u,t}^{(l)}|> M\big\}-\mathbb{E}_{\mathcal{F}_{t-1},\btheta_0}\big[Q_{u,t}^{(l)}I\big\{|Q_{u,t}^{(l)}|>M\big\}\big]
\end{align*}
for some diverging $M>0$ specified later. Notice that $\{\tilde{Q}_{u,t}^{(l)}\}_{t\in[n]\setminus[m]}$ is a martingale difference sequence with $\max_{t\in[n]\setminus[m]}|\tilde{Q}_{u,t}^{(l)}|\leq 2M$. By the Azuma's inequality (\citeauthor{Azuma1967_supp}, \citeyear{Azuma1967_supp}; 
see also Theorem 3.1 of  \cite{Lesigne2001_supp}), 
we have
\begin{align}\label{eq:tailmart}
\max_{l\in[q]}\max_{u\in[K^q]}\mathbb{P}\bigg\{\bigg|\frac{1}{n-m}\sum_{t=m+1}^n\tilde{Q}_{u,t}^{(l)}\bigg|\geq x\bigg\}\leq 2\exp\bigg\{-\frac{(n-m)x^2}{8M^2}\bigg\}
\end{align}
for any $x>0$. By \eqref{eq:tailQ}, 
\begin{align*}
&\max_{u\in[K^q]}\mathbb{P}\bigg[\bigg|\frac{1}{n-m}\sum_{t=m+1}^nQ_{u,t}^{(l)}I\big\{|Q_{u,t}^{(l)}|>M\big\}\bigg|>0\bigg]\leq \max_{u\in[K^q]}\mathbb{P}\bigg\{\max_{t\in[n]\setminus[m]}|Q_{u,t}^{(l)}|>M\bigg\}\\
&~~~~~~~~~\leq 2n\exp\bigg\{-\frac{3|\mathcal{S}_l|^{1/2}M^2}{6|\mathcal{S}_l|^{1/2}(\log C_1^{-1})^2+2M\log C_1^{-1}}\bigg\}\,.
\end{align*}
Together with \eqref{eq:tailmart}, by the Bonferroni inequality, we have
\begin{align*}
%&\mathbb{P}\bigg(\max_{l\in\mathcal{G}}\max_{u\in[K^q]}\bigg|\frac{1}{n-m}\sum_{t=m+1}^n\big[\tilde{Q}_{u,t}^{(l)}+Q_{u,t}^{(l)}I\big\{|Q_{u,t}^{(l)}|>M\big\}\big]\bigg|>x\bigg)\\
%&~~~~~~~~~~~\leq 2|\mathcal{G}|K^q\exp\bigg\{-\frac{(n-m)x^2}{8M^2}\bigg\}\\
%&~~~~~~~~~~~~~~+2|\mathcal{G}|K^qn\max_{l\in\mathcal{G}}\exp\bigg\{-\frac{3|\mathcal{S}_l|^{1/2}M^2}{12|\mathcal{S}_l|^{1/2}(\log C_1^{-1})^2+4M\log C_1^{-1}}\bigg\}\,,\\
&\mathbb{P}\bigg(\max_{l\in\mathcal{H}}\max_{u\in[K^q]}\bigg|\frac{1}{n-m}\sum_{t=m+1}^n\big[\tilde{Q}_{u,t}^{(l)}+Q_{u,t}^{(l)}I\big\{|Q_{u,t}^{(l)}|>M\big\}\big]\bigg|>x\bigg)\\
&~~~~~~~~~~~\leq 2|\mathcal{H}|K^q\exp\bigg\{-\frac{(n-m)x^2}{8M^2}\bigg\}\\
&~~~~~~~~~~~~~~+2|\mathcal{H}|K^qn\max_{l\in\mathcal{H}}\exp\bigg\{-\frac{3|\mathcal{S}_l|^{1/2}M^2}{6|\mathcal{S}_l|^{1/2}(\log C_1^{-1})^2+2M\log C_1^{-1}}\bigg\}
\end{align*}
for any $x>0$. Recall $S_{\mathcal{H},\min}=\min_{l\in\mathcal{H}}|\mathcal{S}_l|$ and $|\mathcal{H}|\leq q$. Selecting 
\begin{align}\label{eq:Msele}
M= C\max\bigg(\sqrt{q\log K},\,\frac{q\log K}{\sqrt{S_{\mathcal{H},\min}}},\,\sqrt{\log n},\,\frac{\log n}{\sqrt{S_{\mathcal{H},\min}}}\bigg)
\end{align}
for some sufficiently large constant $C>0$, it holds that 
\begin{align}\label{eq:uniconrate1}
&\max_{l\in\mathcal{H}}\max_{u\in[K^q]}\bigg|\frac{1}{n-m}\sum_{t=m+1}^n\big[\tilde{Q}_{u,t}^{(l)}+Q_{u,t}^{(l)}I\big\{|Q_{u,t}^{(l)}|>M\big\}\big]\bigg|=O_{\rm p}\bigg(M\sqrt{\frac{q\log K}{n}}\bigg)\notag\\
&~~~~~~~~=O_{\rm p}\bigg(\frac{q\log K}{\sqrt{n}}\bigg)+O_{\rm p}\bigg\{\frac{(q\log K)^{3/2}}{\sqrt{nS_{\mathcal{H},\min}}}\bigg\}\\
&~~~~~~~~~~~+O_{\rm p}\bigg\{\frac{\sqrt{q(\log K)(\log n)}}{\sqrt{n}}\bigg\}+O_{\rm p}\bigg\{\frac{(q\log K)^{1/2}(\log n)}{\sqrt{nS_{\mathcal{H},\min}}}\bigg\}\,.\notag
\end{align}
On the other hand, by \eqref{eq:condtailQ}, it holds that
\begin{align*}
&\mathbb{E}_{\mathcal{F}_{t-1}}\big[|Q_{u,t}^{(l)}|I\big\{|Q_{u,t}^{(l)}|>M\big\}\big]=M\mathbb{P}_{\mathcal{F}_{t-1}}\big\{|Q_{u,t}^{(l)}|>M\big\}+\int_{M}^\infty\mathbb{P}_{\mathcal{F}_{t-1}}(|Q_{u,t}|>x)\,{\rm d}x\\
&~~~~~~\leq 2M\exp\bigg\{-\frac{3|\mathcal{S}_l|^{1/2}M^2}{6|\mathcal{S}_l|^{1/2}(\log C_1^{-1})^2+2M\log C_1^{-1}}\bigg\}\\
&~~~~~~~~~+2\int_M^\infty\exp\bigg\{-\frac{x^2}{4(\log C_1^{-1})^2}\bigg\}\,{\rm d}x+2\int_M^{\infty}\exp\bigg(-\frac{3|\mathcal{S}_l|^{1/2}x}{4\log C_1^{-1}}\bigg)\,{\rm d}x\\
&~~~~~~\leq 2M\exp\bigg\{-\frac{3|\mathcal{S}_l|^{1/2}M^2}{6|\mathcal{S}_l|^{1/2}(\log C_1^{-1})^2+2M\log C_1^{-1}}\bigg\}\\
&~~~~~~~~~+4(\log C_1^{-1})^2M^{-1}\exp\bigg\{-\frac{M^2}{4(\log C_1^{-1})^2}\bigg\}+\frac{8\log C_1^{-1}}{3|\mathcal{S}_l|^{1/2}}\exp\bigg(-\frac{3|\mathcal{S}_l|^{1/2}M}{4\log C_1^{-1}}\bigg)\,,
\end{align*}
which implies
\begin{align}
&\max_{u\in[K^q]}\bigg|\frac{1}{n-m}\sum_{t=m+1}^n\mathbb{E}_{\mathcal{F}_{t-1}}\big[Q_{u,t}^{(l)}I\big\{|Q_{u,t}^{(l)}|>M\big\}\big]\bigg|\notag\\
&~~~~~~~~\leq 2M\exp\bigg\{-\frac{3|\mathcal{S}_l|^{1/2}M^2}{6|\mathcal{S}_l|^{1/2}(\log C_1^{-1})^2+2M\log C_1^{-1}}\bigg\}\label{eq:bdcond1}\\
&~~~~~~~~~~~+4(\log C_1^{-1})^2\exp\bigg\{-\frac{M^2}{4(\log C_1^{-1})^2}\bigg\}+\frac{8\log C_1^{-1}}{3|\mathcal{S}_l|^{1/2}}\exp\bigg(-\frac{3|\mathcal{S}_l|^{1/2}M}{4\log C_1^{-1}}\bigg)\,.\notag
\end{align}
Due to $S_{\mathcal{H},\min}=\min_{l\in\mathcal{H}}|\mathcal{S}_l|\rightarrow\infty$ and $M\rightarrow\infty$ satisfying \eqref{eq:Msele}, by \eqref{eq:bdcond1}, we have
\begin{align*}
&\max_{l\in\mathcal{H}}\max_{u\in[K^q]}\bigg|\frac{1}{n-m}\sum_{t=m+1}^n\mathbb{E}_{\mathcal{F}_{t-1}}\big[Q_{u,t}^{(l)}I\big\{|Q_{u,t}^{(l)}|>M\big\}\big]\bigg|\\
&~~~~~~~~~~~\lesssim \exp(-\tilde{C}_1M^2)+\exp(-\tilde{C}_2MS_{\mathcal{H},\min}^{1/2})=o(n^{-1/2})\,,
\end{align*}
where $\tilde{C}_1>0$ and $\tilde{C}_2>0$ are two universal constants. Due to $Q_{u,t}^{(l)}=\tilde{Q}_{u,t}^{(l)}+Q_{u,t}^{(l)}I\{|Q_{u,t}^{(l)}|>M\}-\mathbb{E}_{\mathcal{F}_{t-1},\btheta_0}[Q_{u,t}^{(l)}I\{|Q_{u,t}^{(l)}|>M\}]$, together with \eqref{eq:uniconrate1}, it holds that 
\begin{align}
\max_{l\in\mathcal{H}}\max_{u\in[K^q]}\bigg|\frac{1}{n-m}\sum_{t=m+1}^nQ_{u,t}^{(l)}\bigg|=&~O_{\rm p}\bigg(\frac{q\log K}{\sqrt{n}}\bigg)+O_{\rm p}\bigg\{\frac{(q\log K)^{3/2}}{\sqrt{nS_{\mathcal{H},\min}}}\bigg\}\label{eq:toprovebd1}\\
&+O_{\rm p}\bigg\{\frac{\sqrt{q(\log K)(\log n)}}{\sqrt{n}}\bigg\}+O_{\rm p}\bigg\{\frac{(q\log K)^{1/2}(\log n)}{\sqrt{nS_{\mathcal{H},\min}}}\bigg\}\,.\notag
\end{align}
By \eqref{eq:ellexpa}, we have 
\begin{align*}
\max_{l\in\mathcal{H}}\max_{u\in[K^q]}|\hat{\ell}_{n,p}^{(l)}(\btheta_u) - \ell_{n,p}^{(l)}(\btheta_u)|=&~O_{\rm p}\bigg(\frac{q\log K}{\sqrt{nS_{\mathcal{H},\min}}}\bigg)+O_{\rm p}\bigg\{\frac{(q\log K)^{3/2}}{\sqrt{n}S_{\mathcal{H},\min}}\bigg\}\\
&+O_{\rm p}\bigg\{\frac{\sqrt{q(\log K)(\log n)}}{\sqrt{nS_{\mathcal{H},\min}}}\bigg\}+O_{\rm p}\bigg\{\frac{(q\log K)^{1/2}(\log n)}{\sqrt{n}S_{\mathcal{H},\min}}\bigg\}\,.
\end{align*}
Recall $K\asymp \epsilon^{-1}$. Together with \eqref{eq:unicov1}, it holds that 
\begin{align*}
\max_{l\in\mathcal{H}}\sup_{\btheta \in \bTheta}|\hat{\ell}_{n,p}^{(l)}(\btheta) - \ell_{n,p}^{(l)}(\btheta)|\leq&~2sC_1^{-1}C_2\epsilon+O_{\rm p}\bigg(\frac{q\log K}{\sqrt{nS_{\mathcal{H},\min}}}\bigg)+O_{\rm p}\bigg\{\frac{(q\log K)^{3/2}}{\sqrt{n}S_{\mathcal{H},\min}}\bigg\}\\
&+O_{\rm p}\bigg\{\frac{\sqrt{q(\log K)(\log n)}}{\sqrt{nS_{\mathcal{H},\min}}}\bigg\}+O_{\rm p}\bigg\{\frac{(q\log K)^{1/2}(\log n)}{\sqrt{n}S_{\mathcal{H},\min}}\bigg\}\,.
\end{align*}
Due to $s\leq q$, with selecting $\epsilon\asymp n^{-1/2}S_{\mathcal{H},\min}^{-1/2}$, we have
\begin{align*}
\max_{l\in\mathcal{H}}\sup_{\btheta \in \bTheta}|\hat{\ell}_{n,p}^{(l)}(\btheta) - \ell_{n,p}^{(l)}(\btheta)|=O_{\rm p}\bigg\{\frac{q\log(nS_{\mathcal{H},\min})}{\sqrt{nS_{\mathcal{H},\min}}}\bigg\}+O_{\rm p}\bigg\{\frac{q^{3/2}\log^{3/2}(nS_{\mathcal{H},\min})}{\sqrt{n}S_{\mathcal{H},\min}}\bigg\}\,.
\end{align*}
We complete the proof of Lemma \ref{la:uniformcov}. $\hfill\Box$

\subsubsection{Proof of Lemma \ref{la:asym1}}\label{pf:laasymp1}

By the Taylor expansion, we have
\begin{align*}
&\frac{1}{n-m}\sum_{t=m+1}^n\big\{\hat{f}_t^{(l)}(\theta_{0,l},\tilde{\btheta}_{-l})-\hat{f}_t^{(l)}(\btheta_0)\big\}\\
&~~~~~~~~~~=\frac{1}{n-m}\sum_{t=m+1}^n\frac{\partial\hat{f}_t^{(l)}(\theta_{0,l},\bar{\btheta}_{-l})}{\partial\btheta_{-l}^{\T}}(\tilde{\btheta}_{-l}-\btheta_{0,-l})\\
&~~~~~~~~~~=\underbrace{\frac{1}{n-m}\sum_{t=m+1}^n\frac{\partial\hat{f}_t^{(l)}(\theta_{0,l},\tilde{\btheta}_{-l})}{\partial\btheta_{-l}^{\T}}(\tilde{\btheta}_{-l}-\btheta_{0,-l})}_{R_{1,l}}\\
&~~~~~~~~~~~~~+\underbrace{\frac{1}{n-m}\sum_{t=m+1}^n\bigg\{\frac{\partial\hat{f}_t^{(l)}(\theta_{0,l},\bar{\btheta}_{-l})}{\partial\btheta_{-l}^{\T}}-\frac{\partial\hat{f}_t^{(l)}(\theta_{0,l},\tilde{\btheta}_{-l})}{\partial\btheta_{-l}^{\T}}\bigg\}(\tilde{\btheta}_{-l}-\btheta_{0,-l})}_{R_{2,l}}\,,
\end{align*}
where $\bar{\btheta}_{-l}$ is on the joint line between $\btheta_{0,-l}$ and $\tilde{\btheta}_{-l}$. %In the sequel, we first consider the convergence rates of $|R_{1,l}|$ and $|R_{2,l}|$.

For $R_{1,l}$, by the Taylor expansion, it holds that
\begin{align*}
R_{1,l}=&~\frac{1}{n-m}\sum_{t=m+1}^n\frac{\partial\hat{f}_t^{(l)}(\tilde{\btheta})}{\partial\btheta_{-l}^{\T}}(\tilde{\btheta}_{-l}-\btheta_{0,-l})\\
&+(\theta_{0,l}-\tilde{\theta}_l)\bigg\{\frac{1}{n-m}\sum_{t=m+1}^n\frac{\partial^2\hat{f}_t^{(l)}(\bar{\theta}_l,\tilde{\btheta}_{-l})}{\partial\theta_l\partial\btheta_{-l}^{\T}}\bigg\}(\tilde{\btheta}_{-l}-\btheta_{0,-l})\,,
\end{align*}
where $\bar{\theta}_l$ is on the joint line between $\theta_{0,l}$ and $\tilde{\theta}_l$. Recall $\hat{f}_t^{(l)}(\btheta)=\hat{\bvarphi}_{l}^{\T}\bg_t^{(l)}(\btheta)$. By the definition of $\hat{\bvarphi}_{l}$ given in \eqref{eq:hatanl}, we have 
\[
\bigg|\frac{1}{n-m}\sum_{t=m+1}^n\frac{\partial\hat{f}_t^{(l)}(\tilde{\btheta})}{\partial\btheta_{-l}}\bigg|_\infty\leq \tau\,,
\]
which implies
\begin{align}
\bigg|\frac{1}{n-m}\sum_{t=m+1}^n\frac{\partial\hat{f}_t^{(l)}(\tilde{\btheta})}{\partial\btheta_{-l}^{\T}}(\tilde{\btheta}_{-l}-\btheta_{0,-l})\bigg|\leq&~ \bigg|\frac{1}{n-m}\sum_{t=m+1}^n\frac{\partial\hat{f}_t^{(l)}(\tilde{\btheta})}{\partial\btheta_{-l}}\bigg|_\infty |\tilde{\btheta}_{-l}-\btheta_{0,-l}|_1\notag\\
\leq&~\tau |\tilde{\btheta}_{-l}-\btheta_{0,-l}|_1\,.\label{eq:asybd1}
\end{align}
For any $k\in[q]$, due to
\begin{align*}
\frac{\partial^2\hat{f}_t^{(l)}(\bar{\theta}_l,\tilde{\btheta}_{-l})}{\partial\theta_l\partial\theta_k}=\hat{\bvarphi}_l^{\T}\frac{\partial^2\bg_t^{(l)}(\bar{\theta}_l,\tilde{\btheta}_{-l})}{\partial\theta_l\partial\theta_k}\,,
\end{align*}
we then have
\begin{align*}
\bigg|\frac{\partial^2\hat{f}_t^{(l)}(\bar{\theta}_l,\tilde{\btheta}_{-l})}{\partial\theta_l\partial\theta_k}\bigg|\leq&~ |\hat{\bvarphi}_l|_1\bigg|\frac{\partial^2\bg_t^{(l)}(\bar{\theta}_l,\tilde{\btheta}_{-l})}{\partial\theta_l\partial\theta_k}\bigg|_\infty\\
\leq&~\frac{|\hat{\bvarphi}_l|_1}{|\mathcal{S}_l|}\sum_{(i,j)\in\mathcal{S}_l}\bigg|\frac{\partial^2}{\partial\theta_l\partial\theta_k}\bigg[\frac{X_{i,j}^t-\gamma_{i,j}^{t-1}(\btheta)}{\gamma_{i,j}^{t-1}(\btheta)\{1-\gamma_{i,j}^{t-1}(\btheta)\}}\frac{\partial\gamma_{i,j}^{t-1}(\btheta)}{\partial\btheta}\bigg]\bigg|_\infty
\end{align*}
Notice that
\begin{align*}
&\frac{\partial^2}{\partial\theta_{l_1}\partial\theta_{l_2}}\bigg[\frac{X_{i,j}^t-\gamma_{i,j}^{t-1}(\btheta)}{\gamma_{i,j}^{t-1}(\btheta)\{1-\gamma_{i,j}^{t-1}(\btheta)\}}\frac{\partial\gamma_{i,j}^{t-1}(\btheta)}{\partial\theta_{l_3}}\bigg]\\
&~~~~~~=2\bigg[\frac{1-2\gamma_{i,j}^{t-1}(\btheta)}{\{\gamma_{i,j}^{t-1}(\btheta)\}^2\{1-\gamma_{i,j}^{t-1}(\btheta)\}^2}+\frac{X_{i,j}^t-\gamma_{i,j}^{t-1}(\btheta)}{\{\gamma_{i,j}^{t-1}(\btheta)\}^2\{1-\gamma_{i,j}^{t-1}(\btheta)\}^2}\notag\\
&~~~~~~~~~~~~~~+\frac{\{X_{i,j}^t-\gamma_{i,j}^{t-1}(\btheta)\}\{1-2\gamma_{i,j}^{t-1}(\btheta)\}^2}{\{\gamma_{i,j}^{t-1}(\btheta)\}^3\{1-\gamma_{i,j}^{t-1}(\btheta)\}^3}\bigg]\frac{\partial\gamma_{i,j}^{t-1}(\btheta)}{\partial\theta_{l_1}}\frac{\partial\gamma_{i,j}^{t-1}(\btheta)}{\partial\theta_{l_2}}\frac{\partial\gamma_{i,j}^{t-1}(\btheta)}{\partial\theta_{l_3}}\notag\\
&~~~~~~~~~-\bigg[\frac{1}{\gamma_{i,j}^{t-1}(\btheta)\{1-\gamma_{i,j}^{t-1}(\btheta)\}}+\frac{\{X_{i,j}^t-\gamma_{i,j}^{t-1}(\btheta)\}\{1-2\gamma_{i,j}^{t-1}(\btheta)\}}{\{\gamma_{i,j}^{t-1}(\btheta)\}^2\{1-\gamma_{i,j}^{t-1}(\btheta)\}^2}\bigg]
% \label{eq:th3der}
\\
&~~~~~~~~~~~~~~~\times\bigg\{\frac{\partial^2\gamma_{i,j}^{t-1}(\btheta)}{\partial\theta_{l_1}\partial\theta_{l_2}}\frac{\partial\gamma_{i,j}^{t-1}(\btheta)}{\partial\theta_{l_3}}+\frac{\partial^2\gamma_{i,j}^{t-1}(\btheta)}{\partial\theta_{l_2}\partial\theta_{l_3}}\frac{\partial\gamma_{i,j}^{t-1}(\btheta)}{\partial\theta_{l_1}}+\frac{\partial^2\gamma_{i,j}^{t-1}(\btheta)}{\partial\theta_{l_3}\partial\theta_{l_1}}\frac{\partial\gamma_{i,j}^{t-1}(\btheta)}{\partial\theta_{l_2}}\bigg\}\notag\\
&~~~~~~~~~+\frac{X_{i,j}^t-\gamma_{i,j}^{t-1}(\btheta)}{\gamma_{i,j}^{t-1}(\btheta)\{1-\gamma_{i,j}^{t-1}(\btheta)\}}\frac{\partial^3\gamma_{i,j}^{t-1}(\btheta)}{\partial\theta_{l_1}\partial\theta_{l_2}\partial\theta_{l_3}}\,.
\end{align*}
By the triangle inequality and Condition \ref{as:gambd},  we know 
\begin{align}\label{eq:asybd2}
\max_{k\in[q]}\bigg|\frac{\partial^2\hat{f}_t^{(l)}(\bar{\theta}_l,\tilde{\btheta}_{-l})}{\partial\theta_l\partial\theta_k}\bigg|\leq C_*|\hat{\bvarphi}_l|_1
\end{align}
for some universal constant $C_*$ specified in \eqref{eq:h3bd}, which implies
\begin{align*}
&\bigg|(\theta_{0,l}-\tilde{\theta}_l)\bigg\{\frac{1}{n-m}\sum_{t=m+1}^n\frac{\partial^2\hat{f}_t^{(l)}(\bar{\theta}_l,\tilde{\btheta}_{-l})}{\partial\theta_l\partial\btheta_{-l}^{\T}}\bigg\}(\tilde{\btheta}_{-l}-\btheta_{0,-l})\bigg|\\
&~~~~~~~~\leq |\theta_{0,l}-\tilde{\theta}_l|\bigg|\frac{1}{n-m}\sum_{t=m+1}^n\frac{\partial^2\hat{f}_t^{(l)}(\bar{\theta}_l,\tilde{\btheta}_{-l})}{\partial\theta_l\partial\btheta_{-l}^{\T}}\bigg|_\infty|\tilde{\btheta}_{-l}-\btheta_{0,-l}|_1\\
&~~~~~~~~\leq C_*|\hat{\bvarphi}_l|_1 |\theta_{0,l}-\tilde{\theta}_l||\tilde{\btheta}_{-l}-\btheta_{0,-l}|_1\,.
\end{align*}
Together with \eqref{eq:asybd1}, it holds that 
$
|R_{1,l}|\leq \tau |\tilde{\btheta}_{-l}-\btheta_{0,-l}|_1+C_*|\hat{\bvarphi}_l|_1 |\theta_{0,l}-\tilde{\theta}_l||\tilde{\btheta}_{-l}-\btheta_{0,-l}|_1$. 

For $R_{2,l}$, by the Taylor expansion, we have
\begin{align*}
R_{2,l}=(\bar{\btheta}_{-l}-\tilde{\btheta}_{-l})^{\T}\bigg\{\frac{1}{n-m}\sum_{t=m+1}^n\frac{\partial^2\hat{f}_t^{(l)}(\theta_{0,l},\dot{\btheta}_{-l})}{\partial\btheta_{-l}\partial\btheta_{-l}^{\T}}\bigg\}(\tilde{\btheta}_{-l}-\btheta_{0,-l})
\end{align*}
for some $\dot{\btheta}_{-l}$ on the joint line between $\bar{\btheta}_{-l}$ and $\tilde{\btheta}_{-l}$, which implies
\[
|R_{2,l}|\leq |\bar{\btheta}_{-l}-\tilde{\btheta}_{-l}|_1|\tilde{\btheta}_{-l}-\btheta_{0,-l}|_1\bigg|\frac{1}{n-m}\sum_{t=m+1}^n\frac{\partial^2\hat{f}_t^{(l)}(\theta_{0,l},\dot{\btheta}_{-l})}{\partial\btheta_{-l}\partial\btheta_{-l}^{\T}}\bigg|_\infty\,.
\]
Since $\bar{\btheta}_{-l}$ is on the joint line between $\btheta_{0,-l}$ and $\tilde{\btheta}_{-l}$, then $|\bar{\btheta}_{-l}-\tilde{\btheta}_{-l}|_1\leq |\tilde{\btheta}_{-l}-\btheta_{0,-l}|_1$. Parallel to \eqref{eq:asybd2}, we can also show
\[
\bigg|\frac{1}{n-m}\sum_{t=m+1}^n\frac{\partial^2\hat{f}_t^{(l)}(\theta_{0,l},\dot{\btheta}_{-l})}{\partial\btheta_{-l}\partial\btheta_{-l}^{\T}}\bigg|_\infty\leq C_*|\hat{\bvarphi}_l|_1\,.
\]
Hence, 
$
|R_{2,l}|\leq C_*|\hat{\bvarphi}_l|_1|\tilde{\btheta}_{-l}-\btheta_{0,-l}|_1^2$. Then
\begin{align*}
&\bigg|\frac{1}{n-m}\sum_{t=m+1}^n\big\{\hat{f}_t^{(l)}(\theta_{0,l},\tilde{\btheta}_{-l})-\hat{f}_t^{(l)}(\btheta_0)\big\}\bigg|\leq |R_{1,l}|+|R_{2,l}|\\
&~~~~~~~~~\leq \tau |\tilde{\btheta}_{-l}-\btheta_{0,-l}|_1+C_*|\hat{\bvarphi}_l|_1 |\theta_{0,l}-\tilde{\theta}_l||\tilde{\btheta}_{-l}-\btheta_{0,-l}|_1+C_*|\hat{\bvarphi}_l|_1|\tilde{\btheta}_{-l}-\btheta_{0,-l}|_1^2\\
&~~~~~~~~~= \tau |\tilde{\btheta}_{-l}-\btheta_{0,-l}|_1+C_*|\hat{\bvarphi}_l|_1 |\tilde{\btheta}-\btheta_{0}|_1|\tilde{\btheta}_{-l}-\btheta_{0,-l}|_1\\
&~~~~~~~~~\leq \tau |\tilde{\btheta}-\btheta_{0}|_1+C_*|\hat{\bvarphi}_l|_1 |\tilde{\btheta}-\btheta_{0}|_1^2\,.
\end{align*}
We complete the proof of Lemma \ref{la:asym1}. $\hfill\Box$

\subsubsection{Proof of Lemma \ref{la:asym2}}\label{pf:laasymp2}

Due to $\hat{f}_t^{(l)}(\btheta)=\hat{\bvarphi}_l^{\T}\bg_t^{(l)}(\btheta)$ and $f_t^{(l)}(\btheta)=\bvarphi_l^{\T}\bg_t^{(l)}(\btheta)$, then
\begin{align*}
\bigg|\frac{1}{n-m}\sum_{t=m+1}^n\big\{\hat{f}_t^{(l)}(\btheta_0)-f_t^{(l)}(\btheta_0)\big\}\bigg|\leq |\hat{\bvarphi}_l-\bvarphi_l|_1\bigg|\frac{1}{n-m}\sum_{t=m+1}^n\bg_t^{(l)}(\btheta_0)\bigg|_\infty\,.
\end{align*}
Write $(G_1^{(l)},\ldots,G_q^{(l)})^{\T}=(n-m)^{-1}\sum_{t=m+1}^n\bg_t^{(l)}(\btheta_0)$. Due to 
\[
\bg_t^{(l)}(\btheta_0)=\frac{1}{|\mathcal{S}_l|}\sum_{(i,j)\in\mathcal{S}_l}\frac{X_{i,j}^t-\gamma_{i,j}^{t-1}(\btheta_0)}{\gamma_{i,j}^{t-1}(\btheta_0)\{1-\gamma_{i,j}^{t-1}(\btheta_0)\}}\frac{\partial\gamma_{i,j}^{t-1}(\btheta_0)}{\partial\btheta}\,,
\]
we then have
\begin{align*}
G_k^{(l)}=&~\frac{1}{(n-m)|\mathcal{S}_l|}\sum_{t=m+1}^n\sum_{(i,j)\in\mathcal{S}_l}\frac{X_{i,j}^t-\gamma_{i,j}^{t-1}(\btheta_0)}{\gamma_{i,j}^{t-1}(\btheta_0)\{1-\gamma_{i,j}^{t-1}(\btheta_0)\}}\frac{\partial\gamma_{i,j}^{t-1}(\btheta_0)}{\partial\theta_k}\\
:=&~\frac{1}{(n-m)|\mathcal{S}_l|^{1/2}}\sum_{t=m+1}^n\check{Q}_{k,t}^{(l)}
\end{align*}
for any $k\in[q]$. Using the same arguments for deriving \eqref{eq:toprovebd1}, we can also show
\begin{align}\label{eq:asymbd3}
\max_{l\in\mathcal{H}}\max_{k\in[q]}|G_k^{(l)}|=O_{\rm p}\bigg\{\frac{\sqrt{(\log q)\log(qn)}}{\sqrt{nS_{\mathcal{H},\min}}}\bigg\}+O_{\rm p}\bigg\{\frac{(\log q)^{1/2}\log(qn)}{\sqrt{n}S_{\mathcal{H},\min}}\bigg\}\,,
\end{align}
where $S_{\mathcal{H},\min}=\min_{l\in\mathcal{H}}|\mathcal{S}_l|$. Together with Condition \ref{as:al}, it holds that
\begin{align}
&\max_{l\in\mathcal{H}}\bigg|\frac{1}{n-m}\sum_{t=m+1}^n\big\{\hat{f}_t^{(l)}(\btheta_0)-f_t^{(l)}(\btheta_0)\big\}\bigg|\notag\\
&~~~~~~~~=O_{\rm p}\bigg\{\frac{\omega_n\sqrt{(\log q)\log(qn)}}{\sqrt{nS_{\mathcal{H},\min}}}\bigg\}+O_{\rm p}\bigg\{\frac{\omega_n(\log q)^{1/2}\log(qn)}{\sqrt{n}S_{\mathcal{H},\min}}\bigg\}\,.\label{eq:asymbd4}
\end{align}
Notice that
\[
\frac{1}{n-m}\sum_{t=m+1}^nf_t^{(l)}(\btheta_0)=\frac{1}{(n-m)|\mathcal{S}_l|}\sum_{t=m+1}^n\sum_{(i,j)\in\mathcal{S}_l}\frac{X_{i,j}^t-\gamma_{i,j}^{t-1}(\btheta_0)}{\gamma_{i,j}^{t-1}(\btheta_0)\{1-\gamma_{i,j}^{t-1}(\btheta_0)\}}\bigg\{\bvarphi_l^{\T}\frac{\partial\gamma_{i,j}^{t-1}(\btheta_0)}{\partial\btheta}\bigg\}\,.
\]
By Conditions \ref{as:gambd} and \ref{as:al}, we know
\[
\max_{l\in[q]}\max_{t\in[n]\setminus[m]}\max_{(i,j)\in\mathcal{S}_l}\bigg|\bvarphi_l^{\T}\frac{\partial\gamma_{i,j}^{t-1}(\btheta_0)}{\partial\btheta}\bigg|\leq C_2C_4\,.
\]
Using the same arguments for deriving \eqref{eq:asymbd3}, we can also show
\begin{align*}
&\max_{l\in\mathcal{H}}\bigg|\frac{1}{n-m}\sum_{t=m+1}^nf_t^{(l)}(\btheta_0)\bigg|\\
&~~~~~~~~~=O_{\rm p}\bigg\{\frac{\sqrt{\log(1+|\mathcal{H}|)\log(n|\mathcal{H}|)}}{\sqrt{nS_{\mathcal{H},\min}}}\bigg\}+O_{\rm p}\bigg\{\frac{\sqrt{\log(1+|\mathcal{H}|)}\log(n|\mathcal{H}|)}{\sqrt{n}S_{\mathcal{H},\min}}\bigg\}\,.
\end{align*}
Together with \eqref{eq:asymbd4}, due to $\omega_n(\log q)^{1/2}\log(qn)=o(1)$, we have Lemma \ref{la:asym2}. $\hfill\Box$

\subsubsection{Proof of Lemma \ref{la:asym3}}\label{pf:laasymp3}

For any $\theta_l\in B(\tilde{\theta}_l,\tilde{r})$, by the Taylor expansion, we have
\begin{align*}
\frac{1}{n-m}\sum_{t=m+1}^n\bigg\{\frac{\partial\hat{f}_t^{(l)}(\theta_l,\tilde{\btheta}_{-l})}{\partial\theta_l}-\frac{\partial \hat{f}_t^{(l)}(\tilde{\btheta})}{\partial\theta_l}\bigg\}=\bigg\{\frac{1}{n-m}\sum_{t=m+1}^n\frac{\partial^2\hat{f}_t^{(l)}(\bar{\theta}_l,\tilde{\btheta}_{-l})}{\partial\theta_l^2}\bigg\}(\theta_l-\tilde{\theta}_l)
\end{align*}
for some $\bar{\theta}_l$ on the joint line between $\theta_l$ and $\tilde{\theta}_l$, which implies
\[
\bigg|\frac{1}{n-m}\sum_{t=m+1}^n\bigg\{\frac{\partial\hat{f}_t^{(l)}(\theta_l,\tilde{\btheta}_{-l})}{\partial\theta_l}-\frac{\partial \hat{f}_t^{(l)}(\tilde{\btheta})}{\partial\theta_l}\bigg\}\bigg|=\bigg|\frac{1}{n-m}\sum_{t=m+1}^n\frac{\partial^2\hat{f}_t^{(l)}(\bar{\theta}_l,\tilde{\btheta}_{-l})}{\partial\theta_l^2}\bigg||\theta_l-\tilde{\theta}_l|\,.
\]
Due to $\theta_l\in B(\tilde{\theta}_l,\tilde{r})$, by \eqref{eq:asybd2}, it holds that
\begin{align*}
\sup_{\theta_l\in B(\tilde{\theta}_l,\tilde{r})}\bigg|\frac{1}{n-m}\sum_{t=m+1}^n\bigg\{\frac{\partial\hat{f}_t^{(l)}(\theta_l,\tilde{\btheta}_{-l})}{\partial\theta_l}-\frac{\partial \hat{f}_t^{(l)}(\tilde{\btheta})}{\partial\theta_l}\bigg\}\bigg|\leq C_*\tilde{r}|\hat{\bvarphi}_l|_1\,.
\end{align*}
Since
\[
\frac{1}{n-m}\sum_{t=m+1}^n\frac{\partial\hat{f}_t^{(l)}(\tilde{\btheta})}{\partial\theta_l}=\frac{\hat{\bvarphi}_l^{\T}}{n-m}\sum_{t=m+1}^n\frac{\partial \bg_t^{(l)}(\tilde{\btheta})}{\partial\theta_l}\,,
\]
by \eqref{eq:hatanl}, we know
\[
\bigg|\frac{1}{n-m}\sum_{t=m+1}^n\frac{\partial\hat{f}_t^{(l)}(\tilde{\btheta})}{\partial\theta_l}-1\bigg|\leq \tau\,,
\]
which implies
\begin{align*}
\sup_{\theta_l\in B(\tilde{\theta}_l,\tilde{r})}\bigg|\frac{1}{n-m}\sum_{t=m+1}^n\frac{\partial\hat{f}_t^{(l)}(\theta_l,\tilde{\btheta}_{-l})}{\partial\theta_l}-1\bigg|\leq \tau+C_*\tilde{r}|\hat{\bvarphi}_l|_1\,.
\end{align*}
%It follows from Condition \ref{as:al} that $\max_{l\in[q]}|\hat{\bvarphi}_l|_1=O_{\rm p}(1)$. Notice that $r=o(1)$ and $\tau=o(1)$. Hence,
%\begin{align*}
%\max_{l\in[q]}\sup_{\theta_l\in B(\tilde{\theta}_l,r)}\bigg|\frac{1}{n-m}\sum_{t=m+1}^n\frac{\partial\hat{f}_t^{(l)}(\theta_l,\tilde{\btheta}_{-l})}{\partial\theta_l}-1\bigg|=o_{\rm p}(1)\,.
%\end{align*}
We complete the proof of Lemma \ref{la:asym3}. $\hfill\Box$

\setstretch{1.2}

\end{document}